\documentclass[12pt,reqno]{amsart}
\usepackage[a4paper,margin=30mm,headheight=15pt]{geometry}
\usepackage{amsmath,amssymb,amsthm,mathtools,mathrsfs,tikz-cd,enumitem,booktabs,array,longtable,microtype,placeins}
\DeclareFontFamilySubstitution{TS1}{cmr}{lmr}
\usetikzlibrary{decorations.pathmorphing,arrows.meta}
\usepackage[hidelinks,bookmarksnumbered]{hyperref}
\definecolor{qiBlue}{HTML}{4E6087}
\definecolor{qiInk}{HTML}{343A48}
\definecolor{qiTeal}{HTML}{286F72}
\definecolor{qiAmber}{HTML}{A95724}

\setlist[enumerate]{leftmargin=*,itemsep=3pt,topsep=5pt}
\allowdisplaybreaks[2]
\numberwithin{equation}{section}
\newtheorem{theorem}{Theorem}[section]
\newtheorem{proposition}[theorem]{Proposition}
\newtheorem{lemma}[theorem]{Lemma}
\newtheorem{corollary}[theorem]{Corollary}
\theoremstyle{definition}
\newtheorem{definition}[theorem]{Definition}
\newtheorem{hypothesis}[theorem]{Hypothesis}

\theoremstyle{remark}
\newtheorem{remark}[theorem]{Remark}
\newcommand{\OO}{\mathcal O}\newcommand{\FF}{\mathcal F}\newcommand{\XX}{\mathcal X}
\newcommand{\MM}{\mathcal M}\newcommand{\EE}{\mathfrak E}\newcommand{\GG}{\mathbb G_m}
\newcommand{\CC}{\mathbb C}\newcommand{\PP}{\mathbb P}\renewcommand{\AA}{\mathbb A}\newcommand{\QQ}{\mathbb Q}
\newcommand{\LL}{\mathbf L}\newcommand{\RR}{\mathbf R}\newcommand{\vir}{\mathrm{vir}}
\newcommand{\rig}{\mathrm{rig}}\newcommand{\red}{\mathrm{red}}
\newcommand{\rel}{\mathrm{rel}}\newcommand{\perf}{\operatorname{Perf}}\newcommand{\Pic}{\operatorname{Pic}}
\newcommand{\Ext}{\operatorname{Ext}}\newcommand{\ext}{\operatorname{ext}}\newcommand{\Hom}{\operatorname{Hom}}\newcommand{\Tor}{\operatorname{Tor}}
\newcommand{\End}{\operatorname{End}}\newcommand{\Aut}{\operatorname{Aut}}\newcommand{\Supp}{\operatorname{Supp}}
\newcommand{\Spec}{\operatorname{Spec}}\newcommand{\ch}{\operatorname{ch}}\newcommand{\rk}{\operatorname{rk}}
\newcommand{\pd}{\operatorname{pd}}\newcommand{\depth}{\operatorname{depth}}\newcommand{\vdim}{\operatorname{vdim}}
\newcommand{\RHom}{\mathbf R\!\operatorname{Hom}}\newcommand{\RG}{\mathbf R\Gamma}
\newcommand{\fib}{\operatorname{fib}}\newcommand{\cofib}{\operatorname{cofib}}\newcommand{\tr}{\operatorname{tr}}
\newcommand{\spc}{\operatorname{sp}}\newcommand{\id}{\operatorname{id}}\newcommand{\coker}{\operatorname{coker}}
\newcommand{\length}{\operatorname{length}}
\newcommand{\dt}{\operatorname{DT}}
\hypersetup{pdftitle={Surface sheaves on the quintic I. Expanded degenerations and virtual gluing},pdfauthor={Artan Sheshmani}}
\newcommand{\companionref}[1]{\textup{\csname companion@#1\endcsname}}
\title[Surface sheaves on the quintic I]{Surface sheaves on the quintic I.\\ Expanded degenerations and virtual gluing}
\author{Artan Sheshmani}
\address{Artan Sheshmani}
\address{Beijing Institute of Mathematical Sciences and Applications, A6, Room 205 No. 544, Hefangkou Village, Huaibei Town, Huairou District, Beijing 101408, China}
\address{Massachusetts Institute of Technology (MIT), IAiFi Institute, 182 Memorial Drive, Cambridge, MA 02139, USA}
\address{National Research University, Higher School of Economics, Russian Federation, Laboratory of Mirror Symmetry, NRU HSE, 6 Usacheva Street, Moscow, Russia, 119048}
\email{artan@mit.edu}
\date{17 November 2026}
\begin{document}
\begin{abstract}
We study the specialization of Donaldson--Thomas virtual cycles for surface sheaves on a Calabi Yau threefold degenerating to a simple normal crossing threefold, composed of two quasi Fano threefolds glued along a common (double) divisor. The compactification of the two relative sheaf moduli spaces involves construction of suitable moduli space of sheaves on expanded degenerations of each quasi Fano variety. On a fixed expanded fiber the matching datum is an equivalence of the derived restrictions of the component sheaves. We identify the corresponding deformation complexes and give valuative and obstruction-theoretic hypotheses under which the specialized virtual class is a sum over marked joining divisors. When a marked space is the full fiber product of proper relative factors over a smooth simple boundary moduli space, its class is a refined diagonal pullback. For the case of sheaves supported on hyperplane sections of the threefold we discuss construction of the proper relative family by descending the chartwise stable loci, and identify both marked pieces, including the one with zero boundary. To provide some examples, we carry out explicit constructions for the moduli space of surface structure sheaf and its twists. We also prove valuative results for specified sheaves with noninvertible hulls. In a sequel, Part II, we use the construction of degeneration technique in this article to compare quintic invariants with the predicted S-duality modularity of \textit{M5-brane modified elliptic genus}  provided by string theorists, Gaiotto Strominger and Yin.
\end{abstract}\maketitle\tableofcontents
\section{Introduction}
Let $X$ be a smooth quintic threefold and let $H$ be its hyperplane class. For sheaves with supporting divisor in $|H|$, S-duality predicts the vector $h^{\mathrm{DT}}_1=(h^{\mathrm{DT}}_{1,\mu})_{\mu\in\mathbb Z/5}$ of generating series of the virtual degrees, indexed by the curve-charge residue $\mu$. To use a degeneration in this problem one needs to know precisely which sheaves and which obstruction spaces survive on the special fiber. A surface sheaf on $X$ can move in two ways. Its supporting divisor moves in a linear system, and the sheaf on that divisor moves even when the divisor is fixed. We ask what happens to the virtual class of this moduli problem when $X$ degenerates to two threefolds meeting along a K3 surface.

Write the central fiber as $Y_1\cup_S Y_2$. On the open set where a supporting surface cuts $S$ in a curve, the two component sheaves can be matched along that curve. In the closure of this open set the supporting divisor can contain $S$. The restriction of a sheaf to $S$ then remembers the kernel of multiplication by a local equation of $S$. This is the reason for keeping the derived restriction in the gluing problem.

There are two questions. On a given union, do the two sheaves together with their identification on $S$ recover the sheaf and its deformations? In a family, does a sheaf on the smooth fiber have a unique stable limit on some expansion, up to the automorphisms of the expansion? Theorem~\ref{thm:intro-general} answers the first question for the perfect pure sheaves in its statement. Theorem~\ref{thm:main} identifies the specialized virtual class under precise hypotheses on the second question and on the relative obstruction theories.

Li and Wu construct proper moduli of admissible quotients on expansions and use them to degenerate the invariants of ideal sheaves and stable pairs \cite[Theorems~4.14--4.15 and \S6]{LW}. A quotient retains its presentation. Here the object is the surface sheaf itself, whose supporting divisor moves through a linear system. In earlier work on a K3 fibration the surfaces supporting the sheaves were fibers of the fibration \cite{GSFiber}. In the quintic degeneration the surface can acquire the joining K3 as a component. On the fixed union the identification of the two derived restrictions recovers the sheaf even in this case. The same normalization triangle, applied to its endomorphisms, recovers its deformation theory. This is the step which carries the marked-divisor construction of Li and Wu to these surface-sheaf virtual classes. For the related question of degeneration for derived Quot schemes and categorified Donaldson--Thomas invariants, see \cite{KryczkaSheshmani2026}.

\textbf{Obstructions to generalization}: Two obstructions determine the scope of this extension. Corollary~\ref{cor:noseparatedall} concerns the sheaf. At $\OO_T^{\oplus r}$, $r\ge2$, quotienting by common scalars leaves $\mathrm{PGL}_r$ as stabilizer. This group cannot be a fiber of a proper diagonal. Keeping only stable sheaves in the same charge also loses limits, as the family in Proposition~\ref{prop:noordinaryproperextension} shows. We therefore formulate the virtual-cycle theorem for charges and stable loci satisfying the valuative criteria. In the primitive hyperplane class every pure sheaf on the smooth quintic already has rank one on an integral surface and is stable. The existence and uniqueness questions then concern what happens on the expansions.

Proposition~\ref{prop:unlabelledpolarization} concerns the target. On the one-bubble chart either invertible smoothing parameter can be deleted. Away from the locus where both vanish, a line bundle descending through those identifications is pulled back from the contraction. The omitted locus has codimension two in the regular chart, so the line bundle is determined there as well. Its degree on every contracted ruling is zero. Our descent datum is instead the open locus of stable sheaves for each chart polarization. Equality of these loci under the target isomorphisms, including in families, is the condition in \eqref{eq:generalstableoverlaps}. Section~\ref{sec:kimexactfull} relates this condition to the polarization and stability assumptions of Kim \cite[Assumptions~3.0.1 and~3.2.2]{Kim}.

\textbf{The remedy}. Fix $v=(0,H,\beta,n)$. On the smooth quintic every pure sheaf of this character is $i_*G$ for a rank-one torsion-free sheaf on an integral normal hyperplane section $i:T\hookrightarrow X$. Such a sheaf is stable. To follow it to the central fiber we keep the sheaves on the individual components of an expansion and identify their derived restrictions at each join. Theorem~\ref{thm:intro-general} recovers the sheaf from these data. Applied to its endomorphisms, the same construction recovers the deformation and obstruction maps. Both assertions hold on a fixed expansion for perfect pure sheaves in every divisor multiplicity.

For primitive $v$ the lower charges $\beta,n$ are arbitrary in the following calculation. With $\Pic(S)=\mathbb Z h$, conservation of the divisor class forces the successive boundary classes of a perfect pure limit to be $b_eh$ with $0\le b_0\le\cdots\le b_q\le1$. At a zero cut the boundary vanishes. At a positive cut the support on the branch to its right meets $S$ properly. A local equation of the join is therefore a non-zero-divisor on that branch sheaf, so the common derived restriction is a sheaf on $S$ (Corollary~\ref{cor:primitiveautomaticnew}). This sheaf can still have point torsion. Proposition~\ref{prop:route19-simple} shows that a positive restriction gives a point of the compactified Jacobian precisely when it is simple, equivalently pure. Thus one can test on the actual limit whether a Jacobian evaluation is available.

To construct such limits we write $R=G^{**}$ and retain the finite quotient $R/G$. On the expanded supports of Proposition~\ref{prop:divisorial-hull27} the divisor class of $R$ extends. Vertical twists place the extension in a stable chamber under the cut inequalities there. Once an admissible Cohen--Macaulay hull satisfying the hypotheses of Proposition~\ref{prop:hull-punctual-DVR} has been extended, any quotient of length $N$ extends after finite base change and insertion of at most $N$ further ruled components. The minimal refinement and the quotient are unique over the given hull. This treats nonreduced point defects as well. Uniform stable reduction of the hull, including uniqueness of the expanded target, remains to be proved for arbitrary $(\beta,n)$.

For the structure-sheaf character $1-e^{-H}$ and its twists by $\OO_X(aH)$, Theorem~\ref{thm:verifiedcompactification} verifies descent of the chartwise stable loci and both valuative criteria for the complete smooth-fiber moduli. The parameter space is the blowup of the hyperplane system times the base curve at the point representing the joining hyperplane in the central fiber. A support approaching that point acquires the limiting equation of a graph on the inserted ruled component. Its two joins explain why the specialized class has both a term over the compactified Jacobian and a term with zero boundary. For $a=0$, Part~II \cite{Companion} computes their degrees as $4$ and $1$, respectively. Their sum gives $[q^{-55/24}]h^{\mathrm{DT}}_{1,0}=4+1=5=[q^{-55/24}]Z^{\mathrm{GSY}}_0$, where $Z^{\mathrm{GSY}}$ is the vector of theta coefficients of the quintic elliptic genus of Gaiotto, Strominger and Yin \cite[arXiv version, (3.10)]{GSYquintic}. For $G=I_{p/T}$, Part~II identifies the two positive marked classes as restrictions of refined diagonal pullbacks and computes the degree $-320-360-120=-800$, including the zero-boundary term. It also treats marked specializations for a line, a conic and a line with one point defect. A second degeneration sharing the quartic branch provides a comparison for further relative classes.

For the primitive S-duality conjecture, the distinction between the two equalities in Theorem~\ref{thm:main} is essential. Marking a joining divisor gives the specialization of a virtual cycle in the Chow group of the complete sheaf moduli. Writing a marked term as an intersection on the K3 requires the further identification \eqref{eq:mainmatching}. The structure-sheaf calculation verifies both steps at its positive marking and constructs the relative obstruction theory at its zero marking. Part~II uses surface-sheaf degrees to determine the full polar part predicted by Gaiotto, Strominger and Yin. The modular form $F^{\mathrm{DT}}$ reconstructed there from four geometric degrees and the support bound satisfies $F^{\mathrm{DT}}=Z^{\mathrm{GSY}}$. In weight $-3/2$ that polar part determines the series once the modular transformation law of the geometric degrees is proved. Part~II constructs a weight-$26$ period series $\mathcal N_X$ for the full hyperplane system and identifies the contributions of specified surface sheaves at every point length. Part~III \cite{CompanionThree} studies the point defects on an individual hyperplane section $T$, including at its singular points. If $G$ is locally free along $Z\subset T$ and $\OO_Z$ has finite projective dimension on $T$, it proves that $i_*G$ and $i_*(I_ZG)$ have the same ambient Behrend weight. This accounts for a factor in the geometric series of Part~II in every point length. 

\subsection{Which invariant is being counted?}
More generally, we fix a pure two-dimensional sheaf $F$ of character
\begin{equation}\label{eq:charge}
       v=\ch(F)=(0,rH,\beta,n),\qquad r\geq 1.
\end{equation}
Purity means that every nonzero subsheaf has two-dimensional support. In particular the divisor counted with multiplicity is
\[
       [F]_{\rm div}=\sum_{T}\length_{\OO_{X,\eta_T}}(F_{\eta_T})[T],\qquad [F]_{\rm div}\sim rH,
\]
where $T$ runs over the irreducible components of the support and $\eta_T$ is its generic point. A divisor in $|rH|$ is cut out on $X$ by a form of degree $r$. For $r=1$ it is integral, its equation annihilates $F$, and $F=i_*G$ for a rank-one torsion-free sheaf $G$ on it (Proposition~\ref{prop:primitive}). We call this the \emph{primitive divisor class}.

Higher $r$ includes sheaves of higher generic rank on a component, sheaves on reducible divisors and sheaves on nonreduced divisors. This covers $|L^{\otimes k}|$ for every ample line bundle $L$ on a smooth quintic, since $\Pic(X)=\mathbb Z H$.

Write $d=\int_X H\beta$, $n=\int_X\ch_3(F)$ and $\chi=\chi(X,F)$. As $H^3=5$ and $c_2(X)=10H^2$, Riemann--Roch gives
\[
 P_F(m)=\chi(X,F(mH))=\frac{5r}{2}m^2+dm+\chi,
 \qquad \chi=n+\frac{25r}{6}.
\]
To impose Gieseker stability we divide the Hilbert polynomial of each sheaf by its leading coefficient. We require $p_A(m)<p_F(m)$ for $m\gg0$ for every nonzero proper subsheaf $A\subset F$. Replacing $<$ by $\leq$ defines semistability.

When every semistable sheaf of character $v$ is stable, the moduli space $M_X(v)$ is projective. At $[F]$ its tangent and obstruction spaces are $\Ext^1_X(F,F)$ and $\Ext^2_X(F,F)$. Serre duality on the quintic identifies one with the dual of the other. The virtual dimension is zero, and the invariant is
\[
 \Omega_X(v)=\deg[M_X(v)]^{\vir}
     =\sum_{j\in\mathbb Z}j\,e_c(\nu_{M_X(v)}^{-1}(j)).
\]
Here $\nu$ is the Behrend function \cite[Theorem~4.18]{Behrend27} and $e_c$ is Euler characteristic with compact support. Thus a family of sheaves can make a nonzero contribution even when none of its points is isolated. If there are strictly semistable sheaves, we use the rational generalized Donaldson--Thomas invariant $\overline\Omega_X(v)$ of Joyce and Song \cite[\S5]{JS}, formed from the stack of $H$-Gieseker-semistable sheaves. In the stable case we set $\overline\Omega_X(v)=\Omega_X(v)$.

Sections~\ref{sec:smoothob} and~\ref{sec:hall} construct the obstruction theory and give the Hall logarithm defining $\overline\Omega_X(v)$.

\subsection{What a surface sheaf remembers on the common K3 surface}
Choose general homogeneous forms $f_i$ of degrees $i=1,4,5$ on $\PP^4$.
The equation $f_1f_4=tf_5$ moves the quintic to a union of a quartic
threefold and a hyperplane. Resolving the total space gives a family
$\XX\to C$ with central fiber
\begin{equation}\label{eq:central}
 X_0=Y_1\cup_S Y_2,
\end{equation}
where, writing $\PP^3=V(f_1)$,
\begin{equation}\label{eq:quinticgeom}
 Y_1=Q=V(f_4),\qquad Y_2=B=\operatorname{Bl}_\Gamma\PP^3,
 \qquad \Gamma=V(f_4,f_5)\subset\PP^3.
\end{equation}
The two components meet along the quartic K3 surface $S=V(f_4)\subset\PP^3$.
The term \emph{double divisor} refers to this surface. It is Cartier
on each component, with opposite normal bundles $\OO_S(h)$ and
$\OO_S(-h)$, where $h=\OO_S(1)$ in additive notation.

Let $i:S\hookrightarrow Y$ denote either inclusion, and let $F_Y$ be
a pure two-dimensional sheaf on $Y$. Its support can contain $S$ as a
component. If $s$ cuts out $S$, the Cartier resolution gives
\[
 i_*\LL i^*F_Y\simeq [F_Y(-S)\xrightarrow{\,s\,} F_Y]
 \qquad\text{in degrees $-1,0$.}
\]
The kernel records sections annihilated by $s$. Replacing this complex by its
ordinary restriction would discard precisely this kernel. A complex is
\emph{perfect} if it is locally quasi-isomorphic to a bounded complex
of finite rank vector bundles. For the singular threefold $X_0$, perfectness
of a sheaf is a hypothesis on its local resolutions. 

\begin{theorem}[Sheaves on a fixed union]\label{thm:intro-general}
A perfect pure two-dimensional sheaf $F$ on $X_0$ is equivalent, in families,
to a pair of pure two-dimensional sheaves $F_i$ on $Y_i$ and an isomorphism
$\LL i_1^*F_1\simeq\LL i_2^*F_2$ of their derived restrictions, taken with
its homotopies in the moduli stack of perfect complexes on $S$. On each geometric fiber at least one branch is nonzero. Write $E$ for the common
complex, $u_i:Y_i\hookrightarrow X_0$ and $j:S\hookrightarrow X_0$. On each
geometric fiber the gluing gives an exact sequence
\[
 0\to j_*\mathcal H^{-1}(E)\to F\to
 u_{1*}F_1\oplus u_{2*}F_2\to j_*\mathcal H^0(E)\to0.
\]
The assertion holds at each join of a chain of smooth components.
\end{theorem}
Near a join the equation is $xy=0$. The two-periodic resolution of
either branch and the finite projective dimension of $F$ give
$\operatorname{Tor}_{k}^{X_0}(F,\OO_{Y_i})=0$ for $k>0$. Conversely one
patches the two perfect objects along $E$ and tensors the normalization
sequence with the resulting complex. Its cohomology gives the displayed
sequence. Theorem~\ref{thm:perfectpurematchingnew} proves purity, flatness in
families and the assertion on morphisms.

Proposition~\ref{prop:route19-two-arcs} exhibits the boundary issue already on the quartic root. Two one-parameter families have the same limiting surface sheaf, but their generic boundary sheaves lie on two different plane sections of $S$. Their limits in the separated compactified Jacobian are distinct. The construction of Theorem~\ref{thm:intro-general} uses the common perfect complex on $S$, including its degree $-1$ cohomology and its extension class. Passage to a Jacobian diagonal in the virtual formula will require a separate assertion on the resulting marked space.

When $\Pic(S)=\mathbb Zh$ and $r=1$,
Corollary~\ref{cor:primitiveautomaticnew} gives $\mathcal H^{-1}(E)=0$ on
the expanded targets under consideration. The remaining sheaf can still have
a zero-dimensional subsheaf.

\subsection{Compactification and specialization of the virtual class}
To retain limits of sheaves, one allows chains of $q$ projective-line bundles
$\PP_S(\OO_S\oplus\OO_S(h))$ between $Q$ and $B$. Adjacent components
meet along copies of $S$. The expansion stack $\EE$ has local chart
\begin{equation}\label{eq:chartintro}
 [\AA^{q+1}/\GG^q],\qquad t=t_0\cdots t_q .
\end{equation}
The $q$ torus factors scale the rulings. The parameter $t_i$ smooths the
$i$th join. Marking a join means selecting one divisor $t_i=0$ and allowing
the other parameters to vary. Cutting there means partially normalizing the
threefold along the selected copy of $S$. The two exposed copies are the
\emph{outer distinguished divisors} of the resulting relative expansions.

Over $T\to C$ an expanded sheaf consists of a contraction
$c:Z\to\XX_T$ and a $T$-flat sheaf $F$ on $Z$. An arrow $(u,\phi):(Z,c,F)\to (Z',c',F')$ over $T$ consists of two isomorphisms
\begin{equation}\label{eq:intro-arrow-explicit}
 u:Z\xrightarrow{\sim} Z',\qquad c'\circ u=c,\qquad
 \phi:F\xrightarrow{\sim}u^*F'.
\end{equation}
An invertible function $a$ on $T$ acts by $(\mathrm{id}_Z,a\,\mathrm{id}_F)$.
Dividing the sheaves of isomorphisms by these scalar automorphisms in the
fppf topology defines scalar rigidification. A dilation $u$ preserving the
isomorphism class of $F$ still acts on the rigidified object.

For a fixed character $v$, Hypothesis~\ref{hyp:compact} imposes the following
conditions on this moduli problem.
\begin{enumerate}[label=(K\arabic*)]
\item Gieseker stability on each chart,
for its specified relatively ample line bundle, defines an open locus of
pure sheaves perfect on fibers, preserved by pullback along every target
isomorphism and by base change.
\item The lengths $q$ are bounded, finitely
many finite-type relative Quot schemes with frames present all families,
and the generic fiber is the required full stable moduli $M_\eta(v)$.
\item Every object over a discretely valued field extends over the valuation
ring after a finite extension, including objects on already expanded targets.
\item For two families over that ring each rigidified generic isomorphism
extends uniquely, and at a geometric point the group of target automorphisms
preserving the sheaf up to isomorphism is finite.
\end{enumerate}
Section~\ref{sec:kimexactfull}
compares these requirements with \cite[Assumptions~3.0.1 and~3.2.2]{Kim}.

\begin{theorem}[Specialization and relative intersection]\label{thm:main}
Assume (K1)--(K4), openness in the flat sheaf locus of the stack of perfect
complexes, and cohomology of the structure sheaf as in Proposition~\ref{prop:structurecohom}. This requires $h^\bullet(\OO_Z)=(1,0,0,1)$ on a complete expansion and $h^\bullet(\OO_W)=(1,0,0,0)$ on each relative piece, in degrees zero through three and compatibly with base change.

The rigidified moduli $\MM$ is proper and separated over $C$ and has a relative
perfect obstruction theory of virtual rank zero over $\EE$. On its special fiber,
\begin{equation}\label{eq:mainfull}
 \spc[M_\eta(v)]^{\vir}=\sum_\gamma (g_\gamma)_*[\MM^\dagger_\gamma]^{\vir}
 \quad\text{in }A_0(\MM_0)_\QQ .
\end{equation}
Here $\MM^\dagger_\gamma$ parametrizes sheaves with one marked joining divisor
and fixed discrete numerical data $\gamma$ on the resulting pieces.
The finite map $g_\gamma$ forgets the marking. Suppose such a space is the
full fiber product of two independently proper relative moduli of simple
sheaves under restriction to $S$,
\begin{equation}\label{eq:mainmatching}
 \MM^\dagger_\gamma\simeq M_{1,\gamma}^{\rel}
 \times_{J_{\lambda_\gamma}} M_{2,\gamma}^{\rel}.
\end{equation}
The relative factors carry the open sheaf loci and obstruction theories of Hypothesis~\ref{hyp:product}. Here $J_{\lambda_\gamma}$ is a smooth separated moduli space of nonzero perfect
boundary objects of class $\lambda_\gamma$ with $\Hom_S(E,E)=\CC$ and
$\Ext^k_S(E,E)=0$ for $k<0$. With the compatible scalar and determinant reductions of
Section~\ref{sec:reduction}, the virtual intersection is
\begin{equation}\label{eq:mainJ}
 [\MM^\dagger_\gamma]^{\vir}=\Delta_{J_{\lambda_\gamma}}^!
 \bigl([M_{1,\gamma}^{\rel}]^{\vir}\times[M_{2,\gamma}^{\rel}]^{\vir}\bigr).
\end{equation}
\end{theorem}
Equation~\eqref{eq:mainfull} identifies zero-cycles modulo rational equivalence
on the same space. On the left one closes a representative of the generic
virtual cycle in the family and intersects with the special fiber. On the
right each marked contribution is placed at its image in $\MM_0$. Taking
degrees gives the numerical DT identity. Before taking degrees the equality
retains where the terms specialize up to rational equivalence.

In the proof,
$t=t_0\cdots t_q$ gives multiplicity one to each marked divisor. Virtual
pullback carries this divisor identity to \eqref{eq:mainfull}
(Section~\ref{sec:specialization}). For \eqref{eq:mainJ} the diagonal imposes
equality of the two boundary objects. The differential on their deformations
is the difference of restriction maps. Theorem~\ref{thm:diagonal} identifies
the resulting obstruction theory with the one for glued sheaves.

When the
boundary is a stable curve sheaf of class $(0,h,\chi)$, the dimension count
is $3+3-(2+h^2)=0$.

\subsection{Main results and plan of the paper}

Section~\ref{sec:smoothmoduli} constructs the canonical Cohen--Macaulay hull of a pure surface sheaf and bounds the length of its quotient by the original sheaf for each fixed semistable character. In primitive class the hull of $i_*G$ is $i_*G^{**}$ and the bound is explicit. The same section gives the polystable limits in higher multiplicity and the Joyce--Song invariant on the smooth fiber. A degeneration identity for the Hall logarithm remains a separate question. Section~\ref{sec:geometry} constructs the two roots of the quintic degeneration and their log Calabi--Yau structures. Section~\ref{sec:generalbranches} proves the fixed-union gluing theorem and computes the index of a simple perfect sheaf $F_Y$ on a relative branch in terms of the Mukai vector of its restriction $E$ by the formula
\[
               2\chi_Y(F_Y,F_Y)=-\langle v(E),v(E)\rangle .
\]
Two families on the quartic root have the same limiting sheaf and distinct limiting points of the compactified Jacobian (Proposition~\ref{prop:route19-two-arcs}). This explains the choice of $\mathbf{Perf}(S)$ for the general matching problem.

Section~\ref{sec:compactification} gives the conditions for descent, boundedness and stable reduction on the expansion charts. When the support determines a stable expanded Hilbert point, a generic isomorphism determines the target isomorphism, and stability on that target gives uniqueness for the sheaf up to scalar (Proposition~\ref{prop:support-separated27}). Sections~\ref{sec:local} and~\ref{sec:classical} characterize purity of a nonzero curve restriction by the Cohen--Macaulay condition at an admissible join, glue flat families, and compute the endomorphism triangle and the stabilizers under gluing.

Section~\ref{sec:divisorial} proves conservation of the component divisor classes in every multiplicity. In primitive class and Picard rank one on $S$, there are three support profiles for a perfect pure limit with finite expansion stabilizers (Theorem~\ref{thm:primitiveprofiles}). In the first profile both root supports restrict to the same curve $C$, and the support on each intervening bubble is the inverse image of $C$. In the second the quartic root is empty and the first bubble carries a graph over $S$. In the third the whole sheaf lies on the blowup of $\mathbb P^3$ away from $S$. Section~\ref{sec:stability} determines when stability is preserved on cutting a chain. For rank-one branches in the setting of Theorem~\ref{thm:cutcriterion}, the Hilbert polynomials give two inequalities $0<D_e<B_e$ at each cut. Vertical twists vary these inequalities independently. On a fixed numerical sector Theorem~\ref{thm:weightedmatching} gives the boundary weights for which the two relative stability tests recover stability of the glued sheaf. The section also finds chambers valid for bounded point ideals on the specified supports and extends prescribed rank-one reflexive hulls.

Sections~\ref{sec:derived}--\ref{sec:POT} prove descent after derived base change, identify the scalar and determinant reductions under matching, and construct the rank-zero obstruction theory for sheaves relative to the stack of complete expansions. Sections~\ref{sec:relative}--\ref{sec:specialization} construct the relative theories at nonzero and zero boundary and prove the two virtual-cycle identities of Theorem~\ref{thm:main}. Marking one factor of $t=t_0\cdots t_q$ gives the Chow specialization under the compactification hypotheses. Identifying a marked space with the full fiber product of independently proper relative factors over a smooth moduli space of simple boundary objects gives the refined diagonal formula. When the common boundary is nonsimple, its automorphisms modulo scalars can already give a nonproper family of matchings of a fixed pair of branches (Proposition~\ref{prop:boundaryisomproper}).

Section~\ref{sec:cartiercase} extends finite quotients of a given hull and computes their effect on the deformation theory. For every $r\ge2$ the determinantal family of Theorem~\ref{thm:matrixexpandedallr} has a stable limit of generic rank one on its occupied components. Its supporting surface on the bubble maps finitely onto $S$ with degree $r$, and the hull requires $r$ generators at a boundary point. Corollary~\ref{cor:detdefect-local-vfc} proves a local virtual diagonal formula for these sheaves and their finite modifications at every positive marking. The same section constructs the projective global family and its relative factors for the hyperplane structure sheaf and its twists. It also shows how a further bubble restores perfectness for a point ideal approaching a node as in Proposition~\ref{prop:point}. Appendix~\ref{app:supportS} gives the equation of the ruled graph when a support specializes to $S$, including the case of a singular outer curve. Appendix~\ref{app:algorithm} gives the cut and twist calculations on the walls of stability.

\subsection*{Acknowledgements}
The author thanks Amin Gholampour who extensively discussed this project with the author over the past 15 years. The idea of double degeneration was Amin's, however due to difficulties in compactifying the relative branch moduli spaces and also having no control on cohomology theory of compactified Jacobian (and even worse, no morphisms to it either) the project halted for many years, it was the latter obstruction, which led the author to develop ideas of intersection theory on moduli space of perfect complexes on double divisor rather than the compactified Jacobian. The author would like to thank Eyal Markman for teaching him about middle cohomology of compactified Jacobian which helped the author to come up with some of numerical calculations in some specific examples. The author expresses sincere gratitude to Sheldon Katz as well as Richard Thomas for teaching him many facts about geometry of torsion sheaf moduli spaces over the years. The author would like to especially thank John Pardon for stimulating conversations which might provide a different approach to the problem discussed in this article, taking advantage of auxiliary divisors and how they induce maps to compactified Jacobians. Finally especial thanks go to Yukinobu Toda, Maxim Kontsevich, Ludmil Katzarkov, and Vladimir Baranovsky for discussions during the earlier stages of this project. He is grateful to Charles Doran, Tony Pantev, Shing-Tung Yau, Martijn Kool for valuable comments and discussions. The author acknowledges support by grants from Beijing Institute of Mathematical Sciences and Applications (BIMSA), the Beijing NSF BJNSF-IS24005, and the China National Science Foundation (NSFC) NSFC-RFIS program W2432008. He would also like to thank China's National Program of Overseas High Level Talent for generous support. Finally, he would like to thank NSF AI Institute for Artificial Intelligence and Fundamental Interactions at Massachusetts Institute of Technology (MIT) which is funded by the US NSF grant under Cooperative Agreement PHY-2019786.

\section{Surface sheaves on the smooth quintic}\label{sec:smoothmoduli}
We define the moduli problem in class \eqref{eq:charge}, construct its obstruction theory and describe the defect of a pure sheaf from being Cohen--Macaulay. 

\subsection{General ample divisor classes and the stable moduli}\label{sec:generalsetup}

\begin{definition}\label{def:generalmoduli}
Fix $r>0$, $\beta\in H^4(X,\QQ)$ and $n\in H^6(X,\QQ)\cong\QQ$ occurring as the Chern character of a coherent sheaf.  The stack $\mathfrak M^{\mathrm{ss}}_H(r,\beta,n)$ has objects $T$-flat coherent sheaves on $X\times T$ whose geometric fibers are pure two-dimensional and Gieseker semistable with character $(0,rH,\beta,n)$.  Its open stable substack is $\mathfrak M^s_H(r,\beta,n)$.  A sheaf is semistable if
\[
 \frac{P_E}{a(E)}\le\frac{P_F}{a(F)}
 \quad(0\ne E\subsetneq F),
\]
where $a(E)$ is the quadratic coefficient, and the inequality is for sufficiently large arguments.  Replace $\le$ by $<$ for stability.  Denote the projective coarse moduli space constructed by Simpson \cite[Theorem~1.21]{Simpson} by $M^{\mathrm{ss}}_H(v)$ and its stable locus by $M^s_H(v)$.
\end{definition}

\begin{proposition}\label{prop:generalpolynomial}
For a pure surface sheaf of class $(0,rH,\beta,n)$, put $d=H\beta$.  Then
\begin{equation}\label{eq:generalpolynomial}
 P_F(m)=\frac{5r}{2}m^2+dm+n+\frac{25r}{6},\qquad
 d+\frac{5r}{2}\in\mathbb Z,\quad n+\frac{25r}{6}\in\mathbb Z.
\end{equation}
Every nonzero effective supporting component has class $a_iH$, $a_i>0$, and its generic length $m_i$ satisfies
\begin{equation}\label{eq:generalcycle}
 \sum_i m_i a_i=r.
\end{equation}
Stability with respect to any ample line bundle on $X$ is equivalent to stability with respect to $\OO_X(1)$.
\end{proposition}
\begin{proof}
The tangent sequence of $X\subset\PP^4$ gives
\[
 c(T_X)=\frac{(1+H)^5}{1+5H},\quad
 c_1(X)=0,\quad c_2(X)=10H^2,\quad \int_XH^3=5.
\]
Multiply $(rH+\beta+n)e^{mH}(1+c_2(X)/12)$ and keep degree six.  The three contributions are $(rH)(m^2H^2/2)$, $\beta mH$, and $n+rHc_2(X)/12$.  This proves the polynomial.  Its value at zero is integral and its first difference at zero is $d+5r/2$, proving integrality.  Equivalently
\[
 P_F(m)=5r\binom m2+(d+5r/2)m+\chi(F).
\]
At a codimension-one generic point the ambient local ring is a discrete valuation ring and the sheaf has finite length.  Its determinant cycle is the sum of those lengths times the prime divisors.  Lefschetz gives $\Pic(X)=\mathbb ZH$; positivity against $H^2$ gives $a_i>0$.  Taking divisor classes gives \eqref{eq:generalcycle}.  Finally an ample line bundle is $\OO_X(\ell H)$, $\ell>0$, and its polynomial is $P_F(\ell m)$.

This substitution preserves the first nonzero sign in every difference of reduced polynomials.
\end{proof}

\begin{proposition}\label{prop:generalcoprime}
Put $e=d+5r/2\in\mathbb Z$.  If $\gcd(r,e,\chi(F))=1$, in particular if $\gcd(r,\chi(F))=1$, every semistable sheaf of class $(0,rH,\beta,n)$ is stable.  In this case $M^{\mathrm{ss}}_H(v)=M^s_H(v)$ is projective and carries the symmetric perfect obstruction theory of sheaves on $X$ of virtual dimension zero.
\end{proposition}
\begin{proof}
Suppose $0\ne E\subsetneq F$ and $p_E=p_F$. Write $\ch_1(E)=r'H$ and $e'=H\ch_2(E)+5r'/2$. If $r'=r$, equality of leading coefficients gives $P_{F/E}=0$, contradicting $F/E\ne0$. Thus $0<r'<r$. Coefficient comparison gives
\[
     re'=r'e,\qquad r\chi(E)=r'\chi(F).
\]
Choose $u,w,z\in\mathbb Z$ with $ur+we+z\chi(F)=1$. Then
$r'=r(ur'+we'+z\chi(E))$, a contradiction. Hence $M_H^{\rm ss}(v)=M_H^s(v)$. The numerical character is locally constant, so this moduli space is an open and closed part of Simpson's projective moduli space \cite{Simpson}.

For stable $F$, Serre duality identifies the dual of the unit $H^0(\OO_X)\xrightarrow{\sim}\End(F)$ with the top trace. The resulting symmetric obstruction theory \cite{HT,Thomas} has virtual rank
$\ext^1_X(F,F)-\ext^2_X(F,F)=0$.
\end{proof}

\subsection{The complete smooth-target charts}
Let $\chi=n+25r/6$, $d=H\beta$, and $P(m)=5rm^2/2+dm+\chi$.  Boundedness of semistable sheaves supplies an integer $N$ for which every such $F$ is $N$-regular.  Put $V=\CC^{P(N)}$.  The fixed-character open-and-closed part of the open subscheme
\[
 U\subset\operatorname{Quot}_X(V\otimes\OO_X(-N),P)
\]
imposes purity, semistability and the condition that $V\to H^0(F(N))$ is an isomorphism.  The latter condition is open by cohomology and base change; $H^{>0}(F(N))=0$ on this locus.  The family groupoid is
\begin{equation}\label{eq:generalQuotgroupoid}
 \mathfrak M_H^{\mathrm{ss}}(v)=[U/\operatorname{GL}(V)].
\end{equation}
Indeed a flat semistable family has locally free $p_*F(N)$ of rank $P(N)$, and its frame bundle gives a map to $U$.  Conversely the universal quotient descends along the frame changes.  These operations are mutually inverse on objects and arrows.  On the stable locus the only automorphisms over a fixed $X$ are scalar, so scalar rigidification is the stable coarse moduli space.

Simpson's determinant-of-cohomology linearization gives the projective good quotient on the semistable locus \cite{Simpson}.  One may use $\det p_*F(m)$ for $m\gg N$ with the $\operatorname{SL}(V)$ action; the determinant character removes the remaining central weight when a $\operatorname{GL}(V)$ presentation is used.  Its GIT semistability is exactly the polynomial inequality of Definition~\ref{def:generalmoduli}.

Thus \eqref{eq:generalQuotgroupoid} also parametrizes sheaves on singular supports with noninvertible hull.

\subsection{The obstruction theory on the smooth quintic}\label{sec:smoothob}
For the fixed smooth threefold $X$, let $p:X\times U\to U$ be a chart of stable sheaves carrying the universal family $\FF$. The deformation complex before removing scalars is $K=\RR p_*\mathcal RHom(\FF,\FF)$, with tangent complex $K[1]$. For a square-zero extension with ideal $I$, the obstruction to lifting $F$ is the product of its Atiyah class with the Kodaira--Spencer class of the extension, in $\Ext_X^2(F,F\otimes I)$.

When it vanishes, lifts up to isomorphism form a torsor under $\Ext_X^1(F,F\otimes I)$ \cite{HT}. Scalar rigidification and the determinant construction of Section~\ref{sec:reduction} give
\begin{gather*}
 C^{\red}=\fib\!\left(
  \cofib(\OO_U\xrightarrow{1\mapsto\id}K)
  \longrightarrow\cofib(\OO_U\longrightarrow\RR p_*\OO_{X\times U})
 \right),\\
 (C^{\red})^\vee[-1]\longrightarrow L_{U}
\end{gather*}
on an \emph{\'etale} chart of the rigidified moduli, with the obstruction theory pulled back along the chart. One can work on the scalar gerbe when a universal sheaf exists only there. The unit removes $\Ext^0_X(F,F)=\CC$. The trace on $\Ext^3_X(F,F)$ is dual to this unit by Serre duality and is an isomorphism. Since $H^1(\OO_X)=H^2(\OO_X)=0$, the two remaining groups are
\[
 T_{[F]}M_H^s(v)=\Ext_X^1(F,F),\qquad
 \operatorname{Ob}_{[F]}=\Ext_X^2(F,F)=\Ext_X^1(F,F)^\vee .
\]
This yields virtual dimension zero on each stable component. Scalar conjugation is trivial on $K$, so these charts carry the descended theory. For a proper stable moduli its virtual degree is the invariant defined in the introduction. The Hall-algebra definition below applies when strictly semistable objects occur.

\subsection{All semistable charges and the Hall algebra}\label{sec:hall}
For fixed $v$, the stack $\mathfrak M_H^{\rm ss}(v)$ has as objects over $T$ the $T$-flat families on $X\times T$ with semistable geometric fibers of character $v$, and arrows the isomorphisms of these families. We use this stack when the character admits strict semistability.  Its numerical invariant is the rational generalized Donaldson--Thomas invariant, denoted $\overline{\dt}_H(v)$ \cite{JS}.

The definition uses the intrinsic Behrend function and the Hall algebra of coherent sheaves. The Hall stack retains the full automorphism group, for example $\mathrm{GL}_r$ at $\OO_T^{\oplus r}$. Its quotient by scalars is the $\mathrm{PGL}_r$ which enters Corollary~\ref{cor:noseparatedall}.

Let $\delta_v$ be the stack function of $H$-semistable sheaves of class $v$.  The Hall product records short exact sequences, with the first factor the subobject and the second the quotient.  The logarithm at a fixed reduced polynomial is
\begin{equation}\label{eq:Hallloggeneral}
 \epsilon_v=
 \sum_{j\ge1}\frac{(-1)^{j-1}}j
 \sum_{\substack{v_1+\cdots+v_j=v\\p(v_i)=p(v)}}
       \delta_{v_1}*\cdots*\delta_{v_j}.
\end{equation}
On the numerical charge lattice of a quintic this is coefficientwise finite.  Indeed every nonzero factor has $\ch_1(v_i)=r_iH$ with $r_i\ge1$, so $j\le r$.  Equality of the reduced polynomials fixes $d_i/r_i=d/r$ and $\chi_i/r_i=\chi/r$.  For each ordered partition of $r$ there is therefore at most one numerical tuple.  Terms that are not integral charges or not represented by sheaves are zero.

The integration morphism of Joyce and Song \cite[\S5]{JS} on the virtual-indecomposable Lie algebra sends $\epsilon_v$ to the coefficient defining $\overline{\dt}_H(v)$, in the normalization in which it equals the ordinary virtual degree when no strictly semistable objects occur.  For example
\[
 \epsilon_{2v_1}=\delta_{2v_1}-\tfrac12\delta_{v_1}*\delta_{v_1}
\]
when these are the only equal-polynomial decompositions.  The second term records the stack of filtrations with two stable factors of character $v_1$.

The Euler form governing the Hall Lie bracket is explicit here:
\begin{equation}\label{eq:generalEulerpairing}
 \chi_X(F_1,F_2)=r_2d_1-r_1d_2.
\end{equation}
To prove it, multiply
\[
 \ch(F_1)^\vee=(0,-r_1H,\beta_1,-n_1),\qquad
 \ch(F_2)=(0,r_2H,\beta_2,n_2)
\]
and retain degree six.  The only surviving terms are $\beta_1r_2H-r_1H\beta_2$, since $c_1(X)=0$.  This proves \eqref{eq:generalEulerpairing}.  Integration in the Hall algebra defines the generalized invariant \cite[\S5]{JS}.

\subsection{The canonical Cohen--Macaulay extension}
On the smooth quintic use its fixed trivialization of $\omega_X$ and write
\[
 F^{\mathrm d}=\mathcal Ext^1_X(F,\OO_X),\qquad
 T_F=\mathcal Ext^2_X(F,\OO_X),\qquad
 F^{\mathrm{CM}}=(F^{\mathrm d})^{\mathrm d}.
\]
The notation $F^{\mathrm{CM}}$ is not the ambient double dual $F^{**}$, which is zero for a torsion sheaf.  It is a duality construction in codimension one.

\begin{theorem}\label{thm:generalCMhull}
For every pure two-dimensional coherent sheaf $F$ on a smooth threefold with trivial canonical bundle, $T_F$ has finite support, $F^{\mathrm d}$ and $F^{\mathrm{CM}}$ are Cohen--Macaulay of dimension two, and there is a canonical exact sequence
\begin{equation}\label{eq:generalCMsequence}
 0\longrightarrow F\longrightarrow F^{\mathrm{CM}}
 \longrightarrow T_F^\dagger\longrightarrow0,
 \qquad T_F^\dagger=\mathcal Ext^3_X(T_F,\OO_X).
\end{equation}
Moreover $\operatorname{length}T_F^\dagger=\operatorname{length}T_F$ and
\[
 \ch(F^{\mathrm{CM}})=\ch(F)+\ch(T_F^\dagger).
\]
On the quintic, in numerical cohomology the last term is $\operatorname{length}(T_F)[\mathrm{pt}]$.
If $F=i_*G$ for a rank-one torsion-free sheaf on a normal Cartier surface, then $F^{\mathrm{CM}}=i_*(G^{**})$.  There is no normality, reducedness, or generic-rank-one assumption in \eqref{eq:generalCMsequence}.
\end{theorem}
\begin{proof}
At a closed point $x$, let $A=\OO_{X,x}$ and $M=F_x$.  Purity excludes a finite-length submodule, so $\depth_A M\ge1$ when $M\ne0$.  The regular local ring $A$ has dimension three; Auslander--Buchsbaum gives $\pd_A M\le2$.  At a codimension-two point, purity similarly gives depth at least one over a regular local ring of dimension two, hence projective dimension at most one.

At a codimension-one point the ring is a discrete valuation ring and the torsion module has projective dimension one.  It follows that $\mathcal Ext^2(F,\OO_X)$ is supported at finitely many closed points and that $\mathcal Ext^j(F,\OO_X)=0$ for $j\ge3$.  Also $\mathcal Hom(F,\OO_X)=0$, since a torsion sheaf has no nonzero map to an invertible sheaf.

Put $\mathbb D=\RR\mathcal Hom_X(-,\OO_X)$.  The cohomology just computed gives the canonical truncation triangle
\begin{equation}\label{eq:generaldualtriangle}
 F^{\mathrm d}[-1]\longrightarrow\mathbb D(F)
       \longrightarrow T_F[-2]\longrightarrow F^{\mathrm d}.
\end{equation}
The sheaf $F$ is perfect on the smooth threefold, so $\mathbb D\mathbb D(F)=F$.  For a finite-length sheaf $T$, local duality over a regular ring of dimension three gives
\[
 \mathbb D(T)=T^\dagger[-3],\qquad
 \mathcal Ext^j(T,\OO_X)=0\quad(j\ne3).
\]
For a residue field this follows from its Koszul resolution; induction on the length, using a composition series and the Ext sequence, proves the assertion and equality of lengths for every $T$.  Dualizing \eqref{eq:generaldualtriangle} therefore gives
\begin{equation}\label{eq:generalbidualtriangle}
 T_F^\dagger[-1]\longrightarrow F
    \longrightarrow\mathbb D(F^{\mathrm d})[1]
    \longrightarrow T_F^\dagger.
\end{equation}
Its cohomology sequence yields, with all arrows fixed by biduality,
\[
 \mathcal Hom(F^{\mathrm d},\OO_X)=0,\quad
 0\to F\to\mathcal Ext^1(F^{\mathrm d},\OO_X)
       \to T_F^\dagger\to0,
\]
and $\mathcal Ext^j(F^{\mathrm d},\OO_X)=0$ for $j\ge2$.

A finite module over a regular local ring whose positive Ext groups above degree one vanish has projective dimension at most one: apply Hom to a minimal free resolution; a nonzero last free term would otherwise give a nonzero last Ext.  Since $F^{\mathrm d}$ is torsion and nonzero along the divisorial support, its projective dimension there is exactly one.  The depth formula shows that it is Cohen--Macaulay of codimension one.

Dualizing a length-one free resolution of this module gives a length-one free resolution of $F^{\mathrm{CM}}$; the initial Hom term is zero.  Thus the latter is Cohen--Macaulay as well.  This proves \eqref{eq:generalCMsequence} and the Chern-character assertion.

For the final assertion, Cartier duality identifies
\[
 \mathcal Ext^1_X(i_*G,\OO_X)
   =i_*\mathcal Hom_T(G,\OO_T(T)).
\]
On a normal surface the reflexive extension of a rank-one torsion-free sheaf agrees with $G$ in codimension one and has finite quotient.  It is Cohen--Macaulay because the surface is normal Gorenstein of dimension two.  Applying the displayed duality twice cancels the invertible twist $\OO_T(T)$ and gives $i_*G^{**}$.  The biduality map is its natural inclusion, proving compatibility with \eqref{eq:generalCMsequence}.
\end{proof}

\begin{proposition}\label{prop:generaldefectbound}
For a fixed semistable character $v=(0,rH,\beta,n)$ there is an integer $N_v$ such that
\[
 \operatorname{length}(F^{\mathrm{CM}}/F)\le N_v
\]
for every sheaf in $\mathfrak M_H^{\mathrm{ss}}(v)$.  There is a finite stratification of the reduced underlying scheme of a finite-type Quot atlas on which $F^{\mathrm d}$, $F^{\mathrm{CM}}$, and $F^{\mathrm{CM}}/F$ are flat and commute with base change.  
\end{proposition}
\begin{proof}
Use the atlas $U$ in \eqref{eq:generalQuotgroupoid} with universal sheaf $\mathcal F$.  The morphism $X\times U\to U$ is smooth of relative dimension three.  The sheaf is $U$-flat, and its geometric fibers have local projective dimension at most two.  Let $A=\OO_{U,u}\to B=\OO_{X\times U,x}$ and choose local presentations
$0\to K_2\to B^{a_1}\to B^{a_0}\to\mathcal F_x\to0$.
Both syzygies are $A$-flat. On the fiber $K_2\otimes_A k(u)$ is
locally free by the projective-dimension bound. Lift a basis to
$B^b\to K_2$.

Its cokernel and then its kernel vanish by Nakayama,
using flatness of $K_2$ for the kernel. Thus $\mathcal F$ is perfect, and $\mathbb D(\mathcal F)$ is a bounded coherent perfect complex, with derived base change.

Apply generic flatness to its finitely many cohomology sheaves on every irreducible component of the reduced base.  Remove the resulting dense open sets and repeat on the smaller-dimensional complements, computing the restricted derived complex there.  Noetherian induction gives finitely many locally closed strata with flat cohomology.  On a stratum $B$, the spectral sequence
\[
 \Tor_p^{\OO_B}\bigl(\mathcal H^q(\mathbb D\mathcal F),\OO_{B'}\bigr)
 \Longrightarrow
 \mathcal H^{q-p}(\LL f^*\mathbb D\mathcal F)
\]
collapses because the cohomology sheaves are $B$-flat.  In particular $T_{\mathcal F}$ is flat and its fibers have finite support.  Its support is proper and quasi-finite over $B$, hence finite.  Its Hilbert polynomial is its length, constant on each connected stratum.  A finite-type stratum has finitely many connected components.  Taking the maximum of these finitely many lengths and using Theorem~\ref{thm:generalCMhull} gives $N_v$.

On each stratum $\mathcal F^{\mathrm d}$ is flat with Cohen--Macaulay fibers of codimension one, so it has a length-one locally free resolution on the total space by the same lifting argument.  Its dual has fiber cohomology only in degree one.  The fiberwise flatness criterion for a bounded complex of base-flat modules makes $\mathcal F^{\mathrm{CM}}$ a flat sheaf commuting with base change.  Relative duality for the finite flat sheaf $T_{\mathcal F}$ similarly makes $T_{\mathcal F}^\dagger$ finite flat.

Triangle \eqref{eq:generalbidualtriangle} gives the exact sequence in families.  Further finite refinement, if needed, makes its numerical labels constant.  This proves the stated stratified construction.
\end{proof}

On such a stratum, finite quotients of a specified Cohen--Macaulay sheaf $R$ are represented by $\operatorname{Quot}^j(R)$, and the kernels recover the surface sheaves.  If $R$ is invertible on a normal supporting surface this is a Hilbert scheme of points; it is generally a Quot scheme, not a Hilbert scheme.  For instance the rank-one reflexive module $(x,z)$ over $\CC[[x,y,z]]/(xy-z^2)$ has a two-dimensional fiber at the singular point.

Its length-one quotients with fixed support form $\PP^1$, whereas those of the local ring form a point.  The equality follows directly by identifying a length-one quotient with a one-dimensional quotient of $R/\mathfrak mR$.  For completeness, the exact sequence
\[
 0\to (x,z)\to A\to A/(x,z)=\CC[[y]]\to0
\]
has middle depth two and right depth one.  The depth lemma gives depth two for $(x,z)$.  The surface $A$ is normal (a hypersurface with an isolated singularity), and $(x,z)$ is torsion-free of rank one, hence reflexive.  Its two independent generators in $(x,z)/\mathfrak m(x,z)$ prove non-freeness by Nakayama.  The quotient projective space therefore retains a genuine non-invertible-hull contribution.

Thus these local module types must be retained in the full count.

\begin{proposition}\label{prop:unavoidablesemistability}
For every $r\ge2$ the class
\[
 v_r=r\ch(\OO_T)=\left(0,rH,-\frac r2H^2,\frac{5r}{6}\right),
 \qquad T\in|H|,
\]
contains strictly semistable sheaves for every ample polarization. The scalar-rigidified stack has a positive-dimensional stabilizer at $\OO_T^{\oplus r}$. Without semistability the pure sheaves in this numerical class form an unbounded family.  
\end{proposition}
\begin{proof}
Choose a smooth integral $T$ by Bertini. For $0\ne A\subsetneq\OO_T$,
the quotient has dimension at most one, so $P_A=P_{\OO_T}-P_{\OO_T/A}$
has the same quadratic coefficient as $P_{\OO_T}$ and $p_A<p_{\OO_T}$
for every ample polarization.  The direct sum $F=\OO_T^{\oplus r}$ is semistable: successive intersections with its summands express every subsheaf as extensions of subsheaves of $\OO_T$, whose reduced polynomials are at most $p_{\OO_T}$.

A summand $\OO_T\subset F$ has exactly that polynomial, so $F$ is not stable.  This calculation is independent of the polarization.  Since $H^0(T,\OO_T)=\CC$, one has
\[
 \Aut_X(F)=\operatorname{GL}_r,\qquad
 \Aut_X(F)/\GG=\operatorname{PGL}_r.
\]
The last group has dimension $r^2-1>0$, proving the stack assertion.

For the last assertion it suffices to take $r=2$ and a smooth $T$.  For every integer $a>0$ choose a finite subscheme $Z_a\subset T$ of length $5a^2$ and set
\[
 F_a=i_*\bigl(\OO_T(aH)\oplus I_{Z_a/T}(-aH)\bigr).
\]
Both summands are torsion-free on $T$, so $F_a$ is pure.  The identity
\[
 (e^{aH}+e^{-aH})(1-e^{-H})-5a^2[\mathrm{pt}]
       =2(1-e^{-H})\quad\text{in degrees }\le6
\]
proves $\ch(F_a)=v_2$.  Moreover, note that the family of such sheaves is indeed unbounded. That is, $F_a$ contains $i_*\OO_T(aH)$ with reduced linear coefficient $2a-1$, unbounded as $a\to\infty$.  If these sheaves were bounded, some fixed twist $N$ would make every $F_a(N)$ globally generated.  For $a>N$, its nonzero direct summand $i_*I_{Z_a/T}((N-a)H)$ has no global sections, since it is contained in a negative ample line bundle on $T$.

It cannot be globally generated.  This contradiction proves unboundedness.  For higher $r$ we add $r-2$ copies of $\OO_T$ and achieve the same conclusion.
\end{proof}

For $r=1$ no strict semistability occurs for any $(\beta,n)$. Proposition~\ref{prop:primitive} proves the support classification in this divisor class.

\subsection{A stable sheaf with a polystable limit}
\label{sec:determultiple}
There is a second family in which stabilizers can be read off before
introducing either degeneration. The intrinsic stabilizer of its special member controls separatedness after scalar rigidification.

\begin{proposition}\label{prop:determultiple}
Fix $r\ge2$ and a smooth hyperplane section $T=(s=0)$ of $X$.
There is a flat family over a discrete valuation ring whose general
member is $i_*M$ for a line bundle $M$ on a smooth integral
$D\in|rH|$ and whose special member is $\mathcal O_T^{\oplus r}$.
Its Chern character is $r(1-e^{-H})$. Writing $h=H|_D$, the class
$c=c_1(M)$ satisfies
\[
 c h=\frac{5r(r-1)}2,\qquad
 c^2=\frac{5r(r-1)(r-2)}6,\qquad
 \frac{(ch)^2}{h^2}-c^2=\frac{5r(r^2-1)}{12}>0 .
\]
In particular the general member is stable and its line bundle is
numerically independent of the hyperplane class. The special member
has stabilizer $\mathrm{GL}_r$ before scalar rigidification.
\end{proposition}
\begin{proof}
Take an $r$ by $r$ matrix $A$ of general hyperplane sections.
At any $x\in X$ evaluation of all the entries at $x$ is surjective
onto $\mathrm{Mat}_r(\mathbb C)$ after trivializing $\mathcal O_X(H)|_x$.
The universal evaluation over $X$ is smooth. The determinantal stratum
of matrices of corank $k$ has codimension $k^2$ and is smooth along
its exact-corank locus. Since $\dim X=3$, the incidence for corank
at least two projects into a proper closed subset of matrices of
sections.

Generic smoothness for the corank-one stratum then makes
$D=(\det A=0)$ smooth for general $A$. It is nonempty and connected:
$D\sim rH$ and $H^1(X,\mathcal O_X(-rH))=0$ by Serre duality and
Kodaira vanishing. Hence $D$ is integral. Elementary row and
column operations reduce $A$ locally along $D$ to
$\mathrm{diag}(I_{r-1},u\det A)$ with $u$ a unit. Thus
\[
 0\longrightarrow\mathcal O_X(-H)^{\oplus r}
 \xrightarrow{A}\mathcal O_X^{\oplus r}\longrightarrow i_*M
 \longrightarrow0
\]
for an invertible $M$ on $D$. Stability follows from the
generic-rank-one argument on an integral support.
Grothendieck--Riemann--Roch gives
\[
 \mathrm{ch}(i_*M)=i_*\left(1+c-rh/2+
 c^2/2-r\,c\cdot h/2+r^2h^2/6\right)=r(1-e^{-H}).
\]
Equating the last two Chern-character entries and using $h^2=5r$
gives the three intersection formulae.

Over $\mathbb C[t]_{(t)}$ replace $A$ by $(1-t)sI_r+tA$.
The locus of smooth integral determinants of corank exactly one along
their zeros is open and contains $t=1$ on the original pencil.
The determinant is not identically zero on either fiber of the
discrete valuation ring, so the matrix is injective on both fibers.
The local flatness criterion applied to its two-term locally free
presentation gives a flat cokernel. At $t=0$ this cokernel is
$\mathcal O_T^{\oplus r}$.

Since $H^0(T,\mathcal O_T)=\mathbb C$,
its endomorphism ring is $\mathrm{Mat}_r(\mathbb C)$, as claimed.
\end{proof}

\begin{corollary}\label{cor:noseparatedall}
For $r\ge2$ there is no separated algebraic stack over $\CC$
whose geometric objects include all semistable sheaves of character
$r(1-e^{-H})$ on the smooth quintic and whose automorphism groups
are their intrinsic automorphism groups modulo scalars. In particular
such a stack cannot occur as the smooth fiber of a separated expanded
sheaf stack.
\end{corollary}
\begin{proof}
At $x=[\OO_T^{\oplus r}]$, Proposition~\ref{prop:determultiple} gives
$\Aut(x)=\mathrm{GL}_r/\GG=\mathrm{PGL}_r$. If the diagonal of the
stack were proper, its pullback by $(x,x):\Spec\CC\to\mathcal M\times
\mathcal M$ would be proper. This pullback is the group scheme
$\mathrm{PGL}_r$, which is affine of positive dimension. A proper
affine scheme of finite type over a field is finite, a contradiction.
Properness of the diagonal is preserved under base change to a fiber.
\end{proof}

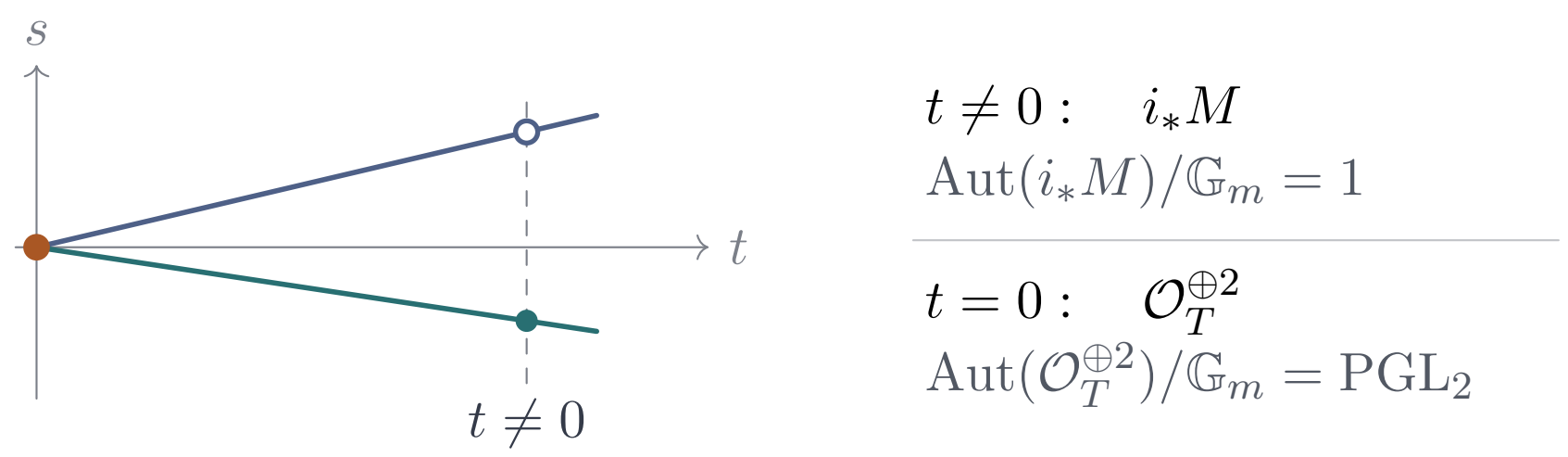
\begin{figure}[htbp]
\centering
\begin{tikzpicture}[x=1cm,y=1cm,font=\small,line cap=round,line join=round]
 \draw[->,qiInk!65] (-.15,0)--(4.8,0) node[right] {$t$};
 \draw[->,qiInk!65] (0,-1.08)--(0,1.30) node[above] {$s$};
 \draw[line width=1.05pt,qiBlue] (0,0)--(4.0,.94);
 \draw[line width=1.05pt,qiTeal] (0,0)--(4.0,-.60);
 \draw[dashed,qiInk!65] (3.5,-.97)--(3.5,1.09);
 \filldraw[fill=white,draw=qiBlue,line width=1pt] (3.5,.8225) circle (2.25pt);
 \fill[qiTeal] (3.5,-.525) circle (2.25pt);
 \fill[qiAmber] (0,0) circle (2.7pt);
 \node[below,qiInk] at (3.5,-.94) {$t\ne 0$};
 \node[anchor=west] at (6.2,.98) {$t\ne0:\quad i_*M$};
 \node[anchor=west,qiInk!85] at (6.2,.47) {$\Aut(i_*M)/\GG=1$};
 \draw[qiInk!30] (6.26,.05)--(10.87,.05);
 \node[anchor=west] at (6.2,-.39) {$t=0:\quad\OO_T^{\oplus2}$};
 \node[anchor=west,qiInk!85] at (6.2,-.91) {$\Aut(\OO_T^{\oplus2})/\GG=\mathrm{PGL}_2$};
\end{tikzpicture}
\caption{A transverse slice of the determinantal family for $r=2$. At a point of $T$ where $A$ has two distinct eigenvalues, the tangent lines to the two local determinant branches have slopes $-\lambda_1$ and $-\lambda_2$ in the $(t,s)$ plane. The drawing illustrates two distinct real slopes. At $t=0$ the matrix is $sI_2$ and its cokernel has rank two on $T$. The groups at right are the scalar-rigidified automorphism groups of the global fibers.}\label{fig:stabilizerslice}
\end{figure}

Figure~\ref{fig:stabilizerslice} illustrates the rank change in Proposition~\ref{prop:determultiple}. The obstruction in Corollary~\ref{cor:noseparatedall} is already visible on the smooth threefold, before a choice of polarization on any expansion.

\subsection{The primitive hyperplane specialization}
\begin{proposition}\label{prop:primitive}
Let $X\subset\PP^4$ be a smooth quintic and let $F$ be a pure two-dimensional sheaf with $\ch_1(F)=H$. Its support is an integral normal hyperplane section $i:T\hookrightarrow X$, and $F=i_*G$ for a rank-one torsion-free sheaf $G$ on $T$. The sheaf $F$ is Gieseker-stable for every ample polarization on $X$.
\end{proposition}
\begin{proof}
By Lefschetz $\Pic(X)=\mathbb ZH$. Write $[T_j]=a_jH$ and $b_j=\length_{\OO_{X,\eta_j}}F_{\eta_j}$ for the components of the divisorial support. Then
\[
   1=\sum_j b_ja_j,\qquad b_j>0,\qquad 5a_j=H^2\cdot[T_j]>0 .
\]
Thus $[F]_{\rm div}=[T]$ for a single integral hyperplane section. Over the discrete valuation ring $\OO_{X,\eta_T}$ the length-one module $F_{\eta_T}$ is annihilated by the equation $s_T$. Hence $\dim\Supp(s_TF)\le1$. Purity gives $s_TF=0$, and $F=i_*G$ with $G$ rank-one torsion-free on $T$.

For the Gauss morphism $\gamma:X\to(\PP^4)^\vee$,
$\gamma^*\OO(1)=\OO_X(4)$. A curve contracted by $\gamma$ would have degree zero against $4H$. Thus $\gamma$ has finite fibers. The singular locus of $T=X\cap P$ lies in $\gamma^{-1}([P])$, so $T$ is regular in codimension one. A Cartier hypersurface is Cohen--Macaulay; Serre's criterion gives normality.

For $0\ne F'\subsetneq F$ the generic ranks agree. If $a$ is their common leading Hilbert coefficient for an ample polarization, then
\[
   p_F-p_{F'}=P_{F/F'}/a>0\qquad(m\gg0).
\]
This proves stability.
\end{proof}

\begin{proposition}\label{prop:hilbertquintic}
With $d=H\beta$ and $n=\int_X\ch_3(F)$, the Hilbert polynomial is
\begin{equation}\label{eq:Hilbertquintic}
P_F(m)=\frac52m^2+dm+n+\frac{25}{6}.
\end{equation}
The stable moduli on a smooth fiber is projective and its virtual dimension is zero.
\end{proposition}
\begin{proof}
In Proposition~\ref{prop:generalpolynomial} set $r=1$. The Riemann--Roch integral becomes
\[
 \chi(X,F(mH))=\int_X\left(\tfrac12m^2H^3+mH\beta+
              \ch_3(F)+\tfrac1{12}Hc_2(X)\right)
       =\tfrac52m^2+dm+n+\tfrac{25}{6}.
\]
Proposition~\ref{prop:primitive} identifies $M_H^s(v)=M_H^{\rm ss}(v)$, which is projective. Serre duality gives $\Ext_X^2(F,F)=\Ext_X^1(F,F)^\vee$, hence virtual dimension zero.
\end{proof}

\subsection{A uniform bound on the punctual defect of the original sheaves}\label{sec:defectbound}
The cokernel of the inclusion of a pure sheaf in its Cohen--Macaulay extension has bounded length in a fixed semistable class. This bounds the lengths of the finite quotients in the Quot construction below.

\begin{proposition}\label{prop:defectbound}
Let $X\subset\PP^4$ be a smooth quintic threefold, $\ch(F)=(0,H,\beta,n)$,
$d=H\beta$, $e=d+5/2$ and $\chi=n+25/6$.  For an integer $b$ put
\[
 B(b)=\sum_{\substack{j\geq0\\5j\leq b}}(b-5j+1).
\]
An empty sum has value zero. Write $F=i_*G$ as in
Proposition~\ref{prop:primitive}. For $R=G^{\vee\vee}$, with duals on $T$, one has
\[
 0\leq\ell(R/G)\leq B(e)+B(5-e)-\chi .
\]
The bound is uniform in $T$ and in the isomorphism class of $G$.
\end{proposition}
\begin{proof}
Set $h=H|_T$, and choose a plane section $D\in|h|$ which is smooth
and misses the singular locus of $T$, the locus where $G$ or $R$
is not locally free, and the support of $Q=R/G$. Each excluded locus
is finite. Indeed, at every point of codimension one the local ring of
the normal surface $T$ is a discrete valuation ring, and a torsion-free
module of rank one over that ring is free. Bertini's theorem supplies
such a $D$.

It is a plane quintic, so $g(D)=6$ and $\deg h|_D=5$.
Put $A=G|_D=R|_D$. Since the defining equation of $D$ is
$G$-regular, for each $m$ the restriction sequence gives
\[
 \begin{aligned}
 \chi(D,A(mh))&=P_F(m)-P_F(m-1)=5m+d-5/2,\\
 \chi(D,A(mh))&=5m+\deg A+1-6.
 \end{aligned}
\]
Hence $\deg A=e$. By adjunction $\omega_T=\OO_T(h)$, so
$(R^\vee\otimes\omega_T)|_D=A^\vee(h)$ has degree $5-e$.

Let $U$ be either $R$ or $R^\vee(h)$ and put $b=\deg(U|_D)$.
It is reflexive of rank one on $T$. A nonzero section of $U(-jh)$
defines an effective Weil divisor in its divisor class on the normal
surface $T$. Its intersection with $h$ equals $b-5j$, so
$H^0(T,U(-jh))=0$ if $5j>b$. For every $j$ the sequence
\[
 0\longrightarrow U(-(j+1)h)\longrightarrow U(-jh)
  \longrightarrow (U|_D)(-jh)\longrightarrow0
\]
is exact. For a line bundle $M$ of degree $c\geq0$ on $D$, evaluation
on a divisor of $c+1$ distinct points embeds $H^0(D,M)$ into a
vector space of dimension $c+1$, since the twisted kernel has degree
$-1$. Consequently
\[
 h^0(T,U)\leq\sum_{0\leq5j\leq b} h^0(D,(U|_D)(-jh))
     \leq B(b).
\]
Serre duality on the projective Gorenstein surface $T$ gives
$h^2(T,R)=h^0(T,R^\vee(h))$. Taking Euler characteristics in
$0\to G\to R\to Q\to0$ now yields
\[
 \ell(Q)=h^0(T,R)-h^1(T,R)+h^0(T,R^\vee(h))-\chi
       \leq B(e)+B(5-e)-\chi.
\]

The link with the ambient defect can also be seen stalkwise. For a
closed point $x\in\operatorname{Supp}Q$, set $A_x=\OO_{T,x}$ and
write $(-)^\dagger=\operatorname{Ext}^2_{A_x}(-,\omega_{A_x})$
on modules of finite length. Since $R_x$ is maximal Cohen--Macaulay,
local duality gives $\operatorname{Ext}_{A_x}^{j}(R_x,\omega_{A_x})=0$
for $j>0$. The connecting homomorphism for the preceding short exact
sequence is therefore an isomorphism
\[
 \operatorname{Ext}_{A_x}^1(G_x,\omega_{A_x})\xrightarrow{\sim}Q_x^\dagger .
\]
The dual of the residue field is the residue field; exactness of
$(-)^\dagger$ and a composition series of $Q_x$ give
$\ell(Q_x^\dagger)=\ell(Q_x)$. By Cartier duality and summation over $x$,
\begin{equation}\label{eq:dualdefectlength}
 \ell\bigl(\mathcal Ext^2_X(F,\omega_X)\bigr)
   =\ell\bigl(\mathcal Ext^1_T(G,\omega_T)\bigr)=\ell(R/G).
\end{equation}
\end{proof}

The proposition bounds the original punctual defect. It is not a bound on all possible bubble chains: non-invertible hulls on singular expanded supports can carry further data. In a fixed-source point-quotient sector every nontrivial bubble must carry a quotient point, so the separate inequality $k\le n$ used in Section~\ref{sec:idealchamber} applies there.

\begin{proposition}[Failure of stable reduction on the smooth quintic]\label{prop:noordinaryproperextension}
For every $r\ge2$ the scalar-rigidified moduli of stable sheaves of
character $v_r=r(1-e^{-H})$ on a smooth quintic fails the existence
part of the valuative criterion for properness, even after arbitrary
finite extension of the discrete valuation ring. Consequently a
proper relative moduli stack over the degeneration base cannot have
these ordinary stable moduli as its entire smooth fibers.
\end{proposition}
\begin{proof}
Let $R,K$ and $\mathcal G$ be the family in
Proposition~\ref{prop:determultiple}, with $\mathcal G_K$ stable and
$\mathcal G_0=\OO_T^{\oplus r}$. Denote the projective Simpson space
of semistable sheaves by $M^{ss}_X(v_r)$; its geometric points classify
$S$-equivalence classes \cite[Theorem~1.21]{Simpson}. For any finite
extension $R\subset R'$ with fraction field $K'$, the given family
induces $f:\Spec R'\to M^{ss}_X(v_r)$. A putative stable extension
would induce $g:\Spec R'\to M^{ss}_X(v_r)$ with $g_{K'}=f_{K'}$.

This also holds for a point of the rigidified stack since its map to
the Simpson space descends through the scalar gerbe. Separatedness
gives $g=f$. Over the geometric closed point, any stable sheaf $G$
giving this point would therefore satisfy
\[
       G=\operatorname{gr}_{\rm JH}(G)
       \simeq\operatorname{gr}_{\rm JH}(\OO_T^{\oplus r})
       =\OO_T^{\oplus r}.
\]
The sheaf on the right has a proper direct summand with its own
reduced Hilbert polynomial, contradicting stability. Since properness
is preserved by base change, restricting a proposed relative
compactification with these smooth fibers to one smooth quintic
gives the same contradiction.
\end{proof}

Together, Corollary~\ref{cor:noseparatedall} and Proposition~\ref{prop:noordinaryproperextension} exclude both attempts to obtain the required proper stable-sheaf stack in these characters. Including the displayed polystable limit violates separatedness with the prescribed inertia. Omitting its entire $S$-equivalence class violates valuative existence. An expansion of the central fiber of the threefold leaves this example on a smooth target untouched.

The proposition exhibits a boundary point in a proper compactification of the smooth-target stable locus. Indeed, take its scheme-theoretic closure
\[
   \overline{M_X^s(v_r)}\ \subset\ M_X^{ss}(v_r).
\]
The right-hand side is projective by \cite[Theorem~1.21]{Simpson}, so the closure is projective. The valuative calculation puts the point $[\OO_T^{\oplus r}]$ in its boundary. Its moduli meaning is the $S$-equivalence class with associated graded $\OO_T^{\oplus r}$.

A compactification by expanded relative sheaves has additional data at a joining K3: an expanded target, a sheaf on it, and its derived restriction to the joining surface. The coarse point above records the associated graded on $X$. Constructing the relative data and its compatible virtual class is the role of Hypothesis~\ref{hyp:compact} and Theorem~\ref{thm:main}.

\section{The first degeneration and its two components}\label{sec:geometry}
We construct a semistable model of the quintic family and determine the canonical bundles, normal bundles, and polarizations used for the two component moduli.

Choose general forms $f_1,f_4,f_5$ of the indicated degrees and degenerate the quintic by
\[
f_1f_4=t f_5.
\]
The usual small semistable resolution of the local equation $xy=tz$ gives \eqref{eq:quinticgeom}. Set $H_B=c_1(\OO_{\PP^3}(1))$ pulled back to $Y_2$, and let $E_B$ be the exceptional divisor of its blowup.

Along $\Gamma$ use coordinates $x=f_1$, $y=f_4$, $z=f_5$ and a fourth coordinate on the ambient fourfold. The total family has local equation $xy=tz$. Blow up $(x,t)$. On $t=xw$ the strict-transform equation is $y=wz$. The special fiber is $xw=0$ in smooth coordinates. Its branches are $w=0$ and $x=0$, and their scheme-theoretic intersection has ideal $(w,x)$. On the other chart $x=tu$ one has $uy=z$. These charts give $Q$ and the blowup of $\PP^3$ along $(f_4,f_5)$.

The gluing surface on the second component is isomorphic to the original quartic K3, since blowing up its Cartier curve $\Gamma$ is an isomorphism.

\begin{figure}[htbp]
\centering
\begingroup
\definecolor{qAmbient}{HTML}{53618A}
\definecolor{qSeam}{HTML}{176568}
\definecolor{qSupport}{HTML}{B15C25}
\definecolor{qInk}{HTML}{272C3D}
\begin{tikzpicture}[x=1cm,y=1cm,font=\small,
  line cap=round,line join=round,
  axisink/.style={draw=qAmbient,line width=.7pt},
  trace/.style={draw=qSupport,line width=1.55pt},
  fleck/.style={draw=qAmbient!38,line width=.35pt}]
\path[use as bounding box] (-3.14,-1.98) rectangle (10.04,2.19);
\begin{scope}
\node[text=qInk] at (0,1.94) {\textsc{Smooth fiber}};
\def\XX{.91*.52*(exp(\z)-exp(-\z))}
\def\YY{.31*.52*(exp(\z)+exp(-\z))}
\path[fill=qAmbient!12] plot[variable=\z,domain=-1.56:1.56,samples=49]
 ({\XX},{1.12-\YY}) --
 plot[variable=\z,domain=1.56:-1.56,samples=49] ({\XX},{-.63-\YY}) -- cycle;
\path[fill=qAmbient!8] plot[variable=\z,domain=-1.56:0,samples=26]
 ({\XX},{1.12-\YY}) --
 plot[variable=\z,domain=0:-1.56,samples=26] ({\XX},{-.63-\YY}) -- cycle;
\foreach \w in {-.54,-.30,-.07,.17,.41,.65,.9} {
  \draw[fleck] plot[variable=\z,domain=-1.56:1.56,samples=43] ({\XX},{\w-\YY});}
\foreach \z in {-1.4,-1.12,-.8,-.44,0,.44,.8,1.12,1.4} {
  \draw[fleck] ({\XX},{1.12-\YY})--({\XX},{-.63-\YY});}
\draw[axisink] plot[variable=\z,domain=-1.56:1.56,samples=49]
 ({\XX},{1.12-\YY}) --
 plot[variable=\z,domain=1.56:-1.56,samples=49] ({\XX},{-.63-\YY}) -- cycle;
\draw[trace] plot[variable=\z,domain=-1.56:1.56,samples=49] ({\XX},{.09-\YY});
\node[text=qAmbient] at (0,1.28) {$X_t$};
\node[text=qSupport,anchor=east,inner sep=1pt] at (-2.42,-.75) {$T_t$};
\draw[qSupport!78,line width=.45pt] (-2.38,-.72)--(-2.16,-.69);
\node[text=qInk] at (0,-1.82) {$xy=t$};
\end{scope}
\draw[-{Stealth[length=2.5mm,width=1.7mm]},qInk!70,line width=.65pt]
 (2.78,-.15)--(4.10,-.15);
\node[above,text=qInk,font=\footnotesize] at (3.44,-.1) {$t\to0$};
\begin{scope}[shift={(6.92,0)}]
\node[text=qInk] at (0,1.94) {\textsc{Central fiber}};
\foreach \d/\v in {-1/8,1/16} {
 \path[fill=qAmbient!\v] (0,1.12)--({\d*2.27},.345)--({\d*2.27},-1.405)--(0,-.63)--cycle;
 \foreach \w in {-.54,-.30,-.07,.17,.41,.65,.9} {
   \draw[fleck] (0,\w)--({\d*2.27},{\w-.775});}
 \foreach \i in {.14,.32,.50,.68,.86} {
   \draw[fleck] ({\d*2.27*\i},{1.12-.775*\i})--
                     ({\d*2.27*\i},{-.63-.775*\i});}
 \draw[axisink] (0,1.12)--({\d*2.27},.345)--({\d*2.27},-1.405)--(0,-.63);
}
\draw[qSeam,line width=1.45pt] (0,-.63)--(0,1.12);
\draw[trace] (-2.27,-.685)--(0,.09)--(2.27,-.685);
\filldraw[fill=white,draw=qSeam,line width=.85pt] (0,.09) circle[radius=2.6pt];
\fill[qSupport] (0,.09) circle[radius=1.23pt];
\node[text=qSeam,inner sep=1.2pt,above] at (0,1.15) {$S$};
\node[text=qAmbient] at (-1.35,1.04) {$Q$};
\node[text=qAmbient] at (1.35,1.04) {$B$};
\node[text=qSupport,anchor=east,inner sep=1pt] at (-2.43,-.73) {$T_1$};
\node[text=qSupport,anchor=west,inner sep=1pt] at (2.43,-.73) {$T_2$};
\draw[qSupport!78,line width=.45pt] (-2.39,-.72)--(-2.29,-.69);
\draw[qSupport!78,line width=.45pt] (2.39,-.72)--(2.29,-.69);
\node[text=qInk] at (0,-1.82) {$xy=0$};
\end{scope}
\end{tikzpicture}
\endgroup
\caption{Transverse surface support in the first degeneration, near a smooth point
of the boundary curve $C$. In local coordinates $(x,y,u,v)$ the threefold
has equation $xy=t$ and the supporting surface has equation $u=0$.
The drawing fixes $v$ and shows the positive real branches. The amber
trace specializes to the traces of $T_1\subset Q$ and $T_2\subset B$.
Their common point in this slice runs through $C=T_1\cap S=T_2\cap S$
as $v$ varies.}
\label{fig:transverse-specialization}
\end{figure}
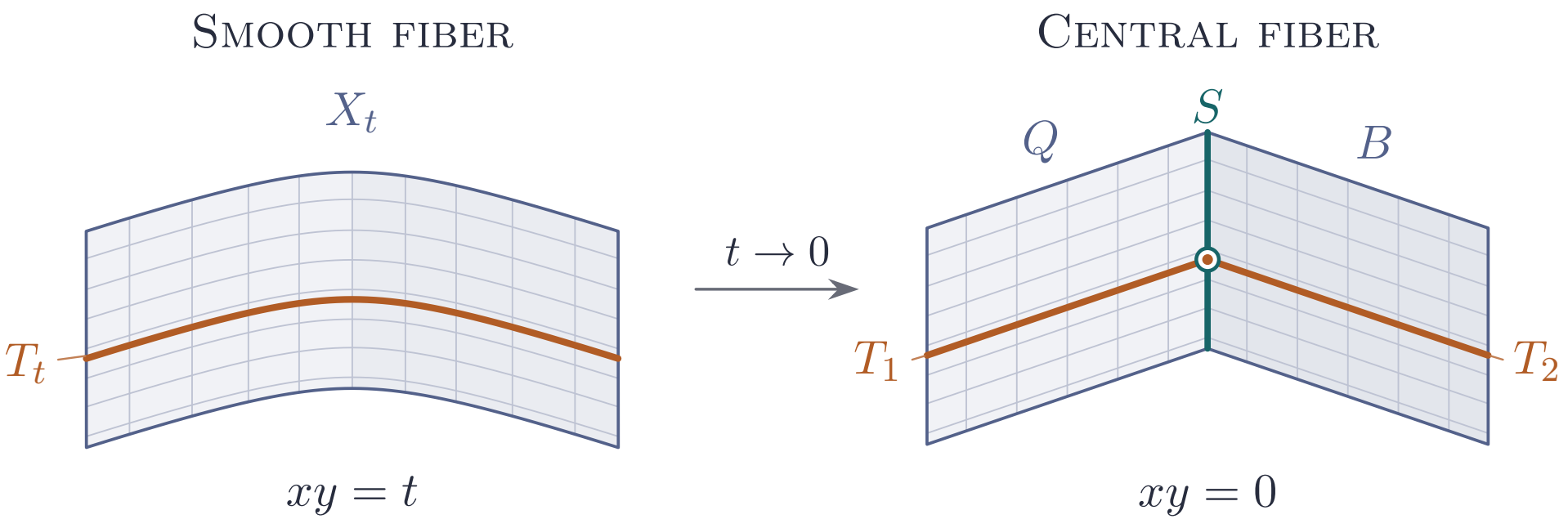

Figure~\ref{fig:transverse-specialization} shows the smoothing near a point of the common curve of two transverse supporting surfaces.

\begin{proposition}\label{prop:logcy}
The double divisor and the normal bundles satisfy
\begin{align}
K_{Y_1}&=-S,& K_{Y_2}&=-S,\\
N_{S/Y_1}&=\OO_S(h),&N_{S/Y_2}&=\OO_S(-h).
\end{align}
Every expanded fiber has trivial dualizing sheaf. A relative chain consisting of one root component and any number of bubbles has dualizing sheaf equal to the negative of its outer relative divisor.
\end{proposition}
\begin{proof}
Adjunction gives $K_Q=-H_Q$ and $S\sim H_Q$. On the blowup,
\[
K_{Y_2}=-4H_B+E_B,\qquad S\sim4H_B-E_B.
\]
The center is a divisor of class $5h$ on $S$, so $E_B|_S=5h$. Thus $S|_S=-h$ on $Y_2$, while $S|_S=h$ on $Y_1$.

For a bubble $p:\Delta=\PP_S(\OO_S\oplus N)\to S$ with distinguished sections $S_-,S_+$, the projective-bundle canonical formula gives
\[
\omega_\Delta(S_-+S_+)=p^*\omega_S=\OO_\Delta.
\]
The dualizing sheaf on a normal-crossing union restricts on a component to its canonical sheaf twisted by the joining divisors. The displayed trivializations glue by the opposite-residue convention along the nodes, inherited from the Calabi--Yau smoothing. If one outer divisor is not glued, precisely its negative remains. This proves both assertions.
\end{proof}

\begin{lemma}\label{lem:ample}
Let $H$ be the pullback of the hyperplane line bundle to the semistable total space. The line bundle $\mathcal A=\OO_{\XX}(2H+Y_1)$ is relatively ample after shrinking the base about zero. Its restrictions are $2H$ on the smooth fibers, $H_Q$ on $Y_1$, and $6H_B-E_B$ on $Y_2$.
\end{lemma}
\begin{proof}
The restrictions follow from $Y_1|_{Y_1}=-S$ and $Y_1|_{Y_2}=S$. The ideal of $\Gamma$ is generated by forms of degrees four and five, so $I_\Gamma(5)$ is globally generated. The associated blowup embeds in a product of projective spaces with the two hyperplane classes $H_B$ and $5H_B-E_B$. Their sum $6H_B-E_B$ is ample. The restriction on the other root is ample as well.

Ampleness on the components implies ampleness on their projective union, and relative ampleness is open on the base.
\end{proof}

\begin{proposition}\label{prop:structurecohom}
Let $Z$ be a complete expanded central fiber, and let $W$ be a relative chain obtained by cutting $Z$ at a node and keeping one root. Then
\begin{align}
&H^0(Z,\OO_Z)=H^3(Z,\OO_Z)=\CC,
&&H^1(Z,\OO_Z)=H^2(Z,\OO_Z)=0,\\
&H^0(W,\OO_W)=\CC,
&&H^j(W,\OO_W)=0\quad(j>0).
\end{align}
These equalities persist in the corresponding flat families of expansions, with locally free cohomology and base change.
\end{proposition}
\begin{proof}
Let $Z_0,\ldots,Z_{k+1}$ be the components of a complete expansion
and $D_e$ their successive intersections. Pushing all sheaves to $Z$
gives the normalization sequence
\[
0\longrightarrow\OO_Z\longrightarrow
  \bigoplus_i\OO_{Z_i}\xrightarrow{\delta}
  \bigoplus_e\OO_{D_e}\longrightarrow0,\qquad
 \delta((f_i))_e=f_i|_{D_e}-f_{i+1}|_{D_e}.
\]
For the blowup $b:B\to\PP^3$, $Rb_*\OO_B=\OO_{\PP^3}$.
The Fano root has $H^{>0}(Q,\OO_Q)=0$. On a bubble $p:\Delta\to S$,
$Rp_*\OO_\Delta=\OO_S$ and each section induces the identity on
$H^2(S,\OO_S)$. Thus $H^0(\delta):\CC^{k+2}\to\CC^{k+1}$ is
surjective with diagonal kernel. In degree two we obtain
\[
 (a_1,\ldots,a_k)\longmapsto
 (-a_1,a_1-a_2,\ldots,a_{k-1}-a_k,a_k),\qquad
 0\to\CC^k\xrightarrow{H^2(\delta)}\CC^{k+1}
                \xrightarrow{\sum}\CC\to0 .
\]
For $k=0$ the exact row is $0\to0\to\CC\xrightarrow{1}\CC\to0$.
For a relative chain oriented from its root, omit the last output
coordinate; the remaining matrix is invertible. The long exact sequence
gives the asserted groups. Partially smoothing the chain gives the
same calculation with fewer bubbles. The hypersurface sequence gives
the stated $h^j$ on smooth quintic fibers.

On the smooth expansion
charts these constant dimensions give locally free $R^jp_*\OO$
by cohomology and base change, and derived base change gives the
claim for pullbacks of those charts. Relative duality on the complete
family identifies $R^3p_*\OO$ with $(p_*\omega_{Z/T})^\vee$;
use the relative Calabi--Yau trivialization.
\end{proof}

\subsection{The quartic root in a second family}\label{sec:partI-pencil-example}
The pair $(Q,S)$ also occurs as one side of a degeneration of a K3 fibration. Choose $s_0,s_1\in H^0(Q,\OO_Q(1))$ with $S=(s_0=0)$ and smooth base curve $C_*=(s_0=s_1=0)\in|h|$ on $S$. The pencil defines
\[
 \kappa:Y=\operatorname{Bl}_{C_*}Q\longrightarrow\PP^1.
\]
Write $A$ for the pullback of $H_Q$ and $E_C$ for the exceptional divisor. Each fiber has class $F_{\rm fib}=A-E_C$. The strict transform of $S$ is a smooth fiber, and $K_Y=-F_{\rm fib}$ by the blowup formula.

Let $s$ denote the degeneration parameter. Blow up $S\times\{0\}$ in $Q\times\AA^1_s$ and then the strict transform of $C_*\times\AA^1_s$. The resulting family $\mathscr Y$ has general fiber $Y$ and special fiber
\[
 \mathscr Y_0=Q\cup_{S_\infty}W,\qquad
 P=\PP_S(\OO_S\oplus\OO_S(h)),\qquad W=\operatorname{Bl}_{C_*\subset S_0}P .
\]
Here $S_0,S_\infty\subset P$ are the disjoint sections with normal line bundles $\OO_S(h),\OO_S(-h)$, respectively. Locally write $S=(x=0)$ and $C_*=(x=u=0)$. In the chart $x=sz$ of the first blowup the second center has equations $z=u=0$ for all $s$. Its special fiber is the curve $C_*$ on $S_0$. Along $S_\infty$ the chart $s=xw$ gives the original normal crossing.

The strict transform $\mathscr B_0$ of $S\times\AA^1_s$ is the specified outer divisor of the family. It meets $W$ along $B_0$, the strict transform of $S_0$. Since the center is Cartier on $S_0$, one has $B_0\simeq S$, $N_{B_0/W}\simeq\OO_S$ and $B_0\cap S_\infty=\varnothing$. The two blowup formulas give
\[
 \omega_{\mathscr Y/\AA^1_s}(\mathscr B_0)\simeq\OO_{\mathscr Y},
 \qquad \omega_W(S_\infty+B_0)\simeq\OO_W .
\]
Thus a sheaf on $W$ has two boundary values, its derived restrictions to $S_\infty$ and $B_0$. Only the former is identified with the restriction from $Q$. The one-boundary cohomology calculation of Proposition~\ref{prop:structurecohom} concerns the roots of the quintic degeneration. Here the projection $W\to S$ gives $H^j(W,\OO_W)=H^j(S,\OO_S)$ by the projective bundle and blowup formulas. Part~II \cite[\S3]{Companion} uses the two-boundary obstruction theory of this second family.
\section{Sheaves on the components and restriction to the K3 surface}\label{sec:generalbranches}
Fix the first central fiber $Z=Y_1\cup_S Y_2$ of \eqref{eq:central}. There are two closed immersions $i_\alpha:S\hookrightarrow Y_\alpha$ and a single immersion $j:S\hookrightarrow Z$. The \emph{double divisor} is the K3 surface $S$. It is Cartier in each smooth threefold $Y_\alpha$ and is the locus where the two branches of $Z$ meet. Locally the inclusions have the form
\[
 Z=(xy=0),\qquad Y_1=(y=0),\quad Y_2=(x=0),\quad S=(x=y=0).
\]
In particular, $S$ has its reduced structure. The equation $xy=0$ does not give $S$ multiplicity two. In an expansion the same description applies at each of the inserted copies of $S$.

On $Y_i$ fix a numerical class $v_i$ of dimension two and the polarization $A_i$ of Lemma~\ref{lem:ample}. Write $\mathfrak M_i(v_i)$ for the stack of flat families of pure two-dimensional sheaves with that class. The Gieseker stable substack is specified by $p_{A_i}(F')<p_{A_i}(F_i)$ for every proper nonzero subsheaf $F'\subset F_i$. Which stable substacks on expansions give proper relative spaces for a degeneration formula is a further question.

Already on these unexpanded components the natural evaluation on the entire stack is
\begin{equation}\label{eq:derivedevaluationnew}
 \mathfrak M_i(v_i)\longrightarrow\mathbf{Perf}(S),\qquad F_i\longmapsto E_i=\LL i_i^*F_i .
\end{equation}
Every coherent sheaf on a smooth threefold is perfect. A flat family of such sheaves is perfect on the product locally on any noetherian parameter scheme: lift a finite free resolution of a fiber using base flatness and Nakayama's lemma as in the local resolution argument immediately below. Thus \eqref{eq:derivedevaluationnew} commutes with derived base change. If $s_i$ is the equation of $S$, resolving $\OO_S$ by $[\OO_{Y_i}(-S)\xrightarrow{s_i}\OO_{Y_i}]$ gives
\begin{align}
 i_{i*}\mathcal H^{-1}(E_i)&=\ker(F_i(-S)\xrightarrow{s_i}F_i),&
 i_{i*}\mathcal H^0(E_i)&=F_i/s_i F_i(-S).\label{eq:twoHboundarynew}
\end{align}
These are the only possibly nonzero cohomology sheaves. For a pure $F_i$, the first is nonzero exactly when $S$ is an irreducible component of its support. Indeed a zero divisor on a finite module belongs to an associated prime, and every associated prime of $F_i$ is the generic point of a surface.

At the generic point of $S$ a nonzero module of finite length over the discrete valuation ring has nonzero kernel under multiplication by its uniformizer. In families it is the complex $E_i$ that commutes with base change; its individual cohomology sheaves can jump.

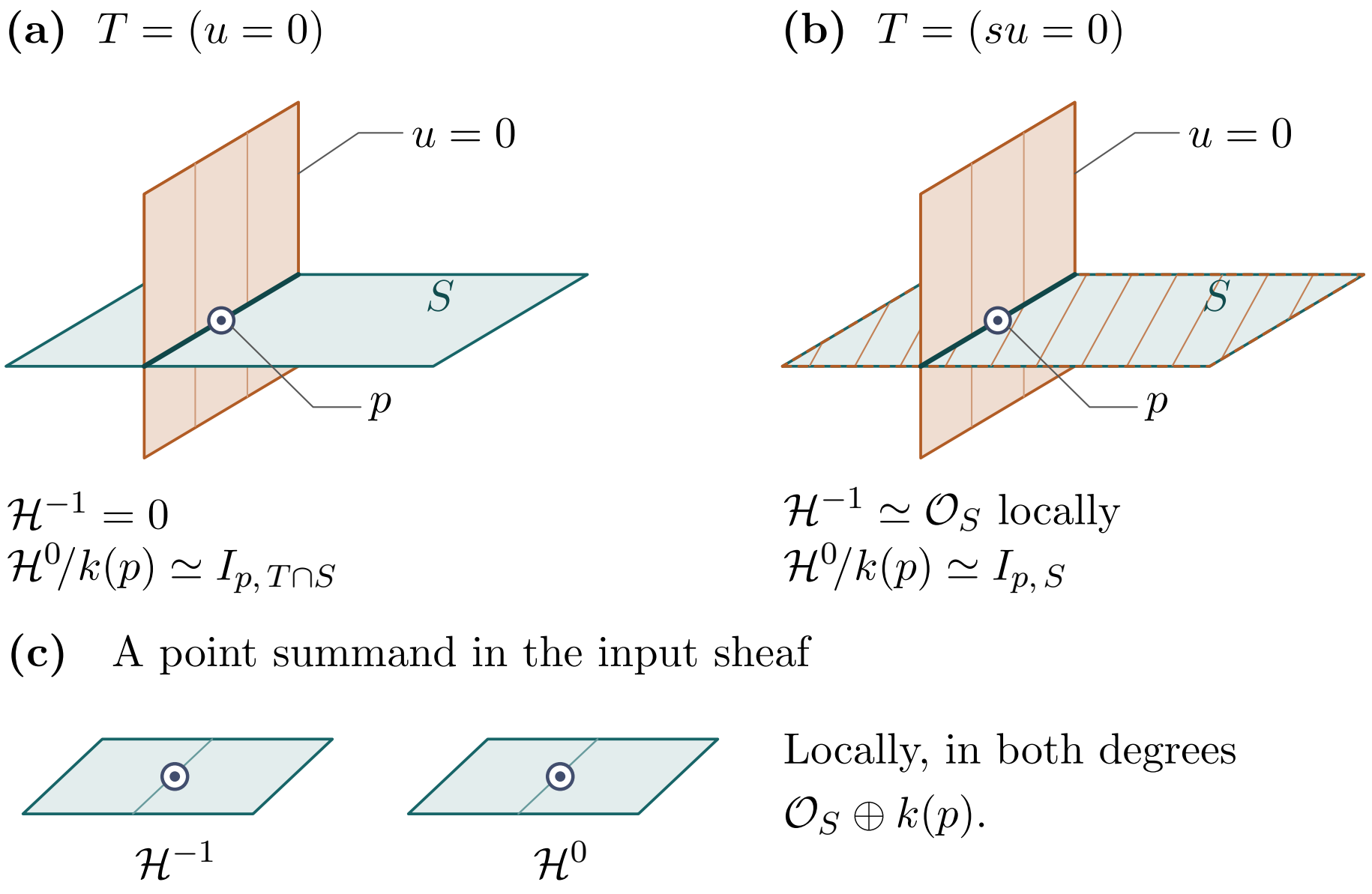
\begin{figure}[htbp]
 \centering
 \begingroup
 \definecolor{brS}{HTML}{176568}
 \definecolor{brT}{HTML}{B15C25}
 \definecolor{brI}{HTML}{53618A}
 \begin{tikzpicture}[x=1cm,y=1cm, line join=round, line cap=round,
   every node/.style={font=\small}, wire/.style={draw=brI!62,line width=.4pt},
   seam/.style={draw=brS!70!black,line width=1.25pt},
   surface/.style={draw=brT,line width=.7pt,fill=brT!21},
   base/.style={draw=brS,line width=.7pt,fill=brS!12},
   leader/.style={draw=black!65,line width=.4pt}]
  \foreach \d/\case in {0/1,6.85/2}{
   \begin{scope}[shift={(\d,0)}]
     \ifnum\case=1
       \node[anchor=west,inner sep=0pt] at (-.02,2.55)
         {\textbf{(a)}\ \ $T=(u=0)$};
     \else
       \node[anchor=west,inner sep=0pt] at (-.02,2.55)
         {\textbf{(b)}\ \ $T=(su=0)$};
     \fi
     \filldraw[surface] (1.22,-1.21)--(2.58,-.40)--(2.58,1.93)--(1.22,1.12)--cycle;
     \draw[brT!62,line width=.38pt] (1.67,-.93)--(1.67,1.39)
                             (2.13,-.67)--(2.13,1.66);
     \filldraw[base] (.00,-.40)--(3.77,-.40)--(5.13,.41)--(1.36,.41)--cycle;
     \ifnum\case=2
       \begin{scope}
         \clip (.00,-.40)--(3.77,-.40)--(5.13,.41)--(1.36,.41)--cycle;
         \foreach \x in {-1,-.6,...,5.8}
            \draw[brT!78,line width=.4pt] (\x,-.47)--++(.56,1.02);
       \end{scope}
       \draw[brT,dash pattern=on 2.1pt off 2pt,line width=.7pt]
          (.00,-.40)--(3.77,-.40)--(5.13,.41)--(1.36,.41)--cycle;
     \fi
     \filldraw[surface] (1.22,-.40)--(2.58,.41)--(2.58,1.93)--(1.22,1.12)--cycle;
     \draw[brT!62,line width=.38pt] (1.67,-.13)--(1.67,1.39)
                                (2.13,.14)--(2.13,1.66);
     \draw[seam] (1.22,-.4)--(2.58,.41);
     \filldraw[fill=white,draw=brI!75!black,line width=.83pt] (1.90,.005) circle (3.15pt);
     \fill[brI!75!black] (1.90,.005) circle (1.15pt);
     \draw[leader] (1.99,-.05)--(2.71,-.76)--(3.13,-.76);
     \node[anchor=west,inner sep=2pt] at (3.13,-.76) {$p$};
     \node[anchor=south east,inner sep=1pt,text=brS!76!black] at (4.00,.05) {$S$};
     \draw[leader] (2.58,1.30)--(3.11,1.66)--(3.50,1.66);
     \node[anchor=west,inner sep=2pt] at (3.50,1.66) {$u=0$};
     \ifnum\case=1
       \node[anchor=west,inner sep=0pt] at (0,-1.67) {$\mathcal H^{-1}=0$};
       \node[anchor=west,inner sep=0pt] at (0,-2.18)
          {$\mathcal H^0\!/k(p)\simeq I_{p,\,T\cap S}$};
     \else
       \node[anchor=west,inner sep=0pt] at (0,-1.67)
          {$\mathcal H^{-1}\simeq\OO_S$ locally};
       \node[anchor=west,inner sep=0pt] at (0,-2.18)
          {$\mathcal H^0\!/k(p)\simeq I_{p,\,S}$};
     \fi
   \end{scope}
  }
 \begin{scope}[shift={(0,-3.05)}]
   \node[anchor=west,inner sep=0pt] at (0,.10) {\textbf{(c)}\quad A point summand in the input sheaf};
   \foreach \dX/\term in {0/{-1},3.40/{0}}{
     \begin{scope}[shift={(\dX,0)}]
       \filldraw[base] (.15,-1.30)--(2.18,-1.30)--(2.88,-.64)--(.85,-.64)--cycle;
       \draw[brS!65,line width=.45pt] (1.12,-1.30)--(1.82,-.64);
       \filldraw[fill=white,draw=brI!80!black,line width=.85pt] (1.49,-.97) circle (3.2pt);
       \fill[brI!80!black] (1.49,-.97) circle (1.3pt);
       \node[anchor=north] at (1.5,-1.40) {$\mathcal H^{\term}$};
     \end{scope}
   }
   \node[anchor=west,align=left,inner sep=0pt] at (6.85,-1.05) {Locally, in both degrees\\[3pt]
                                      $\mathcal O_S\oplus k(p)$.};
 \end{scope}
 \end{tikzpicture}
 \endgroup
 \caption{Derived restriction to the joining surface $S$. Panels (a) and (b) use the ideal $I_{p,T}$.
 In (a) the two surfaces cross along a curve. In (b) the hatched sheet
 is a component of $T$ coinciding with $S$. The other sheet meets it along
 the dark line. The ring marks the support of the skyscraper submodule
 $k(p)\subset\mathcal H^0$. In (c) the support is $S$ itself. The rings mark the point summand in each cohomology sheaf when the input sheaf already has a point summand.}
 \label{fig:cartiercontactprofiles}
\end{figure}

The kernel in \eqref{eq:twoHboundarynew} also gives a purity constraint.
If $F_i$ is pure of dimension two, every nonzero subsheaf of
$i_{i*}\mathcal H^{-1}(E_i)\subset F_i(-S)$ has dimension two.
Thus $\mathcal H^{-1}(E_i)$, whenever nonzero, is torsion-free on the
integral surface $S$. The quotient $\mathcal H^0(E_i)$ can contain a
skyscraper subsheaf in either support configuration of
Figure~\ref{fig:cartiercontactprofiles}.

Work in $Y=\Spec\CC[[s,u,v]]$ with $S=(s=0)$, and write
$B=\CC[[u,v]]$ and $k=\CC$. On a surface $T\subset Y$ take the ideal
$F=\mathfrak m_{p,T}$ of the origin. If $T=(u=0)$ then
$F=(s,v)\subset\CC[[s,v]]$. Multiplication by $s$ is injective and
reduction modulo $s$ gives
\[
 0\longrightarrow k\bar s\longrightarrow F/sF
   \longrightarrow (v)\CC[[v]]\longrightarrow0 .
\]
If $T=(su=0)$ then $A=\OO_T=\CC[[s,u,v]]/(su)$ and
$F=(s,u,v)\subset A$. Both choices give pure surface sheaves since
$\OO_T$ is Cohen--Macaulay of pure dimension two and $F\subset\OO_T$.
In the second case
\[
 \ker(s:F\to F)=(u)A\simeq B,\qquad
 0\longrightarrow k\bar s\longrightarrow F/sF
    \longrightarrow (u,v)B\longrightarrow0 .
\]
Indeed $\operatorname{Ann}_A(s)=(u)$ and
$B\xrightarrow{\,1\mapsto u\,}(u)A$ is an isomorphism.
In both sequences the kernel of reduction is $(s)/sF$; its generator
$\bar s$ survives because $s\notin\mathfrak m_{p,T}^2$ and is killed by
$\mathfrak m_{p,T}$. It accounts for the embedded point in degree zero.

For comparison, one can see the role of input purity in the same chart.
Let $j:S\hookrightarrow Y$ and put $G=j_*\OO_S\oplus k(p)$. The summand
$k(p)\subset G$ makes $G$ impure. Since $sG=0$, the Cartier resolution
computes both boundary cohomology sheaves locally as
\[
 \mathcal H^{-1}(\LL j^*G)\simeq\OO_S\oplus k(p),\qquad
 \mathcal H^0(\LL j^*G)\simeq\OO_S\oplus k(p).
\]
Panel (c) depicts this change when the purity hypothesis is dropped.

\begin{lemma}\label{lem:perfectfamilywithoutnormalnew}
Let $Z_T\to T$ be a flat finitely presented family of the expanded threefold fibers under consideration over a noetherian base. Let $F$ be a finitely presented $T$-flat sheaf whose restriction to each geometric fiber is pure of dimension two and perfect on that threefold. Then $F$ is perfect on $Z_T$, locally with resolution length at most two. In a projective flat family of pure sheaves, the locus of parameters with perfect fibers is open.
\end{lemma}
\begin{proof}
At a closed point on a geometric fiber its local ring is Cohen--Macaulay of dimension three, and a nonzero pure surface module has depth at least one. Perfectness and Auslander--Buchsbaum give projective dimension at most two. At a nonclosed point the same bound follows by localization (or from the dimension of its local ring). Over a local parameter ring $A$ and local ambient ring $B$ take free surjections onto $F$ and its first syzygy.

Writing $M=F_z$ and $K_1,K_2$ for the kernels gives
\[
 0\to K_1\to B^{a_0}\to M\to0,\qquad
 0\to K_2\to B^{a_1}\to K_1\to0 .
\]
Flatness of $M$ and $B$ over $A$ gives $\Tor_i^A(K_j,-)=0$ for $i>0$ in succession. Both sequences thus remain exact over $k=A/\mathfrak m_A$. Where the second fiber syzygy is free, lift a basis and obtain a surjection from a free $B$-module to the second syzygy by Nakayama's lemma. Its kernel is zero by $A$-flatness and Nakayama's lemma again. These arguments commute with extension to the geometric residue field by faithful flatness.

The resulting resolution works on a neighborhood. Conversely such a resolution of a base-flat module restricts to each fiber, so the good locus in $Z_T$ is open. Its complement has closed image in $T$ by properness of the projective family, giving the final assertion.
\end{proof}

\begin{theorem}\label{thm:perfectpurematchingnew}
On a fixed union $Z=Y_1\cup_S Y_2$ as above, perfect pure two-dimensional sheaves are equivalent, in families, to pairs of pure two-dimensional branch sheaves (allowing a zero branch) together with an equivalence $\LL i_1^*F_1\simeq\LL i_2^*F_2$ in the derived stack $\mathbf{Perf}(S)$. The globally zero object is excluded. More generally the same assertion holds for a chain of smooth components with disjoint consecutive double divisors.

No condition on the degree $-1$ cohomology or the purity of the boundary is imposed. If $u_i:Y_i\hookrightarrow Z$ and $E$ is the common complex, the glued sheaf $F$ fits into
\begin{gather}
 0\longrightarrow j_*\mathcal H^{-1}(E)\longrightarrow F\longrightarrow K\longrightarrow0,\label{eq:Torcorrectionnew}\\
 0\longrightarrow K\longrightarrow u_{1*}F_1\oplus u_{2*}F_2
 \longrightarrow j_*\mathcal H^0(E)\longrightarrow0.\nonumber
\end{gather}
The difference map in the second sequence uses the specified gluing equivalence. Both sequences are asserted on geometric fibers. The derived triangle producing them is compatible with arbitrary base change.
\end{theorem}
\begin{proof}
At a join put $B=R[[x,y]]/(xy)$, $B_1=B/(y)$, $B_2=B/(x)$, with $R$ regular. Let $M$ be the stalk of a perfect pure sheaf. The resolution
\[
 \cdots\xrightarrow{y}B\xrightarrow{x}B\xrightarrow{y}B\longrightarrow B_1\longrightarrow0
\]
gives $\Tor^{B}_{j+2}(M,B_1)\cong\Tor^B_j(M,B_1)$ for $j\ge1$. Since $\pd_B M<\infty$, these groups vanish. Interchanging $x,y$ proves the same assertion for $B_2$. In particular $M_i=M\otimes_B^{\LL}B_i=M\otimes_BB_i$.

Localize at any point of the join in $\Supp M$ and choose a minimal finite free resolution $P_\bullet\to M$ there. Then $P_\bullet\otimes_B B_i$ is exact in positive degrees and remains minimal: the entries of each differential lie in the maximal ideal. The top free term has the same positive rank. Thus
\[
 \pd_{B_i}M_i=\pd_B M,\qquad \dim B_i=\dim B .
\]
Auslander--Buchsbaum now gives $\depth_{B_i}M_i=\depth_B M$ at each such localization, since both local rings are Cohen--Macaulay. Faithfully flat completion detects the equalities. Away from the join the branch immersion is open locally. Hence any associated point of $M_i$ is associated to $M$, so each nonzero branch sheaf is pure two-dimensional. Associativity of derived pullback identifies its further restriction to $S$ with the common complex $E$.

Conversely, affine Milnor patching \cite[Theorem~1.3 and Remark~1.7]{Chough} on $B_1\times_R B_2$, followed by Zariski descent, identifies a matched pair, including its matching equivalence and homotopies, with a perfect object $P$ on $Z$. Tensor the normalization sequence with $P$ to obtain
\begin{equation}\label{eq:normaltriangleallnew}
 P\longrightarrow u_{1*}F_1\oplus u_{2*}F_2\longrightarrow j_*E\longrightarrow P[1].
\end{equation}
Since $E$ has cohomology in degrees $-1,0$, the only possible positive cohomology of $P$ is the cokernel of the map to $j_*\mathcal H^0(E)$. This map is onto by the quotient $F_i\twoheadrightarrow F_i/s_iF_i(-S)$ on either branch. There is no negative cohomology. The remaining terms give \eqref{eq:Torcorrectionnew} for $F=\mathcal H^0(P)$.

Its quotient $K$ lies in the direct sum of the pure branch sheaves. Its subobject $j_*\mathcal H^{-1}(E)$ lies in $u_{1*}F_1(-S)$ by \eqref{eq:twoHboundarynew}. For any subsheaf $B'\subset F$ of dimension at most one, its image in $K$ and its intersection with $j_*\mathcal H^{-1}(E)$ both vanish. Hence $B'=0$.

Over a Noetherian parameter scheme $T$, the normalization triangle on $Z\times T$ commutes with derived base change. In both directions the resulting perfect complex has derived parameter fibers in degree zero by the preceding calculation. The fiber criterion for relative Tor amplitude makes it a $T$-flat sheaf. Patching gives the equivalence on mapping spaces, and iteration at the disjoint joins gives the assertion for chains.
\end{proof}

Taking classes and then applying the same normalization triangle to $F^\vee\otimes^{\LL}F$ gives two consequences used in the degeneration formula:
\begin{align}
 [F]&=\sum_i[u_{i*}F_i]-\sum_e[j_{e*}\mathcal H^0(E_e)]
                      +\sum_e[j_{e*}\mathcal H^{-1}(E_e)]\quad\text{in }G_0(Z),\label{eq:generalKnormalnew}\\
 \RHom_Z(F,F)&\simeq\fib\!\left(\bigoplus_i\RHom_{Z_i}(F_i,F_i)
       \longrightarrow\bigoplus_e\RHom_{S_e}(E_e,E_e)\right).\label{eq:generalEndnew}
\end{align}
Here $Z_i$ denotes the $i$th component, including bubbles. Indeed $\LL u_i^*(F^\vee\otimes^{\LL}F)=F_i^\vee\otimes^{\LL}F_i$ and likewise at each $S_e$. Taking derived global sections gives the signed difference of the two restrictions at each double divisor. Formula~\eqref{eq:generalEndnew} governs deformations of the boundary identification itself. In particular deleting $\mathcal H^{-1}(E_e)$ from \eqref{eq:generalKnormalnew} can change the divisorial cycle of the glued sheaf.

\begin{corollary}\label{cor:negextstablegluednew}
If $F$ is simple as a sheaf on the complete fixed expanded target, every matched boundary complex $E_e$ satisfies $\Ext^j_{S_e}(E_e,E_e)=0$ for $j<0$. Consequently a nonzero boundary with negative self-Ext cannot occur in such a matching, even if it occurs in an independently stable branch moduli space.
\end{corollary}
\begin{proof}
The long exact sequence of \eqref{eq:generalEndnew} contains
\[
 0\longrightarrow\bigoplus_e\Ext^{-1}_{S_e}(E_e,E_e)\longrightarrow
 \Hom_Z(F,F)\longrightarrow\bigoplus_i\Hom_{Z_i}(F_i,F_i),
\]
since the branch sheaves have vanishing negative self-Ext.
The map on the right sends $\id_F$ to $(\id_{F_i})_i$, which is nonzero since at least one branch is nonzero. Its one-dimensional source therefore injects, proving the assertion in degree $-1$. In degrees below $-1$ the adjacent negative Ext groups of the sheaves $F$ and $F_i$ both vanish, and exactness gives the remaining assertion.

Higher-dimensional $\Hom_S(E_e,E_e)$ can remain; the argument imposes no purity condition on $\mathcal H^0(E_e)$.
\end{proof}
\begin{theorem}\label{thm:branchindexnew}
Let $Y$ be either root of the first degeneration, possibly expanded into a relative chain with outer divisor $S$. Thus $\omega_Y=\OO_Y(-S)$. Let $F$ be a simple perfect surface sheaf and $E=\LL i^*F$. With $v(E)=\ch(E)\sqrt{\operatorname{td}(S)}$ its Mukai vector, one has
\begin{align}
 2\chi_Y(F,F)&=\chi_S(E,E)=-\langle v(E),v(E)\rangle,\label{eq:generalindexnew}\\
 \Ext^3_Y(F,F)^\vee&\simeq\Ext^{-1}_S(E,E)\qquad(E\ne0).\label{eq:branchthirdnew}
\end{align}
The tangent complex of the derived branch stack at $F$ is $\RHom_Y(F,F)[1]$. On a bounded simple locus with $E\ne0$ and $\Ext^{-1}_S(E,E)=0$, scalar rigidification gives a two-term theory of virtual dimension $1+\langle v(E),v(E)\rangle/2$. This statement allows a boundary complex with cohomology in two degrees.
\end{theorem}
\begin{proof}
Tensor the Cartier sequence $0\to\omega_Y\to\OO_Y\to i_*\OO_S\to0$ with $F^\vee\otimes^{\LL}F$ and apply $\RG$. The result is a triangle with terms $A=\RHom_Y(F,F\otimes\omega_Y)$, $B'=\RHom_Y(F,F)$, $C=\RHom_S(E,E)$ in this order. Serre duality gives $A\simeq (B')^\vee[-3]$, so $\chi(C)=2\chi(B')$. Riemann--Roch on the K3 gives the second equality in \eqref{eq:generalindexnew}. The exact cohomology sequence reads
\[
 0=H^{-1}(B')\longrightarrow\Ext^{-1}_S(E,E)
 \longrightarrow\Hom_Y(F,F\otimes\omega_Y)
 \longrightarrow\CC\,\id_F\longrightarrow\Hom_S(E,E).
\]
The last map sends $\id_F$ to $\id_E$. It is injective for $E\ne0$. Duality gives \eqref{eq:branchthirdnew}. Higher Ext groups vanish by the same duality and the absence of negative Ext between sheaves. The unit removes $H^0(B')$ after scalar rigidification. If \eqref{eq:branchthirdnew} vanishes, the resulting tangent has cohomology only in degrees zero and one. Its virtual rank is $\ext^1_Y(F,F)-\ext^2_Y(F,F)=1-\chi_Y(F,F)$. On this locus $H^{>0}(Y,\OO_Y)=0$ by Proposition~\ref{prop:structurecohom}.

The scalar-rigidified derived moduli has the indicated perfect tangent of amplitude $[0,1]$ since the relative tangent of the scalar gerbe is $\OO[1]$ and its map to $B'[1]$ is the shifted unit. Explicitly the tangent after rigidification is $\operatorname{cofib}(\OO\xrightarrow{1\mapsto\id_F}B')[1]$. Restriction of its dual to the classical truncation supplies the cotangent obstruction morphism.

The fixed-boundary deformation complex before scalar reduction is $A[1]$, the fiber of $B'[1]\to C[1]$. For either smooth unexpanded anticanonical root the restriction morphism on perfect complexes carries the Lagrangian structure of \cite[Corollary~6.2]{BDRelative}; this is consistent with the displayed duality.
\end{proof}

To see why the Ext vanishing matters on the actual components, take $Y=B$, let $V$ be a line bundle on $S$ and put $F=i_*V$. Its support is the whole middle K3 and $F$ is stable for every ample line bundle on $B$: a proper nonzero subsheaf has generic rank one and quotient of dimension at most one. Since $N_{S/B}=\OO_S(-h)$, the two cohomology sheaves of $\LL i^*F$ are $V\otimes\OO_S(h)$ in degree $-1$ and $V$ in degree zero. Directly,
\[
 \Ext^3_B(F,F)^\vee=\Hom_B(i_*V,i_*(V\otimes\OO_S(h)))
 =H^0(S,\OO_S(h))\cong\CC^4 .
\]
Thus the unreduced higher tangent direction of this stable branch really occurs. Simplicity of the branch alone cannot supply a two-term scalar-rigidified relative theory at every boundary object.

There are also matched examples in the first degeneration. For $r\ge4$ take a section of $\OO(rH)$ on the total family with central equation $f_1a_{r-1}+f_4b_{r-4}$, where $a_{r-1}$ and $b_{r-4}$ are general forms (take $b_0\ne0$ when $r=4$). Its restrictions to $Q$ and $B$ are respectively $f_1a_{r-1}$ and the pullback of $f_4b_{r-4}$ from $\PP^3$. Neither vanishes identically on its root, and both vanish along $S$.

On the reduced Cohen--Macaulay union $Z$ this section is a non-zero-divisor. Let $T\subset Z$ be its Cartier zero locus. Near the central fiber the section defines a relatively flat Cartier family, by the fiberwise Cartier criterion. For any vector bundle $V$ on the total family, its restriction $V|_T$ is perfect and pure on $Z$. If $m=\rk V$, the defining two-term resolution gives
\[
 \LL j^*(V|_T)\simeq V|_S\oplus (V|_S\otimes\OO_S(-rh))[1].
\]
Here the differential vanishes because the section vanishes on $S$. The term in degree $-1$ is a rank-$m$ vector bundle on the whole K3. Formula~\eqref{eq:Torcorrectionnew} therefore changes the reconstructed sheaf even at the generic point of $S$. This computation concerns coherent sheaves on the family; a stability assertion for the expanded central fiber would require the choice of polarization.

The same family gives a non-locally-free example with stable generic fiber. On a smooth generic threefold the base of the sections $f_1a_{r-1}+f_4b_{r-4}$ is the curve $\Gamma$. Along $\Gamma$ the conormal directions $df_1,df_4$ are independent, by the complete intersection equations and $d(f_1f_4-tf_5)=-t\,df_5$ there. General $a,b$ have no common zero on $\Gamma$, so the section cuts the generic threefold smoothly along $\Gamma$; away from $\Gamma$ this is Bertini's theorem.

The resulting ample divisor is connected and hence integral. Choose a smooth point of the central surface $T$ outside $S$ on the quartic root and lift it to a section $p$ of the surface family after an etale base change. The ideal of this section is flat over the base by the sequence
\[
 0\longrightarrow I_{p/\mathscr T}\longrightarrow\OO_{\mathscr T}
       \longrightarrow\OO_{p}\longrightarrow0.
\]
Its generic fiber is a rank-one torsion-free sheaf on the integral surface, so it is Gieseker stable. At the special fiber the ideal $I_{p_0/T}$ is pure, perfect on $Z$, and has the same two-term derived restriction as $\OO_T$ at $S$. Perfectness holds near $S$ by the Cartier resolution and away from $S$ because the ambient root is smooth.

This central sheaf is simple. Indeed $T$ is a Cohen--Macaulay Cartier surface and $p_0$ is a smooth point of it. Dualizing the ideal sequence locally gives $\mathcal Hom_T(I_{p_0/T},\OO_T)=\OO_T$, using $\mathcal Ext^j_T(\OO_{p_0},\OO_T)=0$ for $j=0,1$. Multiplication by any function preserves the ideal, giving $\mathcal End_T(I_{p_0/T})=\OO_T$. To compute global sections, contract $Z$ to $\overline Z=(f_1f_4=0)\subset\PP^4$.

The two normalization sequences identify $\RR\pi_*\OO_Z=\OO_{\overline Z}$: the contraction is the identity on $Q$ and $S$, and the blowup $B\to\PP^3$ has $\RR b_*\OO_B=\OO_{\PP^3}$. The hypersurface sequence for $\overline Z$ and the cohomology of $\OO_{\PP^4}(m)$ give $H^0(Z,\OO_Z)=\CC$ and $H^1(Z,\OO_Z(-rH))=0$. The Cartier sequence now gives $H^0(T,\OO_T)=\CC$.

Thus simplicity of the glued special sheaf can coexist with a boundary complex having two nonzero cohomology sheaves of positive rank and two nonzero complementary projection endomorphisms.

Perfectness on the singular union is a genuine restriction on a limit. Locally the pure module $B/(x,y)=R$ supported on $S$ is not perfect over $B$: tensoring the displayed periodic branch resolution with $R$ gives $\Tor_n^B(R,B_1)=R$ for every $n>0$. Establishing that every required stable generic object admits a perfect limit, possibly after inserting components, belongs to the valuative problem.

The fixed-target matching theorem gives no such extension across a discrete valuation ring.

This local module occurs in the unexpanded quintic family itself. Make the base change $t=s^2$ and take the hyperplanes $f_1=s\ell$ for a general linear form $\ell$. The generic hyperplane section is integral by the primitive divisor class on the smooth quintic, so its structure sheaf is stable. Near a general point of $S$ where $f_5\ell\ne0$, trivialize the line bundles and put $x=f_1$, $y=f_4/f_5$. Then $xy=s^2$.

In the flat closure of the generic surface $x=s\ell$ one has $y=s/\ell$; the second equation follows by saturation with respect to $s$. Its special-fiber module there is $B/(x,y)$. In the expansion chart $x=sz$, $zy=s$, its strict transform has $z=\ell$, $y=s/\ell$, so the central support lies in the interior of the bubble on this open set. Thus failure of perfectness after contraction can occur for a flat specialization of a stable primitive sheaf.

This calculation is local along $S$; Proposition~\ref{prop:graphsupport} describes the corresponding global graph in the relative pencil.

\subsection{Derived restriction and the separation of limiting boundary values}
Let $i:S\hookrightarrow Y$ be a smooth Cartier divisor in a smooth threefold, and write $s$ for its defining section. For a coherent surface sheaf $F$ on $Y$ put $E=\LL i^*F$. Resolving $i_*\OO_S$ by $[\OO_Y(-S)\xrightarrow{s}\OO_Y]$ gives
\begin{equation}\label{eq:route19-trunc}
\begin{gathered}
 i_*A=\ker(F(-S)\xrightarrow{s} F),\qquad
 i_*B=\coker(F(-S)\xrightarrow{s} F),\\
 A[1]\longrightarrow E\longrightarrow B
       \xrightarrow{\kappa_E}A[2],\qquad
 \kappa_E\in\Ext^2_S(B,A).
\end{gathered}
\end{equation}
Indeed the two sheaves in the first line are annihilated by the ideal of $S$ and are the pushforwards of the cohomology sheaves in degrees $-1$ and $0$. The second line is the canonical truncation triangle on $S$. It accounts for the extension datum even when the two sheaves are known. Since $S$ is smooth, $E$ is perfect.

For a family flat over a base the same derived restriction commutes with derived base change; flat cohomology sheaves and a fixed ordinary sheaf moduli target require additional conditions. In the following family the general boundary is a stable line bundle on an integral curve, and its value at the special parameter depends on the compactification.

\begin{proposition}\label{prop:route19-two-arcs}
Let $Q$ be the smooth quartic threefold, $S=(s=0)$ a smooth hyperplane section, and $h=H_Q|_S$. Let $D\subset Q$ be a smooth complete-intersection curve meeting $S$ transversely in a nonempty scheme $Z$. Fix $k\in\mathbb Z$. The sheaf $F_0=i_*I_{Z/S}(kh)$ lies in the closure in the fixed-target Simpson space of stable surface sheaves whose restriction to $S$ belongs to the compactified Jacobian $J_{4k-2}$ of vector $(0,h,4k-2)$.

Boundary evaluation on this open locus does not extend to any neighborhood of $[F_0]$ in that closure. Its derived boundary is
\[
 \LL i^*F_0\simeq I_{Z/S}(kh)\oplus I_{Z/S}((k-1)h)[1].
\]
\end{proposition}
\begin{proof}
Choose $\ell\in H^0(Q,\OO_Q(H_Q))$ with $C_\ell=(\ell|_S=0)$ smooth and disjoint from $Z$. Such sections form a nonempty open set by Bertini and point avoidance. Over $A=\CC[[t]]$ let $T_t=(s+t\ell=0)$ and put
\[
 G=I_{D/Q}(kH_Q),\qquad
 F_{\ell}=\coker\big(G_A(-H_Q)\xrightarrow{s+t\ell}G_A\big)
 \quad\text{on }Q_A.
\]
Multiplication by $s$ is injective on $G$, so the displayed arrow is injective and its cokernel is flat over $A$ by the local criterion of flatness. Transversality and properness imply that $D_A\cap T$ is finite flat over $A$ after restriction to this trait. More explicitly, $s+t\ell$ cuts out a relative Cartier divisor on $D_A$ since $s|_D$ is a nonzerodivisor. Restricting the ideal sequence of $D_A$ to $T$ gives $F_{\ell}=j_*I_{D_A\cap T/T}(kH_Q)$ for $j:T\hookrightarrow Q_A$.

The closed and generic surfaces are smooth integral hyperplane sections, by openness at $S$. A rank-one torsion-free sheaf on either surface is Gieseker stable: any proper nonzero subsheaf has rank one and its quotient has positive Hilbert polynomial of lower degree. Thus $F_\ell$ defines a Simpson-space arc with special value $[F_0]$.

Write $V=I_{Z/S}(kh)$. The same ideal sequence and transversality give $\LL i^*G=V$. Hence on the generic boundary $\LL i^*F_{\ell,\eta}$ is represented in degrees $-1,0$ by
\[
 [\,V_\eta(-h)\xrightarrow{t\ell}V_\eta\,]
   \simeq \OO_{C_\ell}(kh)_\eta .
\]
Here $V$ is torsion-free and agrees with $\OO_S(kh)$ near $C_\ell$. At the closed point $F_0$ has the resolution $[G(-H_Q)\xrightarrow{s}G]$ by perfect sheaves on $Q$. Its derived pullback is $[V(-h)\xrightarrow{0}V]$, which gives the asserted splitting. Its Mukai vector is
\[
 v(V)-v(V(-h))=(0,h,(k-\tfrac12)h^2)=(0,h,4k-2).
\]
Its two nonzero summands give two independent scalar automorphisms. In particular it is not a stable complex for any Bridgeland stability condition.

Now take $\ell'$ satisfying the same conditions with $C_{\ell'}\ne C_\ell$. Both arcs have closed value $[F_0]$. An extension of evaluation to a morphism into $J_{4k-2}$ would induce, on the first trait, a morphism agreeing generically with the constant point $[\OO_{C_\ell}(kh)]$. Since $J_{4k-2}$ is separated, it must agree also at the closed point. The second trait forces the same image to be $[\OO_{C_{\ell'}}(kh)]$, a distinct point by its support. This proves the assertion.
\end{proof}

\section{Expanded sheaf moduli and stable reduction}\label{sec:compactification}\label{sec:generalafterexample}
The first degeneration has central fiber $Y_1\cup_S Y_2$. A family of sheaves may have different expanded targets in two different presentations. To test separatedness one must extend both the isomorphism of targets and the isomorphism of sheaves over it. We specify these isomorphisms, then determine which conditions of \cite[Assumptions~3.0.1 and~3.2.2]{Kim} can be checked on the expansions of this family. The same definitions apply to the second degeneration.

\subsection{Expanded targets, sheaves, and their isomorphisms}
Write $\XX\to C$ for either simple degeneration and $S$ for its joining divisor. The two branches of $\XX_0$ will be denoted by $Y_1,Y_2$ here. For the general indexing used in this section we set $D_e=S_{e-1}$ in the notation of the introduction. In an expansion of length $k$ the special fiber is the union
\[
 W_0\cup_{D_1}W_1\cup_{D_2}\cdots\cup_{D_k}W_k\cup_{D_{k+1}}W_{k+1}.
\]
Here $W_0=Y_1$, $W_{k+1}=Y_2$ and $W_a=\PP_S(\OO_S\oplus N_{S/Y_1})$ for $1\le a\le k$. The two disjoint sections of each ruled component (a \emph{bubble}) are identified with the corresponding sections or root divisors on the neighboring components. Each $D_e=W_{e-1}\cap W_e$ is a copy of $S$. The contraction $c$ is the identity on the two roots and the ruling projection on each bubble.

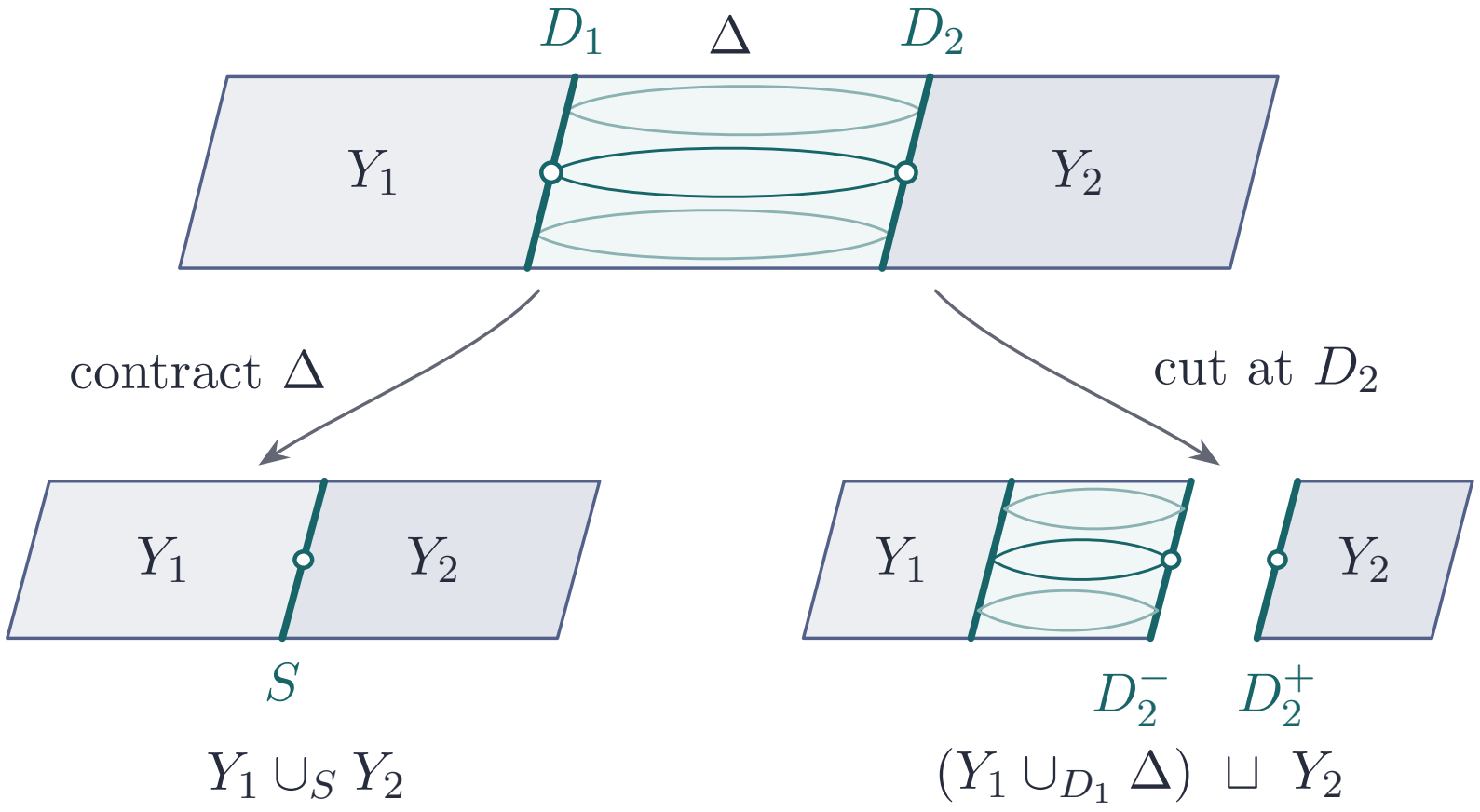
\begin{figure}[htbp]
\centering
\begingroup
\definecolor{qAmbient}{HTML}{53618A}
\definecolor{qSeam}{HTML}{176568}
\definecolor{qInk}{HTML}{272C3D}
\begin{tikzpicture}[x=1cm,y=1cm,font=\small,
  line cap=round,line join=round,
  side/.style={draw=qAmbient,line width=.65pt},
  sectionink/.style={draw=qSeam,line width=1.45pt},
  maparr/.style={-{Stealth[length=2.2mm,width=1.6mm]},draw=qInk!72,line width=.65pt}]
\def\tile#1#2#3{%
 \path[side,fill=#3] (#1,-.68)--({#1+.34},.68)--({#2+.34},.68)--(#2,-.68)--cycle;}
\def\jdiv#1{\draw[sectionink] (#1,-.68)--({#1+.34},.68);}
\def\ruling#1#2#3#4{%
 \pgfmathsetmacro{\offset}{.17+.25*(#3)}
 \draw[draw=#4,line width=.55pt]
 ({#1+\offset},#3)
 ..controls ({#1+\offset+.33},{#3+.23}) and ({#2+\offset-.33},{#3+.23})..
 ({#2+\offset},#3)
 ..controls ({#2+\offset-.33},{#3-.23}) and ({#1+\offset+.33},{#3-.23})..
 cycle;}
\path[use as bounding box] (-6.15,-3.88) rectangle (6.72,2.04);
\begin{scope}[shift={(-.17,.72)}]
\tile{-3.73}{-1.26}{qAmbient!11}
\tile{-1.26}{1.26}{qSeam!6}
\tile{1.26}{3.73}{qAmbient!17}
\ruling{-1.26}{1.26}{-.44}{qSeam!50}
\ruling{-1.26}{1.26}{0}{qSeam}
\ruling{-1.26}{1.26}{.44}{qSeam!50}
\jdiv{-1.26}\jdiv{1.26}
\node[text=qInk] at (-2.35,0) {$Y_1$};
\node[text=qInk] at (2.65,0) {$Y_2$};
\node[text=qInk,above] at (.18,.70) {$\Delta$};
\node[text=qSeam,above] at (-.94,.70) {$D_1$};
\node[text=qSeam,above] at (1.62,.70) {$D_2$};
\foreach \p in {-1.09,1.43}{
 \filldraw[fill=white,draw=qSeam,line width=.85pt] (\p,0) circle[radius=2pt];}
\end{scope}
\draw[maparr] (-1.35,-.12).. controls (-1.75,-.56) and (-2.77,-.98)..(-3.34,-1.36);
\draw[maparr] (1.47,-.12).. controls (1.88,-.56) and (2.92,-.98)..(3.49,-1.36);
\node[anchor=east,text=qInk] at (-2.71,-.68) {contract $\Delta$};
\node[anchor=west,text=qInk] at (2.87,-.68) {cut at $D_2$};
\begin{scope}[shift={(-3.17,-2.03)},xscale=.88,yscale=.82]
\tile{-2.22}{0}{qAmbient!11}
\tile{0}{2.22}{qAmbient!17}
\jdiv{0}
\filldraw[draw=qSeam,line width=.85pt,fill=white] (.17,0) circle[radius=2pt];
\node[text=qInk] at (-.96,0) {$Y_1$};
\node[text=qInk] at (1.22,0) {$Y_2$};
\node[text=qSeam,below] at (0,-.72) {$S$};
\end{scope}
\begin{scope}[shift={(3.25,-2.03)},xscale=.85,yscale=.82]
\tile{-3.20}{-1.80}{qAmbient!11}
\tile{-1.80}{-.30}{qSeam!6}
\tile{.59}{2.05}{qAmbient!17}
\jdiv{-1.80}\jdiv{-.30}\jdiv{.59}
\ruling{-1.80}{-.30}{-.44}{qSeam!50}
\ruling{-1.80}{-.30}{0}{qSeam}
\ruling{-1.80}{-.30}{.44}{qSeam!50}
\foreach \p in {-.13,.76}{
 \filldraw[fill=white,draw=qSeam,line width=.85pt] (\p,0) circle[radius=2pt];}
\node[text=qInk] at (-2.38,0) {$Y_1$};
\node[text=qInk] at (1.49,0) {$Y_2$};
\node[text=qSeam,below] at (-.45,-.73) {$D_2^{-}$};
\node[text=qSeam,below] at (.76,-.73) {$D_2^{+}$};
\end{scope}
\node[text=qInk] at (-3.0,-3.56) {$Y_1\cup_S Y_2$};
\node[text=qInk] at (2.92,-3.56) {$(Y_1\cup_{D_1}\Delta)\ \sqcup\ Y_2$};
\end{tikzpicture}
\endgroup
\caption{An expansion with one ruled component $\Delta$ and a marking at $D_2$.
Each oval depicts a real circle in a ruling $\mathbb{P}^1$ of $\Delta\to S$.
The two open circles on the middle oval lie over the same $s\in S$.
Contraction sends their entire ruling to $s$.
Cutting at $D_2$ keeps the ruling and produces two copies
$D_2^{-},D_2^{+}$ of the marked surface, on $\Delta$ and $Y_2$, respectively.}
\label{fig:expansioncontraction}
\end{figure}

A \emph{marking} chooses one of the $D_e$ in a fiber, or a section of the family of such choices over a base. A \emph{cut at the marking} is partial normalization along $D_e$. It leaves two copies of $D_e$, one on each resulting chain. The two copies must be identified when the chains are glued back. Figure~\ref{fig:expansioncontraction} shows the contraction and a choice of marking for $k=1$.

More generally an expanded pair $(Y,S)$ has a root $Y$, followed by $k$ bubbles and no second root. Its \emph{outer relative divisor} is the free section of its last bubble, or $S\subset Y$ for $k=0$. Normality to this divisor will mean the Tor condition below. The term \emph{singular divisor} in \cite[Definitions~2.4.3 and~2.4.11]{Kim} refers to a joining divisor in the singular locus of the total fiber. The divisor itself is smooth.

Let $U_k=C\times_{\AA^1}\AA^{k+1}$ locally at $0\in C$, where the map on the right sends $(t_0,\ldots,t_k)$ to $\prod t_j$. A standard family $Z[k]\to U_k$ has local equation $xy=t_j$ at its $j$-th crossing, up to the smooth parameters along $S$. Its torus $T_k=\GG^k$ acts by
\[
 t_j\longmapsto\lambda_j^{-1}\lambda_{j+1}t_j
 \quad(0\le j\le k),\qquad \lambda_0=\lambda_{k+1}=1.
\]
Thus the product $t$ is invariant. On $t_j\ne0$ the crossing $xy=t_j$ is smooth. Deleting this parameter, after absorbing its unit into an adjacent parameter, identifies this target with one given by a shorter chart. For $k=1$ the opens $t_0\ne0$ and $t_1\ne0$ both describe the target $\XX_t$ with expansion length zero. The target stack $\EE$ retains these identifications as well as the torus rescalings \cite[Definition~2.4.7 and Theorem~2.4.9]{Kim}.

The charts for pairs use $\AA^k$ and the corresponding identifications \cite[Definitions~2.4.10 and~2.4.15]{Kim}.

Hereafter the \emph{groupoid of expanded targets over $T$} means $\EE(T)$. Its objects are the families $c:Z\to\XX\times_C T$ obtained locally on $T$ by these charts. Its arrows are isomorphisms of expansions preserving $c$. In the case of a pair they must also preserve the outer relative divisor. To specify the arrows and their action on sheaves, put
\[
 U=\coprod_k U_k,\qquad R=U\times_{\EE}U,qquad
 s,t:R\rightrightarrows U .
\]
A $T$-point of $R$ is a triple $(u,u',\sigma)$ with $u,u':T\to U$ and a target isomorphism $\sigma:Z_u\!\xrightarrow{\sim}\!Z_{u'}$; $s$ selects $u$ and $t$ selects $u'$. The universal arrow is an isomorphism $s^*Z\simeq t^*Z$ over $R$.

Fix a numerical type $v$ and a rule specifying the stable locus on each projective chart. In the \emph{groupoid of expanded sheaves before scalar rigidification}, $\widetilde{\MM}(v)(T)$, the objects are $(c:Z\to\XX_T,F)$, with $F$ finitely presented, $T$-flat, pure of dimension two and perfect on each geometric fiber, of the specified numerical type and stable by that rule.

The polynomial recording the original numerical type is computed using $c^*A_0$ for the original relatively ample line bundle $A_0$ on $\XX/C$. This polynomial and the chart polarization used to test stability have different roles. The arrows and their composition are
\[
 \begin{split}
 (\sigma,\phi):(Z,F)\longrightarrow(Z',F'),\quad&
  \sigma\in\operatorname{Isom}_{\EE(T)}(Z,Z'),\qquad
  \phi:F\xrightarrow{\sim}\sigma^*F',\\
 (\tau,\psi)\circ(\sigma,\phi)&=(\tau\sigma,\sigma^*\psi\circ\phi).
 \end{split}
\]
Pullback along $T'\to T$ pulls back both entries. In the scalar rigidification $\MM(v)=\widetilde{\MM}(v)\!\mathbin{\! /\! /}\GG$, one divides the isomorphisms $\phi$ by scalar units, then sheafifies the isomorphism functors and stackifies in the fppf topology \cite[\S5.1]{ACV}. The group of target rescalings preserving a sheaf may still act after this operation.

\subsection{Admissibility, the error polynomial and stable reduction}\label{sec:kimexactfull}
At a joining divisor $D\subset Z$ or at the outer divisor of a pair define normality of $F$ to $D$ by
\[
 \mathcal Tor_1^{\OO_Z}(F,\OO_D)=0 .
\]
Admissibility requires this vanishing for each such divisor on each geometric fiber \cite[Definitions~2.5.1, 2.5.7 and~2.5.8]{Kim}. It is a condition on a sheaf relative to an immersion. Its numerical test uses the following polynomial.

For an ideal sheaf $I$ and coherent $M$ on a Noetherian scheme put
\[
 \Gamma_I(M)=\bigcup_{b\ge1}\ker(M\longrightarrow
          \mathcal Hom(I^b,M)),\qquad m\longmapsto(a\mapsto am).
\]
Thus $\Gamma_I(M)$ consists of sections annihilated locally by a power of $I$. On a Noetherian affine open the kernels form an increasing sequence of submodules of a finite module, hence eventually agree. On our quasi-compact targets a finite cover gives one exponent that suffices. For the ideal $I_e$ of an internal joining divisor $D_e$, set
\[
 U_e=\Gamma_{I_e}(F),\quad
 F'=F\big/\!\bigoplus_{e\ {\rm internal}} U_e,\quad
 V_{i,e}=\Gamma_{I_e\OO_{W_i}}(F'\otimes\OO_{W_i})
       \quad (i=e-1,e).
\]
The sum is direct because the distinct joining divisors are disjoint. For a sheaf $M$ on $Z$ or one of its components, write
$P^c_M(m)=\chi(M\otimes(c^*A_0)^{\otimes m})$, restricting $c^*A_0$ as appropriate. The contribution to the error from $D_e$ is
\begin{equation}\label{eq:kimerror}
 \operatorname{Err}_e(F)=P^c_{U_e}
             +\tfrac12 P^c_{V_{e-1,e}}+\tfrac12 P^c_{V_{e,e}} .
\end{equation}
The first term measures torsion on the union. Even when this term vanishes, ordinary restriction to a branch can create torsion supported on $D_e$. The two half-weighted terms measure it on the two sides \cite[Definition~2.5.9, (2.17)]{Kim}. For a pair there is one more term,
$\operatorname{Err}_{\rm out}(F)=P^c_{\Gamma_{I_{\rm out}}(F)}$.
The total $\operatorname{Err}(F)$ is the sum over the internal divisors, with this additional term in the relative case
\cite[\S2.5.3, Lemmas~2.5.10 and~2.5.13]{Kim},
\cite[Definition~3.15, (3.10) and Lemma~3.19]{LW}.

All these torsion sheaves have support set-theoretically in copies of $S$. The contraction is an isomorphism on each of these reduced divisors. Filtration by powers of its ideal expresses the polynomial of each nonzero term as a sum of Hilbert polynomials with positive top coefficient on $S$. In particular $\operatorname{Err}(F)=0$ exactly when all the torsion terms vanish, equivalently when $F$ is admissible.

A \emph{degree-one error} means a nonzero total polynomial $am+b$ with $a>0$.

Chapter~3 of \cite[p.~59]{Kim} treats generically rank-one sheaves of irreducible effective divisor class. For the fixed numerical type, \cite[Assumption~3.0.1(1)--(3)]{Kim} specifies a bound $N$ and a line bundle $\widetilde H$ on the universal target of expansions of length $\le N$:
\begin{enumerate}
\item Its pullback to each standard projective chart is relatively ample, and its restriction to the unexpanded target is $A_0^{\otimes m_0}$ for a fixed $m_0>0$.
\item Each $\widetilde H$-Gieseker-semistable sheaf of the chosen type is stable.
\item For each $\widetilde H$-stable sheaf of that type, the error polynomial is zero or has degree different from one. This tests the sheaf before imposing admissibility.
\end{enumerate}
For the relative groupoid \cite[Assumption~3.2.2]{Kim} further requires $k<M<N$ for every occurring expansion length $k$. On a multiple $rH$, both the number of support components and the generic ranks can vary.

An object of \cite[Definitions~3.2.1 and~3.2.5]{Kim} is a pure stable admissible sheaf of fixed $c^*A_0$ polynomial admitting a presentation by the prescribed locally free source, with finite stabilizer as a quotient. In the relative case its restriction to the outer divisor is pure of dimension one. For these groupoids, \cite[Theorems~4.1.1--4.1.2 and~4.2.1--4.2.2]{Kim} prove valuative existence after finite base change and extension of generic isomorphisms, respectively.

Our stack $\MM(v)$ allows $\mathcal H^{-1}(\LL i^*F)\ne0$. Its valuative criteria concern the pairs $(Z,F)$ and arrows $(\sigma,\phi)$ of the preceding subsection. The following proposition checks the third numerical condition on the perfect locus.

\begin{proposition}\label{prop:kim-unexpanded-root}
On a geometric fiber of an expansion of either degeneration or of either relative pair, let $F$ be a pure two-dimensional sheaf which is perfect on that fiber. Then $\operatorname{Err}(F)$ is zero or of degree two, in every divisor multiplicity. Thus the numerical exclusion in \cite[Assumption~3.0.1(3)]{Kim} holds on this perfect locus for every choice of chart polarization. On an unexpanded smooth pair the conclusion holds for all pure two-dimensional sheaves.

On the unexpanded quartic root $Q$, every such sheaf with divisor class $H_Q$ is stable for every ample line bundle. Consequently both parts (2) and (3) of that assumption hold on this root in primitive divisor class for all lower charges.
\end{proposition}
\begin{proof}
If $U_e\ne0$, then $U_e\subset F$ has pure dimension two. For its ideal-power filtration, set $L_b=I_e^b U_e/I_e^{b+1}U_e$; only finitely many $L_b$ are nonzero. Their support lies in $D_e\simeq S$, and at least one has positive generic rank on $D_e$. Hence, writing $a=c^*A_0|_{D_e}$, one obtains
\[
 [m^2]P^c_{U_e}(m)=\tfrac12\sum_b c_1(a)^2\cdot[L_b]_2>0 .
\]
Each other term in the total error has degree at most two and nonnegative quadratic coefficient, so the total has degree two. The same calculation applies to nonzero outer torsion.

Suppose all $U_e=0$ and the outer torsion, when present, is zero. Then $F'=F$. Each branch restriction $F_i=F\otimes\OO_{W_i}$ is zero or pure of dimension two. Here is the local reason. At a double divisor the completed local ring is $B=R[[x,y]]/(xy)$ with $R$ regular, and $B_1=B/(y)$. The free resolution of $B_1$ over $B$ is periodic in positive degrees, with alternating differentials $y,x$. For the stalk $M$ of $F$, perfectness and periodicity give
\[
 \Tor^B_j(M,B_1)=0\quad(j\ge1).
\]
A minimal finite free resolution of $M$ at a point of the divisor therefore restricts to a free resolution on the branch. All differential entries remain in the maximal ideal, hence the restricted resolution is minimal of the same length. The two local Cohen--Macaulay rings have equal dimension, so the Auslander--Buchsbaum formula gives
$\depth_{B_1}(M\otimes B_1)=\depth_B M$ at each such localization.

Associated points on the branch are consequently associated points of $F$; off the double divisor the assertion follows from the open immersion. The other branch has the same calculation with $x,y$ interchanged. Each nonzero $V_{i,e}\subset F_i$ thus has dimension two, and the filtration argument gives $[m^2]P^c_{V_{i,e}}>0$. If every $V_{i,e}=0$ then \eqref{eq:kimerror} vanishes identically.

This proves the first assertion without a normality hypothesis on $F$ or purity of its restriction to $S$.

On an unexpanded pair the only contribution is $P_{\Gamma_I(F)}$ with $I=\OO_Y(-S)$. Locally $I=(u)$, and
\[
 \Gamma_I(F)=0\iff
 \ker(u:F\to F)=0\iff\mathcal Tor_1^Y(F,\OO_S)=0 .
\]
If it is nonzero, purity again gives degree two. Finally $\Pic(Q)=\mathbb Z H_Q$, and every effective integral divisor has class $bH_Q$ for $b\in\mathbb Z_{>0}$. Thus $[F]=H_Q$ forces one integral support with generic multiplicity one. For $0\ne G\subsetneq F$ purity gives the same generic multiplicity, so $\dim(F/G)\le1$. The quadratic coefficients of $P_G,P_F$ coincide and $P_F-P_G=P_{F/G}>0$ for $m\gg0$. Dividing by the common leading coefficient gives $p_G<p_F$.
\end{proof}

The perfectness assumption in the first assertion has numerical content on a bubble. For an explicit computation take the one-bubble pair $Z=Q\cup_{D_-}\Delta$, with projection $\rho:\Delta\to S$ and outer section $D_+$. Choose $m\gg0$ and a general member
\[
 T\in |\OO_\Delta(1)\otimes\rho^*\OO_S(mh)|.
\]
This linear system is very ample for $m\gg0$. Thus $T$ is integral and meets both sections in nonempty Cartier curves $C_\pm=T\cap D_\pm$. The divisor sequence and $R\rho_*\OO_\Delta(-1)=0$ give $R\rho_*\OO_T=\OO_S$, so its contracted two-cycle is the primitive root class $[S]=H_Q$. Extend $F=\OO_T$ from $\Delta$ to $Z$ by the closed immersion. Its rank one on the integral $T$ gives stability for every ample line bundle on $Z$.

On $D_-$ the three sheaves in \eqref{eq:kimerror} are respectively $0$, $\OO_{C_-}$ (root), and $0$ (bubble). The outer torsion vanishes. Indeed $F\otimes\OO_Q=\OO_{C_-}$, whereas neither $D_-$ nor $D_+$ contains the generic point of $T$. Hence
\[
 \operatorname{Err}(F)=\tfrac12P_{\OO_{C_-}},\qquad
 [m']\operatorname{Err}(F)(m')=\tfrac12\deg_{A_0} C_->0 .
\]
This sheaf is nonperfect on the union. In the periodic resolution of the branch algebra, one branch parameter annihilates $F$ and the other is regular on $F$ along $C_-$; tensoring gives a copy of $\OO_{C_-}$ in every positive odd degree for the quotient by the annihilating parameter. The example distinguishes the full sheaf category tested by \cite[Assumption~3.0.1(3)]{Kim} from its perfect locus.

\begin{proposition}\label{prop:unlabelledpolarization}
Let $\mathfrak X_{\le N}\to\EE_{\le N}$ be the universal target on the stack of expansions of $\XX/C$ of length at most $N\ge1$. Its arrows include the identifications of charts obtained by deleting invertible parameters \cite[Definition~2.4.7 and Theorem~2.4.9]{Kim}; no joining divisor is marked. Every line bundle on $\mathfrak X_{\le N}$ has degree zero on each ruling of the bubble on the one-bubble chart. In particular this universal target admits no line bundle that is relatively ample on that chart.
\end{proposition}
\begin{proof}
Restrict to length at most one. Let $j:\XX\hookrightarrow\mathfrak X_{\le1}$ be the open target of length zero. For $L\in\Pic(\mathfrak X_{\le1})$ set $A=j^*L$. On the one-bubble chart write $a:Z[1]\to\mathfrak X_{\le1}$ for the chart map, $c:Z[1]\to\XX$ for the contraction, and $s,v$ for the base parameters, with $t=sv$. On $U=Z[1]\setminus V(s,v)$ either $s$ or $v$ can be deleted. Deleting $s$ or $v$ on the overlap gives the same unexpanded target map to $\XX$. Pullback of $L$ along this identification yields
\[
 (a^*L)|_U\simeq(c^*A)|_U.
\]
The one-bubble resolution of $xy=sv$ has charts
\[
 x=sz,\ v=zy;\qquad s=xw,\ y=wv .
\]
Omitting the smooth coordinates of $S$, they have regular coordinates $(s,z,y)$ and $(x,w,v)$ respectively. The complement of $U$ is cut out by $(s,zy)$ and $(xw,v)$, respectively, and has codimension two in both. Elsewhere regularity and the same codimension follow from smoothness away from the joining divisors and flatness over the parameter plane.

Since $Z[1]$ is regular, a line bundle and its inverse extend uniquely across this complement: for $i:U\hookrightarrow Z[1]$ and invertible $B$, the reflexive extension formula is $B=i_*i^*B$. Applying it to $a^*L\otimes(c^*A)^{-1}$ gives $a^*L\simeq c^*A$. On any ruling $\ell\simeq\PP^1$, $c(\ell)$ is a single point, so
\[
 (a^*L)|_\ell\simeq\OO_\ell\otimes A_{c(\ell)},\qquad
 \deg_\ell a^*L=0.
\]
Relative ampleness on the projective chart would give positive degree on this curve in a fiber.
\end{proof}

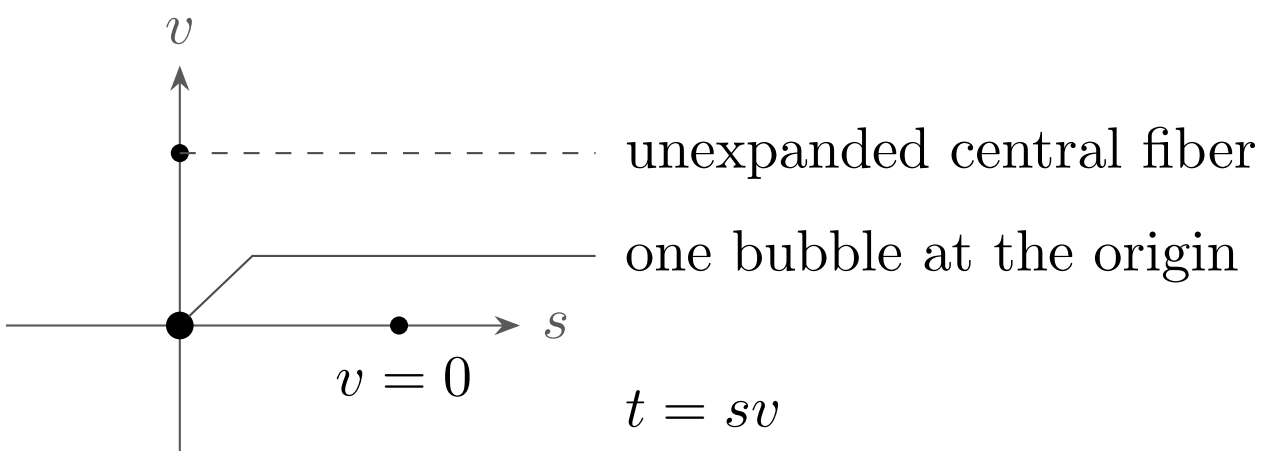
\begin{figure}[htbp]
\centering
\begin{tikzpicture}[x=1cm,y=1cm,>=Stealth,font=\small]
\draw[->,black!65] (-1.15,0)--(2.25,0) node[right]{$s$};
\draw[->,black!65] (0,-.83)--(0,1.72) node[above]{$v$};
\fill (0,0) circle (2.6pt); \fill (1.45,0) circle (1.7pt);\fill (0,1.14) circle (1.7pt);
\node[anchor=north] at (1.48,-.07){$v=0$};
\draw[dashed,black!70] (0,1.14)--(2.75,1.14);
\node[anchor=west] at (2.8,1.14){unexpanded central fiber};
\draw[black!70] (.06,.06)--(.48,.46)--(2.75,.46);
\node[anchor=west] at (2.8,.46){one bubble at the origin};
\node[anchor=west] at (2.8,-.55){$t=sv$};
\end{tikzpicture}
\caption{The one-bubble chart. Both punctured axes parametrize the unexpanded central fiber and the origin retains the ruled component. A line bundle descending through both deletion charts agrees with a contracted line bundle away from the fiber over the origin.}\label{fig:onebubblecoordinates}
\end{figure}

Figure~\ref{fig:onebubblecoordinates} shows where the two presentations of the unexpanded target overlap the one-bubble chart. The obstruction in Proposition~\ref{prop:unlabelledpolarization} applies to descending a polarization. To descend a moduli problem it suffices that the stable open subfunctors for the chosen chart polarizations agree after pullback to each overlap. The tautological pullback isomorphisms then supply the cocycle on triple overlaps. This is a condition to verify for each numerical problem. Section~\ref{sec:verifiedcompactification} checks it on the ordinary, graph and zero-boundary charts of the hyperplane structure-sheaf family. We impose the following conditions on the groupoid of pairs $(c:Z\to\XX_T,F)$ defined above.

\begin{hypothesis}\label{hyp:compact}
Fix $v$ and a neighborhood of $0$ in $C$.
\begin{enumerate}[label=(K\arabic*)]
\item On each expansion chart $U_i$ choose an open subfunctor of the flat families of pure sheaves perfect on geometric fibers. Its stability is Gieseker stability for a specified relatively ample line bundle $A_i$ or a condition proved to define that same open subfunctor. For every target isomorphism $\sigma:Z_u\to Z_{u'}$ over $T$, membership of $F'$ on $Z_{u'}$ is equivalent to membership of $\sigma^*F'$ on $Z_u$. This equivalence holds after arbitrary base change. Admissibility, when imposed, is an additional condition.
\item There is a bound $N(v)$ on the length of every occurring expansion, and finitely many relative Quot schemes of finite type, with prescribed polarizations and Hilbert polynomials, suffice to present all the families locally with a frame. On the smooth generic target their quotient is the full required stable moduli $M_\eta(v)$.
\item If $R$ is a discrete valuation ring with fraction field $K$ and $\Spec R\to C$, every object of $\MM(v)(K)$ extends to an object of $\MM(v)(R')$ after a finite extension of discrete valuation rings $R\subset R'$. This also tests $K$-points represented by already expanded targets.
\item For all $x,y\in\MM(v)(R)$ the restriction map
\[
 \operatorname{Isom}_{\MM(v)(R)}(x,y)\longrightarrow
 \operatorname{Isom}_{\MM(v)(K)}(x_K,y_K)
\]
is bijective. For each geometric $(Z,F)$, the subgroup of target automorphisms carrying $F$ to an isomorphic sheaf is finite.
\end{enumerate}
\end{hypothesis}
The obstruction to condition~(2) of \cite[Assumption~3.0.1]{Kim} in arbitrary divisor multiplicity is already visible on the smooth fiber: for stable $G$, $G^{\oplus r}$ is strictly semistable when $r>1$, as in Proposition~\ref{prop:unavoidablesemistability}. For a fixed charge Proposition~\ref{prop:generalcoprime} tests equality of the stable and semistable loci on that smooth fiber. Conditions (K1)--(K4) concern families on every allowed expansion.

In particular (K4) asks for the extension of the actual pair $(\sigma_K,\phi_K)$ modulo scalars, including its target isomorphism $\sigma_K$.

\subsection{Descent of the general expanded charts}
For a chart $U_i$ write $\mathfrak S_i$ for its stack of flat sheaves perfect on geometric fibers, and $\mathfrak S_i^s$ for the chosen open stable substack. On $R_{ij}=U_i\times_{\EE}U_j$, pullback by the inverse of the universal target isomorphism induces
$\Psi_{ij}:s^*\mathfrak S_i\xrightarrow{\sim}t^*\mathfrak S_j$.
The condition for gluing the chosen open substacks is
\begin{equation}\label{eq:generalstableoverlaps}
 \Psi_{ij}(s^*\mathfrak S_i^s)=t^*\mathfrak S_j^s,\qquad
 \Psi_{jk}\circ\Psi_{ij}\simeq\Psi_{ik}
 \quad\text{on }U_i\times_{\EE}U_j\times_{\EE}U_k .
\end{equation}
The first condition is equality of open subfunctors at all test schemes. The second identification is the canonical one given by composition of pullbacks.

\begin{proposition}\label{prop:generalchartdescent}
Assume \eqref{eq:generalstableoverlaps}, the finite cover by bounded charts of (K2) and finite intrinsic geometric stabilizers as in (K4). The scalar-rigidified stack $\MM(v)$ of simple stable sheaves on these expanded targets is a Deligne--Mumford stack of finite type. The sheaf whose deformations it parametrizes over $\EE$ is $F$ on $Z$.

In particular, before scalar rigidification its first-order deformations on a fixed geometric target have space $\Ext^1_Z(F,F)$, with obstructions in $\Ext^2_Z(F,F)$. No choice of a family of Cohen--Macaulay extensions of $F$ is required.
\end{proposition}
\begin{proof}
On a bounded chart with projective morphism $p:Z_i\to U_i$, choose $m$ so every fiber $F(m)=F\otimes A_i^{\otimes m}$ is globally generated and has vanishing higher cohomology. Treat each of the finitely many chart Hilbert polynomials $P_i$ separately. Cohomology and base change give $p_*F(m)$ locally free of rank $d=P_i(m)$. A frame of this vector bundle yields
\[
 \OO_{Z_i}(-m)^{\oplus d}\twoheadrightarrow F,\qquad
 \OO_{U_i}^{\oplus d}\xrightarrow{\sim}p_*F(m).
\]
The inverse takes the induced basis of sections. With $Q_i$ the framed Quot open and $R_{ij}=U_i\times_{\EE}U_j$, this gives
\[
 \mathfrak S_i^s=[Q_i/\mathrm{GL}_d],\qquad
 \Psi_{ij}:s^*\mathfrak S_i^s\xrightarrow{\sim}t^*\mathfrak S_j^s,
 \qquad \Psi_{jk}\Psi_{ij}\simeq\Psi_{ik}.
\]
The last two identifications are \eqref{eq:generalstableoverlaps}. Effective faithfully flat descent of finitely presented flat sheaves identifies the descended stack with $\widetilde{\MM}(v)$, with its specified target isomorphisms. Its smooth atlas $\coprod_iQ_i\twoheadrightarrow\widetilde{\MM}(v)$ has finite type by (K2).

At a geometric object $x=(Z,F)$ let $H_x$ be the fppf image of $\Aut_{\widetilde{\MM}(v)}(x)\to\Aut_{\EE}(Z)$. On geometric points it consists of the $\sigma$ for which $\sigma^*F\simeq F$. Simplicity gives the following exact sequence of group sheaves:
\begin{equation}\label{eq:intrinsicaut}
 1\longrightarrow\GG\xrightarrow{\lambda\mapsto(1,\lambda\id_F)}
 \Aut_{\widetilde{\MM}(v)}(x)
 \xrightarrow{(\sigma,\phi)\mapsto\sigma}H_x\longrightarrow1 .
\end{equation}
By rigidification \cite[\S5.1]{ACV}, $\Aut_{\MM(v)}(x)=H_x$. For $x,y:T\to\MM(v)$ over $C$ the square
\[
\begin{tikzcd}
 \operatorname{Isom}_T(x,y)\arrow[r]\arrow[d]&\MM(v)\arrow[d,"\Delta_{\MM(v)/C}"]\\
 T\arrow[r,"{(x,y)}"']&\MM(v)\times_C\MM(v)
\end{tikzcd}
\]
is Cartesian. At a geometric isomorphism $g:x_t\simeq y_t$ the fiber is an $H_{x_t}$-torsor. Finiteness and characteristic zero give $T_g\operatorname{Isom}_t(x_t,y_t)=\operatorname{Lie}(H_{x_t})=0$. The locally finitely presented diagonal is therefore unramified, proving the Deligne--Mumford assertion.

For a square-zero thickening with ideal $J$, extensions of a perfect complex $F$ on a prescribed target thickening have obstruction
\[
 (1_F\otimes\mathrm{KS})\circ\mathrm{At}(F)
  \in \Ext^2_Z(F,F\otimes J).
\]
When this class vanishes, the set of isomorphism classes of lifts is a torsor under $\Ext^1_Z(F,F\otimes J)$ \cite[Theorem~3.3]{HT}. Here $J$ is pulled back from the square-zero ideal of the parameter scheme; $\mathrm{At}(F)$ and $\mathrm{KS}$ are the Atiyah class and the Kodaira--Spencer class of the target thickening. The flat-sheaf locus and its stability conditions are open in the families in question.

Equivalences of targets preserve composition in $\RHom_Z(F,F)$, and a scalar acts on this complex by conjugation, which is the identity. This gives the same deformation calculation in every chosen frame. A two-term perfect obstruction theory on the entire rigidified stack requires the additional amplitude and obstruction-morphism verifications in the virtual-class construction.
\end{proof}

\subsection{Algebraicity, automorphisms and the valuative criterion}
\begin{proposition}\label{prop:algebraic}
Under (K1)--(K2), $\widetilde{\MM}(v)$ is algebraic of finite type and has a tautological sheaf on its chart presentations. Its scalar rigidification $\MM(v)$ is Deligne--Mumford if the finite-stabilizer condition in (K4) holds. The endomorphism complex of the tautological sheaf descends through the scalar gerbe.
\end{proposition}
\begin{proof}
The chart construction of Proposition~\ref{prop:generalchartdescent} gives the algebraic presentation before imposing the finiteness of $H_x$; that hypothesis enters only at \eqref{eq:intrinsicaut}. Write $\mathcal F$ for the tautological sheaf on the universal target over $\widetilde{\MM}(v)$. A scalar $\lambda$ acts by $\lambda$ on $\mathcal F$ and $\lambda^{-1}$ on its derived dual. For local transition maps $g_{ij}$ write $g_{ij}g_{jk}=u_{ijk}g_{ik}$ with $u_{ijk}\in\OO^*$. On derived endomorphisms,
\[
 \operatorname{Ad}(g_{ij})\operatorname{Ad}(g_{jk})
  =\operatorname{Ad}(u_{ijk}g_{ik})=\operatorname{Ad}(g_{ik}),
 \qquad \operatorname{Ad}(\lambda\id)=\id .
\] Local trivializations of the scalar gerbe have transition isomorphisms of $\mathcal F$ differing by units, which cancel on this tensor product. The resulting descent data are effective for perfect complexes and identify its pullbacks on each chart.
\end{proof}

\begin{theorem}\label{thm:proper}
Under Hypothesis~\ref{hyp:compact}, the scalar-rigidified stack $\MM(v)\to C$ is separated and proper. The same assertion holds over $\Spec\CC$ for the relative moduli of a pair satisfying the corresponding hypotheses.
\end{theorem}
\begin{proof}
Fix a discrete valuation ring $R$ over $C$, its fraction field $K$, and two objects $x,y\in\MM(v)(R)$. The pullback of the diagonal is
\[
\begin{tikzcd}[column sep=2.8em]
 I=\operatorname{Isom}_R(x,y)\ar[r]\ar[d]&\MM(v)\ar[d,"\Delta_{\MM(v)/C}"]\\
 \Spec R\ar[r,"{(x,y)}"']&\MM(v)\times_C\MM(v).
\end{tikzcd}
\]
The diagonal is representable and of finite type. Indeed, on the bounded target charts the isomorphisms are given by finitely many identifications of surviving parameters and their torus rescalings. Over each such relation, the sheaf isomorphisms form the invertible locus of the relative Hom scheme for two bounded flat families on a projective target. Quotienting by scalar units gives the rigidified Isom spaces.

For finite-type morphisms over the Noetherian bases here the valuative test may be taken over discrete valuation rings. By (K4), $I(R)\xrightarrow{\sim}I(K)$ for each $(x,y)$. Hence $\Delta_{\MM(v)/C}$ is proper. By (K3), for each $\xi_K\in\MM(v)(K)$ there are a finite extension of traits $R\subset R'$ and $\xi_{R'}$ with $\xi_{R'}|_{K'}\simeq \xi_K|_{K'}$. This gives universal closedness. Finite type follows from (K2) by Proposition~\ref{prop:algebraic}. These three properties give properness. The proof over $\Spec\CC$ is identical.
\end{proof}

\subsection{Extension of an isomorphism on a fixed polarized target}
When a target isomorphism over $R$ has already been specified, pullback by it gives two sheaves on the same projective $R$-scheme. Then the remaining isomorphism in (K4) can be extended by the following argument.

\begin{proposition}\label{prop:fixedtargetuniquenessfull}
Let $R$ be a discrete valuation ring with uniformizer $\pi$ and fraction field $K$. On a fixed projective scheme $W\to\Spec R$ with relatively ample line bundle $A$, let $F,G$ be $R$-flat, finitely presented sheaves with stable special fibers of the same reduced Hilbert polynomial. The class of every isomorphism $\alpha:F_K\to G_K$ modulo $K^*$ extends to exactly one class of isomorphisms $F\to G$ modulo $R^*$.
\end{proposition}
\begin{proof}
The finite torsion-free $R$-module $H=\Hom_W(F,G)$ satisfies
$H\otimes_R K=\Hom_{W_K}(F_K,G_K)$ by flat localization. Finiteness follows from properness, and torsion-freeness from $R$-flatness of $G$. Choose $b\in\mathbb Z$ with $\beta=\pi^b\alpha\in H\setminus\pi H$. The exact sequence $0\to G\xrightarrow{\pi}G\to G_0\to0$ gives an injection
\[
 H/\pi H\hookrightarrow\Hom_W(F,G_0)
       =\Hom_{W_0}(F_0,G_0).
\]
Thus $f=\beta_0$ has nonzero image $I$. Put $K_0=\ker f$. If $K_0\ne0$, purity and stability give $p_{K_0}<p_{F_0}$. For $a(-)$ the positive leading Hilbert coefficient, the relation
\[
 a(I)(p_I-p_{F_0})=a(K_0)(p_{F_0}-p_{K_0})
\]
gives $p_I>p_{F_0}$. This contradicts $p_I\le p_{G_0}=p_{F_0}$. Hence $f$ is injective and $p_I=p_{F_0}$. A proper nonzero subsheaf of $G_0$ has strictly smaller reduced polynomial, so $I=G_0$.

Both fibers of $\operatorname{coker}\beta$ vanish, hence $\beta$ is surjective by Nakayama's lemma. Since $G$ is $R$-flat, tensoring its quotient sequence with $R/(\pi)$ identifies $(\ker\beta)_0$ with $\ker f=0$. The generic kernel vanishes as well, so $\ker\beta=0$. If $\gamma$ is a second extension of the same class, $\gamma_K=u\beta_K$ for some $u\in K^*$. Both $\beta$ and $\gamma$ are primitive in the lattice $H$, hence $v_R(u)=0$ and $u\in R^*$. The equality over $R$ follows from injectivity of $H\to H\otimes K$.
\end{proof}

\begin{lemma}\label{lem:openproper}
An open substack of a separated stack is separated. If a proper stack $\mathfrak K\to C$ contains an open substack $\mathfrak U$, then properness of $\mathfrak U\to C$ requires a further valuative existence statement. In particular, if $\mathfrak U$ is stable under every specialization in $\mathfrak K$, it is open and closed and hence proper.
\end{lemma}
\begin{proof}
For the open immersion $j:\mathfrak U\hookrightarrow\mathfrak K$, its relative diagonal is the pullback
\[
 \Delta_{\mathfrak U/C}=(j\times_Cj)^*\Delta_{\mathfrak K/C}.
\]
A proper diagonal remains proper after base change. In the finite-type Noetherian setting $|\mathfrak U|$ is constructible. If it contains all its specializations, it is closed, so $j$ is both open and closed. Then $\mathfrak U\to\mathfrak K\to C$ is proper. In general a separated open substack must still satisfy the existence criterion over valuation rings.
\end{proof}

In particular, applying this test to the open boundary Cohen--Macaulay locus requires checking that its limits satisfy the same depth bound.
\subsection{Numerical labels in families}
On a marked family the branch restrictions are flat sheaves and the boundary is a perfect complex by Theorem~\ref{thm:perfectpurematchingnew}. Its polynomial means $m\mapsto\chi(S,E(m))$, computed by the alternating sum over all its cohomology sheaves. For each twist, perfect proper pushforward and base change give a perfect complex on the parameter space.

Its Euler characteristic is the alternating sum of the ranks in a local finite free representative, hence locally constant. Evaluating at three distinct twists determines all coefficients of each polynomial. Under smoothing of unmarked nodes, the polynomial on a surviving branch is the sum of its component polynomials minus its internal boundary polynomials, by normalization.

This is exactly the continuous weighting convention for expansions in Li--Wu \cite[Section 2.5]{LW}.

\begin{proposition}\label{prop:finitecharges}
If $\MM$ satisfies Hypothesis~\ref{hyp:compact}, then its node-marked space has a finite open-and-closed decomposition by the Hilbert polynomials of the two branches and the common boundary, computed with compatible twisting line bundles on all marked charts. For example, one may pull back a fixed ample line bundle from the contracted target to compute the contracted numerical polynomials.

Any additional discrete numerical class which is locally constant in these flat families refines the decomposition into finitely many nonempty pieces.
\end{proposition}
\begin{proof}
Write $U_k=\AA^{k+1}$ for an expansion chart,
$U_{k,0}=V(t_0\cdots t_k)$ and $H_j=V(t_j)\subset U_k$.
Choosing one unsmoothed joining divisor on a fiber of this chart
gives the morphism
\[
 \coprod_{j=0}^{k}H_j\longrightarrow U_{k,0}.
\]
On coordinate rings it is the homomorphism
\[
 R=k_0[t_0,\ldots,t_k]/(t_0\cdots t_k)
       \longrightarrow\prod_{j=0}^{k}R/(t_j),
       \qquad f\longmapsto(f\bmod t_j)_j ,
\]
where $k_0=\CC$. Its target is a finite $R$-module. Inverting a
smoothing parameter deletes precisely its empty marked summand;
a rescaling preserves each $H_j$. Thus this construction descends
over the stack $\EE_0$ of special expanded targets. If
$\EE_0^\dagger$ denotes targets with a chosen joining divisor, the
marked sheaf stack is
$\MM_0^\dagger=\MM_0\times_{\EE_0}\EE_0^\dagger$.

Finiteness is
local on the target for the smooth topology and preserved by base
change, so $\MM_0^\dagger\to\MM_0$ is finite. Hypothesis
\ref{hyp:compact}(K2) gives quasi-compactness of $\MM_0$ and
therefore of $\MM_0^\dagger$.

Let $F_L,F_R,E$ be the two restrictions and their common boundary
complex on a marked chart over $T$. The boundary value is the
derived restriction. For each $V=L,R,D$, denote the corresponding
perfect object by $F_V$ (so $F_D=E$), the projection by $p_V$,
and the chosen twist by $A_V$. The chartwise complexes
\[
 K_{V,m}=\RR(p_V)_*(F_V\otimes A_V^m)\in\perf(T)
\]
commute with base change. On an open where $K_{V,m}$ is represented
by finite free terms of ranks $b^j$, the value
$P_V(m)=\sum_j(-1)^jb^j$ is constant. These integers agree on
overlaps because the twisting line bundles and restriction complexes
do. The alternating Hilbert polynomials have degree at most two.
Their coefficients are determined by their values at $m=0,1,2$.
Consequently
\[
 |\MM_0^\dagger|\longrightarrow\mathbb Z^9,\qquad
 x\longmapsto\big(P_V^x(m)\big)_{V=L,R,D;\ m=0,1,2}
\]
is locally constant. Each fiber is open and closed. The disjoint
fibers give an open cover of a quasi-compact space, so only finitely
many are nonempty. Appending another locally constant discrete
coordinate gives the same conclusion for every stated refinement.
\end{proof}

\subsection{Recovering the isomorphism of expansions from the support}
One case of (K4) can be tested using only the intrinsic support of the sheaf. Over a discrete valuation ring $R$ write
\[
       D(F)=V(\operatorname{Ann}_{\OO_Z}F)\subset Z
\]
for the scheme-theoretic support of a flat family on an expansion. Say that its support gives a stable expanded Hilbert point if the natural quotient $\OO_Z\twoheadrightarrow\OO_{D(F)}$ is a family of stable admissible quotients with the prescribed Hilbert polynomial in \cite[Theorems~4.14--4.15]{LW}. For this quotient, stability means finiteness of the automorphisms of the expansion preserving the subscheme. We test admissibility and this finiteness for the support itself.

\begin{proposition}[Uniqueness on expansions determined by the support]\label{prop:support-separated27}
Let $(Z,F)$ and $(Z',F')$ be two represented flat families over an excellent discrete valuation ring $R$ as in the valuative test for either the degeneration or a relative pair. Assume their scheme-theoretic supports give stable expanded Hilbert points and commute with restriction to the closed fiber.

Assume that the stability loci are invariant under isomorphisms of expanded targets as in (K1), and that on a common chart the chosen relatively ample line bundle on the target family makes both closed sheaves Gieseker-stable. For every isomorphism of their generic fibers as expanded sheaves there is a unique isomorphism over $R$ extending its class after scalar rigidification.

The residual geometric stabilizers on this locus are finite.
\end{proposition}
\begin{proof}
Put $K=\operatorname{Frac}(R)$, and choose a uniformizer $\pi$. On an affine chart set $J=\operatorname{Ann}(F)$. Flatness of $F$ over $R$ and finite generation give
\[
 (J:\pi)=J,\qquad J[1/\pi]=\operatorname{Ann}(F_K).
\]
Thus $\OO_Z/J$ is $R$-flat, and $D(F)$ is the schematic closure of $D(F_K)$. The closed-fiber identity $D(F)_0=D(F_0)$ is part of the hypothesis.

The generic arrow $(u_K,\phi_K)$ identifies the support quotients. An isomorphism of the canonical quotients $\OO\twoheadrightarrow\OO_D$ is determined by its underlying target isomorphism. The proper diagonal of the expanded Hilbert stack \cite[Theorems~4.14--4.15 and Proposition~5.2]{LW} therefore gives a unique extension
\[
              u:(Z,D(F))\xrightarrow{\sim}(Z',D(F')) .
\]
On a common polarized chart, flatness and $\phi_K$ give $P_{F_0}=P_{u^*F'_0}$. The closed sheaves are stable. Proposition~\ref{prop:fixedtargetuniquenessfull} now extends $[\phi_K]$ uniquely to $[F\xrightarrow{\sim}u^*F']$. Uniqueness identifies the classes on overlaps, giving descent. Any second extension first agrees with $u$, then with $[\phi]$. Finally \eqref{eq:intrinsicaut} embeds the rigidified stabilizer into the finite stabilizer of the support quotient.
\end{proof}

\section{Admissible sheaves and depth at a gluing divisor}\label{sec:local}
The Tor and depth calculations determine when a flat sheaf is perfect on an expanded target $Z$ and when its ordinary restriction to a gluing divisor is pure. They supply local hypotheses for the deformation complex $R\mathcal Hom_Z(F,F)$.

\begin{definition}\label{def:normal}
A coherent sheaf $F$ on a scheme $Z$ is normal to a closed subscheme $D$ if
\[
\Tor_1^{\OO_Z}(F,\OO_D)=0.
\]
Let $W=Y\cup\Delta_1\cup\cdots\cup\Delta_m$ be an expansion of a pair $(Y,D)$.
Write $D_a^-,D_a^+$ for the inner and outer sections of $\Delta_a$.
For $m>0$ the attaching divisors are the identifications
$D=D^-_1$ and $D^+_a=D^-_{a+1}$ for $1\le a<m$.
The \emph{outer distinguished divisor} $D_{\rm out}$ is $D_m^+$ for
$m>0$ and $D\subset Y$ for $m=0$. It is Cartier on $W$.
If $\mathcal J$ denotes the attaching divisors, admissibility means
\[
 \Tor_1^{\OO_W}(F,\OO_{D'})=0\qquad(D'\in\mathcal J\cup\{D_{\rm out}\}).
\]
On a complete expansion of $Y_1\cup_SY_2$ the test runs over its
attaching divisors. In flat families apply it on each geometric fiber.
\end{definition}

\begin{lemma}\label{lem:flatnode}
Set
\[
A_0=\CC[[x,y]]/(xy),\qquad B=A_0[[u,v]],\qquad R=B/(x,y).
\]
For a finite $B$-module $M$, normality to $\Spec R$ is equivalent to flatness over $A_0$. Under these conditions
\[
M\otimes_B^{\LL}R\simeq M\otimes_B R.
\]
\end{lemma}
\begin{proof}
The ring $B$ is flat over the local ring $A_0$, and $R=B\otimes_{A_0}\CC$. Derived tensor associativity gives
\[
M\otimes_B^{\LL}R\simeq M\otimes_{A_0}^{\LL}\CC.
\]
Thus normality is the vanishing of $\Tor_1^{A_0}(M,\CC)$. The local flatness criterion for a finite module over the noetherian local $A_0$-algebra $B$ identifies this condition with $A_0$-flatness. This is the node-normality criterion of Li--Wu \cite[Lemma 3.2]{LW}; it also follows from the local criterion after completion. Flatness gives all the higher Tor vanishings and hence the final assertion.
\end{proof}

\begin{lemma}\label{lem:pdchange}
Let $(A,\mathfrak m)$ be noetherian local, let $B$ be a noetherian local flat $A$-algebra, and let $M$ be a finite $B$-module flat over $A$. Put $R=B/\mathfrak mB$ and $G=M\otimes_A A/\mathfrak m$. If $M\ne0$, then
\[
\pd_BM=\pd_RG
\]
whenever the right side is finite.
\end{lemma}
\begin{proof}
Let $\mathfrak n$ be the maximal ideal of $B$ and let
$k_B=B/\mathfrak n=R/(\mathfrak n/\mathfrak mB)$. Take a
minimal free resolution $P_\bullet\to M$ over $B$. Set $K_0=M$
and use the exact sequences
\[
 0\longrightarrow K_{i+1}\longrightarrow
 P_i=B^{b_i}\longrightarrow K_i\longrightarrow0.
\]
Inductively $K_i$ is $A$-flat: for every $A$-module $N$ the
Tor sequence identifies $\Tor^A_1(K_{i+1},N)$ with
$\Tor^A_2(K_i,N)=0$, since $P_i$ is $A$-flat. Tensoring the
displayed rows with $A/\mathfrak m$ therefore preserves exactness.
Thus $P_\bullet\otimes_B R$ is a free $R$-resolution of $G$.
Minimality means $d_i(P_i)\subseteq\mathfrak n P_{i-1}$ for $i>0$,
so all differentials vanish after the further reduction to $k_B$.
Consequently
\[
 \Tor_i^B(M,k_B)\cong k_B^{b_i}
      \cong\Tor_i^R(G,k_B)\qquad(i\ge0).
\]
Over a noetherian local ring the projective dimension of a nonzero
finite module is the supremum of indices at which its minimal
resolution has positive rank. The two suprema are equal.
\end{proof}

\begin{theorem}\label{thm:nodeperfect}
For the rings of Lemma~\ref{lem:flatnode}, let $M$ be a nonzero finite module normal to $\Spec R$. Then
\begin{equation}\label{eq:pdnode}
\pd_BM=\pd_R(M\otimes_BR)\le2.
\end{equation}
If $M$ is pure of dimension two and its nonzero boundary restriction $G=M\otimes_BR$ has dimension one, the following are equivalent:
\begin{enumerate}[label=(\roman*)]
\item $G$ is pure of dimension one;
\item $\pd_BM=1$;
\item $M$ is Cohen--Macaulay of dimension two.
\end{enumerate}
\end{theorem}
\begin{proof}
Lemmas~\ref{lem:flatnode} and~\ref{lem:pdchange} give $\pd_BM=\pd_RG$ for $G=M\otimes_BR$. Since $R=\CC[[u,v]]$ is regular, $\pd_RG\le2$, proving \eqref{eq:pdnode}. For $G\ne0$ of dimension one, purity is equivalent to $\depth_R G=1$. Auslander--Buchsbaum on the two local rings gives
\[
       \pd_B M=\pd_R G=2-\depth_R G,\qquad
       \depth_B M=3-\pd_B M.
\]
Thus $\depth_R G=1\iff \pd_B M=1\iff \depth_B M=2$. The last equality is the Cohen--Macaulay condition for $\dim M=2$. If $G=0$, Nakayama gives $M=0$ at the node.
\end{proof}

\begin{remark}\label{rem:purityhyp}
Normality and purity of the surface sheaf imply perfectness by \eqref{eq:pdnode}; purity of the boundary is a separate depth condition. For a perfect sheaf without admissibility, the comparison is the full triangle \eqref{eq:normaltriangleallnew}. The compactified-Jacobian formula additionally uses the equivalent conditions of Theorem~\ref{thm:nodeperfect} at the marked divisor. These conditions permit nontransverse, singular, and nonreduced boundary curves.
\end{remark}

\subsection{Nontransverse supports and the meaning of bubbling}\label{sec:transversality}
The local conditions in this section are precisely those used by the matching functor and by the compactification of Section~\ref{sec:verifiedcompactification}.  They concern the sheaf as a module.  They impose neither smoothness of its supporting surface nor geometric transversality of its intersection with the K3.

\begin{proposition}\label{prop:tangentexample}
There are admissible surface sheaves whose two branch supports are smooth, tangent to the double divisor, and meet it in a singular integral curve.  These objects have a length-one resolution at the node and pure boundary.
\end{proposition}
\begin{proof}
Put $B=\CC[[x,y,u,v]]/(xy)$, $R=B/(x,y)$, $g=u^2+v^3$ and $M=B/(x+y+g)$. The element $x+y+g$ avoids the two minimal primes of the Cohen--Macaulay ring $B$, so
\[
 0\to B\xrightarrow{x+y+g}B\to M\to0,\qquad
 M\otimes_B^{\LL}R=[R\xrightarrow{g}R]\simeq R/(g).
\]
Thus $\Tor^B_{>0}(M,R)=0$. The one-dimensional hypersurface domain $R/(g)$ is Cohen--Macaulay. On the two branches the supporting equations are $x+g$ and $y+g$. Their differentials at the origin are $dx$ and $dy$, respectively, equal to those of the joining divisors. The two supports are smooth and tangent there; their common restriction is the cuspidal curve above.
\end{proof}

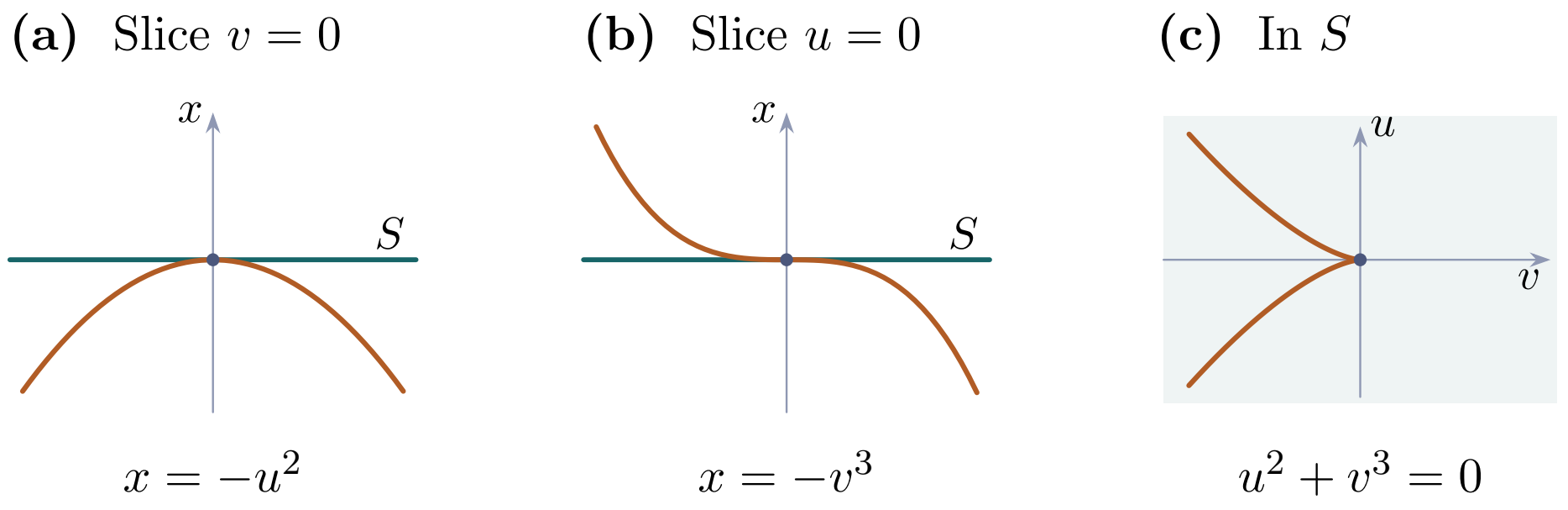
\begin{figure}[htbp]
 \centering\begingroup
 \definecolor{brS}{HTML}{176568}
 \definecolor{brT}{HTML}{B15C25}
 \definecolor{brI}{HTML}{53618A}
 \begin{tikzpicture}[x=1cm,y=1cm,line cap=round,line join=round,
   axis/.style={draw=brI!65,line width=.45pt,-{Stealth[length=1.65mm,width=1.15mm]}},
   zaxis/.style={draw=brS,line width=1.10pt},
   surf/.style={draw=brT,line width=1.15pt},
   shortlabel/.style={font=\footnotesize,inner sep=2pt}]
   \begin{scope}
    \node[font=\small,anchor=west] at (-1.75,1.75) {\textbf{(a)}\ \ Slice $v=0$};
    \draw[axis] (0,-1.20)--(0,1.16) node[shortlabel,anchor=east] {$x$};
    \draw[zaxis] (-1.60,0)--(1.60,0) node[shortlabel,above left] {$S$};
    \draw[surf] plot[domain=-1.5:1.5,samples=50] (\x,{-.46*\x*\x});
    \fill[brI!90!black] (0,0) circle (1.45pt);
    \node[font=\small] at (0,-1.68) {$x=-u^2$};
   \end{scope}
   \begin{scope}[shift={(4.52,0)}]
    \node[font=\small,anchor=west] at (-1.75,1.75) {\textbf{(b)}\ \ Slice $u=0$};
    \draw[axis] (0,-1.20)--(0,1.16) node[shortlabel,anchor=east] {$x$};
    \draw[zaxis] (-1.60,0)--(1.60,0) node[shortlabel,above left] {$S$};
    \draw[surf] plot[domain=-1.5:1.5,samples=50] (\x,{-.31*\x*\x*\x});
    \fill[brI!90!black] (0,0) circle (1.45pt);
    \node[font=\small] at (0,-1.68) {$x=-v^3$};
   \end{scope}
   \begin{scope}[shift={(9.04,0)}]
    \node[font=\small,anchor=west] at (-1.75,1.75) {\textbf{(c)}\ \ In $S$};
    \path[fill=brS!7] (-1.55,-1.13) rectangle (1.55,1.13);
    \draw[axis] (-1.55,0)--(1.50,0) node[shortlabel,below left] {$v$};
    \draw[axis] (0,-1.08)--(0,1.05) node[shortlabel,right] {$u$};
    \draw[surf] plot[domain=-1:1,samples=71] ({-1.35*\x*\x},{.99*\x*\x*\x});
    \fill[brI!90!black] (0,0) circle (1.45pt);
    \node[font=\small] at (0,-1.68) {$u^2+v^3=0$};
   \end{scope}
 \end{tikzpicture}\endgroup
 \caption{Real coordinate slices of the local equations in
 Proposition~\ref{prop:tangentexample}. On the branch with coordinate $x$,
 the smooth support has equation $x=-u^2-v^3$ and the joining K3 has
 equation $x=0$. Their tangent planes at the origin agree. Their
 scheme-theoretic intersection is the integral cusp in (c), with
 parametrization $(v,u)=(-t^2,t^3)$.}
 \label{fig:cuspidalpropercontact}
\end{figure}

The two plane slices in Figure~\ref{fig:cuspidalpropercontact} recover the tangent plane of the support at the cusp.

\begin{proposition}\label{prop:bubblemovesdefect}
In a local family $xy=t$, a point defect approaching the double divisor can be carried into the interior of an inserted bubble.  This operation can restore pure sheaf restriction without making a singular intersection curve geometrically transverse.
\end{proposition}
\begin{proof}
For the explicit transport take the supporting surface $u=0$, make the base change $t=r^2$, and consider the section
\[
 x=y=r,\qquad u=v=0.
\]
Blow up $(x,r)$.  In the chart $x=rz$, the equation is $zy=r$, and the section has $z=1$, $y=r$, $u=v=0$.  At $r=0$ its point is in the bubble interior, away from both distinguished sections.  The ideal of this section in the strict supporting surface is flat over $\CC[[r]]$, by the exact sequence with its flat surface structure sheaf and its finite flat quotient.  Near either gluing divisor that ideal is the surface structure sheaf.

Its restriction is therefore pure and has no positive Tor.  The interior point has finite geometric rescaling stabilizer, in fact trivial in this model.

The same blowup works with the equation $g(u,v)=u^2+v^3=0$ in
place of $u=0$. The point lies at $z=1$, away from both sections.
At either section the point ideal equals the structure sheaf of the
surface, and the boundary ring is $\CC[[u,v]]/(u^2+v^3)$. This
one-dimensional domain has depth one and embedding dimension two.
Its structure sheaf is pure and the boundary curve is singular.
\end{proof}

\begin{proposition}\label{prop:perfectfiber}
An admissible pure two-dimensional sheaf on a standard expanded threefold is perfect, with local resolution length at most two. Near a node at which its boundary is nonzero and pure one-dimensional, the resolution has length one.
\end{proposition}
\begin{proof}
At an internal node Theorem~\ref{thm:nodeperfect} gives $\pd F\le2$, with equality replaced by $\pd F=1$ when the nonzero restriction has pure dimension one. At a closed point $p$ of the smooth locus, purity gives $\mathfrak m_p\notin\operatorname{Ass}(F_p)$ and hence
\[
           \pd_{\OO_{Z,p}}F_p=3-\depth_{\OO_{Z,p}}F_p\le2
           \qquad(F_p\ne0).
\]
Localization gives the same bound throughout the support.
\end{proof}

\begin{proposition}\label{prop:familyperfect}
Let $W\to T$ be a flat finitely presented family of expanded threefolds over a noetherian base, and let $\FF$ be a finitely presented $T$-flat family of admissible pure surface sheaves. Then $\FF$ is perfect on $W$. Locally it has a free resolution of length at most two. For any marked node, its derived restrictions to the two branches and to the gluing divisor are ordinary $T$-flat sheaves and commute with arbitrary base change.
\end{proposition}
\begin{proof}
Fix $w\in W$ over $t\in T$, and write $A=\OO_{T,t}$,
$B=\OO_{W,w}$, $M=\FF_w$, $R=B\otimes_Ak(t)$. After extension
to the geometric residue field, Proposition~\ref{prop:perfectfiber}
bounds the projective dimension of the fiber by two. Faithfully flat
field extension preserves the vanishing of the corresponding Tor
groups, so $\pd_R(M\otimes_A k(t))\le2$ if this module is nonzero.
Lemma~\ref{lem:pdchange} now gives a resolution
\[
 0\longrightarrow P_2\longrightarrow P_1\longrightarrow P_0
     \longrightarrow M\longrightarrow0
\]
by finite free $B$-modules, allowing $P_2=0$. If the fiber module
is zero, Nakayama gives $M=0$. The finitely presented resolution
spreads to a neighborhood of $w$.

For the marked joining divisor or either incident branch, denote the
closed immersion by $j:V\hookrightarrow W$. These families are
flat over the base of the marking. The complex
$L=j^*P_\bullet$ consists locally of three $A$-flat terms in
degrees $-2,-1,0$. At a geometric point $\bar t$ its cohomology
in degree $-i$, $i>0$, is
\[
 \Tor_i^{\OO_{W_{\bar t}}}
             (\FF_{\bar t},\OO_{V_{\bar t}}).
\]
At the joint, complete the local ring as
$B_0=A_0[[u,v]]$, $A_0=k(\bar t)[[x,y]]/(xy)$.
The branch and boundary rings are $B_0\otimes_{A_0}D$ for
$D=A_0/(x)$, $A_0/(y)$ and $k(\bar t)$, respectively. As
$B_0$ and the admissible module $M_0$ are $A_0$-flat by
Lemma~\ref{lem:flatnode}, associativity gives
\[
 M_0\otimes^{\LL}_{B_0}(B_0\otimes_{A_0}D)
      \simeq M_0\otimes^{\LL}_{A_0}D=M_0\otimes_{A_0}D .
\]
Away from the joint a branch restricts by the identity or the zero
map of spaces; at a smooth relative divisor higher Tor vanishes
by the two-term resolution and injectivity of its equation. Thus
$H^{-2}(L_{\bar t})=H^{-1}(L_{\bar t})=0$ at every point.

The local criterion for flatness applied first to
$L^{-2}\to L^{-1}$ and then to
$\operatorname{coker}(L^{-2}\to L^{-1})\to L^0$ makes both maps
injective with $A$-flat cokernels: in each case the source and target
are $A$-flat and the map is injective on the fiber. Hence
$\LL j^*\FF\simeq H^0(L)$ is $T$-flat. For any $A\to A'$,
\[
 (\LL j^*\FF)\otimes_A^{\LL}A'
     =H^0(L)\otimes_A A'
     \simeq\LL j_{A'}^*(\FF_{A'}),
\]
by associativity and flatness of the target families. This gives
the stated base-change isomorphism.
\end{proof}

\begin{corollary}\label{cor:boundaryopen}
On a bounded chart of admissible flat surface sheaves, the condition that the restriction at a specified nonzero one-dimensional boundary be pure is equivalent to the Cohen--Macaulay condition on the surface sheaf along that divisor. It is open. The same is true simultaneously at all the gluing divisors on a bounded expansion chart.
\end{corollary}
\begin{proof}
The equivalence is fiberwise Theorem~\ref{thm:nodeperfect}. The Cohen--Macaulay locus of a flat finitely presented coherent sheaf on a finite-type family is open. On a chart over $T$, write $N$ for the closed complement of the fiberwise Cohen--Macaulay locus, and $p_e:D_e\to T$ for the projective joining divisors. The required open on the chart is
\[
 T^{\rm pure}=T\setminus\bigcup_{e}p_e(N\cap D_e).
\]
Each image is closed by properness of $p_e$, and the index set is finite. Equivalently, one can use openness of purity for the flat boundary family in Proposition~\ref{prop:familyperfect}.
\end{proof}

\section{Gluing coherent sheaves and flat families}\label{sec:classical}
Theorem~\ref{thm:perfectpurematchingnew} treats arbitrary derived restrictions on the perfect locus. This section specializes that construction to an ordinary boundary sheaf and proves explicit formulas for coherent families.

\begin{theorem}\label{thm:classicalnode}
Let $R=\CC[[u,v]]$ and set
\[
B=R[[x,y]]/(xy),\qquad B_x=R[[x]],\qquad B_y=R[[y]].
\]
Finite $B$-modules normal to $R$ are equivalent to triples
\[
(M_x,M_y,\phi),\qquad
\phi:M_x/xM_x\xrightarrow{\sim}M_y/yM_y,
\]
where multiplication by $x$ on $M_x$ and by $y$ on $M_y$ is injective. The inverse functor is
\begin{equation}\label{eq:kernelglue}
M=\ker\left(M_x\oplus M_y\longrightarrow G\right),\qquad
(a,b)\longmapsto\overline a-\phi^{-1}(\overline b),
\end{equation}
with $G=M_x/xM_x$. The equivalence preserves purity of dimension two on the nonzero branches.
\end{theorem}
\begin{proof}
For a normal $M$, Lemma~\ref{lem:flatnode} gives flatness over $A_0$. Tensoring
\[
0\longrightarrow A_0\longrightarrow
\CC[[x]]\oplus\CC[[y]]\longrightarrow\CC\longrightarrow0
\]
with $M$ gives the matching exact sequence. Flat base change to either branch gives injectivity of its normal parameter.

Conversely, construct $M$ by \eqref{eq:kernelglue}. Both projections are surjective. The kernel of the projection to $M_x$ consists of pairs $(0,yb)$ and equals $yM$: choose a lift of the boundary value of $b$ on $M_x$ before multiplying the pair by $y$. Similarly the other kernel is $xM$. Thus
\[
M/yM=M_x,\qquad M/xM=M_y,\qquad M/(x,y)M=G.
\]
To check normality directly, the beginning of the free resolution of $\CC$ over $A_0$ is
\[
A_0^2\xrightarrow{\operatorname{diag}(y,x)}A_0^2
\xrightarrow{(x\ \ y)}A_0\longrightarrow\CC\longrightarrow0.
\]
If $xa+yb=0$ for $a,b\in M$, its image on $M_x$ implies $a\in yM$ and its image on $M_y$ implies $b\in xM$. Hence $(a,b)$ lies in the image of $\operatorname{diag}(y,x)$, proving $\Tor_1^{A_0}(M,\CC)=0$. Lemma~\ref{lem:flatnode} proves normality.

If the branch modules are pure, $M$ embeds in their direct sum and is pure. Conversely, a lower-dimensional submodule $N\subset M_x$ would give the nonzero submodule $(xN,0)\subset M$, of the same lower dimension, since $x$ is injective. This contradicts purity of $M$. The other branch is identical. The restrictions and kernels just constructed are mutually inverse on morphisms as well as on objects, proving an equivalence of categories.
\end{proof}

\begin{proposition}\label{prop:familyglue}
Let $W=W_L\cup_DW_R$ be a fixed node-marked expanded family over $T$. Flat branch sheaves $F_L,F_R$, normal to $D$, with a $T$-flat common boundary sheaf $G$ and an isomorphism of restrictions glue uniquely by the kernel construction. The kernel is $T$-flat, commutes with every base change, and is normal to the marked divisor.

If the branch families are admissible at their other distinguished divisors, the glued family is admissible. Restriction gives the inverse equivalence of fibered categories. The statement iterates on a chain.
\end{proposition}
\begin{proof}
After identifying the two restrictions with $G$, take the difference of the quotient maps. This gives
\begin{equation}\label{eq:familyglue}
 0\longrightarrow F\longrightarrow i_{L*}F_L\oplus i_{R*}F_R
   \xrightarrow{\mathrm{res}_L-\mathrm{res}_R}j_*G\longrightarrow0.
\end{equation}
Write $A$ for the middle term and $B$ for the last term. For a module $N$ on $T$, the two outer terms of
\[
 \Tor_2^T(B,N)\longrightarrow\Tor_1^T(F,N)
                  \longrightarrow\Tor_1^T(A,N)
\]
vanish. Thus $F$ is $T$-flat; $\Tor_1^T(B,N)=0$ gives exactness of \eqref{eq:familyglue} after every base change. Theorem~\ref{thm:classicalnode} gives normality at the marked divisor on each fiber. At any other distinguished divisor the gluing map is locally an isomorphism from the corresponding branch, so the original admissibility test is preserved there.

Conversely, for an admissible family the branch and boundary restrictions are flat by Proposition~\ref{prop:familyperfect}. The restriction and kernel constructions commute with arrows. On a chain their common inverse is the kernel of the sum of the signed restriction maps, since consecutive joining divisors are disjoint.
\end{proof}

\begin{proposition}\label{prop:Endtriangle}
For an admissible sheaf $F$ on an expanded chain $Z$, with branch sheaves $F_i$ and boundary sheaves $G_e$, there is a canonical triangle
\begin{equation}\label{eq:Endtriangle}
\RHom_Z(F,F)\longrightarrow
\bigoplus_i\RHom_{Z_i}(F_i,F_i)
\longrightarrow
\bigoplus_e\RHom_{S_e}(G_e,G_e)\longrightarrow.
\end{equation}
It holds in flat families and is equivariant under all isomorphisms of expanded targets.
\end{proposition}
\begin{proof}
With $D_F=F^\vee\otimes^{\LL}F$, tensor the normalization sequence
\[
 0\longrightarrow\OO_Z\longrightarrow\bigoplus_i(i_i)_*\OO_{Z_i}
       \longrightarrow\bigoplus_e(j_e)_*\OO_{S_e}\longrightarrow0
\]
by $D_F$ in the derived category. By Propositions~\ref{prop:perfectfiber} and~\ref{prop:familyperfect}, $F$ is perfect and $\LL i_i^*F=F_i$, $\LL j_e^*F=G_e$. Thus the projection formula identifies $\LL i_i^*D_F=F_i^\vee\otimes^{\LL}F_i$ and $\LL j_e^*D_F=G_e^\vee\otimes^{\LL}G_e$. Taking $\RG$ gives \eqref{eq:Endtriangle}, whose second map is the signed difference of restrictions. All operations commute with isomorphisms of expanded targets and with derived base change in the stated flat families.
\end{proof}

\subsection{The chain maps and the square-zero diagram chase}\label{sec:chainmapsfull}
We record the comparison at the level of complexes.  This gives an explicit check of the maps used later in the cotangent diagram.

On an affine nodal chart choose a bounded free resolution $P^\bullet\to F$.  It exists by the perfectness theorem.  For a coefficient module $I$ from a square-zero base extension, put
\[
 \mathcal E^n(P;I)=\bigoplus_j\Hom(P^j,P^{j+n}\otimes I),\qquad
 d_{\mathcal E}(a)=d_Pa-(-1)^n a d_P.
\]
Let $P_L,P_R,P_D$ be its restrictions to the two branches and the divisor.  Tensoring the normalization sequence with each finite free matrix module gives an exact sequence in every degree:
\begin{equation}\label{eq:Endshortfull}
0\to\mathcal E^\bullet(P;I)
\to\mathcal E^\bullet(P_L;I)\oplus\mathcal E^\bullet(P_R;I)
\xrightarrow{\mathrm{res}_L-\mathrm{res}_R}
\mathcal E^\bullet(P_D;I)\to0.
\end{equation}
All three arrows commute with $d_{\mathcal E}$, since restriction commutes with the resolution differential.  After a compatible affine cover and totalization of its Cech complexes, this is the global endomorphism triangle, with the connecting map induced by the two restrictions.

For complexes $K_L,K_R,K_D$ representing the three global endomorphism complexes, a model for its homotopy fiber is
\begin{equation}\label{eq:coneexplicitfull}
 K^n=K_L^n\oplus K_R^n\oplus K_D^{n-1},\qquad
 d_K(a_L,a_R,h)=
 \bigl(d a_L,d a_R,\mathrm{res}_La_L-\mathrm{res}_Ra_R-dh\bigr).
\end{equation}
A direct calculation gives $d_K^2=0$: in the third component the two restrictions of $d a_i$ cancel the differential of the two restrictions of $a_i$, and $d^2h=0$.

Over $k[\epsilon]/(\epsilon^2)$ a deformation of a free complex differential is $d_P+\epsilon a$, with $a$ of degree one.  Squaring gives
\[
 (d_P+\epsilon a)^2=\epsilon(d_Pa+a d_P).
\]
It is a differential exactly when $d_{\mathcal E}a=0$.  A change of basis $1+\epsilon b$, with $b$ of degree zero, changes $a$ by a coboundary.  Thus the degree-one cohomology represents first-order deformations, and degree-zero cohomology represents infinitesimal automorphisms.

For matched branch deformations, let $a_L,a_R$ be degree-one cocycles.  An infinitesimal change $1+\epsilon h$ of the boundary identification is a chain isomorphism precisely when
\begin{equation}\label{eq:matchingfirstorderfull}
 \mathrm{res}_La_L-\mathrm{res}_Ra_R=d h.
\end{equation}
Equations $d a_L=d a_R=0$ and \eqref{eq:matchingfirstorderfull} are exactly the cocycle equation in degree one of \eqref{eq:coneexplicitfull}.  Gauge changes are exactly its degree-zero coboundaries.  This proves, on the displayed model, that deformation of the boundary identification is retained; equality of support curves would omit the term $h$.

For a small extension $A'\twoheadrightarrow A$ with square-zero kernel $I$, choose lifts $\widetilde d_L,\widetilde d_R$ and a lift $\Phi$ of the boundary identification, in a common boundary model. Set
\[
 o_i=\widetilde d_i^2,\qquad m=\Phi\widetilde d_L-\widetilde d_R\Phi .
\]
Here and below the differentials on the boundary are restricted, and the reduction of $\Phi$ identifies the two copies of the boundary. Then $o_i$ has degree two and $m$ degree one, with coefficients in $I$. Direct multiplication gives
\[
 [\widetilde d_i,o_i]=0,\qquad
 \widetilde d_Rm+m\widetilde d_L=\Phi o_L-o_R\Phi .
\]
Since $I^2=0$, the latter equality in the boundary complex reads $dm=\mathrm{res}_Lo_L-\mathrm{res}_Ro_R$. Hence $(o_L,o_R,m)$ is a degree-two cocycle of \eqref{eq:coneexplicitfull}. Changing the lifts by degree-one data changes it by a coboundary. Vanishing of its class is equivalent to simultaneous correction of the two differentials and the gluing map.

On a compatible affine cover use the total differential
\[
 D=\delta_{\rm Cech}+(-1)^p d_{\mathcal E}
 \quad\text{in \v Cech degree }p.
\]
In bidegree $(p,q)$ the total degree is $p+q$. Thus corrections use $(0,1),(1,0),(2,-1),\ldots$, whereas their obstructions use $(0,2),(1,1),(2,0),\ldots$. Restriction commutes with both differentials. Totalization gives the same class as the universal square-zero lifting problem under derived matching, and hence the Atiyah-class obstruction map.

\subsection{Units, trace, and the compatible cotangent square}
For a finite free complex the trace is the supertrace
\[
 \operatorname{str}(a)=\sum_j(-1)^j\operatorname{tr}(a^j)
\]
on the degree-zero diagonal matrices, extended as a map of complexes.  The identity
$\operatorname{str}([d_P,b])=0$ proves that it is a chain map.  Restriction commutes with supertrace.  The unit map is $1\mapsto\id_{P^\bullet}$, and
\[
 \operatorname{str}(\id_{P^\bullet})
 =\sum_j(-1)^j\operatorname{rank}P^j=\operatorname{rank}F=0.
\]
This is why no division by the rank can remove scalars in the present problem.

Put $H_V=R\Gamma(\OO_V)$. Normalization gives commuting triangles for endomorphisms and structure sheaves:
\[
\begin{tikzcd}[column sep=large,row sep=large]
 K_Z\arrow[r]\arrow[d,"\operatorname{tr}"']
 &K_L\oplus K_R\arrow[r]\arrow[d,"\operatorname{tr}"]
 &K_D\arrow[d,"\operatorname{tr}"]\arrow[r]&K_Z[1]\arrow[d,"\operatorname{tr}"]\\
 H_Z\arrow[r]
 &H_L\oplus H_R\arrow[r]
 &H_D\arrow[r]&H_Z[1].
\end{tikzcd}
\]
The units give a second commuting diagram.  Removing the common scalar group uses the exact incidence sequence
\[
 0\to k\xrightarrow{a\mapsto(a,a)}k^2
 \xrightarrow{(a,b)\mapsto a-b}k\to0.
\]
Taking the indicated cofibers of units, then the fibers of the determinant-trace maps, gives the reduced endomorphism triangle whenever the compatible determinant construction is used.  This is a three-by-three argument in the stable category: the vertical fibers of a morphism of distinguished triangles form another distinguished triangle.

On a stable pure K3 boundary the reduced boundary complex is $T_J[-1]$.  Dualizing the reduced triangle and shifting gives the obstruction-theory triangle used in the diagonal comparison.  Its maps to cotangent complexes form
\begin{equation}\label{diag:explicitPOTfull}
\begin{tikzcd}[column sep=large,row sep=large]
 \LL j^*(E_1\boxplus E_2)\arrow[r]\arrow[d]
 &E_G\arrow[r]\arrow[d]&e^*\Omega_J[1]\arrow[d]\\
 \LL j^*L_{M_1\times M_2}\arrow[r]
 &L_G\arrow[r]&L_{G/(M_1\times M_2)}.
\end{tikzcd}
\end{equation}
The upper connecting morphism is dual to the signed restriction map in \eqref{eq:coneexplicitfull}.  The right vertical arrow is induced by the Cartesian square with the diagonal of $J$.  Since the square is an actual derived Cartesian square after the stated moduli identifications, these maps commute, including the connecting arrows.  The explicit cocycle calculation above checks that both rows give the same obstruction for every small extension.

Functoriality of virtual pullback applies to this compatible square, not merely to the equality of its virtual ranks.

\subsection{Intrinsic stabilizers under a stable-boundary gluing}\label{sec:stabilizergluefull}
There is one automorphism check which can be completed without a stable-reduction theorem.

\begin{proposition}\label{prop:stabilizerproductfull}
Let $Z=Z_L\cup_D Z_R$ be a marked expansion with the roots fixed.  Let $F$ be the matching of simple nonzero branch sheaves $F_L,F_R$ along a stable nonzero boundary sheaf $G$.  Denote by
\[
 T_i=\{\sigma_i\in\Aut(Z_i/Y_i):\sigma_i^*F_i\simeq F_i\}
\]
the intrinsic geometric stabilizer, interpreted as the fppf sheaf of such target automorphisms.  Then the intrinsic stabilizer of the matched object is
\[
 T_F\simeq T_L\times T_R.
\]
In particular it is finite if and only if both branch stabilizers are finite.  These are stabilizers of the unframed sheaves, not stabilizers of a fixed Quot presentation.
\end{proposition}
\begin{proof}
A target automorphism preserving the marked node is a pair of relative expansion automorphisms.  If it preserves $F$, restriction preserves $F_L$ and $F_R$, giving an injection $T_F\to T_L\times T_R$.

Conversely choose a pair in the latter group and pass to an fppf cover on which sheaf isomorphisms
\[
 u_L:\sigma_L^*F_L\xrightarrow{\sim}F_L,\qquad
 u_R:\sigma_R^*F_R\xrightarrow{\sim}F_R
\]
exist.  Relative expansion automorphisms are the identity on the outer divisor.  Relative to the chosen matching, the two boundary isomorphisms therefore differ by an automorphism of $G$.  Stability gives $\Aut(G)=\GG$, so the difference is a scalar $c$.  Replacing $u_L$ by $c^{-1}u_L$ makes the boundary square commute:
\[
\begin{tikzcd}[column sep=large]
 (\sigma_L^*F_L)|_D\arrow[r,"\sigma^*\phi"]\arrow[d,"c^{-1}u_L"']
 &(\sigma_R^*F_R)|_D\arrow[d,"u_R"]\\
 F_L|_D\arrow[r,"\phi"']&F_R|_D.
\end{tikzcd}
\]
The kernel construction for coherent gluing then supplies an isomorphism $\sigma^*F\simeq F$.  Different choices of $u_i$ change only the common scalar of this isomorphism, not the target automorphism.  This establishes surjectivity fppf-locally and hence equality of the stabilizer sheaves.  The construction respects composition.  The final assertion follows from the product description.
\end{proof}

The stability of $G$ is used exactly in reducing its boundary automorphism to a scalar.  With an impure boundary, the full automorphism matching remains part of the derived boundary stack.  A separate application of this proposition is then not justified.

\section{Divisorial classes on the expanded central fiber}\label{sec:divisorial}
We first determine the possible divisor classes for a perfect pure sheaf with arbitrary derived boundary values. We then specialize to the hyperplane class.

Use the relatively ample line bundle $\mathcal A$ of Lemma~\ref{lem:ample}. On an expansion its pullback is used only to record the original Hilbert polynomial, not as an ample line bundle on the bubbles. Let
\[
 Z=Y_1\cup_{S_0}\Delta_1\cup\cdots\cup_{S_{k-1}}\Delta_k\cup_{S_k}Y_2,
 \qquad p_j:\Delta_j\longrightarrow S.
\]
Assume $\Pic(S)=\mathbb Zh$. Orient every bubble so that its left section $D_{j,-}$ has normal bundle $\OO_S(-h)$ and its right section $D_{j,+}$ has normal bundle $\OO_S(h)$. Then
\begin{equation}\label{eq:sectionrelation}
 D_{j,+}=D_{j,-}+p_j^*h.
\end{equation}
Write $H_B,E_B$ for the classes on $Y_2$ as before. The intersection numbers needed below are
\begin{equation}\label{eq:blowupnumbers}
 H_B^3=1,\quad H_B^2E_B=0,\quad H_BE_B^2=-20,\quad E_B^3=-180.
\end{equation}
The last two identities follow from the degree $20$ and normal bundle $\OO_\Gamma(4)\oplus\OO_\Gamma(5)$ of the blowup center, whose normal degree is $180$.

\begin{lemma}\label{lem:boundarydimension}
Let $F$ be a pure two-dimensional admissible sheaf on an expanded threefold. At every gluing divisor its restriction is either zero or has dimension one. In particular a nonzero restriction cannot consist solely of a finite-length sheaf. No purity of the restriction is asserted here.
\end{lemma}
\begin{proof}
Work at a closed point of a gluing divisor and use
\[
 A_0=\CC[[x,y]]/(xy),\quad B=A_0[[u,v]],\quad M=F^{\wedge}.
\]
Normality makes $M$ flat over $A_0$ by Lemma~\ref{lem:flatnode}. The element $x+y$ is regular on $A_0$, and hence on $M$. If $M\ne0$, purity gives $\dim_BM=2$, so $\dim(M/(x+y)M)=1$. In the ring $B/(x+y)$ the class of $x$ is nilpotent, since $x^2=0$ and $y=-x$. Passing from $M/(x+y)M$ to its quotient by $x$ does not change its support as an $R=\CC[[u,v]]$-module: localization and Nakayama give this assertion at every prime. The latter quotient is $M/(x,y)M$.

It therefore has dimension one. If it vanishes, Nakayama gives $M=0$ near that point. Applying the argument at all closed points proves the claim.
\end{proof}

\subsection{Divisorial conservation for an arbitrary positive multiple}
Retain the original polarization $\mathcal A$ of Lemma~\ref{lem:ample}.  Its restriction to the smooth fiber is $2H$, whereas on the roots it is $H_Q$ and $6H_B-E_B$.  On a bubble it is $p^*h$; this last line bundle is used for the fixed original polynomial, not for testing ampleness.

\begin{theorem}\label{prop:generalprofiles}
Assume $\Pic(S)=\mathbb Zh$. Let $F$ be a pure two-dimensional sheaf, perfect on a fixed expanded quintic central fiber, with quadratic Hilbert coefficient $10r$ for the original polarization. Write $E_e=\LL j_e^*F$, $c_1(E_e)=b_eh$, and $[F_2]_{\mathrm{div}}=aH_B-cE_B$. Then
\begin{align}
 [F_Q]_{\mathrm{div}}&=b_0H_Q,\nonumber\\
 [F_{\Delta_j}]_{\mathrm{div}}
    &=b_jp_j^*h+(b_j-b_{j-1})D_{j,-},\label{eq:generalcomponentclasses}\\
 a+4b_q&=5r,\qquad c=r-b_q,\label{eq:generaldegreeconservation}\\
 0\le b_0\le\cdots\le b_q&\le\lfloor5r/4\rfloor.\nonumber
\end{align}
Moreover $b_e=0$ if and only if $E_e=0$. The statement uses the class $[\mathcal H^0(E_e)]-[\mathcal H^{-1}(E_e)]$ on $S$. Neither cohomology sheaf is required to be a curve sheaf.
\end{theorem}
\begin{proof}
Theorem~\ref{thm:perfectpurematchingnew} gives ordinary pure branch sheaves with effective divisorial cycles. Write the left divisor as $nH_Q$ and the $j$-th bubble divisor as $u_jp_j^*h+s_jD_{j,-}$. Their effectivity gives $n\ge0$ and $s_j\ge0$, the latter by intersection with a general ruling fiber. Restriction of the first Chern class of a perfect complex commutes with pullback. Using $D_{j,-}|_{D_{j,-}}=-h$, $D_{j,-}|_{D_{j,+}}=0$ and $E_B|_S=5h$ gives
\[
 n=b_0,\qquad u_j-s_j=b_{j-1},\quad u_j=b_j,\qquad a-5c=b_q.
\]
Thus $b_j-b_{j-1}=s_j\ge0$. The restriction of a rank-zero perfect complex has rank zero. Consequently the quadratic coefficient of the Hilbert polynomial of $E_e$ on $S$ is zero, even if both its cohomology sheaves have positive rank. Applying the Hilbert polynomial to \eqref{eq:generalKnormalnew} yields
\begin{align*}
 20r&=4b_0+4\sum_{j=1}^q(b_j-b_{j-1})+(6H_B-E_B)^2(aH_B-cE_B)\\
    &=4b_q+16a-60c=4a+16b_q,
\end{align*}
where the middle equality uses $H_B^3=1$, $H_B^2E_B=0$, $H_BE_B^2=-20$, $E_B^3=-180$. The pushforward of the effective right-root cycle to $\PP^3$ is effective or zero, hence $a\ge0$. Substitution now gives \eqref{eq:generaldegreeconservation} and the upper bound. If $b_e=0$, all earlier $b_j$ vanish. The component immediately to the left of that cut has zero divisorial cycle, and its pure branch sheaf is zero. Its derived restriction is therefore zero.

The reverse implication follows by taking $c_1$.
\end{proof}

\begin{corollary}\label{cor:primitiveautomaticnew}
For $r=1$ every perfect pure sheaf of the preceding theorem is admissible at all the double divisors. Its restrictions there can still have a zero-dimensional subsheaf.
\end{corollary}
\begin{proof}
The numerical sequence takes its values in $\{0,1\}$. At a zero cut the boundary vanishes by the theorem. At a positive cut inspect the component immediately to its right. If it is a bubble, monotonicity forces the right-hand coefficient of that bubble to be $1$ as well. Its divisorial class is $p^*h$. Since $p_*\OO_\Delta=\OO_S$, its effective support divisor is $p^{-1}C$ for a curve $C\in|h|$ and cannot contain the attaching section.

If instead the component is $B$, its divisorial class is $H_B$. Its effective support is the pullback of a plane in $\PP^3$. No plane contains the quartic surface $S$, so its strict transform cannot contain $S$. For a pure sheaf on either of these components, the equation of the attaching divisor is a non-zero-divisor by \eqref{eq:twoHboundarynew}. If $s$ cuts out the attaching divisor in this branch $W$, then
\[
 \mathcal H^{-1}(E_e)=\ker\big(F_W(-S_e)\xrightarrow{s}F_W\big)=0,
 \qquad \mathcal H^0(E_e)=F_W/sF_W(-S_e).
\] As $E_e=F\otimes^{\LL}\OO_{S_e}$, this is exactly $\Tor_1^{\OO_Z}(F,\OO_{S_e})=0$. This proves the assertion for every cut, with no assumption of admissibility on the initial sheaf.
\end{proof}

\subsection{The three primitive support profiles}

\begin{theorem}\label{thm:primitiveprofiles}
Assume $\Pic(S)=\mathbb Zh$. Let $F$ be a nonzero pure perfect surface sheaf on $Z$ whose Hilbert polynomial for the pullback of $\mathcal A$ has quadratic coefficient $10$, the coefficient of a primitive hyperplane sheaf on a smooth quintic. Denote its branch restrictions by $F_1,F_{\Delta_j},F_2$, and its boundary restrictions by $G_e$. There are integers $b_e$ with
\begin{equation}\label{eq:binaryprofiles}
 0\le b_0\le b_1\le\cdots\le b_k\le1,
 \qquad c_1(G_e)=b_eh.
\end{equation}
Here $b_e=0$ if and only if $G_e=0$. The divisorial support classes on the roots and bubbles are
\begin{align}
 [F_1]_{\mathrm{div}}&=b_0H_Q,\\
 [F_{\Delta_j}]_{\mathrm{div}}
   &=b_jp_j^*h+(b_j-b_{j-1})D_{j,-},\\
 [F_2]_{\mathrm{div}}&=(5-4b_k)H_B-(1-b_k)E_B.
 \label{eq:profileclasses}
\end{align}
Every nonzero branch has one integral support with generic rank one. If all intrinsic expansion stabilizers are finite, precisely the following support patterns are possible:
\begin{enumerate}[label=(\roman*)]
\item All $b_e=1$. The left root support is a hyperplane section of $Q$, each bubble support is $p_j^{-1}C$ for the common integral $C\in|h|$, and the right root support is the transform of a plane.
\item There is at least one bubble, $b_0=0$, and $b_e=1$ for $e\ge1$. The left root sheaf is zero. The first bubble support is a graph isomorphic to $S$, disjoint from its left section. Every subsequent bubble is supported on $p_j^{-1}C$, and the right root support is the transform of a plane.
\item There are no bubbles and the entire sheaf is supported on $Y_2$, on an integral divisor of class $5H_B-E_B$ disjoint from $S$.
\end{enumerate}
If the nonzero boundary restrictions are additionally pure, all of them have primitive Mukai vector $(0,h,\chi_e)$ and define points of the six-dimensional compactified Jacobians $J_{\chi_e}$.
\end{theorem}
\begin{proof}
By Corollary~\ref{cor:primitiveautomaticnew} the sheaf is admissible. The restrictions to the smooth components are pure, by Theorem~\ref{thm:classicalnode}. Their divisorial cycles are effective. Write the left root cycle as $nH_Q$ with $n\ge0$, and the right root cycle as $aH_B-cE_B$. Its pushforward to $\PP^3$ is effective, so $a\ge0$. Since $\Pic(\Delta_j)=p_j^*\Pic(S)\oplus\mathbb Z[D_{j,-}]$, write its cycle as
\[
 U_jp_j^*h+r_jD_{j,-}.
\]
Effectivity implies $r_j\ge0$: it is the degree of that effective divisor on a general projective-line fiber. Normality excludes the distinguished sections as support components. Derived restriction is ordinary restriction, so the first Chern character of the boundary is the restriction of the branch divisorial class. Using \eqref{eq:sectionrelation} gives
\[
 b_0=n,\qquad U_j-r_j=b_{j-1},\qquad U_j=b_j,\qquad a-5c=b_k.
\]
It follows that $r_j=b_j-b_{j-1}\ge0$. The boundary classes are effective; Lemma~\ref{lem:boundarydimension} shows that a zero class is equivalent to a zero boundary sheaf. Hence every $b_e\ge0$.

We now use the quadratic coefficient, which does not receive any contribution from the curve terms in the normalization sequence. On the roots the polarization restricts to $H_Q$ and $6H_B-E_B$, and on every bubble it is $p_j^*h$. Thus twice that coefficient is
\[
 4n+\sum_j4r_j+(6H_B-E_B)^2(aH_B-cE_B)
 =4b_k+16a-60c.
\]
The equality uses \eqref{eq:blowupnumbers}. Substituting $c=(a-b_k)/5$ simplifies it to $4a+16b_k$. It must equal $20$. Consequently
\begin{equation}\label{eq:totaldegreeprofile}
 a+4b_k=5.
\end{equation}
Since $a\ge0$ and $b_k$ is a nonnegative integer, $b_k\le1$. This proves \eqref{eq:binaryprofiles}, and then $a=5-4b_k$ and $c=1-b_k$, proving all the class formulas.

It remains to identify the effective supports, not just their classes. A divisor of class $H_Q$ is integral because $\Pic(Q)=\mathbb ZH_Q$ and the class is primitive. A divisor of class $H_B$ is the pullback of a plane; no plane contains $\Gamma$, so it has no exceptional divisor component and is integral. On a bubble with $b_{j-1}=b_j=1$ the divisor class is $p_j^*h$. Since $p_{j*}\OO_{\Delta_j}=\OO_S$, all effective divisors in this system are pullbacks of members of $|h|$.

Such a member is integral because $h$ generates $\Pic(S)$. Equality of the boundary sheaves identifies their supporting curves, so all adjacent vertical bubbles use the same curve.

At a transition $b_{j-1}=0,b_j=1$, the divisor class is $D_{j,+}$. Sections of $\OO(D_{j,+})$ are pairs $(c,s)$ with $c\in\CC$ and $s\in H^0(S,\OO_S(h))$. Normality at $D_{j,-}$ forces $c\ne0$: otherwise the divisor contains that section. After rescaling, the divisor is the graph $z=-s/c$ in the normal-bundle chart. Normality at $D_{j,+}$ forces $s\ne0$. The graph is isomorphic to $S$, avoids $D_{j,-}$, and meets $D_{j,+}$ along the integral curve $(s=0)\in|h|$.

If $b_k=0$, the right class is $5H_B-E_B$. Its sections identify with
\[
 H^0(\PP^3,I_\Gamma(5))
 =\CC f_5\oplus f_4H^0(\PP^3,\OO(1)).
\]
A divisor not containing $S$ has equation $cf_5+f_4\ell$ with $c\ne0$. This polynomial is irreducible. Indeed, $\Gamma$ is integral, so a factorization would have a factor vanishing on $\Gamma$. No nonzero form of degree less than four vanishes there, and its degree-four ideal is generated by $f_4$. In degree five a nontrivial factorization would therefore be $f_4$ times a linear form, contrary to $c\ne0$.

There is no exceptional component, since $I_\Gamma^2$ has no element of degree five. Its strict transform is integral. Its restriction to $S$ is a nonzero section of the trivial line bundle $(5H_B-E_B)|_S$, so the support is disjoint from $S$.

The binary nondecreasing sequence has at most one transition. A bubble for which both adjacent $b_e$ vanish has divisor class zero and therefore carries the zero sheaf. Its rescaling fixes the entire sheaf, giving a continuous intrinsic stabilizer. Such bubbles are excluded by the assumed finiteness. Thus a transition must occur at the first bubble, and an all-zero profile must have no bubbles. This proves the three alternatives.

Each identified effective support is integral with multiplicity one. Purity then makes its defining equation annihilate the branch sheaf, exactly as in Proposition~\ref{prop:primitive}; the sheaf is rank-one torsion-free on that support. Finally a pure curve sheaf with class $h$ is rank-one torsion-free on its integral support and is stable. For a stable boundary sheaf $G$, $\tr:\Ext^2_S(G,G)\to H^2(\OO_S)$
is dual to the scalar isomorphism $\CC\to\Hom_S(G,G)$.

The K3
deformation theorem \cite{Mukai} gives smoothness and tangent
dimension $2+\langle v(G),v(G)\rangle=2+h^2=6$.
\end{proof}

\begin{remark}\label{rem:profilemeaning}
The theorem uses the fixed original Hilbert polynomial, not the existence of a universal modified ample line bundle. It classifies divisorial support and boundary curve classes, not all degree-four and degree-six charges. In particular the zero-boundary alternative is part of the primitive quintic problem. It is not replaced by a point of $J_\chi$.
\end{remark}

\begin{corollary}\label{cor:profilesimple}
Under the finite-stabilizer hypotheses of Theorem~\ref{thm:primitiveprofiles}, every such sheaf is simple on its fixed expanded target, independently of a further Gieseker-stability choice.
\end{corollary}
\begin{proof}
A rank-one torsion-free sheaf on an integral projective surface has its global endomorphism algebra embedded in the function field of that surface. This algebra is a finite-dimensional domain over $\CC$, and hence is $\CC$. Each nonzero branch restriction therefore has only scalar endomorphisms.

An endomorphism of the matched sheaf restricts to those scalar maps, and the scalars on adjacent nonzero branches are equal because their common boundary sheaf is nonzero. The nonzero branches form a connected interval in all three profiles. The normalization sequence shows that restriction on endomorphisms is injective. Index the occupied interval by $a,\ldots,b$. Restriction identifies
\[
 \End(F)=\ker\!\left(\CC^{b-a+1}\xrightarrow{\partial}\CC^{b-a}\right),
 \qquad \partial(c)_i=c_i-c_{i+1}.
\]
The kernel is $\CC(1,\ldots,1)$.
\end{proof}

The stability criteria below are applied to the connected interval of nonzero branch sheaves. An empty root has no nonzero subsheaf and contributes no stability inequality. A cut with zero boundary is treated separately in Section~\ref{sec:zerorelative}; it is not forced into an inequality with a zero boundary polynomial.

\begin{figure}[htbp]
\centering
\begin{tikzpicture}[x=1cm,y=.92cm,font=\small,line cap=round,line join=round]
 \foreach \x/\v in {3.4/Q,6.0/\Delta,8.6/B}{\node[qiInk] at (\x,1.02) {$\v$};}
 \foreach \y/\lab in {0/{(i)},-1.05/{(ii)}}{
 \node[anchor=west] at (.3,\y) {\lab};
 \draw[qiBlue!55,line width=.7pt] (2.13,\y)--(9.76,\y);
 \foreach \x in {4.70,7.30}{\draw[qiTeal,line width=.85pt] (\x,\y-.32)--(\x,\y+.32);}
 }
 \draw[qiBlue,line width=1.15pt] (2.33,.11)--(9.67,.11);
 \foreach \x in {4.70,7.30}{\fill[qiAmber] (\x,.11) circle (2.4pt);}
 \draw[qiBlue,line width=1.15pt] (5.25,-1.27).. controls (5.89,-1.36) and (6.5,-.94)..(7.30,-.94)--(9.67,-.94);
 \draw[qiTeal,fill=white,line width=.9pt] (4.70,-1.05) circle (2.6pt);
 \fill[qiAmber] (7.30,-.94) circle (2.4pt);
 \node[above,inner sep=1pt,font=\footnotesize] at (5.94,-.97) {graph};
 \node[anchor=west] at (.3,-2.10) {(iii)};
 \node[above,inner sep=2pt,qiInk] at (3.40,-2.10) {$Q$};
 \node[above,inner sep=2pt,qiInk] at (7.53,-2.10) {$B$};
 \draw[qiBlue!55] (2.13,-2.10)--(9.76,-2.10);
 \draw[qiTeal,line width=.85pt] (5.74,-2.42)--(5.74,-1.78);
 \draw[qiTeal,fill=white,line width=.9pt] (5.74,-2.10) circle (2.6pt);
 \draw[qiBlue,line width=1.15pt] (6.62,-2.22)--(9.62,-2.22);
 \foreach \y/\b in {0/{(1,1)},-1.05/{(0,1)},-2.1/{(0)}}{\node[anchor=west] at (10.1,\y) {$\b$};}
\end{tikzpicture}
\caption{Primitive support profiles. The bars mark joining K3 surfaces and the integers at right are their boundary multiplicities. Filled marks have curve class $h$ and open marks have zero boundary. In (i) a ruled support meets both sections. In (ii) the graph misses the left section. In (iii) there is no bubble and the divisor on $B$ misses $S$. Additional occupied ruled components may occur at the positive end of (i) and (ii).}\label{fig:threecutprofiles}
\end{figure}
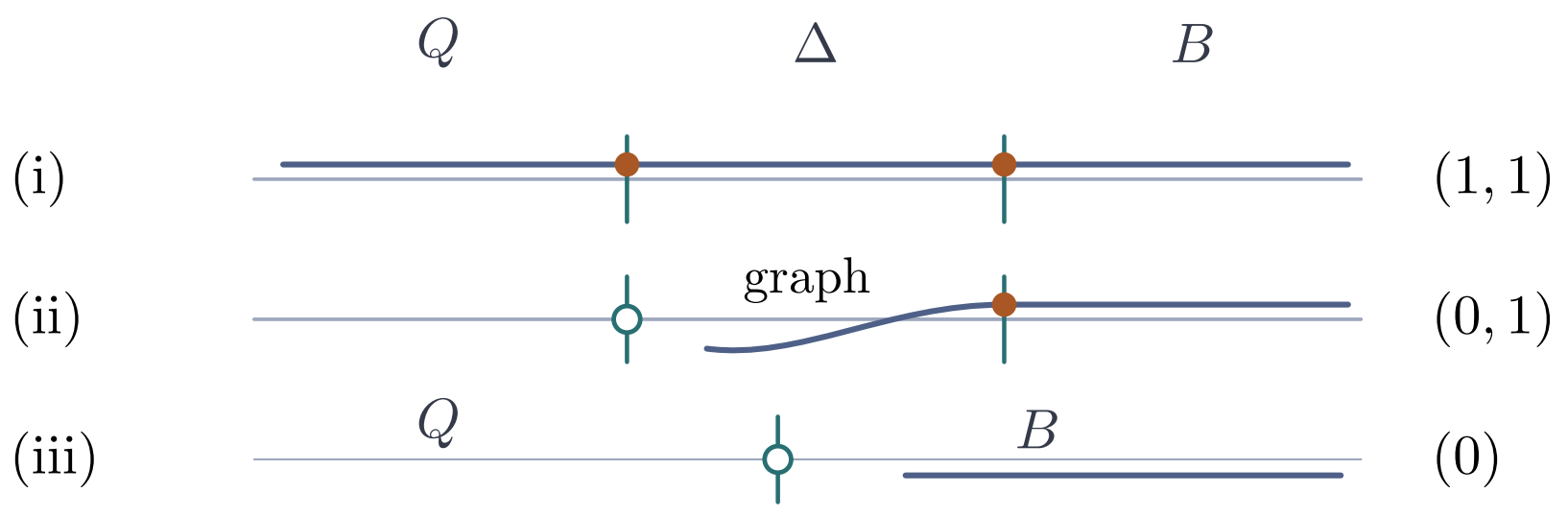

Figure~\ref{fig:threecutprofiles} shows the first occupied component in each of the three cases of Theorem~\ref{thm:primitiveprofiles}. The incidence of the support with the joining divisor determines whether that marking can take a value in a compactified Jacobian.

\subsection{An actual zero-boundary specialization}
Let $P=(f_1=0)\subset\PP^4$ and let $T=(f_5=0)\subset P$ be smooth. In the family $f_1f_4=tf_5$, the fixed surface $T$ is contained in every fiber. It has a lift to the semistable resolution.

\begin{proposition}\label{prop:zeroboundaryfamily}
There is a flat family of embedded surfaces $\mathcal T\cong T\times C$ in the semistable quintic family. Its generic member is a hyperplane section; its special member lies entirely in $Y_2$ and has class $5H_B-E_B$, disjoint from $S$. For any rank-one torsion-free sheaf $G$ on $T$, its pushforward from $\mathcal T$ is a flat family of stable pure surface sheaves. Its special boundary is zero. For $G=\OO_T$ the generic character is $\ch(\OO_H)$.
\end{proposition}
\begin{proof}
The product $T\times C$ embeds in the unresolved family by $f_1=f_5=0$. For $j\colon T\times C\hookrightarrow(f_1 f_4=t f_5)$ the ideal of the center satisfies
\[
 j^{-1}(f_1,t)\cdot\OO_{T\times C}=(t)\simeq\OO_{T\times C}(-T\times\{0\}).
\] The universal property of the blowup therefore lifts the embedding, and its strict transform is still $T\times C$.

On the central right root the lift is the transform of the quintic surface containing $\Gamma$ with multiplicity one, hence has class $5H_B-E_B$. The center $\Gamma$ is a Cartier divisor on $T$, so this transform is isomorphic to $T$. The strict transforms of the two divisors $f_4=0$ and $f_5=0$ are disjoint in their blowup along the complete intersection. This proves disjointness from $S$. Equivalently, the support computation is the last case of Theorem~\ref{thm:primitiveprofiles}.

The constant family $G\boxtimes\OO_C$ is flat, and pushforward by a closed immersion preserves this flatness and purity of its fibers. Every fiber support is integral and the sheaf has generic rank one, so it is stable for every ample polarization by the elementary reduced-polynomial argument. On the generic fiber the divisor is $f_1=0$, yielding the final character assertion.
\end{proof}

\subsection{A bubble bound when the branch hulls are invertible}\label{sec:linehullbound}
The primitive profile theorem also supplies a bound in a larger class than the fixed-support Hilbert scheme. 

\begin{proposition}\label{prop:linehullbubblebound}
Let $F$ satisfy the hypotheses of Theorem~\ref{thm:primitiveprofiles}, with pure nonzero boundary. Suppose that on each nonzero branch support $T_i$ there is an exact sequence
\[
 0\longrightarrow F_i\longrightarrow L_i\longrightarrow Q_i\longrightarrow0,
\]
where $L_i$ is its intrinsic rank-one reflexive hull, $L_i$ is an invertible sheaf, and $Q_i$ is finite. Put $n=\sum_i\operatorname{length}Q_i$. If $k$ is the number of bubbles, then
\begin{equation}\label{eq:linehullbubblebound}
 k\le n+1.
\end{equation}
In the all-nonzero-boundary profile one has the sharper bound $k\le n$. In a flat family with a compatible invertible hull and finite flat quotient of degree $n$, these are uniform bounds for all its admissible stable limits satisfying the stated conditions.
\end{proposition}
\begin{proof}
Every support in Theorem~\ref{thm:primitiveprofiles} is Cartier in a smooth component, hence is Cohen--Macaulay. At a closed point where a finite $Q_i$ is nonzero, the depth lemma applied to the displayed exact sequence gives $\depth F_i=1$: the middle term has depth two and the last has depth zero. At a nonzero pure boundary the sheaf has depth two by Theorem~\ref{thm:nodeperfect}, or its smooth-Cartier relative version. Therefore $Q_i$ is disjoint from the distinguished curves.

Consider a vertical bubble support $T_i=p_i^{-1}C$. If $Q_i=0$, its sheaf is the line bundle $L_i$. Every line bundle on the projective-line bundle over the integral curve $C$ is a tensor product of a pullback line bundle with a power of the tautological line bundle. It is preserved by the fiber-rescaling action, with a linearization whose restrictions on the two sections differ by scalar characters.

These characters can be absorbed by multiplying all sheaves on the left subchain and all sheaves on the right subchain by the corresponding scalars. The normalization kernel then gives an isomorphism between $F$ and every rescaling pullback of $F$. Thus that entire $\GG$ lies in the intrinsic expansion stabilizer, contradicting its finiteness.

Let $J_{\rm vert}$ be the indices of the vertical bubbles. The stabilizer calculation gives
\[
 Q_i=0\ \Longrightarrow\ \GG\subset H_F\quad (i\in J_{\rm vert}),\qquad
 \#J_{\rm vert}\le\sum_{i\in J_{\rm vert}}\length Q_i\le n .
\] The profile theorem allows at most one additional graph bubble; its support already has finite rescaling stabilizer and it need not carry a quotient point. This proves \eqref{eq:linehullbubblebound}, and there is no additional graph in the first profile. In the stated family, flatness makes the total quotient length constant, proving the uniform assertion. No constancy of the canonical hull is asserted for families outside this last hypothesis.
\end{proof}

The hypothesis of invertible branch hulls is important for the rescaling argument. A non-invertible hull on a singular ruled support can itself contain data not preserved by fiber rescaling; its contribution must be bounded separately to apply \cite[Assumption~3.2.2]{Kim}.

\section{Stability on an expanded chain}\label{sec:stability}

\subsection{Subsheaves when the derived restriction has two cohomology sheaves}
Use Theorem~\ref{thm:perfectpurematchingnew} at a marked crossing. Write $D=j_*\mathcal H^{-1}(E)$ and let $K$ be the kernel of the difference map to $j_*\mathcal H^0(E)$. Both are sheaves on the glued threefold $Z$. The gluing data determine a class $\xi\in\Ext^1_Z(K,D)$ by the first sequence in \eqref{eq:Torcorrectionnew}. Fix an ample line bundle on this target, and denote the quadratic coefficient of a Hilbert polynomial by $a(-)$, with $a(0)=0$.

\begin{theorem}\label{thm:extensionstability}
For every $U\subseteq D$ and $V\subseteq K$ let
\[
 \xi_{U,V}\in\Ext^1_Z(V,D/U)
\]
be obtained by first pulling $\xi$ back along $V\hookrightarrow K$ and then pushing out along $D\twoheadrightarrow D/U$. A perfect pure glued sheaf $F$ is Gieseker stable precisely when
\[
 \frac{P_U+P_V}{a(U)+a(V)}<\frac{P_F}{a(F)}
\]
for all $(U,V)$ with $\xi_{U,V}=0$, other than $(0,0)$ and $(D,K)$. Replacing $<$ by $\le$ gives semistability. The test uses the actual extension determined by the equivalence on $S$.
\end{theorem}
\begin{proof}
Let $F_V=F\times_K V$. For $U\subset D$ the pullback and pushout defining $\xi_{U,V}$ give
\[
 0\longrightarrow D/U\longrightarrow F_V/U\xrightarrow{\pi}V\longrightarrow0.
\]
A section $s$ of $\pi$ gives the inverse image $A$ of $s(V)$ in $F_V$. Conversely $A\subset F$ with $A\cap D=U$ and image $V$ induces $s:V\simeq A/U\hookrightarrow F_V/U$. Hence such $A$ exist exactly when $\xi_{U,V}=0$, and each fits into
\[
 0\to U\to A\to V\to0,\qquad
 P_A=P_U+P_V,\quad a(A)=a(U)+a(V).
\]
The pairs $(0,0)$ and $(D,K)$ give $0$ and $F$. Every other $A$ has positive quadratic coefficient by purity, giving the stated test. It includes $D\subset F$ when proper and nonzero, and specializes to the ordinary subobject test for $K=F$ when $D=0$.
\end{proof}

\subsection{The stability test with all subobjects retained}
Let $Z=Z_0\cup_{S_0}\cdots\cup_{S_{q-1}}Z_q$ be a fixed expanded target and choose an ample line bundle on this target.  Let $F_i$ be pure two-dimensional branch sheaves, normal to their adjacent divisors, with specified common boundary sheaves $G_e$.  The $G_e$ may be impure, and the $F_i$ may have arbitrary generic ranks and reducible supports.  Let $F$ be their admissible matching from Theorem~\ref{thm:classicalnode}.  Write $a(E)$ for the quadratic coefficient of a two-dimensional sheaf $E$.

For coherent subsheaves $E_i\subset F_i$, allowing $E_i=0$, use the fixed matching isomorphism to regard both boundary images as subsheaves of $G_e$.  Define
\begin{align}
 I_e(E_\bullet)&=
     \operatorname{im}(E_i\to G_e)+
     \operatorname{im}(E_{i+1}\to G_e),\label{eq:generaledgeimage}\\
 E(E_\bullet)&=
 \ker\left(\bigoplus_i(i_i)_*E_i
        \xrightarrow{\rho_i-\rho_{i+1}}
               \bigoplus_e(j_e)_*I_e(E_\bullet)\right).
 \label{eq:generalmatchedsubobject}
\end{align}
The symbols $i_i,j_e$ are the component and node inclusions.  Thus no equality of boundary supports replaces the matching of boundary modules.

\begin{theorem}\label{thm:generalpolynomialstability}
The sheaf $F$ is Gieseker (semi)stable if and only if, for every tuple $E_\bullet$ with some $E_i\ne0$ and some $E_i\ne F_i$,
\begin{equation}\label{eq:generalstabilitytest}
 \frac{\sum_iP_{E_i}-\sum_eP_{I_e(E_\bullet)}}{\sum_i a(E_i)}
       \le\frac{P_F}{a(F)}.
\end{equation}
Here use $<$ for stability and $\le$ for semistability.  In the denominator zero terms contribute zero.  This criterion is valid for arbitrary branch ranks and arbitrary support multiplicities.
\end{theorem}
\begin{proof}
At each node $e$ the image of the difference map is $I_e(E_\bullet)$ by definition. The joining divisors are disjoint, so \eqref{eq:generalmatchedsubobject} gives a short exact sequence and
\begin{equation}\label{eq:generalEpolynomial}
 P_{E(E_\bullet)}=\sum_iP_{E_i}-\sum_eP_{I_e(E_\bullet)},\qquad
 a(E(E_\bullet))=\sum_i a(E_i).
\end{equation}
Indeed normality of $F_i$ makes each edge image at most one-dimensional. A tuple in the statement has positive quadratic coefficient. Its kernel lies in $F$ and is proper: equality would force $F_i=E_i$ via the surjection $F\to (i_i)_*F_i$. Applying (semi)stability to this kernel proves necessity.

Conversely let $0\ne A\subsetneq F$ and take $E_i=\operatorname{im}(A\to(i_i)_*F_i)$ on $Z_i$. Then
\[
 A\subset E(E_\bullet)\subset F,\qquad
 a(A)=a(E(E_\bullet)),\qquad
 p_A=p_{E(E_\bullet)}-\frac{P_{E(E_\bullet)/A}}{a(A)} .
\]
The quadratic coefficients agree because the two subsheaves have the same stalk at each generic support point, off the nodes. If some $E_i\ne F_i$, the stated inequalities give the required comparison. Otherwise $E(E_\bullet)=F$ and $P_{F/A}(m)>0$ for $m\gg0$, giving $p_A<p_F$ even in the weak test.
\end{proof}

If the branch ranks are one on integral surfaces, the maximal tuples for a prescribed generic-support subset are precisely the interval kernels used in Lemma~\ref{lem:maxsub}.  Their residual polynomials give $D_e-B_e$, $-D_e$, and $(D_m-B_m)-D_{l-1}$, recovering Theorem~\ref{thm:cutcriterion}.  For higher ranks the edge images in \eqref{eq:generaledgeimage} cannot be suppressed.

Even at the generic point of a boundary curve, two one-dimensional subspaces of a two-dimensional boundary fiber have sum of dimension one when they agree and dimension two otherwise.  Thus changing the matching isomorphism can change $P_{I_e}$, while leaving every branch polynomial fixed.  The full criterion retains this information.

\subsection{The numerical criterion}
Let $Z=Z_0\cup_{S_0}\cdots\cup_{S_{r-1}}Z_r$ be a fixed expanded target. Fix an ample line bundle $L$; a rational polarization is treated after clearing denominators. Assume each nonzero $F_i$ is rank-one torsion-free on an integral surface $T_i\subset Z_i$, is normal to its adjacent divisors, and has nonzero pure restrictions identified with $G_e$. We consider a connected chain of these nonzero branches. No smoothness of $T_i$ or of the curves is required.

Write
\[
P_i(z)=a_i z^2+b_i z+c_i,\qquad
B_e(z)=d_ez+e_e,\qquad a_i,d_e>0.
\]
Put
\begin{equation}\label{eq:cutdata}
P=\sum_iP_i-\sum_eB_e,\quad A=\sum_i a_i,\quad p=P/A,
\quad D_e=P_{\le e}-A_ep,
\end{equation}
where
\[
P_{\le e}=\sum_{i\le e}P_i-\sum_{j<e}B_j,\qquad
A_e=\sum_{i\le e}a_i.
\]
All inequalities below are lexicographic inequalities of polynomials for sufficiently large $z$.

\begin{lemma}\label{lem:maxsub}
For a union $U$ of the integral support components, let $K_U$ be the kernel of the restriction of $F$ to the complementary subchains. Every subsheaf of $F$ having generic rank one precisely on $U$ is contained in $K_U$. For an interval $[l,m]$,
\begin{equation}\label{eq:intervalpoly}
P_{K_{[l,m]}}=P_{[l,m]}
-\boldsymbol1_{l>0}B_{l-1}-\boldsymbol1_{m<r}B_m.
\end{equation}
For a disconnected union, the kernel is the direct sum of those for its intervals.
\end{lemma}
\begin{proof}
On a component outside $U$, the image of such a subsheaf has generic rank zero and hence dimension at most one. It is a subsheaf of the pure branch sheaf, so it vanishes. This proves the containment. For an interval, the restriction to either cut boundary is surjective. If there are two cuts, their divisors are disjoint, so the simultaneous restriction to their direct sum is surjective. Its kernel has the polynomial \eqref{eq:intervalpoly}.

Distinct intervals share no gluing divisor, which proves the direct-sum assertion.
\end{proof}

\begin{figure}[htbp]
\centering
\begin{tikzpicture}[x=1cm,y=1cm,font=\small,line cap=round,line join=round]
 \foreach \x in {1.4,3.4,5.4,7.4,9.4}{
  \draw[qiBlue!40] (\x-.97,-.12).. controls (\x-.65,.29) and (\x+.58,.29)..(\x+.97,-.12);
  \draw[qiBlue!35] (\x-.97,-.12).. controls (\x-.65,-.48) and (\x+.58,-.48)..(\x+.97,-.12);}
 \foreach \x in {3.4,5.4,7.4}{
  \path[fill=qiBlue!15,draw=qiBlue,line width=.72pt] (\x-.97,-.12).. controls (\x-.65,.29) and (\x+.58,.29)..(\x+.97,-.12)
  .. controls (\x+.58,-.48) and (\x-.65,-.48)..cycle;}
 \foreach \x in {2.4,8.4}{\draw[qiAmber,line width=1.25pt] (\x,-.47)--(\x,.27);}
 \foreach \x in {4.4,6.4}{\fill[qiTeal] (\x,-.12) circle (2pt);}
 \draw[qiAmber,->] (2.4,.34)--(2.4,.69);
 \draw[qiAmber,->] (8.4,.34)--(8.4,.69);
 \node[above] at (2.4,.67) {$G_{l-1}$}; \node[above] at (8.4,.67) {$G_m$};
 \node[below,qiInk] at (3.4,-.39) {$Z_l$}; \node[below,qiInk] at (5.4,-.39) {$\cdots$};
 \node[below,qiInk] at (7.4,-.39) {$Z_m$};
\end{tikzpicture}
\caption{The maximal subsheaf on an interior interval $[l,m]$. The shaded components retain their internal identifications. A section extends by zero to the complementary components precisely when it vanishes at the two exterior cuts, subtracting $B_{l-1}$ and $B_m$ from the interval polynomial.}\label{fig:intervalends}
\end{figure}
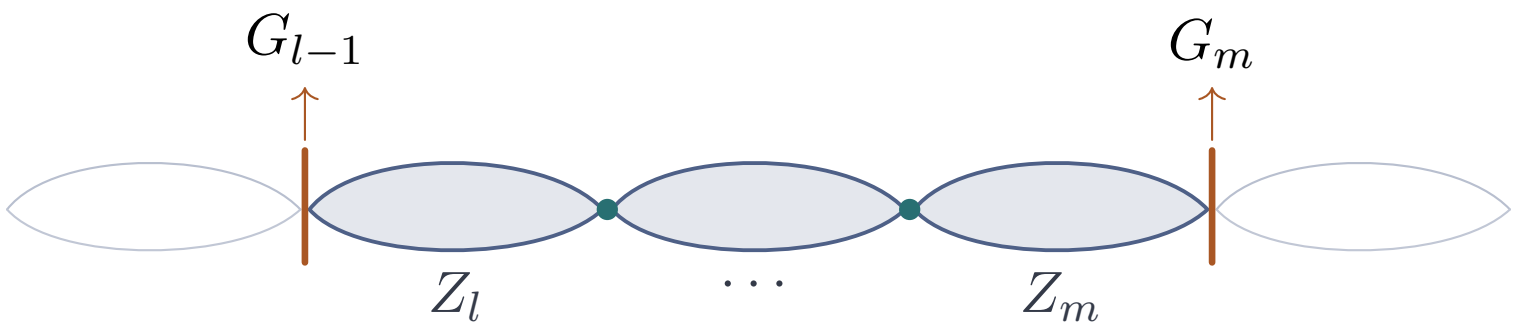

For an interval strictly inside the chain, the two end conditions in Lemma~\ref{lem:maxsub} occur on distinct divisors (Figure~\ref{fig:intervalends}).

\begin{theorem}\label{thm:cutcriterion}
The glued sheaf $F$ is Gieseker-stable if and only if
\begin{equation}\label{eq:cutcriterion}
0<D_e<B_e\qquad(0\le e<r).
\end{equation}
It is semistable if and only if the inequalities are weak.
\end{theorem}
\begin{proof}
The maximal prefix subobject has polynomial $P_{\le e}-B_e$, and the maximal complementary suffix subobject has polynomial $P-P_{\le e}$. Subtracting their leading coefficients times $p$ gives $D_e-B_e$ and $-D_e$. Stability therefore implies \eqref{eq:cutcriterion}.

Conversely, take a nonzero proper subsheaf $E\subset F$. Its generic rank on each integral component is zero or one. If all ranks are one, the quotient is nonzero of dimension at most one and $p_E<p$. Otherwise let $U$ be the proper union where its generic rank is one. By Lemma~\ref{lem:maxsub}, $E\subset K_U$, with quotient of dimension at most one. It is enough to check the maximal kernel.

Prefixes and suffixes were already treated. For an interior interval, normalization gives
\[
P_{[l,m]}=P_{\le m}-P_{\le l-1}+B_{l-1}.
\]
Consequently
\begin{equation}\label{eq:interiorres}
P_{K_{[l,m]}}-(A_m-A_{l-1})p
=(D_m-B_m)-D_{l-1}<0.
\end{equation}
For disconnected $U$ the residual is a sum of negative polynomials. This proves stability. The same argument with weak inequalities proves semistability, since an everywhere-rank-one proper subobject still has a strictly smaller polynomial.
\end{proof}

\begin{corollary}\label{cor:numericstability}
On a fixed numerical stratum, global stability is independent of the matching isomorphisms and of the particular rank-one sheaves with the specified polynomials. The stratum consists entirely of stable matched sheaves or contains none.
\end{corollary}
\begin{proof}
The inequalities in \eqref{eq:cutcriterion} involve only the fixed polynomials and leading coefficients.
\end{proof}

\subsection{Vertical twists}
Assume the fixed target is the central fiber of a regular one-parameter semistable family $\mathscr Z\to\Spec R$, with Cartier components $Z_i$. At a cut $e$ define
\begin{equation}\label{eq:verticaltwist}
\mathcal L_e=\OO_{\mathscr Z}\left(-\sum_{i\le e}Z_i\right)\big|_Z.
\end{equation}
It is trivial generically. It restricts to the positive cut divisor on the prefix and to the negative cut divisor on the suffix, and is trivial on all other components. Its restrictions on the two sides of the selected divisor agree; denote the resulting line bundle on that divisor by $N_e$.

\begin{lemma}\label{lem:twists}
Put $F(k)=F\otimes\mathcal L_e^k$ and $G_e(k)=G_e\otimes N_e^k$. Their total polynomial is independent of $k$. At the selected cut,
\begin{equation}\label{eq:twistrecur}
D_e(k+1)=D_e(k)+B_e(k+1),\qquad B_e(k)=P_{G_e(k)}.
\end{equation}
All $D_f,B_f$ for $f\ne e$ are unchanged. If $\delta_e$ is the linear coefficient of $D_e$, then
\begin{equation}\label{eq:twistslope}
\delta_e(k)=\delta_e(0)+kd_e.
\end{equation}
\end{lemma}
\begin{proof}
Normality gives the exact sequences
\begin{align*}
0&\longrightarrow F_{\le e}(k)\longrightarrow F_{\le e}(k+1)
\longrightarrow G_e(k+1)\longrightarrow0,\\
0&\longrightarrow F_{>e}(k+1)\longrightarrow F_{>e}(k)
\longrightarrow G_e(k)\longrightarrow0.
\end{align*}
Adding their polynomials and subtracting the boundary polynomial as in normalization cancels every change. Surface leading coefficients are unchanged by a line-bundle twist, so the first sequence proves \eqref{eq:twistrecur}. A prefix ending before $e$ contains neither affected component; one ending after $e$ contains both and their boundary, where the same cancellation applies.

The other boundary line restrictions are trivial because distinct distinguished divisors on a bubble are disjoint. Finally a line-bundle twist of a curve sheaf preserves its linear leading coefficient, proving \eqref{eq:twistslope}.
\end{proof}

\begin{theorem}\label{thm:verticalbalance}
Every orbit of the twists \eqref{eq:verticaltwist} contains a semistable representative. It has exactly one stable representative unless at least one cut has the equality wall described below, in which case it has no stable representative.

At a cut choose the unique $k_0\in\mathbb Z$ with $0\le\delta_e(k_0)<d_e$. If $\delta_e(k_0)>0$, the unique stable choice is $k_0$. If $\delta_e(k_0)=0$, write $D_e(k_0)=c_e$. The stable choice is $k_0$ for $c_e>0$ and $k_0+1$ for $c_e<0$. If $c_e=0$, both choices are strictly semistable at that cut and neither is stable.
\end{theorem}
\begin{proof}
Fix a cut and let $d_e=[z]B_e>0$. Lemma~\ref{lem:twists} gives $\delta_e(k)=\delta_e(0)+kd_e$. Put $k_0=-\lfloor\delta_e(0)/d_e\rfloor$. If $0<\delta_e(k_0)<d_e$, then $0<D_e(k_0)<B_e(k_0)$ for $z\gg0$. All other twists fail the weak linear inequalities.

If $\delta_e(k_0)=0$, write $D_e(k_0)=c_e$. Only $k_0$ and $k_0+1$ have linear coefficient in $[0,d_e]$. By \eqref{eq:twistrecur},
\[
  D_e(k_0)=c_e,\qquad
  B_e(k_0+1)-D_e(k_0+1)=-c_e .
\]
At $k_0$ the upper inequality is automatic, and at $k_0+1$ the lower inequality is automatic. Thus $c_e>0$ selects $k_0$, $c_e<0$ selects $k_0+1$, and $c_e=0$ gives equality at both endpoints. The equalities come from the proper prefix or suffix subsheaf of Theorem~\ref{thm:cutcriterion}. The twists at distinct cuts act independently by Lemma~\ref{lem:twists}, so the choices give precisely the stated stable and semistable representatives.
\end{proof}

For a fixed joining divisor $D_e$, the restriction of
$\mathcal L_e$ to the two components adjacent to it is
$\OO(D_e)$ on the left and $\OO(-D_e)$ on the right. This explains
the direction of the step in \eqref{eq:twistslope}. In the diagram
the interior dot is the unique admissible slope when
$\delta_e(0)\notin d_e\mathbb Z$.
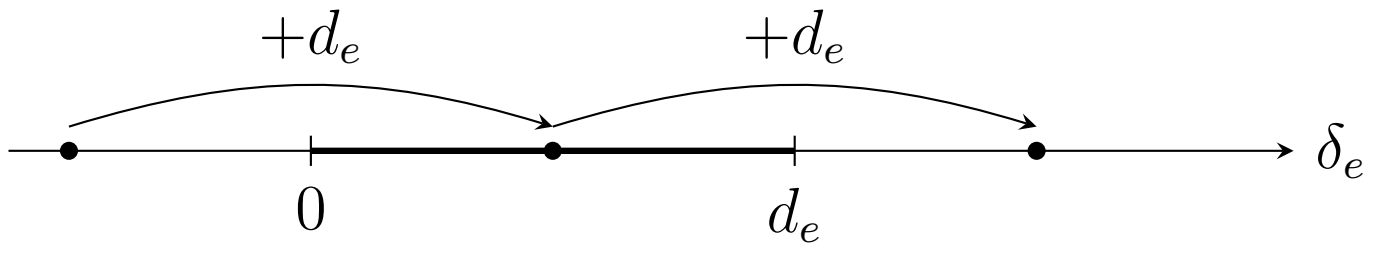
\begin{center}
\begin{tikzpicture}[x=1cm,y=1cm,>=stealth]
 \draw[->] (0,0)--(8.5,0) node[right] {$\delta_e$};
 \draw[very thick] (2,0)--(5.2,0);
 \foreach \x/\lab in {2/0,5.2/d_e}
   {\draw (\x,.10)--(\x,-.10) node[below] {$\lab$};}
 \foreach \x in {.4,3.6,6.8} {\fill (\x,0) circle (1.7pt);}
 \draw[->] (.4,.16) to[bend left=17] node[above] {$+d_e$} (3.6,.16);
 \draw[->] (3.6,.16) to[bend left=17] node[above] {$+d_e$} (6.8,.16);
\end{tikzpicture}
\end{center}

\begin{corollary}\label{cor:twistrepair}
If the no-strictly-semistable condition holds on the entire twist orbit, an admissible rank-one chain sheaf with pure boundary has a unique stable vertical twist. Tensoring preserves normality and the boundary Cohen--Macaulay condition. For an already existing flat one-parameter family, it preserves the generic sheaf and the deformation complex of the sheaf on the expanded threefold.
\end{corollary}
\begin{proof}
The wall alternative is excluded by the hypothesis. Tensoring by a line bundle is exact, commutes with restriction, and preserves depth. It gives the canonical identification
\[
\RR\mathcal Hom(F\otimes\mathcal L,F\otimes\mathcal L)
\simeq\RR\mathcal Hom(F,F).
\]
The line bundles in \eqref{eq:verticaltwist} exist on the one-parameter family and are trivial on its generic fiber. Their canonical equivariant structures also show that tensoring preserves the intrinsic geometric stabilizer.
\end{proof}

The theorem concerns a fixed expanded target and a fixed polarization. Familywise descent of a chosen stability rule remains part of the compactification data in Section~\ref{sec:compactification}. Individual component charges and boundary Euler characteristics may change under the twists, and the charge decomposition must retain those changes.

\subsection{Stability after cutting}
Fix a cut $m$. The two branch polynomials and leading coefficients are
\[
P_L=P_{\le m},\quad a_L=A_m,\qquad
P_R=P-P_{\le m}+B_m,\quad a_R=A-A_m.
\]
It follows that
\begin{equation}\label{eq:branchpolys}
p_L=p+\frac{D_m}{A_m},\qquad
p_R=p+\frac{B_m-D_m}{A-A_m}.
\end{equation}

\begin{proposition}\label{prop:cutting}
For an internal cut $e<m$ in the left branch its stability polynomial is
\begin{equation}\label{eq:leftcut}
D_e^L=D_e-\frac{A_e}{A_m}D_m.
\end{equation}
For an internal cut $e>m$ in the right branch it is
\begin{equation}\label{eq:rightcut}
D_e^R=D_e+\frac{A-A_e}{A-A_m}(B_m-D_m).
\end{equation}
The two branches are stable if and only if their respective polynomials lie strictly between zero and $B_e$ at every internal cut.
\end{proposition}
\begin{proof}
For $e<m$, substitute \eqref{eq:branchpolys} to obtain
\[
 D_e^L=P_{\le e}-A_ep_L
      =(P_{\le e}-A_ep)-\frac{A_e}{A_m}D_m.
\]
For $e>m$, normalization restores the boundary polynomial $B_m$ to the right prefix. Hence
\begin{align*}
 D_e^R
 &=P_{\le e}-P_{\le m}+B_m-(A_e-A_m)p_R\\
 &=D_e-D_m+B_m-\frac{A_e-A_m}{A-A_m}(B_m-D_m)\\
 &=D_e+\frac{A-A_e}{A-A_m}(B_m-D_m).
\end{align*}
Apply Theorem~\ref{thm:cutcriterion} to the left and right chains.
\end{proof}

\begin{corollary}\label{cor:matchingstability}
Suppose that the polynomials are constant on a pair of relative moduli spaces, that their boundary evaluations match, and that \eqref{eq:cutcriterion}, \eqref{eq:leftcut}, and \eqref{eq:rightcut} all hold on the indicated rank-one chain sectors. Then their full matching product consists of globally stable sheaves, and restriction of any globally stable matched object in the sector lies in those relative stability problems.

No further open condition depending on the gluing isomorphism is needed.
\end{corollary}
\begin{proof}
Use the numerical nature of the three sets of inequalities and the family equivalence of Proposition~\ref{prop:familyglue}. The remaining defining conditions of the relative spaces, including numerical class and intrinsic automorphisms, are understood to be the specified ones; the statement establishes their stability compatibility.
\end{proof}

\subsection{A boundary-weighted stability condition compatible with cutting}\label{sec:weighted}
The ordinary relative Gieseker inequalities in Proposition~\ref{prop:cutting} can also be expressed by a stability condition which remembers the cut. The following construction is numerical and intrinsic to the sheaf and its restriction; it does not impose flatness of a reflexive hull or add a freely varying framing map.

Let $W$ be a connected relative chain with outer divisor $D$, and suppose that $F$ has the rank-one integral-component structure of Section~\ref{sec:stability}, with nonzero pure outer restriction $G$. Write $B=P_G$ and $A=a(F)$. For a nonzero subsheaf $E\subset F$, put
\[
 \epsilon_D(E)=\begin{cases}1,&E\longrightarrow G\text{ is nonzero},\\0,&E\longrightarrow G\text{ is zero}.\end{cases}
\]
Fix a rational polynomial $\alpha$ of degree at most one with $0<\alpha<B$. Say that $F$ is $\alpha$-stable at $D$ if
\begin{equation}\label{eq:weighteddefinition}
 \frac{P_E-\epsilon_D(E)\alpha}{a(E)}
 <\frac{P_F-\alpha}{A}
 \qquad(0\ne E\subsetneq F).
\end{equation}
All comparisons are eventual polynomial comparisons. This is a boundary-weighted stability rule; no claim that it is automatically the Gieseker stability condition of \cite[Definitions~3.2.1 and~3.2.5]{Kim} is made.

Orient $W$ from its root towards $D$, and put
\[
 p^\alpha=\frac{P_F-\alpha}{A},\qquad
 D_e^\alpha=P_{\le e}-A_ep^\alpha
\]
at its internal cuts.

\begin{theorem}\label{thm:weightedcriterion}
Under the stated rank-one and normality hypotheses, $F$ is $\alpha$-stable if and only if
\begin{equation}\label{eq:weightedcriterion}
 0<D_e^\alpha<B_e
\end{equation}
at every internal cut.
\end{theorem}
\begin{proof}
A maximal prefix subobject has zero outer image and polynomial $P_{\le e}-B_e$. Its residual in \eqref{eq:weighteddefinition} is $D_e^\alpha-B_e$. A maximal terminal suffix has nonzero outer image and polynomial $P_F-P_{\le e}$. Its residual, after subtracting $\alpha$, is $-D_e^\alpha$. Thus stability implies the displayed conditions.

Conversely, consider a proper subsheaf $E$. If its generic rank is one on all components and its outer image is nonzero, then $F/E$ is nonzero of dimension at most one. Subtracting the same $\alpha$ from both polynomials gives the required strict inequality. If its outer image is zero, then $E\subset F(-D)$. Normality gives $P_{F(-D)}=P_F-B$, so its residual is at most $\alpha-B<0$.

For a proper generic-support union $U$, use the maximal kernel $K_U$ of Lemma~\ref{lem:maxsub}. An interior interval has residual
\[
 (D_m^\alpha-B_m)-D_{l-1}^\alpha<0.
\]
A prefix and a terminal suffix were already computed. If $U$ includes the outer component but $E$ has zero outer image, replace the terminal summand of $K_U$ by its kernel of restriction to $G$. The outer restriction of that summand is surjective, since the other cut divisor is disjoint from $D$. The polynomial decreases by $B$, while its weighted penalty decreases by $\alpha$. Its residual therefore changes by $\alpha-B<0$, preserving strict negativity.

A disconnected support is a sum of intervals and has at most one terminal interval meeting $D$. Adding the corresponding residuals proves the inequality for $K_U$, or for the indicated outer kernel. The further quotient by $E$ has nonnegative Hilbert polynomial of smaller degree and can only decrease the residual. This covers every subsheaf.
\end{proof}

\begin{theorem}\label{thm:weightedmatching}
Fix a cut $m$ in a complete chain and a numerical sector satisfying $0<D_m<B_m$. Give the left and right relative branches the weights
\begin{equation}\label{eq:cutweights}
 \alpha_L=D_m,
 \qquad \alpha_R=B_m-D_m.
\end{equation}
Orient each branch from its own root toward the marked boundary. Then the complete matched sheaf is Gieseker-stable if and only if the left branch is $\alpha_L$-stable and the right branch is $\alpha_R$-stable. Consequently, on this fixed numerical sector the classical globally stable matching functor is the full fiber product of the two boundary-weighted relative functors, provided their other defining conditions agree.
\end{theorem}
\begin{proof}
The branch polynomials satisfy
\[
 P_L-\alpha_L=A_mp,
 \qquad P_R-\alpha_R=(A-A_m)p.
\]
Both corrected reduced polynomials therefore equal the global $p=P/A$. For an internal left cut the corrected defect is exactly $D_e$. For an internal right cut, oriented from the right root, its prefix polynomial is $P-P_{\le e}+B_e$, so its corrected defect is $B_e-D_e$. Conditions \eqref{eq:weightedcriterion} on both branches are therefore exactly the global conditions \eqref{eq:cutcriterion} at all cuts other than $m$.

The condition at $m$ was assumed and also makes both weights strictly between zero and $B_m$.

Classical restriction and gluing are inverse for all flat matched families by Proposition~\ref{prop:familyglue}. The calculation shows that their stability loci agree without a further condition on the matching isomorphism. This proves the functorial statement on the given numerical sector.
\end{proof}

\begin{remark}\label{rem:weightedproperness}
On a fixed projective target the theorem factors the boundary-weighted stability inequalities through the two branches. For expansions, boundedness, descent of the weights, and the valuative criteria must still be established for this particular rule, or it must be proved equivalent to a rule covered by the construction of \cite[Definitions~3.2.1 and~3.2.5]{Kim}.

The ordinary relative Gieseker alternative is Proposition~\ref{prop:cutting}. Either alternative can be used in the virtual diagonal theorem once the independently defined relative factors are proper.
\end{remark}

\subsection{A stable expanded limit with a noninvertible hull in every multiplicity}
\label{sec:matrixexpandedallr}
There are DVR families in every class of divisorial supports $rH$, $r\ge2$, for which the occupied bubble carries a finite cover of the middle K3. The following construction calculates the sheaf on the expansion, its boundary, its polarization and its stabilizer. It uses the first degeneration of Section~\ref{sec:geometry} and gives the two charts of the one-bubble expansion explicitly below.

\begin{theorem}\label{thm:matrixexpandedallr}
Fix $r\ge2$. In a general quintic degeneration there is a flat family of perfect sheaves after the base change $t=s^2$ and a one-bubble expansion. Its generic sheaf is stable of character $r(1-e^{-H})$. The special sheaf is stable for a rational ample polarization on the expanded fiber. It has zero branch on $Q$, an integral degree-$r$ spectral cover on the bubble, and an integral degree-$r$ determinantal surface on $B$.

The contracted bubble cycle is $r[S]$. On an occupied branch the sheaf is maximal Cohen--Macaulay of generic rank one and requires $r$ generators at a boundary point. In particular its Cohen--Macaulay hull is noninvertible and its punctual defect is zero. Its common nonzero boundary is a stable curve sheaf of Mukai vector $(0,rh,-2r)$; it too requires $r$ generators at that point.

The geometric stabilizer on the expanded target consists of $(\mathrm{id},a\,\mathrm{id}_F)$ with $a\in\CC^\times$.
\end{theorem}

\begin{proof}
Fix $p\in S\setminus\Gamma$. Take three independent forms $u,v,w$
on $\PP^3$ vanishing at $p$ and lift them to $\PP^4$, so that
$x=f_1,u,v,w$ have common zero $p$. Put
$A=uA_1+vA_2+wA_3$ with $A_i\in\operatorname{Mat}_{r}(\CC)$.
The coefficient space of the three matrices is irreducible. We first
specify the open conditions on these coefficients.

\emph{The determinantal supports.} For $r\ge3$ the projective determinantal loci satisfy
\[
 \operatorname{codim}\{M:\operatorname{corank}M\ge2\}=4,\qquad
 (\det M=0)_{\mathrm{sm}}=\{M:\operatorname{corank}M=1\}.
\]
Consequently a general $\PP^3\subset\PP(\mathrm{Mat}_r)$ avoids the singular locus and meets the hypersurface transversely; for $r=2$ the determinant quadric in $\PP^3$ is already smooth. Normalize an invertible basis matrix by left multiplication to $I_r$. Bertini for a general plane and the hypersurface sequences give smooth connected sections
\[
 \{\det(zI_r-A)=0\}\subset\PP^3_{[z:u:v:w]},\qquad
 \{\det A=0\}\subset\PP^2_{[u:v:w]} .
\]
Thus the cone over the second curve in the plane root is integral with vertex $p$. Since $\Gamma$ is not contained in the cone, its Cartier pullback to $B$ is its integral strict transform.

On $\operatorname{Tot}_S(\OO_S(h))$ set
$\Sigma=V(\det(zI_r-A|_S))$. Its finite algebra over $S$ is locally
free with basis $1,z,\ldots,z^{r-1}$ in a trivialization of
$\OO_S(h)$. To check integrality for general $A$, take a cyclic
matrix with diagonal $u$ and successive off-diagonal entries
$l_1,\ldots,l_r$ whose plane sections on $S$ are distinct and
integral. Its characteristic polynomial is
$(z-u)^r-\prod l_i$, up to a nonzero scalar multiplying $l_1$.

At the generic point of $(l_1=0)$ it is Eisenstein in $z-u$.
The finite locally free algebra embeds in its generic fiber, a field.
This gives one integral cover. The family of monic covers is proper
and flat over the coefficient space, so geometric integrality is open.

For $C=V(\det A|_S)$ choose the cyclic entries to give
$u^r+v^{r-1}l$, $l\in\langle u,v,w\rangle$. The linear system
$\langle u^r,v^{r-1}u,v^r,v^{r-1}w\rangle$ has only finitely many
base points on $S$ since $S\cap V(u,v)$ is finite for general
$u,v$. The ratios of its last three coordinates give projection
from $p$ to $\PP^2$, which is dominant and generically finite.

Bertini irreducibility in characteristic zero gives an integral general
member (generic smoothness away from the base locus, and a Cartier
curve on $S$ has no embedded points). The chart where the coefficient
of $u^r$ is nonzero gives the displayed cyclic determinants. Thus
integrality of $C$ is another nonempty open condition. Intersect these
opens with $\operatorname{tr}(A|_S)\ne0$.

\emph{The flat sheaf on the expansion.} Put $t=s^2$.
At a crossing $xy=s^2$ blow up $(x,s)$. The two charts and the
matrices in their respective trivializations are
\[
\begin{array}{c|c|c}
 &\text{target chart}&\text{matrix}\\ \hline
 U_z & x=sz,\quad zy=s &M_z=zI_r-A\\
 U_{w'} &s=xw',\quad y=sw' &M_{w'}=I_r-w'A .
\end{array}
\]
On the overlap $w'z=1$ and $M_{w'}=w'M_z$. Globally the entries
are sections of the line bundle
$L=c_X^*\OO_{\XX}(H)\otimes\OO(\mathscr Q)$, where $c_X$ is
contraction and $\mathscr Q$ is the left-root divisor over $s=0$.
Indeed the rational matrix $(x/s)I_r-A$ has only the simple pole
on $\mathscr Q$ removed by the second trivialization. Thus the local
cokernels glue to the cokernel of $(L^{-1})^{\oplus r}\to\OO^{\oplus r}$.
On the generic fiber the matrix is $xI_r-sA$ up to a unit.

Its determinant defines in $\PP^4$ the cone over the preceding
smooth projective determinantal surface with vertex $p$. A general
quintic misses $p$ and cuts the smooth part of this cone transversely.
For fixed $s\ne0$ varying $f_5$ in $f_1f_4=s^2f_5$ varies the
quintic through the entire linear system. Openness in the parameter
therefore gives a smooth complete intersection on the geometric
generic fiber for general data.

Its Koszul resolution over $\overline{\CC((s))}$ gives
$h^0(\OO)=1$ on the geometric generic surface. The
surface is connected and hence integral.

The special-fiber matrices give the following sheaves (bundle twists
are implicit in the local free presentations):
\[
\begin{array}{c|c|c}
 \text{locus} &\text{matrix}&\text{cokernel}\\ \hline
 Q & I_r&0\\
 \Delta\setminus S_\infty &zI_r-A|_S&F_\Sigma\\
 B &-b^*A&F_B\\
 S=\Delta\cap B&-A|_S&E .
\end{array}
\]
Here $b:B\to\PP^3$ is the blowdown and $S_\infty=Q\cap\Delta$.
Along $S_\infty$ the matrix is invertible. On each target fiber
the determinant avoids all generic points of the target components.
The fibers are Cohen--Macaulay, so this determinant is a
nonzerodivisor. The equation $\operatorname{adj}(M)M=\det(M)I_r$
gives fiberwise injectivity.

The local flatness criterion applied
to the two-term presentation gives a sheaf flat over $\CC[[s]]$
and the same presentation after every base change. On the generic
fiber its Chern character is $r(1-e^{-H})$.

On each occupied branch the determinant has multiplicity one
at its generic point. The cokernel therefore has generic rank one
on its support. At a closed point its free resolution on the smooth
threefold gives depth two; on $S$ the analogous resolution gives
depth one.

The determinants are integral, so these are torsion-free
rank-one sheaves on their supports, and are Gieseker-stable there:
a nonzero proper subsheaf has the same generic rank and the nonzero
quotient has smaller dimension. This also proves stability on the
generic target. The cover is finite flat of degree $r$, hence
$(p_\Delta)_*[\Sigma]=r[S]$ as cycles.

At $p$ on $B$, over $p$
at $z=0$ on $\Delta$, and at $p$ on $C$, the presentation tensored
with $k(p)$ is $k(p)^r\xrightarrow{0}k(p)^r$. The minimal number
of generators of each cokernel is $r$. The branch sheaf is already
maximal Cohen--Macaulay on its surface, so its Cohen--Macaulay hull
is the sheaf itself and its punctual defect has length zero. The
boundary on the expanded target is computed by
\[
 0\to\OO_S(-h)^r\xrightarrow{-A|_S}\OO_S^r\to E\to0,
 \qquad v(E)=(0,rh,-rh^2/2)=(0,rh,-2r).
\]
Before expansion the specialization on $Q$ is instead
$\operatorname{coker}(xI_r)=\OO_S^r$, with derived restriction
$[\OO_S(-h)^r\xrightarrow{0}\OO_S^r]$.

\emph{Stability on the union.} Let $P_i$ be the Hilbert polynomial
on a branch, $a_i=[m^2]P_i$, and $P_E$ the polynomial on the common
curve. The boundary presentation has zero kernel, so the normalization
sequence gives
\[
 0\to F_0\to F_\Sigma\oplus F_B\xrightarrow{\mathrm{res}_\Sigma-
 \mathrm{res}_B}E\to0,\quad P=P_\Sigma+P_B-P_E,\quad a=a_\Sigma+a_B .
\]
For a subsheaf of rank one on both components, a proper quotient
lowers the reduced polynomial. For rank $(1,0)$ or $(0,1)$ the
maximal subsheaf has polynomial $P_i-P_E$. Purity follows from
the injection above. The two remaining strict comparisons are
therefore equivalent to
\[
 0<D=P_\Sigma-\frac{a_\Sigma}{a}P<P_E \qquad(m\gg0).
\]
To obtain an ample chamber evaluate first at the nef classes
$(p_\Delta^*h,b^*H_{\PP^3})$. On $\Sigma$ the matrix relation
expresses multiplication by $z$ on $\OO_S^r$ as $A$, so
$(p_\Delta)_*F_\Sigma=\OO_S^r$. On $B$, projection formula,
$\RR b_*\OO_B=\OO_{\PP^3}$ and the pulled-back two-term resolution
give the polynomials
\[
\begin{array}{c|ccc}
 &P_\Sigma&P_B&P_E\\ \hline
 &r(2m^2+2)&r\binom{m+2}{2}&r(4m-2).
\end{array}
\]
Hence $P=r(5m^2/2-5m/2+5)$, $D=2rm-2r$ and $P_E-D=2rm$.
Both linear coefficients are strictly positive. In the order
$B,\Delta,Q$ choose the rational divisor classes
\[
 H_B-\epsilon E_B,\qquad p_\Delta^*(dh)+cS_\infty,\qquad
 (d-c)H_Q,\qquad d=1-5\epsilon,\quad 0<c\ll d .
\]
For $\epsilon>0$ small the first is ample by the blowup criterion.
The middle class has positive ruling degree $c$ and restricts to
$dh$ and $(d-c)h$ on its two sections, so is ample by the split
projective-bundle criterion. The end classes have exactly these
restrictions on the joining surfaces. Clearing denominators gives
line bundles with isomorphic restrictions and thus an ample line bundle
on the union.

Riemann--Roch makes the linear coefficients
of $D,P_E-D$ continuous in $(\epsilon,c)$ near $(0,0)$ since
$a(0,0)>0$. Their strict positivity persists. This proves stability
on the expanded target.

Finally write $P_A(z)=\det(zI_r-A|_S)=z^r-\operatorname{tr}(A|_S)z^{r-1}+\cdots$.
A bubble dilation $z\mapsto\lambda z$ preserving the support satisfies
$\lambda^{-r}P_A(\lambda z)=P_A(z)$. The coefficient of $z^{r-1}$
gives $(1-\lambda^{-1})\operatorname{tr}(A|_S)=0$, hence $\lambda=1$.
Stability gives $\Aut(F_0)=\GG$ from scalar multiplication on $F_0$.
\end{proof}

At this positive cut $\langle v(E),v(E)\rangle=4r^2$. For the simple branches the relative virtual dimension is $1+2r^2$ by Theorem~\ref{thm:branchindexnew}, and the stable boundary has dimension $2+4r^2$. In particular the multiplicity-two calculation is $9+9-18=0$. The strict inequalities above verify fixed-target stability and the construction verifies extension for these DVR families.

A proper virtual comparison for their entire numerical class additionally requires valuative reduction for all families of that class and descent of stability on the expansion stack.

\begin{figure}[htbp]
\centering
\begingroup
\definecolor{coverink}{HTML}{53618A}
\definecolor{coverseam}{HTML}{176568}
\definecolor{coverpt}{HTML}{B15C25}
\begin{tikzpicture}[x=1cm,y=1cm,line cap=round,line join=round,font=\small]
\node[anchor=west] at (-.25,3.62) {\textbf{(a)}\quad A general geometric fiber};
\node[anchor=west] at (6.5,3.62) {\textbf{(b)}\quad The fiber over $p$};
\foreach \a in {0,6.75}{
 \begin{scope}[shift={(\a,0)}]
 \path[draw=coverink!65,fill=coverink!8] (-.12,.63)--(1.94,.03)--(3.46,.43)--(1.40,1.03)--cycle;
 \foreach \i in {.25,.5,.75}{
  \draw[coverink!24] ({-.12+2.06*\i},{.63-.60*\i})--({1.40+2.06*\i},{1.03-.60*\i});
  \draw[coverink!24] ({-.12+1.52*\i},{.63+.40*\i})--({1.94+1.52*\i},{.03+.40*\i});}
 \draw[coverink,line width=.8pt] (1.67,.53)--(1.67,3.06);
 \fill[coverseam] (1.67,.53) circle (1.6pt);
 \node[text=coverseam,anchor=west] at (2.74,.06) {$S$};
 \end{scope}}
\foreach \j/\w in {1.04/1,1.86/2,2.65/3}{
 \fill[coverpt] (1.67,\j) circle (2.1pt);
 \node[anchor=west] at (1.86,\j) {$z_{\w}$};}
\node[anchor=east,text=coverink] at (1.59,.37) {$q$};
\node[anchor=east,text=coverink] at (8.34,.37) {$p$};
\draw[coverpt,line width=.8pt] (8.42,1.53) circle (5pt);
\draw[coverpt,line width=.65pt] (8.42,1.53) circle (3.45pt);
\fill[coverpt] (8.42,1.53) circle (1.95pt);
\node[anchor=west] at (8.81,1.53) {$z=0$};
\node[anchor=west,text=coverink] at (8.81,2.40) {$A(p)=0$};
\node at (2.73,-.43) {$\det(zI_3-A(q))=\prod_{i=1}^3(z-z_i)$};
\node at (9.35,-.43) {$\Sigma_p=\Spec k(p)[z]/(z^3)$};
\end{tikzpicture}
\endgroup
\caption{The degree-three spectral cover on the bubble in Theorem~\ref{thm:matrixexpandedallr}. The dots in (a) label the distinct roots on a general geometric fiber. In (b) the rings indicate a single length-three fiber at $z=0$. At its closed point the presentation matrix vanishes, so the fiber of the cokernel there has dimension three. This accounts for its three local generators although the sheaf has generic rank one on $\Sigma$.}
\label{fig:spectralfiber}
\end{figure}
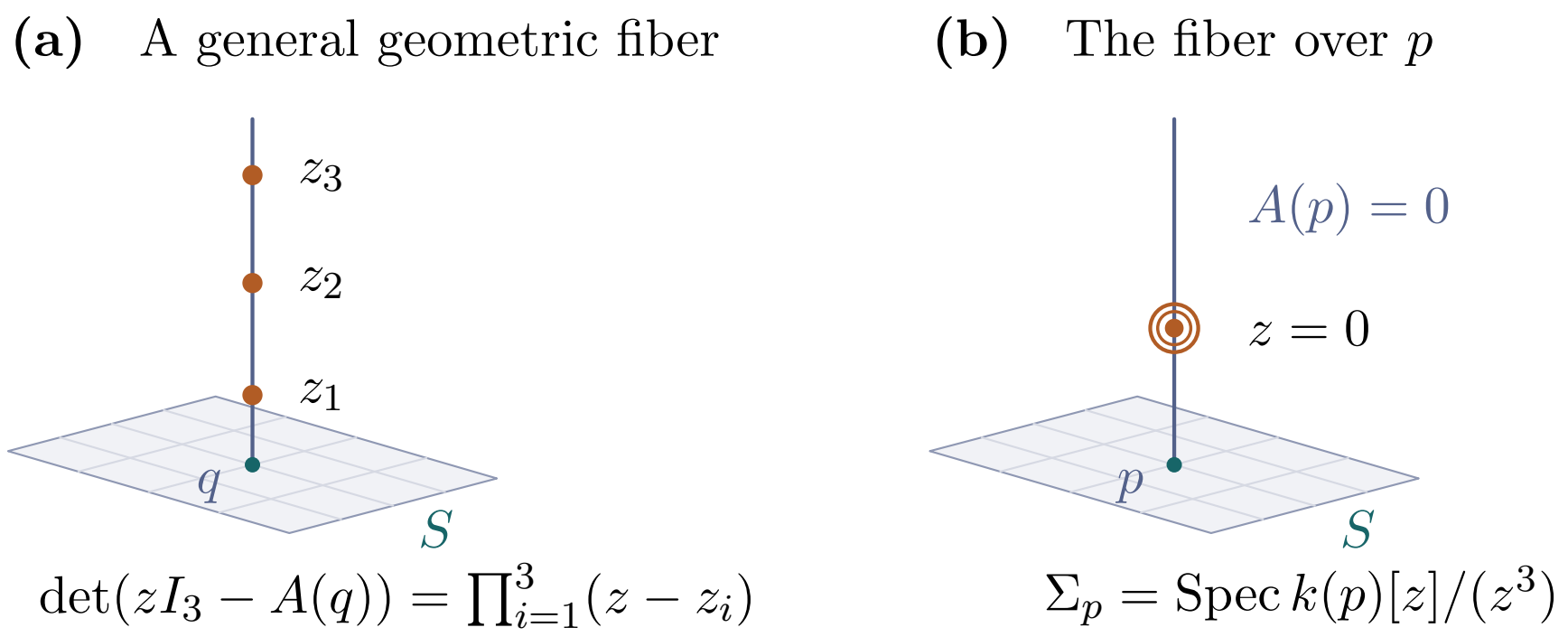

Figure~\ref{fig:spectralfiber} records two fibers of the determinantal example. The fiber over $p$ detects the difference between generic rank one on the support and local freeness there.

\subsection{A uniform chamber for bounded point ideals}\label{sec:idealchamber}
Consider a hyperplane support in the model \eqref{eq:quinticgeom}. Its expanded surface is
\[
T[k]=T_1\cup_C R_1\cup_C\cdots\cup_C R_k\cup_C T_2,
\]
where $T_1$ is a quartic K3 surface, $T_2$ is the blowup of a plane at twenty points on the smooth plane quartic $C$, and $R_j$ is a ruled surface over $C$. The theorem below is a concrete sector to which the abstract stability criterion applies. Its proof uses complete Hilbert polynomials, so it also applies to any other members with the same numerical data and matching hypotheses.

For two inserted ruled surfaces the three copies of $C$ and the
degrees of the polarization on them are as follows. The dots indicate
possible supports of the point subscheme $W$ away from those curves.
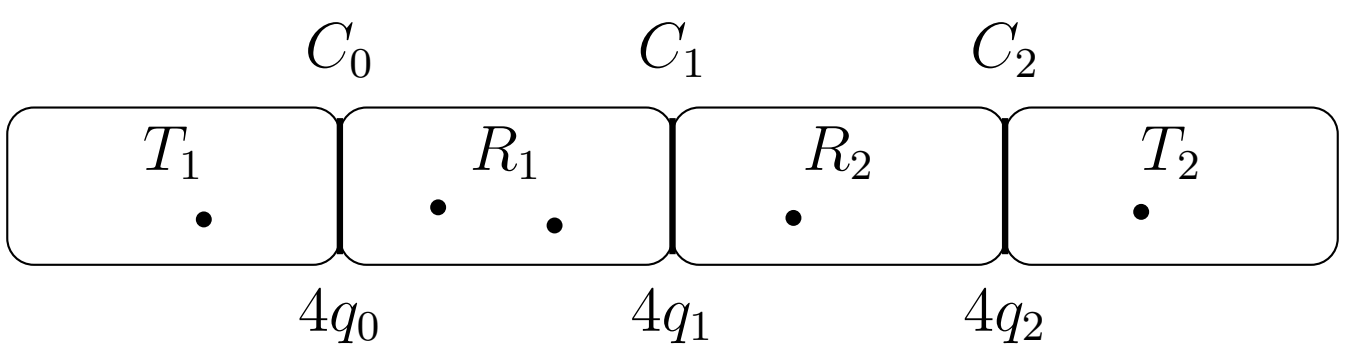
\begin{center}
\begin{tikzpicture}[x=1cm,y=1cm,>=stealth]
 \foreach \x/\lab in {0/T_1,2.2/R_1,4.4/R_2,6.6/T_2}{
   \draw[rounded corners=5pt] (\x,0) rectangle +(2.2,1.04);
   \node at (\x+1.1,.74) {$\lab$};
 }
 \foreach \x/\j in {2.2/0,4.4/1,6.6/2}{
    \draw[very thick] (\x,.07)--(\x,.97);
    \node[above] at (\x,1.10) {$C_{\j}$};
    \node[below] at (\x,-.02) {$4q_{\j}$};
 }
 \foreach \x/\y in {1.3/.30,2.85/.38,3.62/.26,5.2/.31,7.5/.35}
     {\fill (\x,\y) circle (1.5pt);}
\end{tikzpicture}
\end{center}
Here $C_0=T_1\cap R_1$, $C_1=R_1\cap R_2$, and
$C_2=R_2\cap T_2$. Cutting at $C_m$ means partially normalizing
along $C_m$; each of the two resulting surfaces retains one copy
of $C_m$ as its boundary. The stability inequalities below are
tested on each of these surfaces as well as on the complete chain.

Write $h_1$ for the hyperplane on $T_1$, and $\ell,E_1,\ldots,E_{20}$ for the classes on $T_2$. On $R_j$ let $s_j$ be the negative section and $f_j$ a fiber. Then
\[
s_j^2=-4,\quad s_jf_j=1,\quad K_{R_j}=-2s_j,\quad\chi(\OO_{R_j})=-2.
\]
Choose rational numbers
\begin{equation}\label{eq:chamberq}
0<q_0<\cdots<q_k<M
\end{equation}
and use the polarization restrictions
\begin{align}
L|_{T_1}&=q_0h_1,\\
L|_{R_j}&=(q_j-q_{j-1})s_j+4q_jf_j,\\
L|_{T_2}&=(5M-4q_k)\ell-(M-q_k)\sum E_i.
\end{align}
The boundary degrees are $4q_j$. The line bundles glue because the twenty points cut out a divisor of class $5h$ on $C$. They are ample: the bubble expression is a positive combination of the positive section and an ample base pullback; on $T_2$ it is $q_k\ell+(M-q_k)(5\ell-\sum E_i)$, the restriction of the positive product polarization on the blowup embedding used in Lemma~\ref{lem:ample}.

Surface Riemann--Roch gives the following coefficients for structure sheaves:
\begin{equation}\label{eq:polytable}
\begin{array}{c|c|c|c}
 & a&b&c\\ \hline
T_1&2q_0^2&0&2\\
R_j&2(q_j^2-q_{j-1}^2)&4q_{j-1}&-2\\
T_2&\frac52M^2-2q_k^2&4q_k-\frac52M&1.
\end{array}
\end{equation}
A boundary curve has polynomial $4q_jz-2$. If $W$ is a finite subscheme of length $n$, disjoint from the nodes, then $F=I_W$ has
\begin{equation}\label{eq:idealpoly}
P_F=\frac52M^2z^2-\frac52Mz+5-n.
\end{equation}

\begin{theorem}\label{thm:idealchamber}
Assume $n\le N$ and every bubble carries at least one point of $W$, counted with scheme length. In addition to \eqref{eq:chamberq}, when $N>2$ assume
\begin{equation}\label{eq:smallchamber}
1-\frac{q_0^2}{q_k^2}<\frac1{N-2}.
\end{equation}
For $N\le2$ there is no extra inequality. Then $F$ is stable on the full chain and its restriction to both sides of every cut is stable for the induced polarization.
\end{theorem}
\begin{proof}
For the whole chain the slope $b/a$ is $-1/M$. Normalization and subtraction of the cut boundaries give the maximal-subobject slopes
\[
-\frac2{q_j},\qquad
-\frac{2q_j}{q_j^2-q_{i-1}^2},\qquad
\frac{-5M/2}{5M^2/2-2q_{i-1}^2}
\]
for a prefix, an interior bubble interval, and a suffix containing $T_2$, respectively. Each is strictly below $-1/M$. A finite ideal changes only the constant coefficient, so this proves full stability for all point distributions.

A left branch ending at $q_m$ has polynomial
\[
P_L=2q_m^2z^2+2-n_L.
\]
Its proper prefixes and interior intervals have negative linear coefficient. The only remaining patterns are terminal bubble intervals, whose maximal subobject has
\[
a_U=2(q_m^2-q_{i-1}^2),\quad b_U=0,\quad c_U=-n_U.
\]
Here $n_U\ge1$. When $n_L\le2$, the total constant ratio is nonnegative while the subobject ratio is negative. When $n_L>2$, the required inequality is
\[
\frac{a_U}{2q_m^2}<\frac{n_U}{n_L-2}.
\]
The left side is at most $1-q_0^2/q_k^2$, and the right side is at least $1/(N-2)$, so \eqref{eq:smallchamber} proves it strictly.

For a right branch beginning at $q_m$, let
\[
a_R=\frac52M^2-2q_m^2,\qquad b_R=4q_m-\frac52M.
\]
Since $Mb_R+a_R=4Mq_m-2q_m^2>0$, its slope exceeds $-1/M$. A proper bubble prefix has slope $-2/(q_j+q_m)<-1/M$; an interior bubble interval and a suffix containing $T_2$ have slopes below $-1/M$ by the preceding calculation. A branch consisting solely of $T_2$ is stable by rank one and integral support. Disconnected subobjects are handled by adding the residual polynomials. These checks exhaust all generic-rank patterns by Lemma~\ref{lem:maxsub}.
\end{proof}

In the relative point-quotient compactification on a fixed support, an unoccupied bubble has a continuous rescaling stabilizer. Thus every allowed bubble is occupied and $k\le n\le N$. Li--Wu's proper expanded Quot construction applies to this fixed-source problem \cite{LW}; Theorem~\ref{thm:idealchamber} supplies the simultaneous cutting inequalities for its kernels. This is a bounded sector, not an identification of all surface sheaves with point ideals.

\subsection{Extending the rank-one hull on an expanded supporting surface}
We now start with the sheaf on the generic surface. Properness of the expanded Hilbert space in \cite[Theorem~4.14]{LW} gives a support limit after finite base change. The following criterion constructs the sheaf on that limit from a divisor of the generic surface and a computation at the singular points of the total support.

\begin{proposition}[Divisorial extension on a prescribed support]\label{prop:divisorial-hull27}
Work over $R=\CC[[\pi]]$ with fraction field $K$. On a regular semistable expanded family $\mathscr Z/R$ with Cartier components, let $i:\mathscr T\hookrightarrow\mathscr Z$ be a flat, fiberwise effective Cartier divisor and a stable expanded Hilbert point. Assume $\mathscr T$ is integral and normal, with geometrically integral normal generic fiber.

Its closed support is a reduced connected chain $T_0\cup\cdots\cup T_m$, each $T_j$ the integral Cartier surface cut out on the corresponding occupied component of the expansion. Assume successive $T_j$ meet along integral Cartier curves $C_e$ in their joining divisors, and all other components are disjoint from the support.

Let $G_K$ be any rank-one reflexive sheaf on $\mathscr T_K$. Choose a Weil divisor $D_K$ representing it and take its coefficientwise closure $D$ in $\mathscr T$. Assume $\mathscr T$ is regular at the points of its closed fiber, except possibly at finitely many points off the joining curves. At each exceptional point suppose that the completed local ring of $\mathscr T$ and the local divisor class of $D$ satisfy
\begin{gather*}
 A=\CC[[u,v,z,\pi]]/(uv-z(z+\alpha)),\qquad \alpha\in\{0,\pi\},\\
 [D]_A\in\{0,-[V(u,z)],-[V(u,z+\alpha)]\}\subset\mathrm{Cl}(A).
\end{gather*}
Thus when the total support is regular the construction applies to every $G_K$.

For a fixed polarization on the closed expansion use $D_e,B_e$ of \eqref{eq:cutdata} for the extension $i_*\OO_{\mathscr T}(D)$ on its special fiber. Set $\delta_e=[z]D_e$ and $d_e=[z]B_e$. If $\delta_e/d_e\notin\mathbb Z$ at every occupied cut, there is a vertical twist giving an $R$-flat extension of $i_{K*}G_K$ with stable, pure, perfect special fiber. Its branch sheaves are Cohen--Macaulay and their nonzero boundary restrictions are line bundles on $C_e$.

For polarizations satisfying (K1) any two extensions in this class are isomorphic, after scalar rigidification, by the unique extension of their given generic isomorphism.
\end{proposition}
\begin{proof}
Trivialize $G_K$ rationally. The orders at the codimension-one points give $D_K$ and $G_K\simeq\OO_{\mathscr T_K}(D_K)$ by the divisor-class description of rank-one reflexive sheaves \cite[Lemma~31.30.2]{StacksWeil27}. Define $M=\OO_{\mathscr T}(D)$ by the divisorial inequalities at the codimension-one points. Since $\mathscr T$ is normal, this is a coherent reflexive sheaf, and localization at $\pi$ gives $M_K\simeq G_K$. It is $R$-flat because it is torsion-free.

At a regular point of the total support $M$ is invertible: a regular local ring is factorial. At an exceptional point, the two displayed primes of $A$ have ideals $I=(u,z)$ and $I'=(u,z+\alpha)$. For either ideal $J$,
\[
   0\longrightarrow J\longrightarrow A\longrightarrow
     A/J\simeq\CC[[v,\pi]]\longrightarrow0 .
\]
The hypersurface $A$ is Cohen--Macaulay of dimension three. The depth lemma gives $\depth J=3$, and localization of the same sequence gives the $S_2$ condition. Hence $J$ is reflexive, with class the negative of the corresponding prime divisor. By the class hypothesis the completed stalk of $M$ is $A$, $I$, or $I'$. This uses only the divisor $D$ obtained from $D_K$.

Since $A/J$ is $R$-flat, reduction modulo $\pi$ gives
\[
 J/\pi J\ =\ (u,z)\ \subset A/\pi A=\CC[[u,v,z]]/(uv-z^2).
\]
This ideal is rank one and has depth two. The images of $u,z$ are independent modulo the maximal ideal times $(u,z)$, so its minimal number of generators is two. The construction therefore includes a noninvertible special hull. These calculations descend from completion by faithful flatness.

Near a join $M$ is invertible on $\mathscr T$. On each branch its restriction is an invertible sheaf on the integral Cartier surface $T_j$ near that join. The equation of the attaching divisor is regular on $T_j$ and cuts out $C_e$, so the derived restriction of the branch sheaf is a line bundle on $C_e$ in degree zero. On the special expanded fiber the pushforward $i_{0*}M_0$ is perfect there by the Cartier resolution for $\mathscr T_0$.

Away from the joins the expanded fiber is smooth. Its branch sheaf there has depth two at a closed point by the preceding computation or by local invertibility on the Cartier surface. This proves the assertions on depth, purity and perfectness.

For each cut let $k_e=-\lfloor\delta_e/d_e\rfloor$. Tensor $i_*M$ by the product of the line bundles $\OO_{\mathscr Z}(-\sum_{j\le e}Z_j)^{k_e}$ of \eqref{eq:verticaltwist}, with indices taken in the ordered ambient chain. Each is trivial on the generic fiber. Lemma~\ref{lem:twists} gives independent changes at the different cuts. The new linear coefficients satisfy
\[
                 0<\delta_e+k_ed_e<d_e .
\]
Consequently the two polynomial inequalities of Theorem~\ref{thm:cutcriterion} are strict at every cut. They prove stability. The generic sheaf is stable by rank one on its geometrically integral support. For later use, removing a finite-length quotient away from the joining divisors changes only the constant terms of the branch polynomials.

The intrinsic support quotient of this extension is $\OO_{\mathscr Z}\twoheadrightarrow\OO_{\mathscr T}$. Indeed $M$ is faithful on the integral total support. On the closed surface $M_0$ has rank one at every generic point. Its annihilator lies in the intersection of the minimal primes, which is zero since the closed surface is reduced. The generic fiber has the same property.

Proposition~\ref{prop:support-separated27} therefore applies to any two such extensions, using (K1) for the polarized sheaf comparison.
\end{proof}

\section{Derived descent and the matching functor}\label{sec:derived}
\subsection{Derived algebraic charts}
For a projective flat family $W\to T$, let $\mathbf P_{W/T}$ be the derived moduli of perfect complexes, restricted to the bounded loci whose classical fibers are the sheaves under consideration. The sheaf and admissibility loci are the derived open subfunctors determined by their classical open conditions. Their classical truncations are the usual sheaf stacks. The universal complex exists on the unrigidified stack.

\begin{proposition}\label{prop:derivedcharts}
On the bounded charts of Proposition~\ref{prop:algebraic}, these derived sheaf functors are locally geometric and locally of finite presentation. Their relative tangent complex at a family $\FF$ is
\begin{equation}\label{eq:derivedtangent}
\RR\pi_*\mathcal RHom(\FF,\FF)[1].
\end{equation}
The construction is compatible with derived base change and with isomorphisms of the projective target families.
\end{proposition}
\begin{proof}
On an affine chart $T=\Spec R$ put $W_A=W\times_T^{\mathbf R}\Spec A$. For a connective derived $R$-algebra $A$, the open subfunctor of $\mathbf P_{W/T}$ has values
\[
 \{F\in\perf(W_A):F|_{W_{\pi_0 A}}\text{ is a flat sheaf with the prescribed fibers}\}^{\simeq}.
\]
For discrete $A$ its sheaf-flatness condition is equivalently
\[
 \mathcal H^j(F\otimes_A^{\LL} B)=0
 \qquad(j\ne0,\ B\text{ a discrete }A\text{-algebra}).
\]
This identifies the classical truncation. Apply \cite[Theorem~4.12 and Example~4.15]{Pridham} to the projective flat morphism $W\to T$. For a test family $p:W_B\to\Spec B$, fiberwise free resolutions lift by flatness and Nakayama as in Lemma~\ref{lem:perfectfamilywithoutnormalnew}. The failure locus $D\subset W_B$ is closed, so its perfect-fiber locus on the base is $\Spec B\setminus p(D)$, open by properness. Intersect it with the fiber-condition opens of Proposition~\ref{prop:algebraic}.

For the finiteness and the deformation complex, set
\[
 C_A=\RR p_{A*}(F^\vee\otimes^{\LL}F).
\]
The morphism $p_A$ is the base change of a proper flat finitely
presented morphism of finite Tor dimension; proper perfect pushforward
gives $C_A\in\perf(A)$. For an $A$-module complex $I$ and a
derived $R$-algebra map $A\to B$, projection formula and derived
base change give
\[
 \begin{aligned}
 \RHom_{W_A}(F,F\otimes^{\LL}p_A^*I)&\simeq C_A\otimes_A^{\LL}I,\\
 C_A\otimes_A^{\LL}B
 &\simeq \RR p_{B*}(F_B^\vee\otimes^{\LL}F_B)=C_B .
 \end{aligned}
\]
Over a discrete noetherian base the cohomology modules of $C_A$ are
thus finite, and their versions with variable coefficients commute
with filtered colimits. The classical stack is locally of finite
presentation on the bounded charts. Formal effectiveness for sheaves
is Grothendieck existence on the projective $W$ over a complete
noetherian local ring. For endomorphisms the corresponding
completion identity is
\[
 C_A\simeq\mathbf R\!\varprojlim_n
     (C_A\otimes_A^{\LL}A/\mathfrak m^{n+1}),
\]
since $C_A$ is perfect and $A$ complete. These are the finiteness and
effectiveness conditions of the cited representability criterion.

The space of lifts through a split square-zero extension $A\oplus I$,
with specified identification of the reduction, is the mapping space
associated to $C_A\otimes_A^{\LL}I[1]$ (its homotopy group in degree
$j\ge0$ is $H^{1-j}(C_A\otimes_A^{\LL}I)$). Thus the relative
tangent complex is $C_A[1]$, which is \eqref{eq:derivedtangent}.
In particular on an ordinary geometric fiber the infinitesimal
automorphisms, deformations and obstruction group are respectively
$\Ext^0(F,F),\Ext^1(F,F),\Ext^2(F,F)$.

The second base-change
identity identifies the tangent complexes on every chart overlap.
For a target isomorphism $g$ one has
$(g^*F)^\vee\otimes g^*F\simeq g^*(F^\vee\otimes F)$; these
identifications compose under a second isomorphism and hence descend.
\end{proof}

\subsection{Descent at a node}
We use the following derived excision theorem, in its symmetric-monoidal form.

\begin{theorem}\label{thm:Milnorinput}
Let
\[
\begin{tikzcd}
A\arrow[r]\arrow[d]&A_1\arrow[d]\\
A_2\arrow[r]&A_{12}
\end{tikzcd}
\]
be a homotopy pullback of connective commutative derived rings, with $\pi_0(A_2)\to\pi_0(A_{12})$ surjective. Then
\[
\perf(A)\simeq\perf(A_1)\times^h_{\perf(A_{12})}\perf(A_2).
\]
\end{theorem}
This is the perfect-complex consequence of Chough's descent theorem \cite[Theorem 1.3 and Remark 1.7]{Chough}. One can restrict its symmetric-monoidal equivalence on almost perfect complexes to dualizable objects. We use the equivalence itself, not a K-theoretic excision statement.

\begin{theorem}\label{thm:derivedglue}
Let $W=W_L\cup_DW_R$ be a node-marked expanded family over $T$. For every connective derived $T$-scheme $U$ there is a natural equivalence
\begin{equation}\label{eq:Perfglue}
\perf(W\times_T^{\RR}U)
\simeq
\perf(W_L\times_T^{\RR}U)
\times^h_{\perf(D\times_T^{\RR}U)}
\perf(W_R\times_T^{\RR}U).
\end{equation}
It is compatible with further cuts, all derived base changes, and the expansion groupoid. It therefore identifies the actual derived perfect-complex moduli functors on the two sides.
\end{theorem}
\begin{proof}
On a compatible affine chart, write $A$ for the parameter ring and $B=B_L\times_RB_R$ for the nodal ring. The normalization sequence and its base change to a connective derived $A$-algebra $A'$ give fiber sequences of modules
\[
 B\longrightarrow B_L\oplus B_R\longrightarrow R,\qquad
 B'\longrightarrow B'_L\oplus B'_R\longrightarrow R',
 \qquad (-)'=(-)\otimes_A^{\LL}A'.
\]
The ring square is a homotopy pullback since the forgetful functor on derived rings creates limits. The maps $\pi_0B'_i\to\pi_0R'$ are surjective: on $\pi_0$ they are obtained by tensoring the surjections $B_i\to R$ with $\pi_0A'$. Theorem~\ref{thm:Milnorinput} gives \eqref{eq:Perfglue} on the chart.

Taking the limit over an affine hypercover gives \eqref{eq:Perfglue} on $W$, since perfect complexes satisfy descent and limits commute. All arrows are restrictions, hence commute with further base changes and with target isomorphisms. Iterating at several joins takes the limit of a single branch-and-boundary diagram, independently of the order of the cuts. The maximal $\infty$-groupoid functor preserves these limits, giving the assertion for moduli functors.
\end{proof}

\begin{corollary}\label{cor:cotangentglue}
On the sheaf loci of Proposition~\ref{prop:familyglue}, the derived equivalence induces the classical matching equivalence and identifies the cotangent complexes, universal Atiyah classes, and obstruction maps for arbitrary square-zero derived extensions. Its tangent triangle is \eqref{eq:Endtriangle}.
\end{corollary}
\begin{proof}
For normal branch sheaves, $\LL i^*F=i^*F$, and the inverse to restriction on classical families is the kernel in \eqref{eq:familyglue}. Write $\mathbf M$ for either presentation of the derived stack. For a square-zero extension $A'\to A$ with ideal $J$ and $x\in\mathbf M(A)$, the space of lifts is
\[
       \fib_x\bigl(\mathbf M(A')\longrightarrow\mathbf M(A)\bigr).
\]
The equivalence \eqref{eq:Perfglue} identifies these spaces functorially in the extension and in $x$. Their representing cotangent complexes and the universal square-zero obstruction maps therefore agree. Under the equivalence the universal complexes correspond by restriction, hence so do their Atiyah classes by functoriality. Differentiating gives \eqref{eq:Endtriangle}.
\end{proof}

In particular the full perfect-complex boundary covers restrictions with nonzero degree $-1$ cohomology, impure ordinary restrictions, and sectors with zero branch objects. The passage to an ordinary smooth boundary moduli space is a further construction, carried out next.

\section{Scalar rigidification and the determinant construction}\label{sec:reduction}
For rank-zero sheaves, scalar automorphisms and the higher trace directions must be treated separately. This section fixes the conventions used by every obstruction theory in the paper.

\subsection{The moduli construction}
Let $W\to T$ be a complete expanded Calabi--Yau family, a relative root chain, or the constant K3 boundary family. On the locus of nonzero objects with no negative self-$\Ext$ on geometric fibers, let $\mathbf P_{W/T}^{\rig}$ denote rigidification by the central scalar group. This includes the flat sheaf loci. For an arbitrary boundary complex one can first glue and then rigidify the resulting sheaf, as in Section~\ref{sec:diagonal}. Let $\mathbf{Pic}_{W/T}^{\rig}$ be the derived moduli of line bundles modulo their scalar automorphisms, and let $\Pic_{W/T}^{\rig}=t_0\mathbf{Pic}_{W/T}^{\rig}$ denote its classical truncation.

The determinant defines a morphism to the derived Picard stack, and hence to its rigidification. Define
\begin{equation}\label{eq:reducedmoduli}
\mathbf M_{W/T}
=
\mathbf P_{W/T}^{\rig}
\times^{\RR}_{\mathbf{Pic}_{W/T}^{\rig}}
\Pic_{W/T}^{\rig},
\end{equation}
restricted to the desired numerical determinant classes and sheaf loci.

This fixes the higher derived variation of the determinant but allows the stated classical determinant classes. It does not add a chosen trivialization of $\det F$. In particular, for rank zero it does not introduce an additional determinant-framing $\GG$. The classical truncation of \eqref{eq:reducedmoduli} is the original scalar-rigidified sheaf moduli.

When testing against a derived parameter, the classical Picard locus in \eqref{eq:reducedmoduli} is pulled back as an etale space over the original parameter base; it is not retruncated after derived base change.

For the complete expansions and relative root chains of the first degeneration, $H^1(\OO)=H^2(\OO)=0$ by Proposition~\ref{prop:structurecohom}. Their classical rigidified Picard spaces are \'etale over the parameter base. For the constant K3 boundary, the relative space is the constant $\Pic(S)$ and $H^1(S,\OO_S)=0$. Thus an actual determinant line bundle gives an \'etale-local section of the classical Picard locus in \eqref{eq:reducedmoduli}. For the component $W$ of the second degeneration, the separate two-boundary calculation applies; see Section~\ref{sec:partI-pencil-example}.

\begin{lemma}\label{lem:tangentreduction}
Let $\pi:W\times_T\mathbf P\to\mathbf P$ be the universal projection, and put
\begin{align}
K_W&=\RR\pi_*\mathcal RHom(\FF,\FF),\\
C_W&=\cofib(\OO\xrightarrow{\id}K_W),\\
Q_W&=\cofib(\OO\longrightarrow\RR\pi_*\OO_W).
\end{align}
The tangent complex of \eqref{eq:reducedmoduli}, relative to $T$, is
\begin{equation}\label{eq:redtangent}
C_W^{\red}[1],\qquad
C_W^{\red}=\fib(C_W\xrightarrow{\tr}Q_W).
\end{equation}
\end{lemma}
\begin{proof}
On a chart carrying $F$ the scalar units and determinant fit into a
morphism of cofiber sequences
\[
\begin{tikzcd}
 \OO[1]\arrow[r,"{\id_F}"]\arrow[d,"{\rk F}"']&
 K_W[1]\arrow[r]\arrow[d,"\tr"']& C_W[1]\arrow[d]\\
 \OO[1]\arrow[r,"1"']& (R\pi_*\OO_W)[1]\arrow[r]&Q_W[1].
\end{tikzcd}
\]
The tangent of the perfect-complex moduli is $K_W[1]$. The smooth central scalar gerbe contributes $\OO[1]$ with differential the identity endomorphism, so rigidification gives $C_W[1]$. The lower cofiber is the tangent of the rigidified derived Picard stack, and the middle vertical arrow is trace by \cite[Proposition~3.2]{STV}. Taking the tangent of the derived pullback \eqref{eq:reducedmoduli}, and using the zero relative tangent of the classical Picard space, gives \eqref{eq:redtangent}.

These assertions can be made on a cover carrying a universal sheaf. Rigidification is the stackification of quotienting the central scalar subgroup; locally its relative tangent is $\OO[1]$. The complexes and the trace maps are equivariant on overlaps and hence descend. No scalar-trace idempotent is used. Indeed $\tr(\id_F)=\rk(F)=0$ for these sheaves.
\end{proof}

\subsection{Compatibility with matching}
\begin{proposition}\label{prop:redglue}
Suppose a marked family $W=W_L\cup_DW_R$ has connected branches and boundary with $H^0(\OO)=\CC$. On each geometric fiber assume that the two branch sheaves and their common boundary complex $E$ are nonzero and that $\Ext_D^j(E,E)=0$ for $j<0$. Assume the classical rigidified Picard spaces are \'etale over the parameter base, as above. Then derived descent is compatible with \eqref{eq:reducedmoduli}:
\begin{equation}\label{eq:reducedglue}
\mathbf M_W\simeq
\mathbf M_{W_L}\times^{\RR}_{\mathbf M_D}\mathbf M_{W_R}
\end{equation}
on the corresponding matched sheaf loci and compatible numerical determinant classes.
\end{proposition}
\begin{proof}
First retain all scalar automorphisms. The equivalence of Theorem~\ref{thm:derivedglue} is symmetric monoidal, so it also gives descent for invertible objects and determinants.

To check scalars on a family of boundary complexes, let $K_D=\RR\pi_*\mathcal RHom(E,E)$. On a geometric fiber $H^{<0}(K_D)=0$ and $1\mapsto\id_E$ is injective in degree zero. Hence the perfect complex $\cofib(\OO\to K_D)$ has Tor amplitude bounded below by zero. After any classical base change the unit still injects in degree zero. Scalar units therefore act faithfully on families, also when $E$ has two cohomology sheaves. For a simple glued sheaf the required negative-$\Ext$ vanishing follows from Corollary~\ref{cor:negextstablegluednew}.

Choose representatives fppf-locally for the scalar gerbes. On a marked union the change of matching by scalar units is governed by the exact sequence
\[
 1\longrightarrow\GG\xrightarrow{\mathrm{diag}}\GG^2
       \xrightarrow{(a,b)\mapsto ab^{-1}}\GG\longrightarrow1 .
\]
Consequently the branch representatives exhaust all scalar changes of the matching, leaving precisely the common scalar on the glued object. A chain of $m$ joins uses the same exact sequence with middle and right terms $\GG^{m+1},\GG^m$ and incidence map $(a_i)\mapsto (a_{i-1}a_i^{-1})_i$. These group-sheaf identities descend over the parameter.

The same calculation applies to line bundles. For the classical rigidified Picard loci the ordinary matching product computes their derived product by the assumed \'etaleness. In the three compatible determinant diagrams (global, left and right), take the pullbacks defining \eqref{eq:reducedmoduli}. Commutativity of finite limits with the matching product now yields \eqref{eq:reducedglue}, compatibly with base change.

\end{proof}

\begin{remark}\label{rem:zerosectors}
If a branch or boundary object is zero, separate scalar rigidifications of that zero object are not used. In the full matching theorem one first forms the unrigidified derived matching functor, imposes the determinant condition for the glued object, restricts to globally simple sheaves, and rigidifies its single scalar group. Theorem~\ref{thm:derivedglue} still identifies that construction with the global moduli.

Proposition~\ref{prop:redglue} is the additional product description on the nonzero boundary sectors needed for the Jacobian formula.
\end{remark}

\section{The obstruction theory of sheaves on the universal expanded threefold}\label{sec:POT}
On the sheaf stack of Hypothesis~\ref{hyp:compact}, we construct the obstruction theory relative to the stack of expansions. The universal family carries the sheaf and its relative Atiyah class even where the components of the fibers have no common indexing.

\subsection{Perfect pushforward and duality}
\begin{lemma}\label{lem:universalEnd}
Let $T\to\EE$ be a finite-type expansion chart and let $V\to T$ be a bounded chart of flat pure sheaves with perfect fibers. On a cover carrying its universal sheaf, write
\[
\pi_V:\XX_T\times_TV\longrightarrow V,
\qquad K_V=\RR\pi_{V*}\mathcal RHom(\FF,\FF).
\]
Then $K_V$ is perfect and commutes with derived base change. At a simple sheaf on a complete expanded fiber, its cohomology vanishes outside degrees $0,1,2,3$, and
\[
\Ext^0(F,F)=\CC,
\qquad \Ext^3(F,F)=\CC.
\]
The degree-three trace to $H^3(\OO)$ is an isomorphism.
\end{lemma}
\begin{proof}
By Proposition~\ref{prop:familyperfect} (or the fiberwise resolution argument of Lemma~\ref{lem:perfectfamilywithoutnormalnew} in the absence of admissibility), flatness and fiberwise perfectness give $\FF\in\perf(\XX_T\times_TV)$. The map $\pi_V$ is proper and perfect: it is projective and flat of finite presentation. Thus
\[
 K_V=R\pi_{V*}(\FF^\vee\otimes^{\LL}\FF)\in\perf(V),\qquad
 Lg^*K_V\simeq K_{V'}\quad(g:V'\to V)
\]
by perfect pushforward and derived base change.

On a complete expanded fiber $Z$, $\omega_Z\simeq\OO_Z$. Duality yields
$\Ext^i_Z(F,F)\simeq\Ext^{3-i}_Z(F,F)^\vee$. Negative Ext groups vanish for sheaves, so $H^i(K_V|_{[F]})=0$ for $i\notin[0,3]$. Moreover
\[
    \bigl(\Ext^3_Z(F,F)\xrightarrow{\tr}H^3(Z,\OO_Z)\bigr)^\vee
      =\bigl(H^0(Z,\OO_Z)\xrightarrow{1\mapsto\id_F}\End_Z(F)\bigr).
\]
The right side is an isomorphism of lines by simplicity. Relative duality gives the same identification on families.
\end{proof}

\begin{theorem}\label{thm:ambientPOT}
For the complete expanded Calabi--Yau family of Section~\ref{sec:geometry} with the cohomology of Proposition~\ref{prop:structurecohom}, let $p:\MM\to\EE$ be the scalar-rigidified moduli of simple perfect pure sheaves satisfying the algebraic hypotheses of Section~\ref{sec:compactification}. Its reduced derived enhancement of Section~\ref{sec:reduction} is quasi-smooth relative to $\EE$. Its canonical cotangent morphism gives a relative perfect obstruction theory
\begin{equation}\label{eq:ambientPOT}
E_{\MM/\EE}=\bigl(C_{\XX_{\EE}}^{\red}\bigr)^\vee[-1]
\longrightarrow L_{\MM/\EE}
\end{equation}
of amplitude $[-1,0]$ and virtual rank zero. At a sheaf $F$ its tangent and obstruction spaces are, respectively,
\begin{equation}\label{eq:ambientspaces}
\Ext^1(F,F),\qquad\Ext^2(F,F).
\end{equation}
It agrees with the stable-sheaf Donaldson--Thomas obstruction theory on the smooth quintic fiber.
\end{theorem}
\begin{proof}
Let $K,C,Q$ be as in Lemma~\ref{lem:tangentreduction}. At a geometric point $x=[F]$, Proposition~\ref{prop:structurecohom} and Lemma~\ref{lem:universalEnd} give
\[
 Q_x\simeq H^3(\OO_Z)[-3],\qquad
 H^0(K_x)=\CC\id_F,\qquad H^3(C_x)\xrightarrow[\tr]{\sim}H^3(Q_x).
\]
In the triangles
$\OO\to K\to C$ and $C^{\red}\to C\to Q$ the scalar unit and the top trace are therefore isomorphisms at the respective ends. The long exact sequences give
\[
 H^i(C^{\red}_x)=
 \begin{cases}\Ext_Z^i(F,F),&i=1,2,\\0,&\text{otherwise}.\end{cases}
\]
Perfectness and derived base change convert the fiberwise interval into Tor amplitude $[1,2]$. The reduced derived moduli has relative tangent $C^{\red}[1]$ by Lemma~\ref{lem:tangentreduction}, hence is quasi-smooth over $\EE$. Its cotangent map to the classical truncation is \eqref{eq:ambientPOT}, of amplitude $[-1,0]$.

On a universal chart this morphism is the Atiyah-class obstruction morphism. For a square-zero extension with ideal $J$ it gives the obstruction $(1_F\otimes\mathrm{KS})\circ\mathrm{At}(F)$ in $\Ext^2(F,F\otimes J)$ and, when that class vanishes, the torsor of lifts under $\Ext^1(F,F\otimes J)$ \cite{HT,STV}. Restriction preserves the unit, determinant and trace. Scalar conjugation acts identically on $\FF^\vee\otimes^{\LL}\FF$, so the morphism descends through the gerbe and the target groupoid.

On the smooth quintic the same unit and trace reductions give the usual stable-sheaf obstruction morphism, since $H^1(\OO_X)=H^2(\OO_X)=0$. Duality identifies $\Ext^2(F,F)=\Ext^1(F,F)^\vee$ both there and on a complete expanded Calabi--Yau fiber, giving virtual rank zero.
\end{proof}

\begin{remark}\label{rem:relativebase}
The relative theory is over $\EE$, not simply over the original degeneration curve. This distinction is essential at an expanded fiber. The base stack is smooth of dimension one, and its tangent includes both smoothing directions and infinitesimal rescaling of the bubbles.

The complete cotangent construction, followed by the classical truncation, supplies the absolute obstruction theory on the finite-stabilizer moduli; its virtual dimension is one. One must not obtain an absolute theory by retaining all smoothing coordinates while omitting the Lie algebra of the rescaling torus.

More explicitly, on a chart $[U/G]$ pull back the derived sheaf moduli to $U$. Its tangent has amplitude $[0,1]$ relative to $U$, and adjoining $T_U$ retains that amplitude. The tangent of the quotient is the cofiber of $\operatorname{Lie}(G)\otimes\OO$ mapping to this chart tangent. Its possible degree $-1$ cohomology is the infinitesimal intrinsic stabilizer. It vanishes under the finite-stabilizer condition in characteristic zero.

The absolute tangent therefore still has amplitude $[0,1]$. Its rank is the relative rank plus $\dim U-\dim G$. This proves the asserted absolute amplitude and dimension, including the rescaling contribution.
\end{remark}

\subsection{The virtual class over the expansion stack}
The expansion stack has a stratification by quotient stacks, and the bounded parts used here have the intersection theory of \cite{Kresch}. The morphism $p$ is of Deligne--Mumford type after scalar rigidification: over a fixed target a stable sheaf has no remaining automorphisms. The relative perfect obstruction theory defines the bivariant virtual pullback of \cite{Mano}, which we denote
\begin{equation}\label{eq:pvirtual}
p^!:A_*(\EE)_{\QQ}\longrightarrow A_*(\MM)_{\QQ}.
\end{equation}
Its degree is zero. Define
\begin{equation}\label{eq:totalvirtual}
[\MM]^{\vir}=p^![\EE]\in A_1(\MM)_{\QQ}.
\end{equation}
Since $\EE$ is smooth, this is the usual absolute virtual class induced from the relative theory and the smooth base. Its formation does not require flatness of $\MM\to\EE$ or of $\MM\to C$.

\begin{proposition}\label{prop:POTbasechange}
The obstruction theory \eqref{eq:ambientPOT} is compatible with expansion-chart changes and with derived base change to a marked smoothing face. On such a face it is the obstruction theory of the derived matching functor of Theorem~\ref{thm:derivedglue}, with the global determinant and scalar conventions of Section~\ref{sec:reduction}.
\end{proposition}
\begin{proof}
The derived sheaf functor is constructed on the universal family and is compatible with derived pullback. The perfect endomorphism complex and the determinant map have the same property. The classical Picard locus used in \eqref{eq:reducedmoduli} is etale over the parameter base; its derived pullback is the corresponding locus in the base-changed construction. Thus the reduced derived stack and its cotangent morphism commute with base change.

For $g:T'\to T$, write $p:Z\to T$ for the expanded family,
$p':Z'\to T'$ for its pullback and $h:Z'\to Z$. If $F'=Lh^*F$,
trace gives the commutative square
\[
\begin{tikzcd}
 Lg^*Rp_*R\mathcal Hom(F,F) \arrow[r,"Lg^*\tr"]\arrow[d,"\sim"']&
 Lg^*Rp_*\OO_Z \arrow[d,"\sim"]\\
 Rp'_*R\mathcal Hom(F',F')\arrow[r,"\tr"']&Rp'_*\OO_{Z'} .
\end{tikzcd}
\]
Pullback takes $1\mapsto\id_F$ to $1\mapsto\id_{F'}$. Taking
cofibers of the scalar units and then the fiber of trace yields
$Lg^*C_Z^{\red}\simeq C_{Z'}^{\red}$.

On a marked face normalize only the marked node. The two resulting branches are relative expanded families; their remaining nodes are allowed to smooth. Theorem~\ref{thm:derivedglue} applies before specializing any of those other parameters, and identifies this entire family of deformation problems.

The determinant and scalar operations commute with the equivalence as proved in Proposition~\ref{prop:redglue}, or are made globally after matching in the zero-object sectors of Remark~\ref{rem:zerosectors}. Applying the cotangent functor identifies the obstruction theories and all connecting morphisms. Equivariance follows from naturality under isomorphisms of expanded families.
\end{proof}

\begin{corollary}\label{cor:twistPOT}
On a family where a vertical invertible sheaf $\mathcal L$ is defined, tensoring by $\mathcal L$ identifies the relative obstruction theories of the corresponding stable-sheaf loci. For a twist trivial on the generic fiber, the generic Donaldson--Thomas theory is unchanged.
\end{corollary}
\begin{proof}
Tensoring by $\mathcal L$ is an equivalence of the derived deformation functors, with inverse tensoring by $\mathcal L^{-1}$. In a local trivialization this identification is
\[
 \RHom(F,F)\xrightarrow{a\mapsto a\otimes\id_{\mathcal L}}
    \RHom(F\otimes\mathcal L,F\otimes\mathcal L),\qquad
 \tr(a\otimes\id_{\mathcal L})=\tr(a).
\]
These maps preserve the identity endomorphism and agree on overlaps. For the rank-zero objects the determinant shift is $\mathcal L^{\rk F}=\OO$, and in general the compatible prescribed determinant class is translated by that shift. The reduced derived constructions and their cotangent morphisms are therefore identified. The family statement applies over each base carrying this line bundle.
\end{proof}

\section{Relative branches and the K3 boundary}\label{sec:relative}
The component invariants in the diagonal formula require independently defined relative obstruction theories.  Theorem~\ref{thm:branchindexnew} treats a general derived boundary. Here we construct its two-term theory where negative self-Ext vanishes and treat zero boundary using the relative determinant trivialization.

Let $W$ be an expansion of one of the pairs $(Y_i,S)$, and let $D\cong S$ be its outer relative divisor. Thus $W$ consists of the root together with any finite chain of bubbles, and
\begin{equation}\label{eq:pairCY}
\omega_W\cong\OO_W(-D),\qquad
H^0(W,\OO_W)=\CC,\qquad H^{>0}(W,\OO_W)=0.
\end{equation}
The first assertion follows by logarithmic adjunction, and the cohomology statement is Proposition~\ref{prop:structurecohom}. The relative deformation complex is $\RHom_W(F,F)[1]$ on the whole chain.

\begin{lemma}\label{lem:relativeExt3}
Let $F$ be a simple perfect sheaf on $W$, normal to $D$, with nonzero restriction $G=F|_D$. Then
\[
\Hom_W(F,F(-D))=0,\qquad\Ext_W^3(F,F)=0.
\]
\end{lemma}
\begin{proof}
Let $j:F(-D)\hookrightarrow F$ be multiplication by the Cartier equation; it is injective by normality. For $a:F\to F(-D)$ simplicity gives $ja=\lambda\id_F$. Restriction to $D$ gives
\[
     0=(j|_D)(a|_D)=\lambda\id_G .
\]
Since $G\ne0$, $\lambda=0$ and $j$ injective implies $a=0$. By \eqref{eq:pairCY}, $\Ext_W^3(F,F)^\vee=\Hom_W(F,F(-D))=0$.
\end{proof}

\begin{proposition}\label{prop:relativePOT}
On an open locus of simple perfect sheaves with $E=\LL i^*F\ne0$ and $\Ext^{-1}_S(E,E)=0$, the bounded relative sheaf moduli has a canonical relative two-term obstruction theory over the stack of relative expansions. Its tangent complex is $C_W[1]$, where $C_W=\cofib(\OO\to K_W)$. No additional higher determinant-trace term is present. The construction is compatible with outer-boundary evaluation and with smoothing the internal nodes.
\end{proposition}
\begin{proof}
Perfectness and base change follow from Lemma~\ref{lem:perfectfamilywithoutnormalnew}. Equation~\eqref{eq:pairCY} makes $Q_W=0$ in Lemma~\ref{lem:tangentreduction}. Simplicity removes degree zero, and \eqref{eq:branchthirdnew} gives $\Ext^3_W(F,F)=0$. Gorenstein duality and the absence of negative self-Ext give no higher cohomology. On every geometric fiber the reduction therefore gives
\[
 H^j(C_W[1])=
 \begin{cases}\Ext^1_W(F,F)&j=0,\\ \Ext^2_W(F,F)&j=1,\\0&\text{otherwise}.\end{cases}
\]
Perfectness gives $C_W[1]$ Tor-amplitude $[0,1]$. Its dual maps to the cotangent complex by the reduced derived construction, giving the two-term obstruction theory. Compatibility follows from the natural universal restriction functor and Proposition~\ref{prop:POTbasechange}.
\end{proof}

\begin{proposition}\label{prop:relativevdim}
Under the hypotheses of Proposition~\ref{prop:relativePOT}, assume that $G=\LL i^*F$ is a stable pure curve sheaf on $D=S$ of Mukai vector $\lambda$. Then
\begin{equation}\label{eq:relativevdim}
\chi_W(F,F)=\frac12\chi_S(G,G)
=-\frac12\langle\lambda,\lambda\rangle,
\qquad
\vdim M_W^{\rel}=1+\frac12\langle\lambda,\lambda\rangle.
\end{equation}
In particular the relative virtual dimension is three for $\lambda=(0,h,\chi)$, independently of the number of bubbles in $W$.
\end{proposition}
\begin{proof}
Tensor the Cartier sequence of $D$ with $F^\vee\otimes^{\LL}F$. Since $\omega_W=\OO_W(-D)$ and $\LL i^*F=G$, applying $\RG$ gives
\begin{equation}\label{eq:relativeEnd}
 \RHom_W(F,F\otimes\omega_W)\longrightarrow\RHom_W(F,F)
       \longrightarrow\RHom_S(G,G)\longrightarrow .
\end{equation}
Serre duality on $W$ and Riemann--Roch on $S$ now give
\[
 \chi_W(F,F\otimes\omega_W)=-\chi_W(F,F),\qquad
 2\chi_W(F,F)=\chi_S(G,G)=-\langle\lambda,\lambda\rangle .
\]
Since $\ext^0_W(F,F)=1$ and $\ext^3_W(F,F)=0$, subtracting the middle Ext dimensions gives
$\ext^1_W-\ext^2_W=1-\chi_W(F,F)$. The base chart $[\AA^k/\GG^k]$ has dimension zero, so this is also the absolute virtual dimension. For $(0,h,\chi)$ it equals $1+h^2/2=3$.
\end{proof}

\subsection{General stable boundary spaces}
\begin{proposition}\label{prop:generalK3boundary}
Let $\Pic(S)=\mathbb Zh$, $h^2=4$, and let $G$ be a pure $h$-semistable one-dimensional sheaf of Mukai vector $\lambda=(0,bh,\chi_S)$, $b>0$.  If $\gcd(b,\chi_S)=1$, then $G$ is stable.  The corresponding moduli $J_{b,\chi_S}$, when nonempty, is a smooth projective variety of dimension $4b^2+2$.  For a simple admissible relative branch sheaf $F$ with this boundary on an anticanonical root chain, its virtual dimension is $1+2b^2$.  Therefore
\[
 (1+2b^2)+(1+2b^2)-(4b^2+2)=0.
\]
Purity alone does not imply the semistability assumption when $b>1$.
\end{proposition}
\begin{proof}
For a nonzero subsheaf of class $(0,b'h,\chi')$ equality of reduced polynomials reads
\[
       \frac{\chi'}{4b'}=\frac{\chi_S}{4b}
       \quad\Longrightarrow\quad b\mid b'
       \quad\Longrightarrow\quad b'=b
\]
since $\gcd(b,\chi_S)=1$ and $0<b'\le b$. A proper subsheaf with $b'=b$ has a positive-length quotient, so $\chi'<\chi_S$. Thus every semistable sheaf is stable and the Simpson moduli is projective.

The dual of $\Ext^2_S(G,G)\xrightarrow{\tr}H^2(S,\OO_S)$ is the scalar unit, an isomorphism by stability. Hence the trace-free obstruction vanishes, and K3 Riemann--Roch gives
\[
 \dim J_{b,\chi_S}=\ext^1_S(G,G)=2-\chi_S(G,G)
                 =2+\langle\lambda,\lambda\rangle=2+4b^2
\]
at every point \cite{Mukai,HT}. For an admissible branch the boundary is this nonzero sheaf, so $\Ext^{-1}_S(G,G)=0$. Proposition~\ref{prop:relativevdim} gives $\vdim M_W^{\rel}=1+2b^2$, proving the intersection dimension formula.
\end{proof}

If $G$ is strictly semistable, its coarse moduli generally does not have this smooth diagonal and its scalar rigidification can still have positive-dimensional stabilizers.  The matching problem must then retain the boundary stack.  In a proper simple stable expanded moduli the full matching formula \eqref{eq:mainfull} remains the applicable conditional conclusion; it does not become an ordinary Gysin product of arbitrary coarse semistable spaces.

\subsection{The ordinary boundary moduli space}
Assume now that $\Pic(S)=\mathbb Zh$. Fix an integer $\chi$ and put
\[
J_\chi=M_h(S,(0,h,\chi)).
\]
This notation denotes the moduli space of stable pure one-dimensional sheaves, with its classical scheme structure.

\begin{proposition}\label{prop:Jacobian}
Every pure one-dimensional sheaf on $S$ with Mukai vector $(0,h,\chi)$ is stable. The space $J_\chi$ is a smooth projective sixfold, and its support morphism
\begin{equation}\label{eq:Jacobianfibration}
J_\chi\longrightarrow|h|\cong\PP^3
\end{equation}
has fibers the compactified Jacobians of the integral curves in $|h|$. Its derived enhancement obtained from \eqref{eq:reducedmoduli} is the ordinary smooth space $J_\chi$.
\end{proposition}
\begin{proof}
For an effective cycle of class $h$, writing its irreducible components as $a_jh$ with multiplicities $m_j$ gives $\sum_jm_ja_j=1$ and $m_j,a_j>0$. Hence the cycle is an integral curve $C\in|h|$ with multiplicity one. If $s_C$ cuts out $C$, the image $s_CG$ vanishes generically and is zero by purity. Thus $G$ is rank-one torsion-free on $C$.

Any nonzero proper subsheaf $A\subset G$ has finite-length quotient, hence $p_G-p_A=\length(G/A)/4>0$. Taking the support gives \eqref{eq:Jacobianfibration}; its fiber parametrizes rank-one torsion-free sheaves on $C$ with $\chi(G)=\chi$, i.e. its compactified Jacobian of degree $\chi+2$ since $p_a(C)=3$. For every $\chi$ there are line bundles of this degree on smooth members. Proposition~\ref{prop:generalK3boundary} with $b=1$ gives projectivity, smoothness and dimension six.

At each such $G$ the scalar unit and the top trace in Lemma~\ref{lem:tangentreduction} are isomorphisms of lines, so
\[
 T_{\mathbf M_S,[G]}\simeq\Ext^1_S(G,G)[0],\qquad \ext^1_S(G,G)=2+h^2 .
\]
The locally finitely presented derived charts are smooth over $\CC$ by the cotangent criterion. They agree with their smooth classical truncations.
\end{proof}

\begin{proposition}\label{prop:evaluation}
A flat family of relative expanded sheaves, normal to its outer divisor and with pure restriction of vector $(0,h,\chi)$, has a boundary evaluation morphism
\[
\operatorname{ev}:M^{\rel}\longrightarrow J_\chi.
\]
This construction commutes with arbitrary base change on classical parameter schemes. It does not require smoothness or transversality of the supporting surface or of its intersection curve.
\end{proposition}
\begin{proof}
Let $F_T$ be a flat family on $W_T$ and let $D_T$ be the outer Cartier divisor. In
\[
 0\to F_T(-D_T)\xrightarrow{s}F_T\to i_*G_T\to0,qquad
 G_T=i^*F_T ,
\]
multiplication by $s$ stays injective on every geometric fiber by the normality hypothesis. The fiberwise criterion for flatness thus gives exactness on $W_T$ and flatness of $G_T$ over $T$. At internal marked nodes the same assertion is Proposition~\ref{prop:familyperfect}. The restrictions are stable by Proposition~\ref{prop:Jacobian}, so their classifying morphism to the stable-sheaf stack induces the displayed morphism to its scalar-rigidified moduli space.

Universal sheaves can be treated on the scalar gerbe; no choice of a global untwisted universal sheaf is required. Naturality of restriction proves base-change compatibility.
\end{proof}

\begin{remark}\label{rem:boundaryclass}
The vector, not only the dimension of the support, is needed for Proposition~\ref{prop:Jacobian}. In the primitive quintic setting Theorem~\ref{thm:primitiveprofiles} proves that every nonzero boundary has curve class $h$, but purity remains a separate depth condition. Zero-boundary sectors have the relative theory constructed below. Outside the primitive original polynomial, higher-multiplicity and nonprimitive boundary sectors require their actual moduli targets.
\end{remark}

\subsection{Relative objects with zero boundary}\label{sec:zerorelative}
The zero-boundary alternative of Theorem~\ref{thm:primitiveprofiles} has virtual dimension zero. It has its own natural relative determinant construction. It is necessary even in the primitive class, by Proposition~\ref{prop:zeroboundaryfamily}.

Let $W$ be a relative root chain with outer divisor $D$ and the properties \eqref{eq:pairCY}. Let $F$ be simple, perfect and admissible, with $\LL i_D^*F=0$. The multiplication morphism
\[
 F(-D)\longrightarrow F
\]
is an isomorphism, by its Cartier restriction triangle. In particular the support is disjoint from $D$.

The determinant of $F$ has a canonical trivialization on $D$: it is the determinant of the zero restricted complex. Write $\mathbf{Pic}_{W,D}$ for the derived Picard stack of line bundles equipped with a trivialization on $D$. Its tangent is
\begin{equation}\label{eq:relativePictrace}
 \fib\bigl(\RG(\OO_W)[1]\longrightarrow\RG(\OO_D)[1]\bigr)
 =\RG(\OO_W(-D))[1].
\end{equation}
The determinant morphism of the zero-boundary perfect-complex moduli naturally has this target. Since the ambient rank of $F$ is zero, scalar automorphisms act trivially on its determinant and its boundary trivialization. Hence this morphism descends after removing the scalar automorphisms of $F$.

\begin{theorem}\label{thm:zeroboundaryPOT}
On a bounded algebraic chart of simple admissible relative sheaves with zero boundary, impose the compatible classical class of the determinant with its canonical boundary trivialization by derived pullback from $t_0\mathbf{Pic}_{W,D}$. After scalar rigidification this gives a two-term obstruction theory for sheaves on $W$ with tangent $\Ext^1_W(F,F)$, obstruction $\Ext^2_W(F,F)$, and virtual dimension zero.

For fixed targets, gluing with the zero sheaf on the other root branch identifies this theory with the restriction of the complete expanded Calabi--Yau obstruction theory of Theorem~\ref{thm:ambientPOT}. On the finite-stabilizer locus the comparison descends to the marked expansion charts.
\end{theorem}
\begin{proof}
The family perfectness and derived representability arguments of Sections~\ref{sec:local} and \ref{sec:derived} apply on the zero-restriction open locus. Serre duality and $F(-D)\cong F$ give
\[
 \Ext^j_W(F,F)\cong\Ext^{3-j}_W(F,F)^\vee.
\]
Thus $\Ext^0=\Ext^3=\CC$, no negative or above-three groups occur, and the middle dimensions agree.

The relative determinant differential is trace with values in the fiber of the structure-sheaf restriction map. At the complex level this is the natural map
\begin{equation}\label{eq:zerocompacttrace}
 K_W(F)\longrightarrow\RG(\OO_W(-D)).
\end{equation}
To see its lift explicitly, use the isomorphism
$\RHom(F,F(-D))\cong\RHom(F,F)$ induced by multiplication and then take the twisted trace. Equivalently, the boundary endomorphism complex is zero, so the commuting restriction-and-trace square canonically gives the map to the homotopy fiber. These two constructions coincide by the Cartier triangle.

Duality and the structure-sheaf cohomology of $W$ imply
\[
 H^j(\OO_W(-D))=0\ (j\ne3),\qquad H^3(\OO_W(-D))=\CC.
\]
The degree-three map in \eqref{eq:zerocompacttrace} is dual to the unit $H^0(\OO_W)\to\Hom(F,F)$ and is an isomorphism. The classical relative Picard space is locally discrete: its tangent and obstruction groups are $H^1(\OO_W(-D))=H^2(\OO_W(-D))=0$, and its infinitesimal automorphisms vanish as well. Consequently the tangent complex of the stated derived determinant constraint is
\[
 \fib\left(
  \cofib(\OO\xrightarrow{\id}K_W)
     \longrightarrow\RG(\OO_W(-D))
 \right)[1].
\]
It has cohomology only in degrees zero and one, equal to $\Ext^1$ and $\Ext^2$. Perfectness and base change give the same Tor-amplitude assertion in families. Its cotangent morphism is the canonical two-term obstruction theory, as in Theorem~\ref{thm:ambientPOT}. The virtual rank is zero by the displayed duality.

Finally, attach the other relative root $W'$ along $D$ and put the zero object on $W'$. In the full derived gluing problem the boundary object is zero, so the endomorphism complex of the glued object is $K_W$. Its determinant is trivial on $W'$ with the canonical matching trivialization along $D$. Thus the remaining determinant data are exactly a line bundle on $W$ trivialized on $D$.

The global determinant construction, restricted to this matching problem, is the relative determinant construction above. Derived descent therefore identifies the cotangent morphisms, not just their tangent cohomology. On the finite-stabilizer locus the empty branch is an unexpanded root. Indeed, each empty ruled component would admit its entire dilation group as automorphisms of the zero sheaf. Thus the expansion automorphisms come from the occupied branch and the identification descends to the marked quotient charts.
\end{proof}

\begin{figure}[htbp]
\centering
\begin{tikzpicture}[x=1cm,y=.76cm,font=\small,line cap=round]
 \node[anchor=east] at (3,2.36) {degree}; \foreach \x/\v in {3.75/0,5.50/1,7.25/2,9.00/3}{\node at (\x,2.36) {$\v$};}
 \foreach \y/\lab in {1.62/{H^\bullet(W,\OO_W)},.66/{H^\bullet(D,\OO_D)},-.30/{H^\bullet(W,\OO_W(-D))}}{
    \node[anchor=east] at (3,\y) {$\lab$};
    \draw[qiInk!18] (3.48,\y)--(9.29,\y);}
 \foreach \x/\y in {3.75/1.62,3.75/.66,7.25/.66,9.00/-.30}{
    \filldraw[fill=qiTeal,draw=white,line width=.7pt] (\x,\y) circle (3.25pt);}
 \draw[-{Stealth[length=1.7mm]},qiBlue,line width=.75pt] (3.75,1.45)--(3.75,.83);
 \draw[-{Stealth[length=1.7mm]},qiAmber,line width=.85pt] (7.42,.58)--(8.83,-.21);
 \node[above,qiAmber] at (8.13,.31) {$\partial$};
\end{tikzpicture}
\caption{The cohomology sequence for the outer K3 divisor $D$. A dot denotes a copy of $\CC$. Restriction on constants and the connecting map from $H^2(D,\OO_D)$ to $H^3(W,\OO_W(-D))$ are isomorphisms. The latter group receives the top trace in the zero-boundary obstruction theory.}\label{fig:relativecohom}
\end{figure}
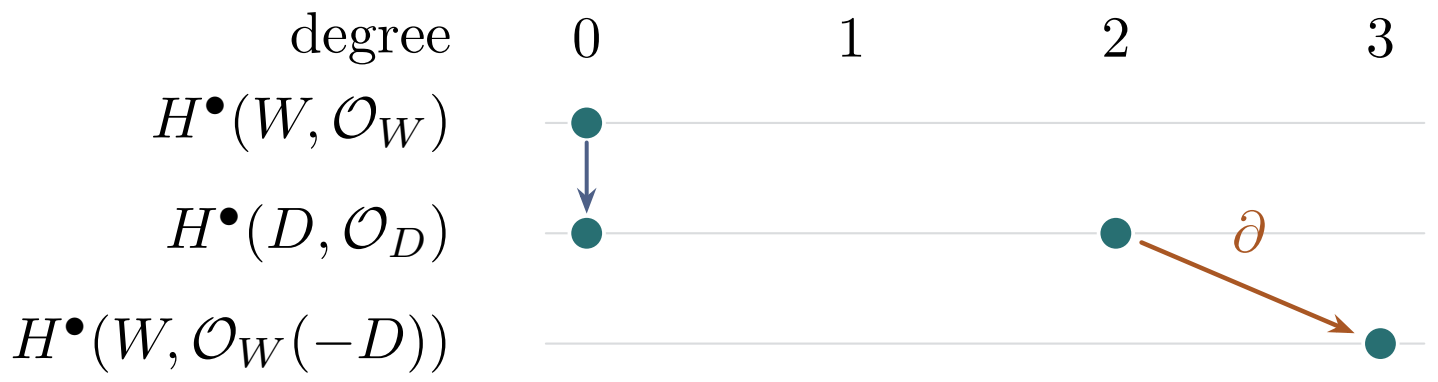

Figure~\ref{fig:relativecohom} locates the class removed by the determinant with its outer trivialization.

\begin{remark}\label{rem:zerotheory}
Ordinary unframed relative Picard reduction on $W$ would not remove this top Ext group, because $H^{>0}(\OO_W)=0$. The relative determinant trivialization on $D$ is essential. It is canonical for a zero-boundary perfect complex and is not an additional choice of framing of the sheaf. Identifying this moduli space with a proper marked sector gives the properness needed to take the degree of the virtual class of Theorem~\ref{thm:zeroboundaryPOT}.
\end{remark}

\section{The virtual diagonal comparison}\label{sec:diagonal}
The purpose of this section is to identify the actual obstruction morphism on a marked moduli space with the pullback along the smooth boundary diagonal.  The input is an equivalence of derived moduli functors.

Let $\mathfrak A_i$ be the stacks of relative expansions of $(Y_i,S)$. Both are smooth of dimension zero. Cutting a marked node in a complete expansion gives two relative expansions, and conversely they glue uniquely as expanded targets. On the standard charts this is the identity
\[
[\AA^{a}/\GG^{a}]\times[\AA^{b}/\GG^{b}]
\cong[\{t_m=0\}\subset\AA^{a+b+1}/\GG^{a+b}],
\]
with the ordered root components held fixed. These chart identifications respect the expansion groupoid. Thus the stack of expansions with one marked node is the product of the two relative expansion stacks; the remaining smoothing coordinates lie entirely on the two sides.

Fix a numerical sector $\gamma$. Let $M_{i,\gamma}^{\rel}\to\mathfrak A_i$ be actual relative sheaf moduli with their independent stability conditions and obstruction theories. Suppose that the nonzero boundary lies in a separated smooth space of simple complexes without negative self-Ext, or in the stable curve-sheaf moduli space, denoted $J_\lambda$. Let
\begin{equation}\label{eq:Gsector}
P_\gamma=M_{1,\gamma}^{\rel}\times M_{2,\gamma}^{\rel},
\qquad
G_\gamma=M_{1,\gamma}^{\rel}\times_{J_\lambda}M_{2,\gamma}^{\rel}.
\end{equation}
The fiber product is interpreted with the boundary-sheaf matching and the scalar conventions of Section~\ref{sec:reduction}. In particular, it is not a fiber product of support-linear-system maps.

\begin{hypothesis}\label{hyp:product}
In the sector $\gamma$, classical restriction and matching identify the node-marked stable moduli with all of $G_\gamma$. Both relative spaces are proper. Their universal boundary complexes take values in the chosen open subspace of the reduced boundary stack of Mukai vector $\lambda$. The two relative stability loci are open in the corresponding perfect-sheaf stacks. Their determinant and scalar conventions are those in \eqref{eq:reducedmoduli}.
\end{hypothesis}

The stability part of this hypothesis is verified on an integral rank-one chain by Theorem~\ref{thm:cutcriterion} and Proposition~\ref{prop:cutting}, provided the full set of cut inequalities holds on every object in the two proper relative moduli. Theorem~\ref{thm:idealchamber} supplies one uniform chamber for this stability test.

\begin{theorem}\label{thm:diagonal}
Under Hypothesis~\ref{hyp:product}, the obstruction theory on the marked moduli induced from Theorem~\ref{thm:ambientPOT} agrees with the obstruction theory of the derived fiber product over $J_\lambda$. Consequently
\begin{equation}\label{eq:diagonalformula}
[G_\gamma]^{\vir}
=
\Delta_{J_\lambda}^{!}
\left([M_{1,\gamma}^{\rel}]^{\vir}
\times[M_{2,\gamma}^{\rel}]^{\vir}\right).
\end{equation}
For $\lambda=(0,h,\chi)$ this is an equality of zero-cycles.
\end{theorem}
\begin{proof}
Over $\mathfrak A_1\times\mathfrak A_2$, write $\mathfrak M_L,\mathfrak M_R$ for the reduced derived relative factors. For a boundary value $E$, K3 duality identifies the trace and the unit:
\[
 \bigl(\Ext^2_S(E,E)\xrightarrow{\tr}H^2(S,\OO_S)\bigr)^\vee
 =\bigl(\CC\xrightarrow{1\mapsto\id_E}\Hom_S(E,E)\bigr).
\]
Simplicity makes this an isomorphism. With negative self-Ext absent, the reduced boundary tangent is $\Ext_S^1(E,E)[0]$. By Theorem~\ref{thm:derivedglue}, Proposition~\ref{prop:redglue} and Hypothesis~\ref{hyp:product},
\[
 \mathfrak G_\gamma\simeq\mathfrak M_L\times^{\mathbf R}_{J_\lambda}\mathfrak M_R,
 \qquad t_0\mathfrak G_\gamma=G_\gamma .
\]
On this product the restriction difference gives
\[
 T_{\mathfrak G_\gamma}\longrightarrow
 T_{\mathfrak M_L}\oplus T_{\mathfrak M_R}
 \xrightarrow{\mathrm d e_L-\mathrm d e_R}e^*T_{J_\lambda}.
\]
All tangent bundles or complexes here are pulled back to the product and taken relative to the two expansion bases. The middle term has amplitude $[0,1]$ by Theorem~\ref{thm:branchindexnew}, and the last term has degree zero; the fiber has amplitude $[0,1]$.

Let $j:G_\gamma\to P_\gamma$ and $e:G_\gamma\to J_\lambda$ be the natural maps. Applying the cotangent functor to this derived Cartesian square gives the commutative diagram of distinguished triangles
\begin{equation}\label{eq:POTdiagram}
\begin{tikzcd}[column sep=large,row sep=large]
\LL j^*E_{P_\gamma}\arrow[r]\arrow[d]
 &E_{G_\gamma}\arrow[r]\arrow[d]
 &e^*\Omega_{J_\lambda}[1]\arrow[d]\\
\LL j^*L_{P_\gamma}\arrow[r]
 &L_{G_\gamma}\arrow[r]
 &L_{G_\gamma/P_\gamma}.
\end{tikzcd}
\end{equation}
Here the first two terms are taken relative to $\mathfrak A_1\times\mathfrak A_2$ initially; the same diagram for the absolute theories follows by adjoining the common smooth-base cotangent complex. The right vertical morphism is induced by base change of the regular diagonal of $J_\lambda$. The connecting arrows in the diagram come from the actual cotangent triangles.

For this compatible triple of obstruction theories, virtual pullback functoriality \cite[Theorem~4.8 and Corollary~4.9]{Mano} gives
\[
 j^![P_\gamma]^{\vir}
 =\Delta_{J_\lambda}^!\bigl([M_{1,\gamma}^{\rel}]^{\vir}\times[M_{2,\gamma}^{\rel}]^{\vir}\bigr)
 =[G_\gamma]^{\vir}.
\]
Restriction is functorial under target isomorphisms, hence this identity descends on the expansion charts. The diagonal of the separated space $J_\lambda$ is closed; $G_\gamma$ is therefore closed in the proper product $P_\gamma$. For $\lambda=(0,h,\chi)$ the virtual dimension is $3+3-(2+h^2)=0$.
\end{proof}

\begin{remark}\label{rem:bubblesincluded}
In \eqref{eq:diagonalformula} each unmarked bubble belongs to one of the two relative expansions. Its sheaf and its smoothing directions enter the corresponding relative factor before the restrictions at the marked node are identified.
\end{remark}

\begin{corollary}\label{cor:complexdiagonalnew}
Let $U$ be a separated open moduli space of simple perfect complexes of fixed Mukai vector $\lambda$ on $S$ with $\Ext^j_S(E,E)=0$ for $j<0$, with scalars rigidified and determinant reduced as in \eqref{eq:reducedmoduli}. It is smooth of dimension $2+\langle\lambda,\lambda\rangle$.

If a marked component of a proper expanded stable problem is the full matching of two independently proper simple relative factors over $U$, with compatible determinant choices, its virtual zero-cycle is the refined diagonal pullback of the product of their virtual classes. In particular this comparison applies when such an open space parametrizes two-term complexes.
\end{corollary}
\begin{proof}
Duality gives $\Ext^j(E,E)=0$ for $j>2$, and the degree-two trace is dual to the isomorphism $\CC\to\Hom(E,E)$. Removing the scalar unit and that trace leaves only $\Ext^1(E,E)$, of dimension $2-\chi_S(E,E)$. The reduced derived boundary is smooth by its cotangent criterion. On each factor \eqref{eq:branchthirdnew} gives $\Ext^3_Y(F,F)=0$.

Derived descent in Theorem~\ref{thm:derivedglue} and the determinant and scalar comparison of Proposition~\ref{prop:redglue} identify the reduced derived marked locus with the derived fiber product over $U$. Pulling back the regular diagonal of $U$ gives its cotangent triangle, hence its virtual pullback as in Theorem~\ref{thm:diagonal}. At a matching pair the triangle on tangent complexes is
\[
 T_{\rm match}\longrightarrow T_L\oplus T_R
     \xrightarrow{\mathrm de_L-\mathrm de_R}\Ext^1_S(E,E).
\]
Its Euler characteristic is $2(1+\langle\lambda,\lambda\rangle/2)-(2+\langle\lambda,\lambda\rangle)=0$.
\end{proof}

\subsection{The general matching class}
When the boundary includes complexes with more than one cohomology sheaf, use Theorem~\ref{thm:perfectpurematchingnew}. Define its matching derived stack by taking the homotopy fiber product of the unrigidified branch perfect-complex moduli over the full boundary perfect-complex moduli. Impose the determinant condition of the glued object, its global stability condition, and its single scalar rigidification; quotient also by the expansion automorphisms.

Theorems~\ref{thm:derivedglue} and \ref{thm:perfectpurematchingnew} identify this functor, on the prescribed perfect pure sheaf loci, with the actual marked expanded-sheaf moduli. Its obstruction theory for sheaves on the glued threefold is consequently the pullback of \eqref{eq:ambientPOT}.

This construction supplies a class $[\MM_\gamma^\dagger]^{\vir}$ even when there is no ordinary Jacobian diagonal to use. Theorem~\ref{thm:diagonal} is the additional identification that makes this class a product of independently defined relative invariants on the smooth-boundary sectors. Keeping the two statements distinct prevents deletion of a boundary charge or a sheaf-cohomology contribution.

\subsection{A properness obstruction from the boundary isomorphisms}
\label{sec:boundaryisomproper}
Even on two fixed unexpanded roots the choice of an isomorphism on the K3 can affect properness. Here we determine the fiber of the forgetful morphism which retains the two branch sheaves and forgets their identification. This also describes the amplitude of the derived diagonal used for the virtual pullback.

\begin{proposition}\label{prop:boundaryisomproper}
Fix nonzero simple pure perfect sheaves $F_L,F_R$ on two fixed roots, with identified derived restrictions $E\ne0$ to $S$. Suppose $\Ext^{<0}_S(E,E)=0$. Put $D=\End_S(E)$ and $d=\dim_\CC D$. After rigidifying the scalar automorphisms of the sheaves, the family of all matchings above the two fixed rigidified branch points is the smooth affine variety
\[
 D^\times/\GG
   =\{[u]\in\PP(D):\det(L_u)\ne0\},
 \qquad L_u:D\longrightarrow D,\quad v\longmapsto uv.
\]
It has dimension $d-1$. If $d>1$, any nonempty open stable locus in this fiber is nonproper. Consequently a proper marked moduli space mapping to separated rigidified branch spaces cannot have such an open as its entire fiber over these two points. At the same $E$ the natural derived diagonal of the scalar-rigidified, determinant-reduced K3 moduli has degree-two relative tangent cohomology $(D/\CC\,\id_E)^\vee$.
\end{proposition}
\begin{proof}
Fix an identification of the restrictions with $E$. For a $\CC$-algebra $A$, perfect base change and $\Ext^{<0}(E,E)=0$ give
\[
 \operatorname{Isom}(E_A,E_A)=(D\otimes A)^\times,\qquad
 (a,b)\cdot u=bua^{-1}\quad ((a,b)\in\GG^2).
\]
The subgroup fixing $u$ is the diagonal scalar group. Thus the rigidified matching fiber is
\[
 D^\times/\GG= D_+(\det L_u)\subset \PP(D).
\]
Indeed $uv=1$ is solvable uniquely iff $L_u$ is invertible, and $\det L_u$ is homogeneous of degree $d$, nonzero at $u=1$. This principal open is affine and irreducible of dimension $d-1$. For $d>1$ it cannot be proper over $\CC$, since global functions on a proper integral variety are constant.

Over separated branch spaces properness of the matching space would force properness of each geometric fiber of the forgetful morphism (closed graph followed by proper projection).

For the diagonal calculation, use $C_S^{\red}$ of \eqref{eq:redtangent}. The scalar unit removes $\CC\,\id_E$ from $\Ext^0_S(E,E)$. Since $H^1(S,\OO_S)=0$, the degree-two term after the determinant condition is
\[
 H^2(C_S^{\red})=
 \ker\bigl(\Ext^2_S(E,E)\xrightarrow{\tr}H^2(S,\OO_S)\bigr)
       \simeq(D/\CC\,\id_E)^\vee
\]
by Serre duality. The relative tangent of a derived diagonal is the fiber of the diagonal map $T_E\to T_E\oplus T_E$, hence is $T_E[-1]=C_S^{\red}$. For $d>1$ this natural relative tangent therefore exceeds the amplitude $[0,1]$ of the dual of an ordinary two-term relative perfect obstruction theory. The global theory in \eqref{eq:ambientPOT} uses the matched object before taking its virtual class.
\end{proof}

For an instance on the quintic degeneration, take the Cartier surface of Section~\ref{sec:generalbranches} with equation $f_1a_{r-1}+f_4b_{r-4}$, $r\ge4$, and the ideal of a smooth point of that surface away from $S$. Its two branch sheaves are simple. In fact each root divisor $T_i$ has $H^0(\OO_{T_i})=\CC$ by the divisor sequence: $H^1(Q,\OO_Q(-rH_Q))=0$ by duality and Kodaira vanishing, and $H^1(B,\OO_B(-rH_B))=H^1(\PP^3,\OO(-r))=0$ by the blowup projection.

On the branch containing the point, its ideal has endomorphism sheaf $\OO_{T_i}$ by depth two there. On the K3 the resolution reduces to
\[
 E=\OO_S\oplus\OO_S(-rh)[1],\qquad
 \Ext^{-1}_S(E,E)=H^0(S,\OO_S(-rh))=0 .
\]
The only potentially nonzero off-diagonal entry in $\End(E)$ is $H^1(S,\OO_S(-rh))$, which vanishes by duality and Kodaira vanishing. Thus $D=\CC\oplus\CC$, and the gluing parameter is the ratio of the two diagonal scalars. The fiber is $\GG$. The arc represented by $(1,t)$ has no limit in this fiber over $\CC[[t]]$, even after finite base change: its scalar ratio has nonzero valuation.

Each gluing gives a simple global sheaf, since in \eqref{eq:generalEndnew} the map $\CC\oplus\CC\to D$ on branch endomorphisms has diagonal kernel. This example concerns the fixed pair of branches; it makes no stability assertion about those gluings. If an open set of them is selected by stability, compactifying it requires limits beyond those isomorphisms of the two fixed branches.

Figure~\ref{fig:matchingratio} displays the missing limits of the matching parameter for this pair of simple branches.

\begin{figure}[htbp]
\centering
\begingroup
\definecolor{ratioaxis}{HTML}{53618A}
\definecolor{ratioseam}{HTML}{176568}
\definecolor{ratiopt}{HTML}{B15C25}
\begin{tikzpicture}[x=1cm,y=1cm,font=\small,line cap=round,line join=round]
\path[draw=ratioaxis,fill=ratioaxis!5,line width=.9pt] (0,0) circle (1.47);
\draw[ratioaxis!55,dashed] (-1.47,0) arc[start angle=180,end angle=0,x radius=1.47,y radius=.37];
\draw[ratioaxis!65] (-1.47,0) arc[start angle=180,end angle=360,x radius=1.47,y radius=.37];
\draw[ratioaxis!40] (0,1.47) arc[start angle=90,end angle=-90,x radius=.63,y radius=1.47];
\draw[ratioaxis!40] (0,-1.47) arc[start angle=270,end angle=90,x radius=.63,y radius=1.47];
\draw[ratioseam,line width=1.1pt,-{Stealth[length=2mm]}] plot[variable=\b,domain=-49:69,samples=40] ({.63*cos(\b)},{1.47*sin(\b)});
\foreach \y in {-1.47,1.47}{\filldraw[fill=white,draw=ratiopt,line width=1.15pt] (0,\y) circle (3.1pt);}
\node[anchor=east] at (-.25,1.47) {$[1:0]$};
\node[anchor=east] at (-.25,-1.47) {$[0:1]$};
\node[anchor=west,text=ratioseam!80!black] at (1.79,.84) {$[1:t],\quad t\longrightarrow0$};
\draw[ratioseam!75,line width=.55pt] (1.75,.84)--(1.34,.84)--(.62,.63);
\node[anchor=west] at (4.92,.74) {$D=\CC\oplus\CC$};
\node[anchor=west] at (4.92,0) {$\varphi=(a,b),\quad ab\ne0$};
\node[anchor=west] at (4.92,-.74) {$[a:b]\in\PP(D)$};
\end{tikzpicture}
\endgroup
\caption{For the split boundary complex in the example, a matching acts by a scalar on each summand. Quotienting by the branch scalars leaves the ratio $[a:b]$. The two open circles on the projective line are the two points where multiplication by $(a,b)$ on $D$ has rank one. The arc ends at a nonunit, explaining the failure of the valuative criterion over the fixed branch points.}
\label{fig:matchingratio}
\end{figure}

\subsection{A simplicity criterion for primitive boundary evaluation}
We consider the positive boundary of a perfect pure sheaf of primitive divisor class on a complete expanded quintic fiber. Corollary~\ref{cor:primitiveautomaticnew} shows that its boundary $E$ is a sheaf of dimension one and Mukai vector $(0,h,\chi)$. Purity is still to be tested.
\begin{proposition}\label{prop:route19-simple}
Assume $\Pic(S)=\mathbb Zh$. For a coherent sheaf $E$ of vector $(0,h,\chi)$ and dimension one on $S$ the following are equivalent:
\[
 \End_S(E)=\CC,\qquad E\text{ is pure},\qquad [E]\in J_\chi .
\]
Consequently a family of such positive perfect boundary restrictions on a classical base, all of whose geometric fibers are simple, defines a morphism to $J_\chi$.
\end{proposition}
\begin{proof}
Let $T_0(E)$ be the maximal zero-dimensional subsheaf. If $T_0(E)\ne0$, choose a geometric point $p$ in its support. The socle of the finite-length stalk gives an injection $b:k(p)\hookrightarrow E$. Nakayama's lemma gives a surjection $a:E\twoheadrightarrow k(p)$ since $E_p\ne0$. Thus $ba\ne0$, while $ba$ is zero at the generic point of the divisorial support. It is linearly independent of $\id_E$. Hence simplicity implies $T_0(E)=0$. For a sheaf of dimension one this is purity.

By Proposition~\ref{prop:Jacobian}, purity and the primitive curve class imply stability, and a stable sheaf over an algebraically closed field is simple. For the family statement the perfect derived restriction has all geometric fibers concentrated in degree zero. The fiberwise flatness criterion identifies it with a flat sheaf compatible with base change.

Its stable fibers give the required classifying map on a chart; tensoring the family by a line bundle from the base leaves the map to the coarse space unchanged, so the maps descend.
\end{proof}
This yields a test for a proposed compactification. Suppose its positive marked matching over a fixed pair of simple branches is a nonempty open set of the scalar-rigidified isomorphisms of their common nonzero boundary $E$, with no residual target automorphisms. Suppose further that this open set itself is proper.

Writing $D=\End_S(E)$, Proposition~\ref{prop:boundaryisomproper} identifies the isomorphisms modulo scalars with the irreducible affine variety $D^\times/\mathbb G_m$ of dimension $\dim D-1$. A nonempty open subset proper over $\CC$ has closed image under its immersion into this separated variety, hence equals the whole variety. A proper irreducible affine variety over $\CC$ has dimension zero. Thus $\dim D=1$, and in the primitive situation the preceding proposition gives a point of $J_\chi$.

Properness of the total sheaf space by itself does not imply the stated description of its fixed-branches fiber: further bubbles or changes of branch at the limit can compactify that fiber.

For a stable-curve-sheaf open $U$ in a proper separated compactification $M$, the closure of the graph of $e:U\to J_\chi$ in $M\times J_\chi$ has a proper projection to $M$ and a boundary evaluation to $J_\chi$. Over $U$ it is the original graph by separatedness of $J_\chi$. Proposition~\ref{prop:route19-two-arcs} shows how different limits of $e$ can occur over one point of $M$.

Using this construction to compute a DT invariant requires a comparison between the obstruction theory on the graph closure and the theory of sheaves on the expanded threefold, including the virtual contributions over the exceptional fibers.

\section{Global specialization through the marked expansion boundary}\label{sec:specialization}
Under Hypothesis~\ref{hyp:compact}, the normalization of the central divisor of the expansion stack records a choice of joining K3. Pulling the virtual class back to this marked divisor and pushing it forward gives the specialization formula. The obstruction theories needed to express its terms as relative intersections were constructed on the loci specified in Sections~\ref{sec:relative} and~\ref{sec:diagonal}.

\subsection{The marked boundary as a finite normalization}
\begin{proposition}\label{prop:normalboundary}
On the bounded part of $\EE$ used by $\MM$, the special fiber $\EE_0$ is a reduced normal-crossing Cartier divisor. Its normalization
\begin{equation}\label{eq:nunorm}
\nu:\EE_0^\dagger\longrightarrow\EE_0
\end{equation}
is finite, representable, and proper. The stack $\EE_0^\dagger$ is smooth of dimension zero. It parametrizes one marked node, and
\begin{equation}\label{eq:basecycle}
[\EE_0]=\nu_*[\EE_0^\dagger]
\quad\text{in }A_0(\EE_0)_{\QQ}.
\end{equation}
\end{proposition}
\begin{proof}
On \eqref{eq:chartintro} put $A=k[t_0,\ldots,t_q]$ and $f=\prod_{i=0}^q t_i$. The normalization map is
\[
 A/(f)\hookrightarrow \prod_{i=0}^{q}A/(t_i),\qquad
 g\longmapsto(g\bmod t_i)_i .
\]
The right-hand side is finite over $A/(f)$, normal, and birational over each component. The torus action preserves every hyperplane. Normalization commutes with smooth base change in characteristic zero, so these maps descend to the finite representable morphism $\nu$ of \eqref{eq:nunorm}. Each quotient chart of the source is smooth of dimension $q-q=0$.

A point on the $i$th summand records the unsmoothed node $t_i=0$; all $t_j$, $j\ne i$, remain free. Deleting an invertible parameter on an overlapping chart preserves this choice. Over the dense open of each component with exactly one node, $\nu$ is an isomorphism, including on stabilizers. Since each factor of $f$ occurs once, proper pushforward of fundamental cycles gives \eqref{eq:basecycle}.
\end{proof}

\begin{figure}[htbp]
\centering
\begin{tikzpicture}[x=1cm,y=1cm,font=\small,line cap=round,line join=round]
 \begin{scope}[shift={(8.65,-.55)}]
 \path[fill=qiBlue!24,draw=qiBlue] (0,0)--(1.72,-.49)--(1.72,1.05)--(0,1.54)--cycle;
 \path[fill=qiTeal!22,draw=qiTeal] (0,0)--(-1.31,-.72)--(-1.31,.82)--(0,1.54)--cycle;
 \path[fill=qiAmber!22,draw=qiAmber] (0,0)--(1.72,-.49)--(.41,-1.21)--(-1.31,-.72)--cycle;
 \draw[qiInk,line width=.8pt] (0,1.54)--(0,0)--(-1.31,-.72);
 \draw[qiInk,line width=.8pt] (0,0)--(1.72,-.49);
 \fill[qiInk] (0,0) circle (1.9pt);
 \end{scope}
 \foreach \x/\s/\colorName/\lab in {1.10/1.05/qiBlue/0,1.10/.05/qiTeal/1,1.10/-.95/qiAmber/2}{
  \path[draw=\colorName,fill=\colorName!20] (\x,\s)--(\x+1.53,\s+.20)--(\x+2.14,\s-.19)--(\x+.61,\s-.39)--cycle;
  \fill[\colorName] (\x,\s) circle (1.8pt);
  \node[anchor=east] at (\x-.12,\s-.01) {$t_{\lab}=0$};}
 \draw[-{Stealth[length=2mm]},qiInk,line width=.8pt] (4.13,-.10)--(6.24,-.10);
 \node[above] at (5.2,-.09) {$\nu$};
 \node[anchor=west] at (10.62,-.34) {$t_0t_1t_2=0$};
\end{tikzpicture}
\caption{Three marked faces above the two-bubble chart. The right panel shows the part of each coordinate plane in one real octant; the full chart is complex. Each face has two free smoothing parameters and occurs once in the special divisor. The three marked origins map to the same target with three unsmoothed joins.}\label{fig:normalizationplanes}
\end{figure}
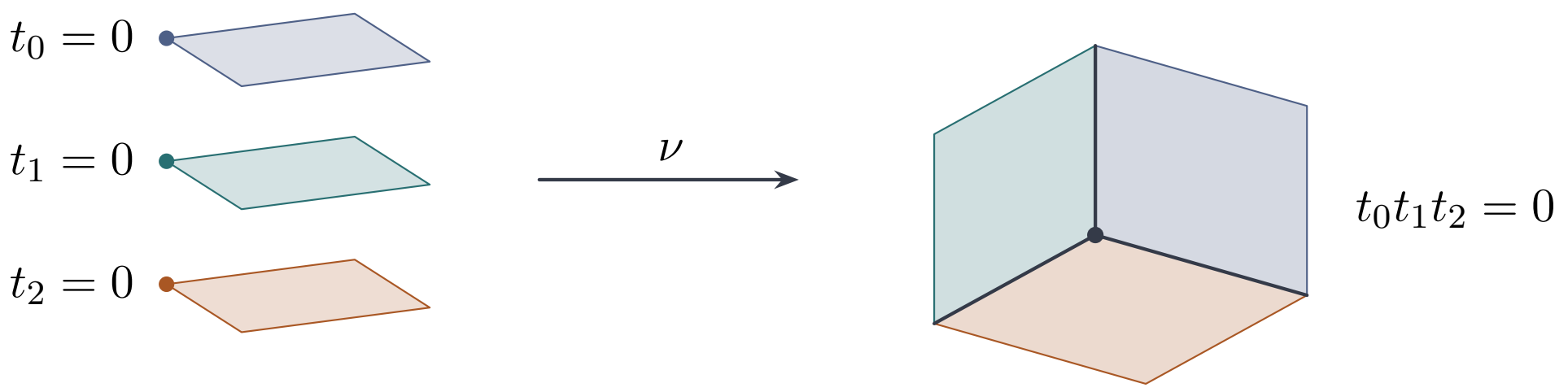

Figure~\ref{fig:normalizationplanes} depicts the finite map before taking the quotient by the rescaling torus.

\begin{lemma}\label{lem:virtualnormalization}
Let $p:\MM\to\EE$ have the relative perfect obstruction theory of Theorem~\ref{thm:ambientPOT}. Form the Cartesian diagram
\[
\begin{tikzcd}[column sep=large,row sep=large]
\MM^\dagger\arrow[r,"\nu_{\MM}"]\arrow[d,"p^\dagger"']
 &\MM_0\arrow[r]\arrow[d,"p_0"']
 &\MM\arrow[d,"p"]\\
\EE_0^\dagger\arrow[r,"\nu"']
 &\EE_0\arrow[r,"i_0"']
 &\EE .
\end{tikzcd}
\]
Give all base changes their induced relative obstruction theories. Then
\begin{equation}\label{eq:globalvirtualnormalization}
i_0^![\MM]^{\vir}
=(\nu_{\MM})_*[\MM^\dagger]^{\vir}.
\end{equation}
The class on the right is the virtual class from the actual base-changed sheaf deformation problem over the smooth stack $\EE_0^\dagger$.
\end{lemma}
\begin{proof}
The notation $i_0^!$ on the upper row means refined pullback along the base Cartier divisor $\EE_0\hookrightarrow\EE$. Bivariance of virtual pullback and its compatibility with refined regular-embedding pullback give
\begin{equation}\label{eq:firstbasechange}
i_0^!p^![\EE]=p_0^!i_0^![\EE]=p_0^![\EE_0].
\end{equation}
Apply \eqref{eq:basecycle} and proper-pushforward compatibility of virtual pullback:
\begin{equation}\label{eq:pushvirtual}
p_0^![\EE_0]
=p_0^!\nu_*[\EE_0^\dagger]
=(\nu_{\MM})_*(p^\dagger)^![\EE_0^\dagger].
\end{equation}
These are the bivariant and proper-pushforward properties proved in \cite[Theorems 4.1, 4.3--4.4]{Mano}. The maps are of Deligne--Mumford type, the relative obstruction complexes are perfect, and the finite normalization is projective, so the proper-pushforward statement applies on the quotient strata as well.

The smooth stack $\EE_0^\dagger$ has dimension zero. Thus $(p^\dagger)^![\EE_0^\dagger]$ is the ordinary virtual class of the base-changed moduli with the obstruction theory pulled back from the universal expanded threefold. Proposition~\ref{prop:POTbasechange} identifies that theory with the independently constructed marked matching theory. Substituting \eqref{eq:totalvirtual} into \eqref{eq:firstbasechange}--\eqref{eq:pushvirtual} proves \eqref{eq:globalvirtualnormalization}. No flatness of $\MM\to C$ is used.
\end{proof}

\subsection{Specialization of the generic class}
\begin{proposition}\label{prop:specialization}
With the proper moduli of Hypothesis~\ref{hyp:compact},
\begin{equation}\label{eq:specialization}
\spc[M_\eta(v)]^{\vir}=i_0^![\MM]^{\vir}
\quad\text{in }A_0(\MM_0)_{\QQ}.
\end{equation}
Here specialization is Chow specialization over the local ring of the smooth base curve at $0$.
\end{proposition}
\begin{proof}
Pass to the trait $C'=\Spec\OO_{C,0}$, with generic point $\eta$,
closed point $0$ and uniformizer $t$. Write $i:\MM_0\to\MM_{C'}$
and $j:\MM_\eta\to\MM_{C'}$. In the dimension convention on the
total space over the trait, localization gives
\[
 A_1(\MM_0)\xrightarrow{i_*}A_1(\MM_{C'})
           \xrightarrow{j^*}A_0(\MM_\eta)\to0 .
\]
If $j^*a=j^*a'$, write $a-a'=i_*z$. Refined intersection with
$0\hookrightarrow C'$ satisfies
\[
 0^!(a-a')=0^!i_*z
   =c_1\bigl(\OO_{C'}(0)|_0\otimes\OO_{\MM_0}\bigr)\cap z=0 .
\]
The normal line bundle of the base point is trivialized by $t$.
Thus $\spc(j^*a)=0^!a$ is independent of the lift. For
$a=[\MM_{C'}]^{\vir}$ one has $j^*a=[M_\eta(v)]^{\vir}$ by virtual
base change. This proves the formula, including when
$\MM_{C'}\to C'$ is not flat.
\end{proof}

\begin{proof}[Proof of Theorem~\ref{thm:main}]
By Theorems~\ref{thm:proper} and~\ref{thm:ambientPOT}, $\MM\to C$ is proper and $p:\MM\to\EE$ has a virtual pullback of degree zero. By Proposition~\ref{prop:finitecharges}, $\MM^\dagger=\coprod_\gamma\MM^\dagger_\gamma$ is a finite decomposition into open and closed substacks. Applying Proposition~\ref{prop:specialization} and Lemma~\ref{lem:virtualnormalization} gives
\begin{align*}
 \spc[M_\eta(v)]^{\vir}
 &=i_0^!p^![\EE]
   =(\nu_{\MM})_*[\MM^\dagger]^{\vir}\\
 &=\sum_\gamma(g_\gamma)_*[\MM^\dagger_\gamma]^{\vir}
       \quad\text{in }A_0(\MM_0)_\QQ.
\end{align*}
The class on each summand is the class of its derived matching problem by Theorem~\ref{thm:derivedglue} and Proposition~\ref{prop:POTbasechange}, with the scalar and determinant operations on the glued object.

For each identification \eqref{eq:mainmatching}, Theorem~\ref{thm:diagonal} evaluates that term by \eqref{eq:mainJ}. On the boundary of vector $(0,h,\chi)$, Propositions~\ref{prop:relativevdim} and~\ref{prop:Jacobian} give virtual dimension $3+3-6=0$.

\end{proof}

\subsection{Node markings, stabilizers, and coefficients}
\begin{remark}\label{rem:multiplicities}
Equations \eqref{eq:mainfull} and \eqref{eq:mainJ} use the reduced semistable base, ordered roots, the actual normalization $\EE_0^\dagger$, and the corresponding stack-theoretic moduli. In this convention the fundamental-cycle identity \eqref{eq:basecycle} has coefficient one on every normalization branch. Finite intrinsic stabilizers are already included in rational stack pushforward; they are not appended as additional ad hoc factors.

For a finite \'etale cover $q:\widetilde{\MM}^\dagger_\gamma\to\MM^\dagger_\gamma$ of constant degree $d_\gamma$ with the pulled-back obstruction theory,
\[
 q_*[\widetilde{\MM}^\dagger_\gamma]^{\vir}
     =d_\gamma[\MM^\dagger_\gamma]^{\vir}.
\]
For a base change with $t=u\prod_i t_i^{a_i}$ and $u$ a unit, $\operatorname{div}(t)=\sum_i a_i\operatorname{div}(t_i)$. The reduced ordered chart \eqref{eq:chartintro} has $a_i=1$ for all $i$.

\end{remark}

\begin{proposition}\label{prop:chartcompatibility}
On a one-bubble chart $[\AA^2/\GG]$, with $t=ab$ and weights $(1,-1)$, the two comparisons restrict to the same iterated matching functor after both nodes have been marked. The two markings remain distinct over the origin. On either normalized face the unmarked node can be smoothed. The virtual formula is the sum of the two normalized-face classes, not a class on the deepest stratum counted twice.
\end{proposition}
\begin{proof}
The normalization algebra on the one-bubble chart is
\[
 k[a,b]/(ab)\hookrightarrow k[a]\times k[b],
 \qquad f\longmapsto (f(a,0),f(0,b)).
\]
On $a=0$ the node indexed by $a$ persists while $b$ smooths the other node. Theorem~\ref{thm:derivedglue} applies over that entire axis. Interchanging $a$ and $b$ gives the other marking. At $a=b=0$ the two iterated matching functors are equivalent by associativity of homotopy pullbacks. Target rescaling acts on all restriction maps.

Each face $[\AA^1/\GG]$ has dimension zero, whereas their common closed stratum $[\{0\}/\GG]$ has dimension $-1$. Thus Lemma~\ref{lem:virtualnormalization} pulls back the sum of the fundamental classes of the two axes, retaining the unmarked smoothing parameter in each summand.
\end{proof}

The same argument applies to any number of bubbles: one node is marked, and every other smoothing coordinate remains part of one of the two relative families. The normalization formulation is invariant under changes of charts and hence does not require a global labeling of all nodes at once.

\section{Expanded sheaves with finite defects and the structure-sheaf example}\label{sec:cartiercase}
We extend finite quotients of a fixed Cohen--Macaulay hull and apply the local virtual comparison to determinantal surface sheaves. For the hyperplane structure-sheaf character we also construct the projective expanded moduli.

\subsection{Valuative extension of finite defects of an expanded hull}
\label{sec:lifting-punctual-defects}

Both the support and the module of the finite quotient may vary. The extension of its hull is fixed, including for the noninvertible hulls of Theorem~\ref{thm:matrixexpandedallr}.

\begin{proposition}\label{prop:hull-punctual-DVR}
Let $R$ be an excellent discrete valuation ring over $\mathbb C$ with algebraically closed residue field, with fraction field
$K$. Let $Z_R$ be a projective expanded family with smooth generic
threefold and a chain as closed fiber, with ordered branches and smoothed
joining divisors. Fix an $R$-flat perfect sheaf $H$ with
Cohen--Macaulay surface fibers. Assume its generic fiber has rank one
on a geometrically integral support.

On the closed fiber assume that the occupied branches
form a connected chain, that $H_i$ is rank-one torsion-free and
Cohen--Macaulay on each integral branch support, and that $H$ is
admissible with pure rank-one restriction $E_e$ on each integral
joining curve of the occupied chain. Assume an ample line bundle $A$
on the closed target satisfies
\[
     0<[z]D_e(z)<[z]P_{E_e}(z)
\]
at these cuts, with $[z]$ denoting the linear coefficient in \eqref{eq:cutdata}. Zero branches outside
the chain are allowed. Assume the branch and boundary conditions hold over the geometric closed point. For any zero-dimensional quotient $H_K\twoheadrightarrow Q_K$
of degree $n=\dim_K H^0(Q_K)$, after finite base change there exists a refinement
$c:\widetilde Z_R\to Z_R$ by at most $n$ further bubbles and an exact
sequence of $R$-flat sheaves
\[
 0\longrightarrow F\longrightarrow \widetilde H\longrightarrow Q
    \longrightarrow0,\qquad
 \widetilde H=\LL c^*H,\qquad \length Q_0=n .
\]
The closed quotient is supported away from all joining divisors. The
kernel is a fiberwise perfect pure surface sheaf, stable on the generic
target and for a polarization on the closed expansion. Its closed
Cohen--Macaulay hull is $\widetilde H_0$. Every added bubble meets
the support of the defect.

Over the given expanded hull the minimal
refinement placing its generic support points off the joining divisors
is unique up to isomorphism, and the extending quotient on that
refinement is unique. Here $Q_K$ is an arbitrary finite module, including
nonreduced punctual quotients.
\end{proposition}

\begin{proof}
Pass to a finite extension splitting the distinct geometric points
$p_1,\ldots,p_l$ in $\Supp Q_K$; $l\le n$. Projectivity gives sections
of $Z_R$ for their closures. At a joining divisor the equation is
$xy=u\pi^b$ with $\pi$ a uniformizer, $u$ a unit and $b>0$.
Absorb $u$ into $y$. A section specializing to this joint satisfies
\[
       a=\operatorname{ord}_\pi x(p_i)\in\{1,\ldots,b-1\},\qquad
       \operatorname{ord}_\pi y(p_i)=b-a .
\]
Insert one level for each distinct such $a$. For levels
$0=a_0<a_1<\cdots<a_k<a_{k+1}=b$ the smoothing parameters are
$\pi^{a_{j+1}-a_j}$. On the interior of level $a_j$ the chart reads
\[
  x=\pi^{a_j}z,\qquad y=\pi^{b-a_j}z^{-1},\qquad z\in\GG .
\]
The adjacent charts have coordinates $X,Y$ with
$x=\pi^{a_j}X$, $y=\pi^{b-a_{j+1}}Y$ and
$XY=\pi^{a_{j+1}-a_j}$. A point of intermediate valuation specializes
to $X=Y=0$ there. Thus every allowable refinement must contain
every selected level. Omitting unoccupied levels gives the indicated
unique list. Units in the factorization are changes of expansion
coordinates.

Changing a local trivialization of the branch normal
bundle multiplies $x$ by a unit, so the lists agree along each
joining divisor and give the global expansion by the standard
expansion charts of Section~\ref{sec:compactification}. Equivalently the points give a stable expanded configuration at each refined joint in the sense of \cite[Definitions~1.2.1--1.2.2]{AF}.

No point can approach a zero-boundary joint: the restriction of the admissible closed hull to that divisor vanishes, so its support misses the divisor by Nakayama, whereas the closure of each generic support point lies in the support of the hull. Sections specializing away from the joins retain their original limits. The total number of added levels is at most $l$.

The perfect complex $\LL c^*H$ on an old component of the closed
fiber restricts in the derived sense to $H_i$. On a new bubble
$p:\Delta\to D_e$ it restricts to $p^*E_e$, and at a new joint to
$E_e$. These identities follow by composition of derived pullbacks,
flatness of $p$, and admissibility of the fixed hull. The normalization
triangle on the refined fiber now has the form
\[
 \widetilde H_0\longrightarrow
   \bigoplus_i H_i\oplus\bigoplus_{\Delta}p^*E_e
      \xrightarrow{\mathrm{res}_{-}-\mathrm{res}_{+}}
   \bigoplus_{\rm joints} E_e\longrightarrow .
\]
Closed pushforward from the strata is understood, and untouched
joints retain their original sheaves. The difference map is locally
surjective, since either branch surjects onto its restriction. Thus
the closed derived fiber is a sheaf, and the same holds generically.
The fiber criterion for relative Tor amplitude gives an $R$-flat
sheaf $\widetilde H$.

At a closed point of a joint, the nonzero branch
terms have depth two; for $p^*E_e$ this follows by smooth pullback
of the Cohen--Macaulay curve sheaf. The joint term has depth one.
The depth lemma gives depth at least two for the kernel. Away from
the joints the branch description proves the same Cohen--Macaulay
assertion, and localization gives it on the surface support.

Choose a relatively very ample line bundle on $\widetilde Z_R$.
The relative Quot scheme of length $n$ quotients of $\widetilde H$
is projective over $R$ by \cite[Theorem~5.1]{NitsureQuot}.
Its valuative criterion extends $Q_K$ uniquely as a quotient,
up to its unique quotient isomorphism, giving $Q$. Its support
is the closure, as a set, of its generic support: every associated
prime of $Q$ avoids $\pi$, by $R$-flatness.

Conversely each of the
$l$ closed sections of the generic support lies in $\Supp Q$ by
closedness. Thus every new bubble is occupied, and $Q$ meets no
joint in the closed fiber. Both $Q$ and $\widetilde H$ are flat,
so their kernel is flat and the sequence stays exact on fibers.
At a closed defect point the target fiber is regular; all finite
modules there are perfect. At all other points $F$ agrees locally
with $\widetilde H$.

Purity follows from the injection into the pure
hull. For $j$ the open immersion of the complement of the finite defect in the closed
fiber, depth two gives
\[
 \widetilde H_0=j_*j^*\widetilde H_0=j_*j^*F_0 .
\]
This also identifies the defect and its occupied levels intrinsically
from the kernel once the expanded target is fixed.

For stability, start with the nef polarization $c^*A$ on the closed
fiber. On a new bubble, $\RR p_*p^*E_e=E_e$, since
$\RR p_*\OO_\Delta=\OO_{D_e}$. Its Hilbert polynomial for the
nef polarization is $P_{E_e}$, of degree one. Inserting this bubble
adds one copy of $P_{E_e}$ to the branch sum and subtracts exactly
one more copy in the joint sum. The total polynomial and quadratic
coefficient are unchanged.

Every cut replacing the old cut has the
same prefix polynomial $P_{\le e}$ and the same prefix quadratic
coefficient $A_e$ as before, by this same cancellation. The strict
linear inequalities therefore hold at the nef limit. Choose an
ample line bundle $B$ on the refined fiber. For small rational
$\epsilon>0$, $c^*A+\epsilon B$ is ample in additive notation.

The coefficients, including $A_e/A_{\rm total}$, vary continuously
at $\epsilon=0$ by the Hilbert polynomial in the two tensor powers,
with $A_{\rm total}>0$ there. All the finitely many strict
inequalities persist; clear denominators. Removing $Q_0$ away
from the joints subtracts only constants from the branch polynomials:
if $n_i$ is its length on the $i$th branch, then
\[
  D_e(F_0)-D_e(\widetilde H_0)
   =-\sum_{i\le e} n_i+ (A_e/A_{\rm total})n .
\]
All branch supports are still integral of generic rank one, and
the boundary sheaves are unchanged. Theorem~\ref{thm:cutcriterion}
now gives stability. On the geometrically integral generic support $F_K$ is
torsion-free of rank one, hence stable.
\end{proof}

Proposition~\ref{prop:divisorial-hull27} supplies the input $H$ of Proposition~\ref{prop:hull-punctual-DVR} whenever the limit of the supporting divisor and the local divisor classes satisfy its hypotheses. In this case the two propositions construct a limit starting from the generic sheaf. Indeed on a normal integral generic surface a rank-one torsion-free sheaf $G_K$ determines
\[
  0\longrightarrow G_K\longrightarrow G_K^{**}
    \longrightarrow Q_K\longrightarrow0,\qquad N=\length(Q_K).
\]
First extend $G_K^{**}$ by closing a representing divisor and balancing its twists. Then extend $Q_K$ over the minimal refinement of that expanded hull, with at most $N$ new bubbles. The uniqueness for the hull uses the stable support quotient. The uniqueness for the quotient in Proposition~\ref{prop:hull-punctual-DVR} concerns refinements of this extended hull.

A comparison with further sheaf limits for which this support construction does not apply would require a contraction back to the same hull.

In particular, choose the expanded determinantal hull of
Theorem~\ref{thm:matrixexpandedallr}. Its strict linear inequalities
and rank-one Cohen--Macaulay branch and boundary conditions were
verified there. For every $r\ge2$ and $n\ge0$ the proposition
therefore supplies stable limits of every length-$n$ modification of
that generic hull, with character $r(1-e^{-H})-n[\mathrm{pt}]$.
Here the $H$ in the character denotes the hyperplane class on the
quintic.

The noninvertible branch hull restricts to the original hull on the old occupied
branches. The new bubble dilations preserving the kernel permute its
finite defect support on each new bubble. A dilation fixing one
interior point is the identity, hence this dilation subgroup is
finite. The original determinantal bubble has trivial dilation
stabilizer by Theorem~\ref{thm:matrixexpandedallr}.

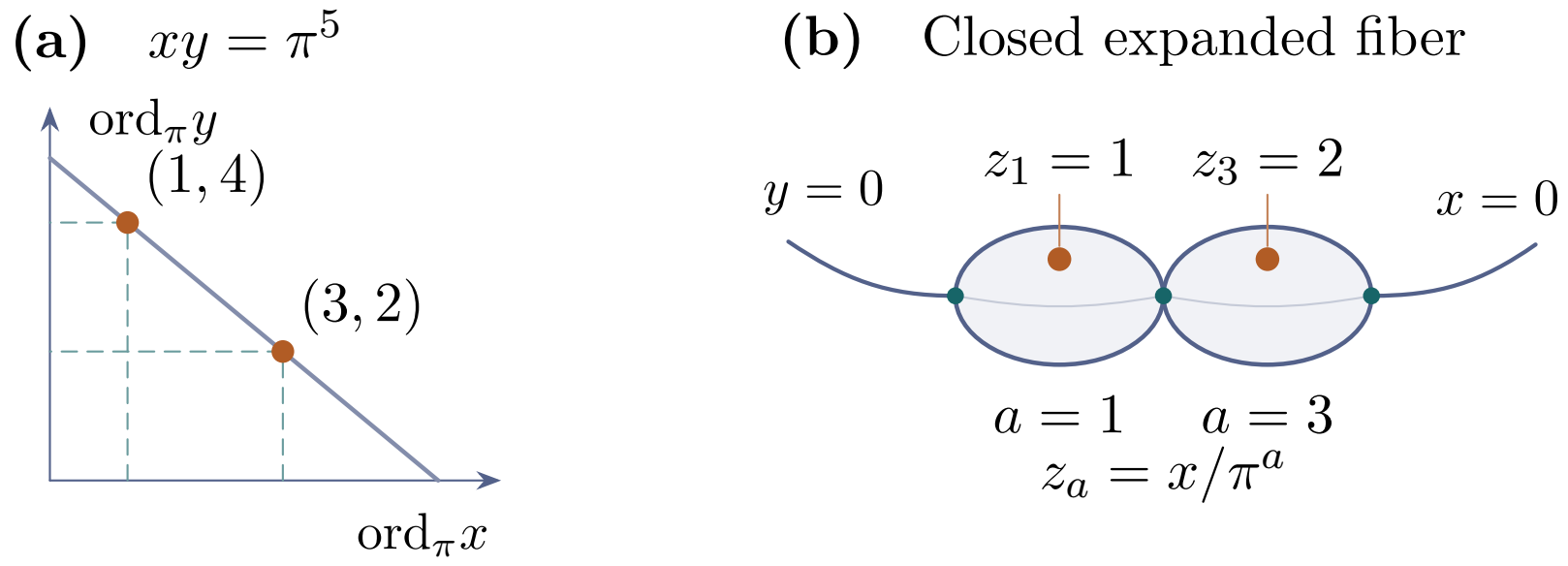
\begin{figure}[htbp]
\centering
\begingroup
\definecolor{levelaxis}{HTML}{53618A}
\definecolor{levelseam}{HTML}{176568}
\definecolor{levelpt}{HTML}{B15C25}
\begin{tikzpicture}[x=1cm,y=1cm,font=\small,line cap=round,line join=round]
\node[anchor=west] at (-.42,3.18) {\textbf{(a)}\quad $xy=\pi^5$};
\node[anchor=west] at (4.83,3.18) {\textbf{(b)}\quad Closed expanded fiber};
\draw[-{Stealth[length=1.7mm]},levelaxis] (0,.18)--(0,2.73);
\draw[-{Stealth[length=1.7mm]},levelaxis] (0,.18)--(3.08,.18);
\node[anchor=west,font=\footnotesize] at (.12,2.63) {$\mathrm{ord}_\pi y$};
\node[anchor=east,font=\footnotesize] at (3.14,-.19) {$\mathrm{ord}_\pi x$};
\draw[levelaxis!72,line width=.85pt] (0,2.38)--(2.65,.18);
\foreach \a/\b in {1/4,3/2}{
 \draw[levelseam!64,densely dashed] ({.53*\a},.18)--({.53*\a},{.18+.44*\b})--(0,{.18+.44*\b});
 \fill[levelpt] ({.53*\a},{.18+.44*\b}) circle (2.2pt);
 \node[anchor=south west,inner sep=3pt] at ({.53*\a},{.18+.44*\b}) {$(\a,\b)$};}
\begin{scope}[shift={(5.28,0)}]
 \draw[levelaxis,line width=.8pt] (-.24,1.81) to[out=-34,in=180] (.9,1.44);
 \draw[levelaxis,line width=.8pt] (3.74,1.44) to[out=0,in=215] (4.86,1.79);
 \foreach \a/\lev/\zl in {1.61/1/1,3.03/3/2}{
   \path[draw=levelaxis,fill=levelaxis!8,line width=.85pt] (\a,1.44) ellipse (.71 and .47);
   \draw[levelaxis!33] ({\a-.71},1.44) to[out=-10,in=190] ({\a+.71},1.44);
   \fill[levelpt] (\a,1.69) circle (2.3pt);
   \draw[levelpt!76,thin] (\a,1.78)--(\a,2.13);
   \node[anchor=south,inner sep=1.8pt] at (\a,2.13) {$z_{\lev}=\zl$};
   \node[anchor=north] at (\a,.90) {$a=\lev$};}
 \foreach \q in {.9,2.32,3.74}{\fill[levelseam] (\q,1.44) circle (1.65pt);}
 \node[anchor=south,font=\footnotesize] at (.00,1.85) {$y=0$};
 \node[anchor=south,font=\footnotesize] at (4.60,1.85) {$x=0$};
 \node at (2.32,.22) {$z_a=x/\pi^a$};
\end{scope}
\end{tikzpicture}
\endgroup
\caption{Two point defects approaching the same point of a joining K3, with $x=\pi$ and $x=2\pi^3$. Their orders determine the distinct occupied levels $a=1,3$, giving smoothing parameters $\pi,\pi^2,\pi^2$. The drawing at right shows the two ruling fibers over that K3 point. The reductions of the units in $x/\pi^a$ determine the marked interior points. Several points with the same valuation would use the same level.}
\label{fig:defectlevels}
\end{figure}

The valuation levels in the proof are illustrated in Figure~\ref{fig:defectlevels} for two distinct generic support points.

\begin{proposition}\label{prop:defect-sequence-defob}
On the refined geometric fiber $Z$ of
Proposition~\ref{prop:hull-punctual-DVR}, write the modification as
a triangle $F\xrightarrow{i}H'\to Q\to F[1]$ and put $A=\RHom_Z(F,F)$.
For the unrigidified problem retaining the exact sequence, define
\[
 L=\fib\!\left(A\oplus\RHom_Z(H',H')
       \xrightarrow{(a,b)\mapsto i a-b i}\RHom_Z(F,H')\right).
\]
The differentials of the two forgetful morphisms on the derived stack
of exact sequences fit into triangles
\begin{align*}
 L[1]&\longrightarrow\RHom_Z(H',H')[1]\longrightarrow\RHom_Z(F,Q)[1]
       \longrightarrow L[2],\\
 L[1]&\longrightarrow\RHom_Z(F,F)[1]\longrightarrow\RHom_Z(Q,H')[2]
       \longrightarrow L[2].
\end{align*}
In the first triangle the tangent complex with fixed $H'$ is
$\RHom_Z(F,Q)$. In the second triangle the tangent complex
with fixed kernel is $\RHom_Z(Q,H')[1]$. The two correction
complexes have cohomology supported at the defect before taking global sections.
Restriction to any joining divisor identifies the arrow $i$ with an
equivalence of the two boundary complexes.

These statements hold
under derived base change for perfect families whose finite defect
is disjoint from the relative joining divisors.
\end{proposition}
\begin{proof}
An endomorphism of the arrow $i$ consists of endomorphisms of its
source and target and a homotopy between the two composites to
$\RHom_Z(F,H')$. Thus $L$ is its derived endomorphism complex
and $L[1]$ its tangent complex, by the square-zero calculation of
Section~\ref{sec:derived} applied to arrows. Since the cofiber
is functorial up to a contractible choice, this also calculates the
tangent of exact triangles.

Project $L$ onto the target endomorphism
summand. Its fiber is
\[
 \fib\bigl(\RHom_Z(F,F)\longrightarrow\RHom_Z(F,H')\bigr)
       =\RHom_Z(F,Q)[-1].
\]
Projection onto the source endomorphism summand instead has fiber
\[
 \fib\bigl(\RHom_Z(H',H')\longrightarrow\RHom_Z(F,H')\bigr)
       =\RHom_Z(Q,H').
\]
Shifting gives the two triangles. Both correction complexes are
formed by internal derived Hom with one argument $Q$, so they
vanish near every joint. All arrows are induced by composition of
maps of perfect complexes; derived base change and proper perfect
pushforward give the stated compatibility.

The second correction can be nonzero even on a single smooth
branch. At a smooth point $p$ of an integral support where $H'$
is invertible on that support, let $Q=k(p)$. Complete the local
threefold there and choose parameters with
$B=\mathbb C[[u,v,w]]$ and $H'_p=B/(w)$.
The Koszul resolution of $k(p)$ for $(u,v,w)$ gives
\[
 \Ext_B^j(k(p),B/(w))=
 \begin{cases}
    \mathbb C,&j=2,3,\\
    0,&j\ne2,3.
 \end{cases}
\]
Indeed $u,v$ is a regular sequence on $B/(w)$ and $w$ acts by zero.
Thus forgetting the exact sequence has a nonzero relative derived tangent
even here. To identify the virtual classes of the original surface sheaves,
one must use their obstruction morphism; replacing it by the obstruction
theory of the length-$n$ Quot scheme would discard the displayed term.
\end{proof}

\begin{corollary}\label{cor:detdefect-local-vfc}
For the determinantal family of Theorem~\ref{thm:matrixexpandedallr}
and every finite-length modification constructed in
Proposition~\ref{prop:hull-punctual-DVR}, mark any joining divisor
with nonzero boundary. Write $F_L,F_R$ for the sheaves on the two
resulting relative chains and $\lambda=(0,rh,-2r)$.

There exist
bounded scalar-rigidified open neighborhoods $U_L,U_R$ of these
objects in the respective stacks of simple relative sheaves, and
a smooth open space $J$ of stable K3 sheaves of vector $\lambda$,
with boundary evaluation maps $e_i:U_i\to J$. On a neighborhood
$V$ of the matched point in $U_L\times_J U_R$ the obstruction
theory of the marked surface sheaf gives
\[
 [V]^{\vir}
  =\left.\Delta_J^!\bigl([U_L]^{\vir}\times[U_R]^{\vir}\bigr)
      \right|_V \quad\text{in } A_0(V)_{\QQ}.
\]
Here the obstruction theories are those of the sheaves on the
expanded threefolds. The assertion applies at every $r\ge2$ and
$n\ge0$, and imposes no coprimality condition on $\lambda$.
\end{corollary}
\begin{proof}
At a positive cut the boundary is the sheaf $E$ of
Theorem~\ref{thm:matrixexpandedallr}; the inserted ruled surfaces
have restrictions $p^*E$, and all defects avoid the joints.
The sheaf $E$ is torsion-free of rank one on its integral support
curve, hence stable. On every occupied integral surface of either
half, write $B_i$ for the corresponding rank-one torsion-free sheaf.
Localization at its generic point embeds
$\End(B_i)$ into its function field.

This finite-dimensional unital
domain over $\mathbb C$ equals $\mathbb C$. The normalization
triangle on either half injects its endomorphisms into
$\bigoplus_i\End(B_i)$ since its boundary sheaves have no negative
self-Ext. For the occupied vertices $I$ and internal edges $J$ of either half the normalization triangle gives
\[
 0\to\End(F_{\rm half})\to\bigoplus_{i\in I}\CC
 \xrightarrow{\partial}\bigoplus_{e\in J}\End(E_e),
 \qquad \partial(a)_e=(a_{i}-a_{i+1})\id_{E_e}.
\] The occupied components
of each half form a connected chain; zero components add no
endomorphisms. Thus $\End(F_L)=\End(F_R)=\mathbb C$. The full
target stabilizer is finite by the preceding defect computation, so
both half stabilizers are finite by Proposition~\ref{prop:stabilizerproductfull}.

Purity, admissibility and stable boundary restriction are open in
the flat perfect loci under consideration. Restrict the independent
quasi-compact charts of Section~\ref{sec:derived} around $F_L,F_R$
to these conditions, simplicity, and finite geometric inertia;
the last condition here follows locally by upper semicontinuity of
the infinitesimal stabilizer dimension in characteristic zero.

These are finite-type Deligne--Mumford neighborhoods after scalar
rigidification. Each $F_i$ is on an anticanonical relative chain,
so \eqref{eq:branchthirdnew} gives
$\Ext^3(F_i,F_i)^\vee=\Ext^{-1}_S(E,E)=0$.
The relative obstruction theories therefore have amplitude
$[-1,0]$ by Proposition~\ref{prop:relativePOT}. On the stable
locus of vector $\lambda$ on $S$, the trace in degree two is
dual to $\mathbb C\xrightarrow{1\mapsto\id_E}\End(E)$.

Consequently the reduced boundary stack is smooth there; take
$J$ to be this stable locus, shrinking the branch neighborhoods
as necessary.

Take $V$ to be an open neighborhood in the matching on which the
glued sheaf is simple. Theorem~\ref{thm:derivedglue} and
Proposition~\ref{prop:redglue} identify its derived enhancement
with the corresponding open in
$\mathbf U_L\times^{\RR}_J\mathbf U_R$. In particular the virtual tangent complexes on $V$ satisfy
\[
 T_V^{\vir}\to T_{U_L}^{\vir}|_V\oplus T_{U_R}^{\vir}|_V
       \xrightarrow{\mathrm de_L-\mathrm de_R} T_J|_V\to T_V^{\vir}[1].
\]
Dualizing gives the two rows of
\eqref{eq:POTdiagram} on this open. The diagonal of $J$ is regular;
virtual pullback for this compatible triangle gives the displayed
equality \cite[Theorem~4.8 and Corollary~4.9]{Mano}. Since
$\langle\lambda,\lambda\rangle=4r^2$, the virtual dimension is
$2(1+2r^2)-(2+4r^2)=0$. The comparison restricts further to any
open stable chart through the point. It is an equality in the Chow
group of the indicated neighborhood.

Pushing classes to a point
or summing the marked pieces over a special fiber requires the proper
moduli morphisms and the specialization of Section~\ref{sec:specialization}.
\end{proof}

\subsection{Verification in the structure-sheaf charge}\label{sec:verifiedcompactification}
For the hyperplane structure-sheaf character the two obstructions of the introduction have a concrete resolution. On the smooth quintic Lemma~\ref{lem:vacuumgenericfull} identifies the entire moduli space with $|H|$. The residual stabilizer in Corollary~\ref{cor:noseparatedall} does not arise in this character. On the one-bubble chart we choose a polarization with positive ruling degree and verify stability by the strict cut inequalities. On the overlap the same sheaf is stable for the induced unexpanded polarization. This verifies descent of the stability locus despite Proposition~\ref{prop:unlabelledpolarization}. The explicit inverse to the divisor construction below then gives (K3) and (K4), including for families whose generic target is already expanded.

Put
\[
 v_0=\ch(\OO_H)=1-e^{-H},\qquad
 P_0(m)=10m^2-5m+5,
\]
where the polynomial is computed with the original line bundle $\mathcal A$ of Lemma~\ref{lem:ample}.  Let $V=H^0(\PP^4,\OO(1))$, write $P=\PP(V)$, and let $p_0=[f_1]\in P$.  The convention for $P$ is that its points are one-dimensional spaces of equations.  Denote the smooth degeneration base by $C$ as before.

\begin{definition}\label{def:Cartierfunctor}
Let $\mathfrak D\to C$ be the following scalar-rigidified fibered groupoid.  An object over $T\to C$ is an expansion $Z_T\to T$ and a flat sheaf which, locally on $T$, is
\[
 \OO_{\mathcal T}\otimes p_T^*A_T,
\]
where $A_T$ is invertible and $\mathcal T\subset Z_T$ is a relative effective Cartier divisor.  Its fibers are admissible, have original Hilbert polynomial $P_0$, and have finite intrinsic expansion stabilizer.  Isomorphisms include isomorphisms of targets over the original family.  The scalar rigidification identifies the ambiguity of tensoring by a line bundle from the base.

The relative groupoids $\mathfrak D_Q$, $\mathfrak D_B^+$ and $\mathfrak D_B^0$ are defined in the same manner, with, respectively, contracted divisor classes $H_Q$, $H_B$ and $5H_B-E_B$; the first two have boundary $\OO_C$, and the last has zero outer restriction.  The original generic sheaf class is retained in each flat family, so the entire divisor class, and not merely a leading coefficient, is fixed.
\end{definition}

\begin{theorem}\label{thm:verifiedcompactification}
Assume that the semistable quintic is the model of Section~\ref{sec:geometry} and that $\Pic(S)=\mathbb Zh$.  Then
\begin{equation}\label{eq:vacuumblowup}
 \mathfrak D\simeq\mathscr P:=\operatorname{Bl}_{(p_0,0)}(P\times C).
\end{equation}
In particular $\mathfrak D$ is a smooth scheme, projective and separated over $C$.  Its general fiber is the \emph{complete} stable-sheaf moduli of character $v_0$, and its central fiber is
\begin{equation}\label{eq:vacuumcentral}
 \mathscr P_0=P_+\cup_{E_P}P_-,\qquad
 P_+=\operatorname{Bl}_{p_0}\PP^4,\quad P_-=\PP^4,\quad E_P=\PP^3.
\end{equation}
All its objects are Gieseker-stable for the chart polarizations specified below.  These stable loci agree under the expansion groupoid.  Every nonzero boundary restriction is pure, no support transversality is imposed, at most one bubble occurs, and every intrinsic geometric stabilizer is trivial.

For the twists in the conclusion the chart polarizations are chosen for that twist. Moreover the independent relative functors are represented by
\begin{equation}\label{eq:vacuumrelativefactors}
 \mathfrak D_Q=\operatorname{Bl}_{[S]}|H_Q|,\qquad
 \mathfrak D_B^+=|H_B|,\qquad
 \mathfrak D_B^0=|5H_B-E_B|.
\end{equation}
The last linear system is interpreted with the relative expansion at its hyperplane of divisors containing $S$.  Thus the three spaces in \eqref{eq:vacuumrelativefactors} are projective schemes.  Restriction at a marked node identifies the two numerical marked pieces with $P_+$ and $P_-$, respectively.  Consequently (K1)--(K4), including the assertions for the relative factors, hold for this groupoid and for its twists by $\OO(aH)$.
\end{theorem}

We construct the charts in \eqref{eq:vacuumblowup} and their inverse on arbitrary parameter schemes.

\begin{lemma}\label{lem:vacuumgenericfull}
For every smooth quintic $X$, the scalar-rigidified moduli of stable sheaves of character $v_0$ is $|H|=\PP^4$, with universal sheaf the structure sheaf of the universal hyperplane surface, up to a line bundle from the parameter space.
\end{lemma}
\begin{proof}
By Proposition~\ref{prop:primitive}, write $F=i_*G$ on an integral normal hyperplane surface $T$.  Put $R=G^{**}$ and $Q=R/G$.  The sheaf $R$ is reflexive of rank one, and $Q$ is finite.  On a general smooth plane quintic $D\subset T$, avoiding the finite singular and nonfree loci, the Hilbert-polynomial difference is
\[
 P_G(m)-P_G(m-1)=5m-5.
\]
Since $g(D)=6$, the restriction formula $\chi(R|_D(m))=\deg(R|_D)+5m+1-6$ gives $\deg(R|_D)=0$.  Also $\chi(G)=5$.

A nonzero section of a reflexive rank-one sheaf of degree zero has no divisorial zero locus: every nonzero effective Weil divisor has positive degree against $h$.  It identifies that sheaf with $\OO_T$ away from finitely many points, and reflexivity extends the identification.  Therefore if $R\not\simeq\OO_T$, then
\[
 H^0(T,R)=H^0(T,R^*)=0.
\]
Adjunction and Serre duality give $h^2(R)=h^0(R^*(h))$.  In the sequence
\[
 0\longrightarrow R^*\longrightarrow R^*(h)
   \longrightarrow R^*(h)|_D\longrightarrow0
\]
the last sheaf is a line bundle of degree $5$ on the smooth genus-six curve $D$.  Riemann--Roch gives Euler characteristic zero.  If it has a section, it is special, and Clifford's inequality gives $h^0\le3$; if it has no section the same bound holds.  Thus
\[
 \chi(R)=-h^1(R)+h^2(R)\le3,
 \qquad
 5=\chi(G)=\chi(R)-\operatorname{length}Q\le3,
\]
a contradiction.  We obtain $R=\OO_T$, then $Q=0$.

We also check the family assertion.  All geometric fibers have $(h^0,h^1,h^2)=(1,0,4)$.  Cohomology and base change therefore make $p_*F$ a line bundle, with arbitrary base change in degree zero.  The evaluation $p^*p_*F\to F$ is fiberwise surjective, hence surjective.  After removing $p_*F$, its kernel is flat over the base and is an invertible ideal on every fiber.  Locally choose a lift of a fiber generator of that ideal.

Nakayama gives a surjection from a free rank-one module; base-flatness makes its kernel reduce to zero, and Nakayama again makes the kernel zero.  Thus the ideal is invertible on the total family and defines a flat Cartier divisor.  The divisor class is $H$: $H^1(X,\OO_X)=H^2(X,\OO_X)=0$ makes its relative Picard class locally constant, with only the ambiguity of tensoring by a line bundle from the base.

Its equation is consequently a relative hyperplane equation.  These constructions are inverse and commute with base change, proving the functorial identification.
\end{proof}

\begin{lemma}\label{lem:vacuumprofiles}
An object of $\mathfrak D$ has either no bubble, or one graph bubble and no sheaf on the quartic root.  Its supporting divisor is respectively a matched hyperplane and plane transform, a divisor of class $5H_B-E_B$ disjoint from $S$, or the union of the graph and a plane transform.  All these supports are connected.  Their structure sheaves have
\[
 h^0=1,\qquad h^1=0,\qquad h^2=4.
\]
\end{lemma}
\begin{proof}
Apply Theorem~\ref{thm:primitiveprofiles}.  A vertical bubble supporting $\OO_{p^{-1}C}$ is fixed by the whole fiber-rescaling torus.  Its identity restrictions on the two boundary structure sheaves extend that automorphism to the complete matched sheaf.  This contradicts finite intrinsic stabilizer.  Hence there is no vertical bubble.  The classification leaves precisely the three alternatives in the statement.

A graph $z=-\ell$ with $\ell|_S\ne0$ has trivial rescaling stabilizer, because $\lambda\ell=\ell$ implies $\lambda=1$.

For a hyperplane surface in $Q$, the divisor sequence and $K_Q=-H_Q$ give cohomology $(1,0,1)$.  A plane transform in $B$ has cohomology $(1,0,0)$, since $R\rho_*\OO_B=\OO_{\PP^3}$ and $\OO_{\PP^3}(-1)$ is acyclic.  The joining plane quartic has cohomology $(1,3)$.  In the normalization sequence the map $\CC^2\to\CC$ is surjective, so the union has cohomology $(1,0,4)$.  The same computation applies to the graph, which is isomorphic to $S$.

For a divisor $T\in|D_B|$, $D_B=5H_B-E_B$, use
\[
 0\to\OO_B(-D_B)\to\OO_B\to\OO_T\to0,
 \qquad K_B+D_B=H_B.
\]
Serre duality and $H^{>0}(B,\OO_B(H_B))=0$ give $H^3(B,\OO_B(-D_B))=\CC^4$ and all its other cohomology groups zero.  The asserted cohomology of $\OO_T$ follows.  The normalization argument uses the genus and Cartier property of the joining curve, not its smoothness.
\end{proof}

\subsubsection{The charts and their target families}
Choose a complement $V'$ to $\CC f_1$ in $V$, with basis $\ell_1,\ldots,\ell_4$.  In the affine neighborhood of $p_0$, a hyperplane has equation
\[
 f_1+\sum_{i=1}^4b_i\ell_i=0.
\]
The center of \eqref{eq:vacuumblowup} is the ideal $(t,b_1,b_2,b_3,b_4)$.

On its $t$-chart write $b_i=tc_i$, and put $\ell_c=\sum c_i\ell_i$.  No target expansion is required.  The supporting surface in the unresolved family is specified by
\begin{equation}\label{eq:vacuumtchart}
 f_1+t\ell_c=0,\qquad f_5+f_4\ell_c=0.
\end{equation}
These equations imply $f_1f_4=tf_5$.  The center ideal $(f_1,t)$ of the small resolution restricts to the invertible ideal $(t)$ on this surface family.  Thus the surface lifts to the semistable target.  Intrinsically, its divisor is the section
\[
 f_1/t+\ell_c\quad\hbox{of}\quad \OO_{\XX}(H+Y_1).
\]
The expression is a rational section of $\OO(H)$ whose only allowed pole is on $Y_1$; it is a regular section of the displayed invertible sheaf.  At $t=0$, its value on $Y_1$ is nonzero, and its divisor on $B$ is the transform of $f_5+f_4\ell_c=0$.  This divisor is integral and disjoint from $S$, by the argument in Theorem~\ref{thm:primitiveprofiles}.

On the $b_j$-chart write
\begin{equation}\label{eq:vacuumbchart}
 b_j=s,\quad b_i=sc_i\ (i\ne j),\quad c_j=1,\qquad t=sv.
\end{equation}
Let $\ell_c=\sum c_i\ell_i$.  Use the one-bubble expansion obtained, locally at the double divisor, by blowing up $(x,s)$ in $xy=sv$.  Its two charts are
\begin{equation}\label{eq:vacuumtargetcharts}
 x=sz,\quad zy=v;
 \qquad s=xw,\quad y=wv.
\end{equation}
Here $x$ is the local equation of the right root and may be chosen as the pullback of $f_1$ near the double divisor.  The construction is globalized as the blowup of the sum of the ideal of that root and $(s)$.  Let $\mathscr Q$ be the Cartier divisor formed by the left root over $s=0$.  The supporting divisor is the section
\begin{equation}\label{eq:vacuumsection}
 f_1/s+\ell_c\quad\hbox{of}\quad
 c_X^*\OO_{\XX}(H)\otimes\OO(\mathscr Q).
\end{equation}
In \eqref{eq:vacuumtargetcharts}, its equations are respectively
\begin{equation}\label{eq:vacuumstrict}
 z+\ell_c=0,
 \qquad 1+w\ell_c=0.
\end{equation}
The second equation excludes the left root at the deepest fiber.  The first gives the graph $z=-\ell_c|_S$ on the bubble and the plane $\ell_c=0$ on the right root.  Their joining curve is the Cartier curve $\ell_c|_S=0$.  It can be singular; the equations still give flat Cartier surfaces and pure curve restrictions.

On $s\ne0$, contraction identifies the divisor with the original hyperplane $f_1+s\ell_c=0$.  On $v\ne0$, contraction identifies it with \eqref{eq:vacuumtchart} with coefficients $c_i/v$.  On the overlap with the $b_i$-chart, $c_i$ is invertible and the transition is
\begin{equation}\label{eq:vacuumoverlaps}
 s_i=sc_i,\qquad v_i=v/c_i,\qquad
 \ell^{(i)}=\ell_c/c_i.
\end{equation}
The target isomorphism is the corresponding rescaling of its bubble.  For a third index these ratios multiply, so the target and sheaf isomorphisms satisfy the cocycle condition.  Off the center the construction is the original hyperplane family.  The universal sheaf has therefore been constructed on all of $\mathscr P$.

\begin{lemma}\label{lem:vacuumflatness}
The constructed support is a relative effective Cartier divisor, flat over $\mathscr P$.  Its structure sheaf is admissible.  At every positive boundary its restriction is the structure sheaf of an integral Cartier curve in $|h|$, and therefore is pure.  The construction commutes with arbitrary base change.
\end{lemma}
\begin{proof}
The expanded target is flat over every chart.  In \eqref{eq:vacuumtchart} the section has no zero component on $Y_1$ and defines an effective Cartier divisor on $B$.  In \eqref{eq:vacuumstrict} the section is a unit along the unused root and is a regular equation on both nonzero branches.  Its restriction to their common $S$ is $\ell_c|_S$, a nonzero section of $\OO_S(h)$.  At a node the resulting two-term free resolution reduces to
\[
 0\to\OO_S(-h)\xrightarrow{\ell_c|_S}\OO_S\to\OO_{C_c}\to0,
\]
which is exact on the left.  Thus the quotient has vanishing $\Tor_1$ with the boundary; by Lemma~\ref{lem:flatnode} this is admissibility.  A nonzero section of $\OO_S(h)$ defines an integral Cartier curve because $\Pic(S)=\mathbb Zh$; its structure sheaf is Cohen--Macaulay of dimension one.

The defining section is regular on every target fiber.  The local flatness criterion applied to its length-one presentation proves flatness of its cokernel over the chart.  Its presentation also proves arbitrary base-change compatibility.  The checked overlap equations identify these presentations, so the assertions are global.  No differential of $\ell_c|_S$ is required to be nonzero.
\end{proof}

\subsubsection{Recovering the parameters from a family}
We give the inverse construction, including infinitesimal parameter schemes.  This is the step which proves that the charts describe the sheaf groupoid and not merely its geometric points.

\begin{lemma}\label{lem:vacuumlinearsystems}
On each of the three target types, the determinant line bundle of the support has $h^0=5$ and no higher cohomology.  On the corresponding expansion charts these groups form a rank-five vector bundle and commute with arbitrary base change.  Locally on the parameter, every line bundle of the indicated support class is the displayed one up to a line bundle from the parameter.
\end{lemma}
\begin{proof}
For the unexpanded ordinary type, the line bundle restricts to $\OO_Q(H_Q)$ and $\OO_B(H_B)$, with common restriction $\OO_S(h)$.  Normalization gives
\[
 h^0=5+4-4=5,\qquad H^{>0}=0.
\]
For the unexpanded zero-boundary type the restrictions are $\OO_Q$ and $\OO_B(5H_B-E_B)$, with common restriction $\OO_S$.  The blowup formula and the Koszul resolution of $I_\Gamma(5)$ give $h^0(B,\OO_B(5H_B-E_B))=5$, $H^2(B,\OO_B(5H_B-E_B))=\CC$, and no other positive cohomology.  In the exact sequence
\[
 0\to\OO_B(H_B)\to\OO_B(5H_B-E_B)\to\OO_S\to0
\]
the restriction on $H^2$ is an isomorphism, because $\OO_B(H_B)$ has no positive cohomology.  The degree-zero restriction is surjective.  In normalization both boundary groups are therefore accounted for, giving again $h^0=5$ and $H^{>0}=0$ on the complete union.

For the one-bubble type the restrictions are
\[
 \OO_Q,\qquad\OO_\Delta(D_+),\qquad\OO_B(H_B).
\]
Here $p_*\OO_\Delta(D_+)=\OO_S\oplus\OO_S(h)$ and $R^jp_*=0$ for $j>0$.  The boundary line bundles are $\OO_S$ on the left and $\OO_S(h)$ on the right.  The degree-zero normalization map is surjective, and
\[
 h^0=1+5+4-1-4=5.
\]
The only possible higher group on a component is the copy of $H^2(S,\OO_S)$ from the bubble.  Its restriction to the left boundary is the identity.  The higher cohomology of the line bundle on the union therefore vanishes.  The same calculations apply when one node is smoothed.  Cohomology and base change give the asserted vector bundle.

Finally, Proposition~\ref{prop:structurecohom} gives $H^1(Z,\OO_Z)=H^2(Z,\OO_Z)=0$.  In a square-zero extension with ideal $I$, the exact sequence
\[
 1\to 1+I\OO_Z\to\OO_{Z'}^*\to\OO_Z^*\to1,
 \qquad 1+I\OO_Z\simeq I\OO_Z
\]
shows that the line bundle has no infinitesimal deformations and no obstruction.  The relative Picard functor is consequently etale in these classes.  After an etale localization its class is the one specified by \eqref{eq:vacuumtchart} or \eqref{eq:vacuumsection}; the remaining ambiguity is an invertible sheaf pulled back from the base.
\end{proof}

\begin{proof}[Proof of \eqref{eq:vacuumblowup}]
In a flat sheaf family in Definition~\ref{def:Cartierfunctor}, Lemma~\ref{lem:vacuumprofiles} and base change recover the quotient
\[
 \OO_{Z_T}\longrightarrow\OO_{\mathcal T}
 \quad\hbox{from}\quad
 p_T^*p_{T*}F\longrightarrow F.
\]
Its kernel is an invertible ideal.  Thus the scalar-rigidified sheaf determines its divisor and conversely.  This recovery also proves that the structure-sheaf condition is an open condition near these objects in the bounded ambient flat-sheaf charts: vanishing of $H^1,H^3$, constancy of $\chi=5$, and upper semicontinuity force $h^0=1,h^2=4$ near such a point; surjectivity of evaluation and invertibility of its kernel are open.

By Lemma~\ref{lem:vacuumprofiles}, the target has at most one bubble.  On an unexpanded ordinary chart its divisor equation is a nonzero vector in $V$ up to scale, different from $f_1$ on the central fiber.  These are the usual charts of $P\times C$ away from $(p_0,0)$.  On the zero-boundary chart the coefficient of $f_1/t$ is nonzero; dividing by it gives uniquely $f_1/t+\ell_c$.  This is the $t$-chart.

On a chart containing a graph bubble, Lemma~\ref{lem:vacuumlinearsystems} gives the basis
\[
 f_1/s,\quad\ell_1,\ldots,\ell_4
\]
for the line bundle defining the support.  The coefficient of $f_1/s$ is a unit because the divisor avoids the left root.  Normalize it to one.  At least one coefficient of $\ell$ is a unit near the chosen point, because normality excludes the zero graph.  The expansion action is
\begin{equation}\label{eq:vacuumaction}
 (s,v,\ell)\longmapsto
 (\lambda s,\lambda^{-1}v,\lambda^{-1}\ell).
\end{equation}
Its invariants are $t=sv$ and $b_i=s\ell_i$.  On $\ell_j\ne0$ the unique gauge $\ell_j=1$ produces exactly \eqref{eq:vacuumbchart}.  These are statements about units, free modules of sections, and equations over the test base; they remain valid on nonreduced bases.  On overlaps the inverse transformations are precisely \eqref{eq:vacuumoverlaps}.  Thus the inverse morphisms of fibered groupoids glue and are mutually inverse to the universal-family construction.

The action in \eqref{eq:vacuumaction} is free on the graph locus.  On the unexpanded charts there is no expansion torus, and the only automorphisms of a connected structure sheaf are the scalars already rigidified.  The resulting groupoid is therefore the scheme $\mathscr P$.  Blowing up a smooth point in the smooth fivefold $P\times C$ proves smoothness. The blowup and the projection $P\times C\to C$ are projective, giving projectivity over $C$.

Its reduced central fiber is \eqref{eq:vacuumcentral}: the proper transform of $P\times\{0\}$ and the exceptional $\PP^4$ intersect in the exceptional $\PP^3$ of the former.
\end{proof}

\subsubsection{Stability, boundedness, and both valuative criteria}
It remains to check that this representation verifies the stability and compactification properties used later.  A scheme representing a different auxiliary functor would not suffice.

On an unexpanded ordinary support, the polarization $\mathcal A$ gives
\begin{equation}\label{eq:vacuumordinarypolys}
 \begin{split}
 P_Q&=2m^2+2,\qquad P_B=8m^2-m+1,\qquad B_S=4m-2,\\
 P&=10m^2-5m+5,\qquad D=P_Q-\tfrac15P=m+1.
 \end{split}
\end{equation}
Thus $0<D<B_S$ as eventual polynomials.  Theorem~\ref{thm:cutcriterion} proves stability.  A zero-boundary support is integral, so its rank-one structure sheaf is stable for every ample polarization.

For a one-bubble chart choose integers $0<q_0<q_1<M$.  If $\mathscr Q$ and $\mathscr B$ are its left and right root divisors, the line bundle
\begin{equation}\label{eq:vacuumlocalpolarization}
 L=M H+(M-q_0)\mathscr Q+(q_1-M)\mathscr B
\end{equation}
has restrictions
\[
 q_0H_Q,\qquad
 (q_1-q_0)D_-+q_1p^*h,\qquad
 (5M-4q_1)H_B-(M-q_1)E_B.
\]
They are ample by the projective-bundle and blowup calculations of Section~\ref{sec:idealchamber}; hence $L$ is relatively ample after shrinking the chart.  Its graph restriction is $q_1h$, because the graph misses $D_-$.  With $q=q_1$, the two active support components have total polynomial and cut defect
\begin{equation}\label{eq:vacuumgraphpolys}
 \begin{split}
 P&=\tfrac52M^2m^2-\tfrac52Mm+5,\qquad B_S=4qm-2,\\
 D&=\frac{2q^2}{M}m+2-\frac{4q^2}{M^2}.
 \end{split}
\end{equation}
Since $0<q<M$, its linear coefficient lies strictly between $0$ and $4q$.  Stability follows again.  On a punctured axis the induced unexpanded polarization has the corresponding parameter $q_0$ or $q_1$; the same strict inequality, or integral-support stability, applies.  Thus every object of the functor is stable for these chart polarizations, and the stable loci agree on all overlaps.

The line bundles themselves are not required to descend to one line bundle on the stack of expansions.

There is no strictly semistable object in this chosen functor.  Lemma~\ref{lem:vacuumflatness} proves that the admissibility error is identically zero on its constructed universal families.  The number of bubbles is at most one and the stabilizer is trivial.  A finite affine cover of $\mathscr P$, with the projective target charts just constructed, has a uniform sufficiently positive presentation for the universal sheaf by relative Serre vanishing.

Thus a finite collection of finite-type Quot charts suffices.  This verifies (K1)--(K2) and the intrinsic stabilizer assertion, without an unproved bound on non-invertible hulls.

For completeness the valuative extension is explicit.  A map from a fraction field to $P$ extends uniquely to $P$ over a discrete valuation ring $R$.  Near $p_0$ let the resulting coordinates be $t,b_1,\ldots,b_4$.  Unless they all vanish identically, set
\[
 r=\min\{\operatorname{ord}(t),\operatorname{ord}(b_1),\ldots,
                  \operatorname{ord}(b_4)\}.
\]
Divide these coordinates by a uniformizer to the power $r$.  At least one becomes a unit; its chart is exactly the corresponding blowup chart, and gives the unique extension to $\mathscr P$.  When the generic point maps to the center, its extra datum is already a point of the exceptional $\PP^4$; its homogeneous coordinates extend uniquely by the same minimum-valuation procedure.

This also treats families whose generic target is expanded and whose degeneration parameter is identically zero.  Pulling back the universal target and sheaf gives existence (K3).

Two extensions with a common generic isomorphism define two maps to the separated blowup $\mathscr P$.  They coincide on the generic point and hence on $\Spec R$.  The universal-family inverse proves that the resulting sheaf isomorphism extends in the scalar-rigidified groupoid.  The trivial intrinsic stabilizer makes it unique.  This proves (K4).

For a fixed twist $a\in\mathbb Z$, tensoring by $\OO(aH)$ transports
the Cartier families, their automorphisms and their valuative extensions.
Stability of the transported sheaves on the reducible supports requires a
choice of polarization. Use \eqref{eq:vacuumlocalpolarization} with
$q_0=M-2$, $q_1=M-1$ and $M>2+4|a-1/2|$. Choose the same $M$ and the
induced root parameter on each unexpanded chart.

On an ordinary or graph
cut, let $q$ denote the induced degree parameter on the quartic or graph
branch. The linear coefficient of $P_Q$ (respectively that of the graph)
increases by $4aq$. The linear coefficient of $P$ increases by $5Ma$
and the proportion of the leading term on this branch is
$2q^2/(5M^2/2)$. Consequently the cut defect $D_a$ satisfies
\[
 \frac{[m]D_a}{4q}
  =\frac{1}{4q}\left(\frac{2q^2}{M}
                    +4aq-\frac{4q^2}{5M^2}\,5Ma\right)
  =\frac12+(a-\tfrac12)(1-q/M).
\]
For $q=q_0,q_1$ its distance from $1/2$ is at most
$2|a-1/2|/M<1/2$. Hence $0<[m]D_a<4q=[m]B_S$.
The cut criterion gives stability, and an integral zero-boundary support
has a stable rank-one invertible sheaf for every polarization.
This proves the assertion on twists with these charge-dependent
chart polarizations.

\subsubsection{The independent relative factors and the marked matching}
It remains to construct the independent relative factors in Theorem~\ref{thm:verifiedcompactification}.

For $(Q,S)$ the parameter of hyperplane supports is $|H_Q|$.  Away from $[S]$ the support is its usual hyperplane surface.  Blow up $[S]$.  On a chart write its equation $f_1+s\ell_c=0$, with $\ell_c|_S\ne0$ modulo scaling.  Expand the pair by blowing up the ideal $(f_1,s)$.  On the support this ideal is $(s)$, so its strict transform is isomorphic to the original support family.  The central graph is $z=-\ell_c|_S$.

The inverse chart calculation is the restriction of \eqref{eq:vacuumaction} with the unused smoothing parameter removed.  This proves
\[
 \mathfrak D_Q=\operatorname{Bl}_{[S]}\PP^4
 \quad\hbox{over}\quad |h|=\PP^3.
\]
All its nonzero supports are integral; stability and pure boundary are automatic from the equations.  This is a projective, separated moduli functor, not just a resolution of a rational support map.

For $(B,S)$, the class $H_B$ is the pullback of a plane.  Every plane transform is integral and no such support contains $S$.  An additional vertical bubble would be fixed by rescaling, and a graph bubble would change the contracted divisor class.  Thus no bubble occurs and $\mathfrak D_B^+=|H_B|=\PP^3$.

Finally $H^0(B,\OO_B(5H_B-E_B))$ is generated by $f_5,f_4\ell_1,\ldots,f_4\ell_4$.  Write its projective parameters as $[c:\ell]$.  For $c\ne0$ its support is disjoint from $S$.  For $c=0$ the unexpanded divisor contains $S$ and a plane transform.  Expand the pair at this parameter hyperplane.  In a trivialization near $S$, the universal divisor has equation
\begin{equation}\label{eq:vacuumBzeroeq}
 c+u\ell=0,
\end{equation}
where $u$ is the equation of $S$ and $f_5$ supplies the nonzero trivialization of $(5H_B-E_B)|_S$.  The expansion is the blowup of $(u,c)$.  On \eqref{eq:vacuumBzeroeq} this ideal is $(u)$, an invertible ideal.  The strict support is consequently the same abstract family as the universal divisor on $B\times\PP^4$.  On $u=cz$ the support equation is $1+z\ell=0$, and on $c=uw$ it is $w+\ell=0$.  It avoids the outer boundary of the expansion.

At $c=0$ it is the graph together with the plane transform.  Its rescaling stabilizer is trivial.  The coefficient $[c:\ell]$ is recovered from these equations, including on arbitrary bases; this gives the third identification in \eqref{eq:vacuumrelativefactors}.

Let $J=J_{-2}$ and let $\sigma:|h|\to J$ be $C\mapsto\OO_C$.  The evaluations of $\mathfrak D_Q$ and $\mathfrak D_B^+$ are $\sigma\rho$ and $\sigma$, where $\rho:\operatorname{Bl}_{[S]}\PP^4\to\PP^3$ is projection.  Since the support map $J\to|h|$ composed with $\sigma$ is the identity, $\sigma$ is a section and hence a monomorphism. Thus
\begin{equation}\label{eq:vacuumpositiveproduct}
 \mathfrak D_Q\times_J\mathfrak D_B^+\simeq\mathfrak D_Q=P_+.
\end{equation}
The zero-boundary marked problem glues the zero sheaf on $Q$ to $\mathfrak D_B^0$ and is $P_-$.  At a graph-bubble object there are two nodes; marking the node adjacent to the empty quartic root gives the second piece, and marking the other node gives the first.  Therefore the full marked normalization is $P_+\sqcup P_-$.  All arrows, evaluations and identifications were obtained from universal divisors, so they are identifications of family groupoids.

This completes the proof of Theorem~\ref{thm:verifiedcompactification}.

\begin{remark}\label{rem:verifiedrange}
The theorem proves stable reduction and the independent proper matching factors for the whole character $\ch(\OO_H)$ and its line-bundle twists.  Its valuative proof includes singular hyperplane surfaces, singular boundary curves, a support which becomes the K3 and moves into a graph bubble, and a zero-boundary limit.
\end{remark}

\begin{figure}[htbp]
\centering
\begingroup
\definecolor{markaxis}{HTML}{53618A}
\definecolor{markseam}{HTML}{176568}
\definecolor{markpt}{HTML}{B15C25}
\begin{tikzpicture}[x=1cm,y=1cm,font=\small,line cap=round,line join=round]
\node[anchor=west] at (-.32,3.23) {\textbf{(a)}\quad A slice across $E_P$};
\node[anchor=west] at (5.48,3.23) {\textbf{(b)}\quad Marking the expanded target};
\draw[markseam,line width=1.15pt] (1.38,.17)--(1.38,2.57);
\draw[markaxis,line width=1.15pt] (-.23,1.16)--(3.62,1.16);
\fill[markpt] (1.38,1.16) circle (2.15pt);
\node[anchor=north west] at (1.47,1.09) {$E_P$};
\node[anchor=south east,text=markseam!78!black] at (1.31,2.51) {$P_-$};
\node[anchor=south west,text=markaxis!85!black] at (2.83,1.19) {$P_+$};
\node[anchor=west,font=\footnotesize,text=markseam!78!black] at (1.45,2.28) {$s=0$};
\node[anchor=north,font=\footnotesize,text=markaxis!85!black] at (3.01,1.15) {$v=0$};
\begin{scope}[shift={(5.72,0)}]
 \foreach \yy in {.81,2.14}{
  \draw[markaxis,line width=.8pt] (.07,\yy)--(4.14,\yy);
  \foreach \xx/\lab in {0/Q,2.12/\Delta,4.23/B}{
    \node[draw=markaxis!82,fill=white,rounded corners=2pt,inner xsep=5pt,inner ysep=2.7pt] at (\xx,\yy) {$\lab$};}
  \foreach \xx in {1.03,3.22}{
    \filldraw[fill=white,draw=markseam,line width=.85pt] (\xx,\yy) circle (2pt);}
 }
 \draw[markpt,line width=1.35pt] (1.03,2.14) circle (4.2pt);
 \draw[markpt,line width=1.35pt] (3.22,.81) circle (4.2pt);
 \node[anchor=south,font=\footnotesize] at (1.03,2.39) {$E=0$};
 \node[anchor=south,font=\footnotesize] at (3.22,1.06) {$E=\OO_C$};
 \node[anchor=west,text=markseam!78!black] at (4.72,2.14) {$P_-$};
 \node[anchor=west,text=markaxis!85!black] at (4.72,.81) {$P_+$};
\end{scope}
\end{tikzpicture}
\endgroup
\caption{The structure-sheaf character. Panel (a) fixes the direction $[\ell]$ in $E_P\simeq\PP^3$ and shows the two branches of $sv=0$. At their crossing the quartic root carries zero, the bubble carries the graph and the right root carries the plane transform. Panel (b) circles the chosen joining divisor on the very same target. The left mark has zero boundary and lies on the $P_-$ summand of the marked normalization. The right mark has boundary $\OO_C$ for $C=(\ell|_S=0)$ and lies on the $P_+$ summand.}
\label{fig:vacuumtwomarkings}
\end{figure}
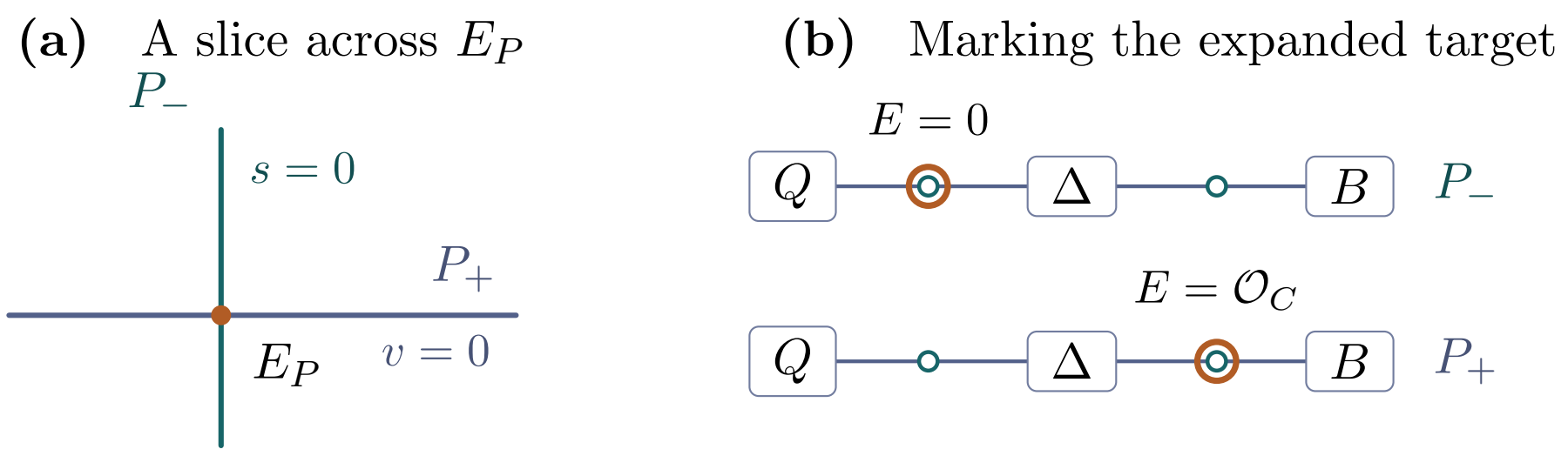

Figure~\ref{fig:vacuumtwomarkings} identifies the two markings over a point of $E_P$ in the chart $t=sv$.

\subsection{A point approaching a node of the expanded support} Suppose a family of relative Cartier support surfaces $T$ has already been extended, including a support that contracted to the middle K3. At a node on this expansion assume the Cartier equation of $T$ cuts a Cartier curve on the double K3. The curve may be singular. A point defect approaching that node gives the following local step for the sheaf itself.

\begin{proposition}\label{prop:point}
On a central expanded threefold $Z$, let $T$ be Cartier and $p\in T$ a closed point. Then $I_{p/T}$ is perfect as an $\OO_Z$-module exactly when $Z$ is regular at $p$. For a DVR section approaching a smoothed node with $xy=s^m$ and $0<a=\mathrm{ord}_s x(p)<m$, a single further bubble places the limit of the section at a regular point of the new fiber. The point is intrinsic to the resulting ideal and kills the new bubble dilation.
\end{proposition}
\begin{proof}
In the triangle associated to $0\to I_{p/T}\to\OO_T\to k(p)\to0$, the middle term has a Cartier resolution. Thus
\[
 I_{p/T}\in\perf(Z)\iff k(p)\in\perf(Z)
               \iff \OO_{Z,p}\text{ is regular}.
\]

Write $x(p)=s^au$, $y(p)=s^{m-a}u^{-1}$ with $u$ a unit. For the expansion with parameters $s^a,s^{m-a}$ the bubble interior is given by
\[
 x=s^az,\qquad y=s^{m-a}z^{-1},\qquad z(p)=u.
\]
The specialization $z(p)=\bar u\in\GG$ is regular on the expanded threefold fiber. The supporting equation restricts to the pullback of the assumed Cartier curve on the joining K3, and hence cuts a relative Cartier divisor near the new bubble. It contains the section by specialization. Flatness of $\OO_{\mathcal T}$ and of the section gives flatness of the kernel. Its special fiber is perfect at the mark by the equivalence above; elsewhere it is the Cartier structure sheaf.

Since $\OO_{T,p}$ has depth two, $\Ext^j_T(k(p),\OO_T)=0$ locally for $j<2$. Applying $\mathcal Hom_T(-,\OO_T)$ to the ideal sequence gives
\[
 I_{p/T}^{\vee}\simeq\OO_T,\qquad
 I_{p/T}^{\vee\vee}/I_{p/T}\simeq k(p).
\]
The ideal is faithful on $T$ (a section vanishing on the ideal would be supported at $p$ and hence zero by depth). Its annihilator in $\OO_Z$ therefore recovers $T$, and the double-dual quotient recovers $p$. A dilation preserving the ideal fixes $\bar u$, so it is the identity.
\end{proof}
\appendix
\section{A support specializing to the K3 divisor}\label{app:supportS}
When the entire supporting surface specializes to the relative K3 divisor, the following family exhibits the support and its sheaf on an expanded target.

Let $Y=Q$, let $S=(s_0=0)$ be the anticanonical K3, and choose $s_1\in H^0(Y,\OO_Y(S))$ with nonzero restriction $a=s_1|_S$. Consider the Cartier family
\begin{equation}\label{eq:pencilsupport}
\mathcal T=(s_0+\tau s_1=0)\subset Y\times\AA^1_\tau.
\end{equation}
Its central support is $S$. Blow up $Y\times\AA^1$ along $S\times\{0\}$ and denote the strict transform by $\widetilde{\mathcal T}$.

\begin{proposition}\label{prop:graphsupport}
The morphism $\widetilde{\mathcal T}\to\mathcal T$ is an isomorphism. In the central expanded relative target, the support is the graph of $-a$ in the normal-bundle chart of the bubble. It avoids the section glued to the root and meets the outer relative section in the Cartier curve
\[
C_a=(a=0)\subset S.
\]
The boundary curve need not be smooth. The geometric rescaling stabilizer of this graph is trivial.
\end{proposition}
\begin{proof}
On $\mathcal T$ the center ideal $(s_0,\tau)$ becomes $(\tau s_1,\tau)=(\tau)$. The family \eqref{eq:pencilsupport} is flat near zero: locally the special equation $s_0$ is a non-zero-divisor, so its Cartier deformation is flat. Thus $(\tau)$ is invertible on $\mathcal T$. By the universal property of the blowup, the embedding of $\mathcal T$ lifts to the blowup. Its image is the closure of the generic support and is isomorphic to $\mathcal T$, proving the first assertion.

Locally write $x=s_0$ in a trivialization. In the bubble chart $x=\tau z$, the strict-transform equation is
\[
z+s_1(\tau z,\ldots)=0.
\]
At $\tau=0$ it becomes $z=-a$. These local equations glue as the graph of a section of the normal line bundle $N_{S/Y}=\OO_S(h)$. It lies in the finite normal-bundle chart, so avoids the infinity section that meets the root. Its intersection with the zero section is $a=0$. A nonzero section on the smooth integral surface $S$ defines a Cartier curve, possibly singular or nonreduced; no smoothness is used.

A fiber rescaling sends $z$ to $\lambda z$. Preserving the graph requires $\lambda a=a$. Since $a\ne0$, this implies $\lambda=1$.
\end{proof}

\begin{figure}[tbp]
\centering
\begingroup
\definecolor{gSeam}{HTML}{176568}
\definecolor{gSurface}{HTML}{B15C25}
\definecolor{gRoot}{HTML}{53618A}
\begin{tikzpicture}[x=1cm,y=1cm,>=Stealth,line join=round,line cap=round,
  every node/.style={font=\small},guide/.style={gRoot!44,thin},
  graph/.style={gSurface,line width=1.45pt}]
\begin{scope}
\node[anchor=west] at (-2.60,2.21) {\textbf{(a)}\quad The support in the root};
\node[anchor=west] at (-2.60,1.68) {$s_0=-\tau u$};
\draw[guide,->] (0,-1.20)--(0,1.35) node[right,text=black] {$s_0$};
\draw[gSeam,very thick,->] (-2.34,0)--(2.42,0) node[right,text=black] {$u$};
\draw[graph] (-2.04,1.02)--(2.04,-1.02);
\draw[gSurface,line width=.81pt,dashed] (-2.16,.41)--(2.16,-.41);
\fill[gSeam!75!black] (0,0) circle (1.7pt);
\node[anchor=north west,inner sep=2pt,text=gSeam!72!black] at (-2.30,-.02) {$S$};
\node[anchor=north] at (0,-1.63) {$\tau\longrightarrow 0$};
\end{scope}
\draw[gRoot!80,->,line width=.8pt] (2.86,.25)--(3.68,.25);
\begin{scope}[shift={(6.85,0)}]
\node[anchor=west] at (-2.60,2.21) {\textbf{(b)}\quad The support on the bubble};
\node[anchor=west] at (-2.60,1.68) {$z=s_0/\tau$,\qquad $z=-u$};
\foreach \i in {-2,-1,0,1,2} {
  \draw[guide] (\i,-1.12)--(\i,.86);
  \draw[guide] (\i,1.01)--(\i,1.17);
}
\draw[gRoot,densely dashed,line width=.8pt] (-2.27,1.17)--(2.27,1.17);
\node[anchor=west,font=\footnotesize,inner sep=2pt,text=gRoot!85!black] at (2.27,1.17) {$\infty$};
\draw[gSeam,line width=1.1pt] (-2.27,0)--(2.27,0);
\node[anchor=west,font=\footnotesize,inner sep=2pt,text=gSeam!75!black] at (2.27,0) {$0$};
\draw[graph] (-1.58,.80)--(2.04,-1.04);
\draw[gRoot!78!black,fill=white,line width=.80pt] (0,0) circle (3.2pt);
\fill[gSurface] (0,0) circle (1.3pt);
\draw[gRoot!75,thin] (.11,.10)--(.53,.57)--(.80,.57);
\node[anchor=west,inner sep=2pt] at (.80,.57) {$u=0$};
\node[anchor=north] at (0,-1.63) {$\tau=0$};
\end{scope}
\end{tikzpicture}
\endgroup
\caption{A supporting surface specializing to the entire K3, shown near a
smooth point of $C_a$. Use $s_0$ and $u=s_1$ as local coordinates in a trivialization
and fix the coordinate along $C_a$. In (a) the two sloping lines are fibers
of the pencil for two positive values of $\tau$. Both tend to $s_0=0$.
In (b) the strict transform on the bubble is the graph $z=-u$.
It meets the zero section along $C_a$ (the circled point in the slice)
and is disjoint from the infinity section attached to the root.}
\label{fig:support-becomes-graph}
\end{figure}
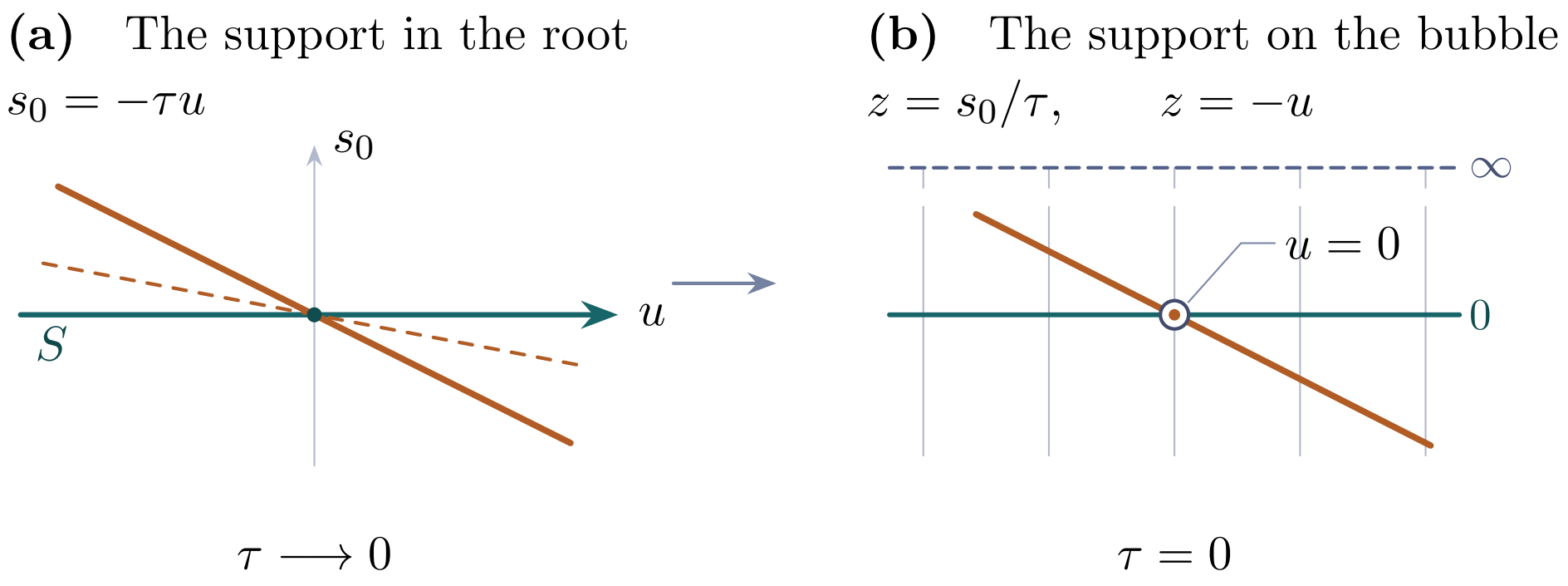

Figure~\ref{fig:support-becomes-graph} depicts the change of coordinates in the proof at a smooth point of the zero curve of $a$.

\begin{corollary}\label{cor:graphsheaf}
A line bundle on $\mathcal T$, or its tensor product with the ideal of a finite flat subscheme whose central support avoids $C_a$, transports to a flat sheaf on $\widetilde{\mathcal T}$. As a surface sheaf on the expanded relative threefold, it is normal to the distinguished divisors and has pure one-dimensional outer boundary. No assumption that $C_a$ is smooth is required.
\end{corollary}
\begin{proof}
The total support is unchanged by Proposition~\ref{prop:graphsupport}, so flatness is preserved. Near the outer boundary the sheaf is a line bundle on the graph. The boundary equation restricts to the non-zero-divisor $a$ on the graph surface. Its ordinary restriction is an invertible sheaf on the Cartier curve $C_a$, hence pure of dimension one. The support avoids the root node. Normality and boundary purity follow. Finite defects away from the boundary do not alter this local calculation.
\end{proof}

This verifies the stated extension for line bundles and finite defects away from the boundary. Section~\ref{sec:compactification} gives the stable-reduction hypotheses for general sheaves.

\section{The numerical stability calculation in explicit form}\label{app:algorithm}
For completeness, the stability test can be carried out over the rational numbers without a numerical approximation. Write
\[
P_i=a_iz^2+b_iz+c_i,\qquad B_e=d_ez+e_e,
\]
form $P=\sum P_i-\sum B_e$ and its coefficient $A$, and compute all $D_e=P_{\le e}-(A_e/A)P$. A degree-at-most-one polynomial is positive for large $z$ if its linear coefficient is positive, or if its linear coefficient is zero and its constant coefficient is positive. Apply this test to both $D_e$ and $B_e-D_e$.

For a vertical twist at $e$, let $\delta_e$ be the linear coefficient of $D_e$. Select the unique integer $k_0$ such that $0\le\delta_e+k_0d_e<d_e$. If this coefficient is nonzero the strict choice is determined. If it is zero, compute the constant coefficient using the recurrence \eqref{eq:twistrecur}; its sign gives the two alternatives of Theorem~\ref{thm:verticalbalance}. An equality is a genuine strict-semistability wall and must not be resolved by a floating-point sign choice.
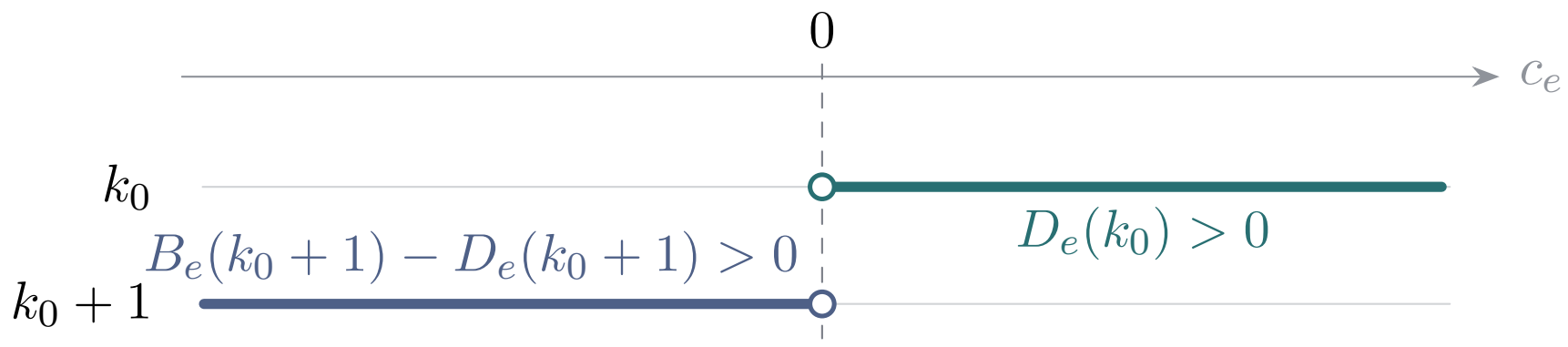
\begin{figure}[htbp]
\centering
\begin{tikzpicture}[x=1cm,y=1cm,font=\small,line cap=round]
 \draw[qiInk!55,-{Stealth[length=2mm]}] (1.15,.8)--(10.72,.8) node[right] {$c_e$};
 \foreach \y/\lab in {0/{k_0},-.85/{k_0+1}}{
  \node[anchor=east] at (1.08,\y) {$\lab$};\draw[qiInk!22] (1.3,\y)--(10.36,\y);}
 \draw[dashed,qiInk!65] (5.8,.89)--(5.8,-1.1);
 \node[above] at (5.8,.87) {$0$};
 \draw[qiTeal,line width=2.1pt] (5.8,0)--(10.30,0);
 \draw[qiBlue,line width=2.1pt] (1.31,-.85)--(5.8,-.85);
 \foreach \y/\cc in {0/qiTeal,-.85/qiBlue}{\draw[fill=white,draw=\cc,line width=.95pt] (5.8,\y) circle (2.5pt);}
 \node[below,qiTeal] at (8.13,0) {$D_e(k_0)>0$};
 \node[above,qiBlue] at (3.25,-.85) {$B_e(k_0+1)-D_e(k_0+1)>0$};
\end{tikzpicture}
\caption{The constant-coefficient step of the twist calculation. On a colored interval the indicated twist is stable at this cut. At the open endpoints both twists are strictly semistable there. The other cuts are tested independently.}\label{fig:constantwall}
\end{figure}

For a cut $m$, form the polynomials in \eqref{eq:leftcut} and \eqref{eq:rightcut} and repeat the test. This gives three separate assertions: stability of the whole sheaf, of the left relative sheaf, and of the right relative sheaf. The tests are finite on each numerical stratum, but familywise compatibility across changing expansion types is still part of the geometric moduli construction.

For example, at a cut with $\delta_e(k_0)=0$ the two potentially semistable twists are distinguished by a single constant. Figure~\ref{fig:constantwall} shows the strict choices on either side of this wall.

\subsection*{AI use declaration}
AI assisted with  proofreading and with Figures.

\end{document}